\documentclass[10pt]{amsart}
\usepackage{amsmath,amsfonts,amssymb,amsthm,amscd,amsxtra,latexsym}
\usepackage{mathrsfs,enumitem,fancyhdr,mathtools}
\usepackage{graphicx,needspace}
\usepackage{tikz}
\usepackage[plainpages=false]{hyperref}
\usepackage{accsupp}
\definecolor{wingInk}{HTML}{203638}
\definecolor{wingGrey}{HTML}{637879}
\definecolor{wingPale}{HTML}{F1F6F5}
\definecolor{wingAxis}{HTML}{A04B29}
\definecolor{wingMirror}{HTML}{225D99}
\definecolor{wingPath}{HTML}{087563}
\definecolor{wingOrbit}{HTML}{843E80}
\newcommand{\WingOrbitSphere}[3]{%
  \begin{scope}[shift={#1},scale=#2,rotate=-13]
    \shade[left color=#3!6,right color=#3!38,shading angle=20]
      (0,0) circle[radius=1];
    \begin{scope}
      \clip (0,0) circle[radius=1];
      \foreach \latitude in {-38,0,38}{
        \pgfmathsetmacro{\startangle}{asin(-tan(\latitude)*tan(16))}
        \pgfmathsetmacro{\endangle}{180-\startangle}
        \draw[#3!80,line width=.31pt]
          plot[domain=\startangle:\endangle,samples=32,variable=\ang]
          ({cos(\latitude)*cos(\ang)},
           {sin(\latitude)*cos(16)-cos(\latitude)*sin(\ang)*sin(16)});
      }
      \foreach \longitude in {35,90,145}{
        \pgfmathsetmacro{\firstlat}{-atan(sin(\longitude)/tan(16))}
        \draw[#3!85,line width=.36pt]
          plot[domain=\firstlat:90,samples=27,variable=\ang]
          ({cos(\ang)*cos(\longitude)},
           {sin(\ang)*cos(16)-cos(\ang)*sin(\longitude)*sin(16)});
      }
    \end{scope}
    \draw[#3,line width=.67pt] (0,0) circle[radius=1];
  \end{scope}%
}

\makeatletter
\renewcommand{\section}{\@startsection{section}{1}%
  \z@{.7\linespacing\@plus\linespacing}{.5\linespacing}%
  {\normalfont\bfseries\centering}}
\renewcommand{\@secnumfont}{\bfseries}
\makeatother

\numberwithin{equation}{section}
\newcommand{\beq}{\begin{equation}}
\newcommand{\eeq}{\end{equation}}
\newcommand{\beqs}{\begin{eqnarray*}}
\newcommand{\eeqs}{\end{eqnarray*}}
\newcommand{\beqn}{\begin{eqnarray}}
\newcommand{\eeqn}{\end{eqnarray}}
\newcommand{\beqa}{\begin{array}}
\newcommand{\eeqa}{\end{array}}

\def\lra{\longrightarrow}

\def\na{\nabla}

\def\lo{\vskip .1cm}

\def\p{\prime}

\def\bc{\begin{center}}
\def\ec{\end{center}}

\def\p{\partial}
\def\b{\bar}
\def\a{\alpha}
\def\g{\gamma}

\def\t{\tau}

\def\R{{\mathbb R}}

\def\tr{\rm tr}

\def\g{\mathfrak g}

\def\R{{\mathfrak R}}

\def\begeq{\begin{equation}}
\def\endeq{\end{equation}}
\def\and{\quad{\rm and}\quad}

\let\lra=\longrightarrow

\def\mapright\#1{\,\smash{\mathop{\lra}\limits^{\#1}}\,}

\def\lab{\ }
\def\lab{\label}

\def\Rm{\mathrm{Rm}}
\def\Ric{\mathrm{Ric}}

\def\tr{\mathrm{tr}}
\def\Scal{\mathrm{Scal}}
\newtheorem{prop}{Proposition}[section]

\newtheorem{lem}[prop]{Lemma}

\newtheorem{rem}[prop]{Remark}

\newtheorem{defi}[prop]{Definition}
\newtheorem{conj}[prop]{Conjecture}

\usepackage{amsthm}

\allowdisplaybreaks[1]
\numberwithin{equation}{section}
\newtheorem{theorem}{Theorem}[section]
\newtheorem{lemma}[theorem]{Lemma}
\newtheorem{proposition}[theorem]{Proposition}
\newtheorem{corollary}[theorem]{Corollary}

\newcommand{\dif}{\,d}
\newcommand{\Sec}{\operatorname{sec}}

\title[Classification of four-dimensional steady Ricci solitons]
  {Classification of four-dimensional $\kappa$-noncollapsed steady gradient Ricci solitons with nonnegative sectional curvature}
\author{Yuxing Deng$^{*}$}
\address{School of Mathematics and Statistics, Beijing Institute of Technology, Beijing, 100081, China}
\email{6120180026@bit.edu.cn}

\thanks {\raggedright * Supported by National Key R$\&$D Program of China 2022YFA1007600 and NSFC No.12571059.}
\date{}

\begin{document}
\begin{abstract}
\emergencystretch=1.7em
In this paper, we classify four-dimensional complete $\kappa$-noncollapsed  steady gradient Ricci solitons
with nonnegative sectional curvature.
Up to scaling, such a Ricci soliton  is  isometric to the Euclidean space
$\mathbb{R}^4$, the product of the three-dimensional Bryant soliton with a line,
the four-dimensional Bryant soliton, or a member of Lai's $O(3)\times\mathbb{Z}_2$-symmetric steady gradient Ricci solitons.
\end{abstract}
\maketitle

\section{Introduction}\label{sec:intro}
A complete Riemannian manifold $(M^n,g)$ is called a steady gradient Ricci soliton if there is a smooth function $f$ such that
\begin{align}
    \Ric=\na^2 f.\notag
\end{align}

Steady gradient Ricci solitons are important since they may arise as the type II singularity models of the Ricci flow. The classification of steady gradient Ricci solitons is also useful for the classification of the $\kappa$-solutions to the Ricci flow.  In dimension three, Perelman conjectured that any noncompact $\kappa$-solution with positive sectional curvature is a Ricci flow generated by the Bryant soliton \cite{Per02}. He also claimed that the three-dimensional  steady gradient Ricci soliton with positive sectional curvature is unique. In \cite{Bre13},
Brendle proved that every complete nonflat $\kappa$-noncollapsed steady gradient Ricci soliton is isometric to the Bryant soliton.  Brendle proved rotational symmetry for steady solitons with positive sectional curvature and asymptotically cylindrical condition\cite{Bre14}.  Zhu and the author generalized Brendle's work and classified steady solitons with linear scalar curvature decay\cite{DZ19, DZ2020, DZ20b}. They also obtained the rigidity of $\kappa$-noncollapsed steady K\"{a}hler-Ricci soliton with nonnegative bisectional curvature\cite{DZKR}. Based on these classification results, Brendle classified the noncompact three-dimensional $\kappa$-solutions to the Ricci flow. Up to scaling and isometry, they are the Bryant soliton, or a family of shrinking round cylinders $S^2\times\mathbb{R}$(or a quotient thereof)\cite{clBrendle}. Brendle, Daskalopoulos and \v{S}e\v{s}um classified the compact three-dimensional $\kappa$-solutions on $S^3$ to the Ricci flow. They are a family of shrinking round spheres $S^3$, or Perelman's ancient solution\cite{Per02,BDS21}. Zhu and the  author proved that every compact $\kappa$-solution of the K\"{a}hler-Ricci flow is a finite quotient of a product of irreducible compact Hermitian symmetric spaces\cite{DZCompact}. Li proved that every noncompact $\kappa$-solution of the K\"{a}hler-Ricci flow of complex dimension $m$ is isometrically biholomorphic to a finite quotient of $N^k\times\mathbb{C}^{m-k}$, where $N^k$ is a simply connected compact Hermitian symmetric space\cite{Li24}.

In dimension four, the classification of $\kappa$-noncollapsed steady gradient Ricci soliton with positive sectional curvature is still open. The canonical Ricci flow of a four-dimensional steady gradient Ricci soliton is similar to the  translating three-dimensional
hypersurfaces in $\mathbb R^4$.  Choi, Haslhofer and Hershkovits~\cite{CHH} classified
the complete nonflat noncollapsed translating three-dimensional
hypersurfaces in $\mathbb R^4$. Up to scaling and rigid motion, they
are the rotationally symmetric three-dimensional bowl, the
one-parameter family of oval bowls constructed by Hoffman, Ilmanen,
Mart{\'i}n and White~\cite{HIMW19}, and the product of the two-dimensional bowl
with a line. In the strictly convex case only the round bowl and the
oval bowls occur.

Haslhofer gives a conjecture for four-dimensional $\kappa$-solutions,
which have nonnegative curvature operator and positive scalar curvature\cite[Conjecture~1.3]{Has}. According to Haslhofer's conjecture, up to scaling and finite quotients, such steady Ricci solitons 
consist of the four-dimensional Bryant soliton, the three-dimensional
Bryant soliton times a line, and Lai's family. For steady Ricci solitons with positive curvature operator, it is natural to make the following conjecture.

\begin{conj}\label{conj:bryant-lai}
Every complete four-dimensional $\kappa$-noncollapsed steady gradient Ricci soliton
with positive curvature operator is, up to scaling and isometry, the four-dimensional Bryant
soliton or a member of Lai's $O(3)\times\mathbb Z_2$-symmetric family.
\end{conj}

Very recently, Fang proved the existence of $O(3)$ actions preserving the potential function\cite{F}. It is an important progress to Conjecture \ref{conj:bryant-lai}. In this paper, we will prove Conjecture \ref{conj:bryant-lai} and classify four-dimensional $\kappa$-noncollapsed steady gradient Ricci solitons with nonnegative sectional curvature, see Theorem \ref{cl:main} and Theorem \ref{th:classification-sec-nonnegative}.

Consider a complete four-dimensional $\kappa$-noncollapsed steady gradient Ricci solitons
$(M^4,g,f)$ with positive curvature operator. By Fang's theorem \cite[Theorem~1.1]{F},  the potential attains its minimum at a
unique point $o$. By scaling of the metric and translation of $f$, we may assume
\begin{align}\label{eq:normalization}
    \qquad R+|\nabla f|^2=1,\qquad \min f=0.
\end{align}

%\begin{defi}
%Let $(M^n,g)$ be a complete Riemannian manifold
%and let $\kappa>0$.
%We say that $(M^n,g)$ is $\kappa$-noncollapsed
%at all scales if, for every $p\in M$ and $r>0$,
%the ball satisfies the implication
%\[
 %\sup_{B_g(p,r)} |\Rm|_g\le r^{-2}
 %\quad\Longrightarrow\quad
 %\Vol_g(B_g(p,r))\ge\kappa r^n .
%\]
%\end{defi}

We first show that the steady Ricci soliton under consideration has uniform curvature decay.

\begin{theorem}\label{th:ray}
Let \((M^4,g,f)\) be a complete $\kappa$-noncollapsed  steady gradient Ricci soliton with positive curvature operator. For every \(p\in M\),
\[
 (M,\lambda g,p)\xrightarrow[\lambda\to0]{\mathrm{pGH}}
 ([0,\infty),d_{\mathrm{Eucl}},0),
\]
and \(R(x)\to0\) uniformly at infinity.
\end{theorem}

We can also give a sharp curvature decay estimate.

\begin{corollary}\label{cor:sphere-decay}
Under the assumption of Theorem \ref{th:ray}. Let $o$ be the minimum point of $f$. There are constants $C>0$ and
$r_0>1$ such that
\[
 \sup_{\partial B_g(o,r)}R\le C\frac{\log r}{r}\qquad(r\ge r_0).
\]
On the non-Bryant soliton with the normalization $R+|\nabla f|^2=1$, we have
\[
 \lim_{r\to\infty}\frac{r}{\log r}\sup_{\partial B_g(o,r)}R=1 .
\]
\end{corollary}

\begin{figure}[t]
\centering
\IfFileExists{figures/lai-action.pdf}{%
  \includegraphics{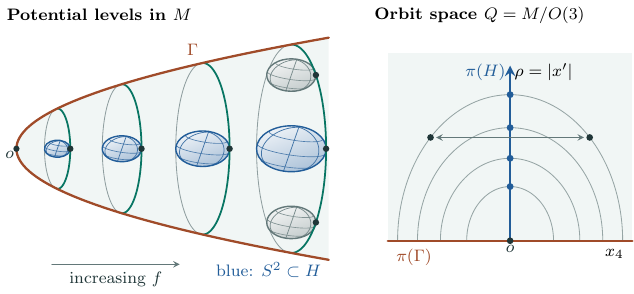}%
}{%
\begin{tikzpicture}[x=.95cm,y=.95cm,>=stealth,line cap=round,line join=round,
                    font=\footnotesize]
  \node[anchor=west] at (0,2.38) {\textbf{Potential levels in $M$}};
  \node[anchor=west] at (6.55,2.38) {\textbf{Orbit space $Q=M/O(3)$}};
  % A level is indexed by its horizontal drawing position. In the equivariant
  % model its points have x_4=b cos(t), x'=a sin(t) theta. Draw an orbit by
  % (X,Y)=(xc+x'_1,x_4+.66 x'_2). The marked representative has theta=e_1.
  \def\levels{1.03/.23,2.18/.35,3.62/.48,5.20/.62}
  \fill[wingPale] plot[domain=-1.98:1.98,samples=81] ({.30+1.42*\x*\x},\x)
       -- cycle;
  \foreach \xc/\a in \levels{
    \pgfmathsetmacro{\b}{sqrt((\xc-.30)/1.42)}
    \draw[wingGrey!72] plot[domain=0:180,samples=81]
       ({\xc-\a*sin(\x)},{\b*cos(\x)});
    \begin{scope}[shift={(\xc,0)},yscale=.66]
      \WingOrbitSphere{(0,0)}{\a}{wingMirror}
    \end{scope}
    \ifdim\xc pt>5pt
      \foreach \t in {45,135}{
        \pgfmathsetmacro{\rad}{\a*sin(\t)}
        \pgfmathsetmacro{\yc}{\b*cos(\t)}
        \begin{scope}[shift={(\xc,\yc)},yscale=.66]
          \WingOrbitSphere{(0,0)}{\rad}{wingGrey}
        \end{scope}
      }
    \fi
    \draw[wingPath,line width=.85pt] plot[domain=0:180,samples=181]
       ({\xc+\a*sin(\x)},{\b*cos(\x)});
    \fill[wingInk] ({\xc+\a},0) circle(1.55pt);
    \ifdim\xc pt>5pt
      \foreach \t in {45,135}
        \fill[wingInk] ({\xc+\a*sin(\t)},{\b*cos(\t)}) circle(1.55pt);
    \fi
  }
  \draw[wingAxis,line width=1.1pt] plot[domain=-1.98:1.98,samples=81]
       ({.30+1.42*\x*\x},\x);
  \fill[wingInk] (.30,0) circle(1.5pt) node[below left,inner sep=1pt] {$o$};
  \node[wingAxis,above] at (3.44,1.53) {$\Gamma$};
  \node[wingMirror,anchor=west] at (3.73,-2.16) {blue: $S^2\subset H$};
  \draw[->,wingGrey] (.93,-2.06) -- (3.21,-2.06)
      node[midway,below,text=wingInk] {increasing $f$};
  % The same samples are projected by (x',x_4) -> (x_4,|x'|) in plot order.
  % The scales of the two coordinate axes are independent.
  \begin{scope}[shift={(9.10,-1.64)}]
    \fill[wingPale] (-2.18,0) rectangle (2.18,3.35);
    \foreach \xc/\a in \levels{
      \pgfmathsetmacro{\b}{sqrt((\xc-.30)/1.42)}
      \draw[wingGrey!67,thin] plot[domain=0:180,samples=81]
           ({1.08*\b*cos(\x)},{4.2*\a*sin(\x)});
    }
    \draw[wingAxis,line width=1.1pt] (-2.18,0)--(2.18,0);
    \draw[wingMirror,line width=1.2pt,->] (0,0)--(0,3.12);
    \foreach \xc/\a in \levels
      \fill[wingMirror] (0,{4.2*\a}) circle(1.65pt);
    \pgfmathsetmacro{\pairx}{1.08*sqrt((5.20-.30)/1.42)*cos(45)}
    \pgfmathsetmacro{\pairy}{4.2*.62*sin(45)}
    \draw[<->,wingGrey] ({-\pairx+.10},\pairy)--({\pairx-.10},\pairy);
    \foreach \sgn in {-1,1}
       \fill[wingInk] ({\sgn*\pairx},\pairy) circle(1.6pt);
    \fill[wingInk] (0,0) circle(1.6pt) node[below,inner sep=2pt] {$o$};
    \node[wingAxis,anchor=north west] at (-2.16,-.02) {$\pi(\Gamma)$};
    \node[below] at (1.86,-.02) {$x_4$};
    \node[wingMirror,left,inner sep=2pt] at (0,3.00) {$\pi(H)$};
    \node[right,inner sep=2pt] at (0,3.00) {$\rho=|x'|$};
  \end{scope}
\end{tikzpicture}%
}
\caption{Sphere orbits in a four-dimensional soliton with $O(3)$-symmetry. At each level set,
the green meridian meets each $S^2$ at the black point, with
the same angular coordinate throughout. The poles lie on the
$O(3)$-fixed axis $\Gamma$. Suppose the soliton is $O(3)\times\mathbb{Z}_2$-symmetric. The reflection fixes each blue sphere
pointwise and interchanges the two gray spheres on the outer level.
The blue spheres over all positive levels, together with $o$, form $H$.
In $Q$ each level projects to an arc, each blue sphere to a point on
$\pi(H)$, and the gray pair to the endpoints of the double arrow.
The projection shows only the incidence relations, while the intrinsic angles are measured by $g$.}
\label{fig:lai-action}
\end{figure}

By \cite[Theorem~1.1]{F},  the soliton is $O(3)$-symmetric, see Figure \ref{fig:lai-action}.  It suffices to consider the non-Bryant soliton. If the uniform scalar curvature decay fails, there exists a sequence of points $p_i$  such that $R(p_i)$ has a uniform positive lower bound. Rescaling the soliton by $R(p_i)$ and taking limit, the limit must be the Bryant soliton times a line. Following the tip line of this limit, we see that the scalar curvature has a positive limit along one half of the fixed axis of the original soliton as in Figure \ref{fig:lai-action}.   Let $U(f)$ be the largest orbit radius at the height of $f$. On one hand, by integration of the level set equation, we can prove that
\begin{align*}
    U(f)^2\le 2f+Cf^{3/4}.
\end{align*}
On the other hand, such an upper bound is ruled out  in Proposition \ref{prop:noenv} due to  Lemma \ref{lem:shift} and Lemma \ref{lem:heat}. Hence, we get the  uniform curvature decay.
 Fang's sector theorem~\cite[Theorem~1.3]{F} shows that the opening angle
is zero if and only if the scalar curvature tends uniformly to zero at
infinity.  This gives Theorem \ref{th:ray}.

 For Corollary \ref{cor:sphere-decay}, it suffices to work on the potential level set $\Sigma_c=\{f=c\}$ of a non-Bryant soliton. By Proposition \ref{prop:noenv}, the cylindrical profile belongs to the neutral alternative. Many asymptotic estimates in \cite{MMS} hold without fast curvature decay assumption.  From the asymptotics, we know that the curvature scale at both poles is $k_c=\log c/c$. After rescaling by this scale,
the caps converge to a three-dimensional Bryant model with tip scalar curvature one. Hence, we have $R/k_c\to1$ at each pole. Away from the caps, the curvature is bounded by the inverse square of the orbit radius. Hence, the supremum on the whole level is asymptotic to $k_c$. Let $d=d_g(o,p)$. It's easy to check that $f/d\to1$ uniformly, so the upper bound also holds
on distance spheres. Finally, we note that $\log d/d=o(d^{-\eta}) $ for all $0<\eta<1$.

To give a complete classification, we also need to deal with the $\mathbb{Z}_2$-symmetry of the positively curved steady gradient Ricci solitons.
\begin{theorem}\label{th:ref}
Let \((M^4,g,f)\) be a complete $\kappa$-noncollapsed  steady gradient Ricci soliton with positive curvature operator. The \(O(3)\)-action extends to an
effective isometric \(O(3)\times\mathbb Z_2\)-action preserving \(f\). 
\end{theorem}

On a Bryant soliton, Theorem \ref{th:ref} is satisfied since  it is $O(4)$-symmetric. On a non-Bryant soliton, the maximal sphere divides each  level set into two sheets. Let's consider
points of the same height $f$ and orbit radius $r$ on the two sheets. To compare the quotient metrics at these points, we introduce
a weighted radial integral of the scalar curvature on each sheet. In this way we get two functions $b_+$ and
$b_-$, and the quotient metrics are determined by these functions together with their first derivatives.
We show that the difference of $b_+$ and $b_-$ satisfies a linear elliptic equation. By constructing a strict supersolution, we prove
that this difference decays exponentially at equal height and radius. Next we fix $(r,b)$ and compare the heights
at which the two sheets attain that value. If the two heights are different somewhere, then by the maximum principle a discrepancy
also occurs at arbitrarily high truncation boundaries. On the other hand, since $b_f=-K$, we can integrate
a lower bound for the quotient curvature $K$ along a sequence of high slabs, and this contradicts
the exponential decay. Thus, the profiles of the two sheets agree exactly. It follows that the map exchanging the two sheets at the same
angular coordinate preserves the metric. In a signed coordinate across the maximal sphere, this map is a smooth reflection.
It extends over the axis by metric completion, and we get
a global isometric involution.

By Theorem \ref{th:ref}, it remains to check whether a non-Bryant soliton with \(O(3)\times\mathbb Z_2\)-symmetry is isometric to one of the solitons constructed by Lai. Lai~\cite{Lai} constructed complete steady gradient Ricci solitons
with positive curvature operator and $O(n-1)\times\mathbb Z_2$
symmetry in dimension $n\geq3$. Her examples are not rotationally symmetric. In dimension three,
they realize Hamilton's conjectured \emph{flying wings} and the solitons are collapsed. In 
dimension $n\geq4$,  her examples are $\kappa$-noncollapsed. In dimension four, Lai prescribes
the scale-invariant tip ratio
\[
             \alpha(g)=
             \frac{\lambda_{\min}(\Ric_g(o))}
                  {\lambda_{\max}(\Ric_g(o))}.
\]
At the tip $o$, the Ricci tensor has an axial eigenvalue and a transverse
eigenvalue of multiplicity three. In the proof of
\cite[Theorem~1.1]{Lai}, for each $0<\alpha<1$, Lai chooses members of
successive expander paths with axial-to-transverse ratio $\alpha$ and
takes a pointed smooth subsequential limit. The steady
soliton obtained in this way has the prescribed tip ratio, which gives an existence result
for each $\alpha$. We fix the sequence of four-dimensional expanding gradient Ricci soliton paths used in the proof
of \cite[Theorem~1.1]{Lai}. For each $i$, denote the $i$th path as
$\{E_{i,\mu}:0\leq\mu\leq1\}$, where each $E_{i,\mu}$ is based at
its tip and scaled to have scalar curvature one there.

Using these paths, we define
\[
\mathcal L :=
\left\{
[M,g,f,o]\;\middle|\;
\begin{array}{l}
E_{i_j,\mu_j}\xrightarrow{\mathrm{pCG}} (M,g,o)
  \text{ for some } i_j\to\infty,\ \mu_j\in[0,1],\\g \text{ has positive 
  curvature operator,}\\
  g\text{ is nonrotationally symmetric.}
\end{array}
\right\},
\]
Here $[M,g,f,o]$ denotes the point isometry class of the 
complete and connected limiting four-dimensional steady gradient Ricci soliton,
with tip $o$ and potential $f$ normalized by\eqref{eq:normalization}. Two such limits $(M,g,f,o)$ and $(M',g',f',o')$ represent the same
element of $\mathcal L$ if and only if there is an isometry
$F:(M,g)\to(M',g')$ such that $F(o)=o'$ and $f'\circ F=f$. We regard $\mathcal L$ as the collection of Lai's examples.

\begin{figure}[!htb]
\centering
\IfFileExists{figures/equatorial-phase.pdf}{%
  \includegraphics[width=\linewidth]{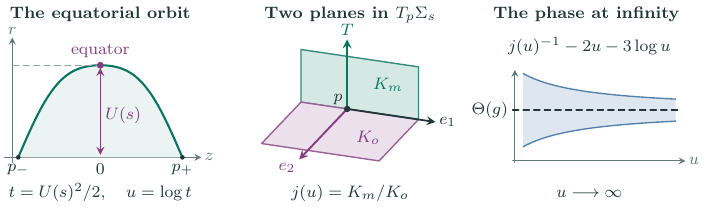}%
}{%
\resizebox{\linewidth}{!}{%
\begin{tikzpicture}[x=1bp,y=1bp,line cap=round,line join=round,
                    font=\footnotesize,
                    lab/.style={anchor=base west,inner sep=0pt}]
  \useasboundingbox (0,0) rectangle (338.99,100.19);
  % Left panel: the profile of an orbit sphere and its equator.
  \draw[wingGrey,line width=.478bp] (2.23,24.66)--(94.51,24.66);
  \fill[wingGrey] (96.63,24.66)--(93.25,26.35)--(94.51,24.66)--(93.25,22.97)--cycle;
  \draw[wingGrey,line width=.478bp] (5.92,22.1)--(5.92,80.09);
  \fill[wingGrey] (5.92,82.2)--(4.23,78.82)--(5.92,80.09)--(7.61,78.82)--cycle;
  \filldraw[fill=wingPath!8,draw=wingPath,line width=.996bp]
    (8.75,24.66)--(9.69,26.83)--(10.63,29)--(11.57,31.15)--(12.5,33.28)--(13.44,35.38)--(14.38,37.44)--(15.32,39.47)--(16.25,41.44)--(17.19,43.36)--(18.13,45.22)--(19.07,47.02)--(20.01,48.75)--(20.94,50.41)--(21.88,52)--(22.82,53.5)--(23.76,54.93)--(24.7,56.28)--(25.63,57.54)--(26.57,58.73)--(27.51,59.83)--(28.45,60.85)--(29.39,61.8)--(30.32,62.67)--(31.26,63.47)--(32.2,64.19)--(33.14,64.85)--(34.07,65.44)--(35.01,65.97)--(35.95,66.44)--(36.89,66.86)--(37.83,67.23)--(38.76,67.55)--(39.7,67.83)--(40.64,68.07)--(41.58,68.27)--(42.52,68.44)--(43.45,68.58)--(44.39,68.69)--(45.33,68.77)--(46.27,68.83)--(47.21,68.86)--(48.14,68.87)--(49.08,68.86)--(50.02,68.83)--(50.96,68.77)--(51.89,68.69)--(52.83,68.58)--(53.77,68.44)--(54.71,68.27)--(55.65,68.07)--(56.58,67.83)--(57.52,67.56)--(58.46,67.24)--(59.4,66.87)--(60.34,66.45)--(61.27,65.98)--(62.21,65.45)--(63.15,64.86)--(64.09,64.21)--(65.03,63.48)--(65.96,62.69)--(66.9,61.82)--(67.84,60.88)--(68.78,59.85)--(69.72,58.75)--(70.65,57.57)--(71.59,56.3)--(72.53,54.96)--(73.47,53.53)--(74.4,52.03)--(75.34,50.45)--(76.28,48.79)--(77.22,47.06)--(78.16,45.26)--(79.09,43.4)--(80.03,41.48)--(80.97,39.51)--(81.91,37.49)--(82.85,35.42)--(83.78,33.32)--(84.72,31.19)--(85.66,29.04)--(86.6,26.87)--(87.54,24.7);
  \draw[wingGrey!70,line width=.3985bp,dash pattern=on 2.989bp off 1.993bp]
    (6.2,68.88)--(48.15,68.88);
  \draw[wingOrbit,line width=.717bp] (48.15,28.32)--(48.15,65.11);
  \fill[wingOrbit] (48.15,25.85)--(50.13,29.8)--(48.15,28.32)--(46.18,29.8)--cycle;
  \fill[wingOrbit] (48.15,67.58)--(46.18,63.63)--(48.15,65.11)--(50.13,63.63)--cycle;
  \fill[wingOrbit] (48.15,68.88) circle[radius=1.645];
  \fill[wingInk] (8.75,24.66) circle[radius=1.115];
  \fill[wingInk] (87.56,24.66) circle[radius=1.115];
  \node[lab,wingInk] at (4.79,91.27) {\textbf{The equatorial orbit}};
  \node[lab,wingGrey] at (98.26,22.94) {$z$};
  \node[lab,wingGrey] at (3.89,83.83) {$r$};
  \node[lab,wingOrbit] at (50.6,43.07) {$U(s)$};
  \node[lab,wingOrbit] at (34.15,74) {equator};
  \node[lab,wingInk] at (3.49,17.93) {$p_-$};
  \node[lab,wingInk] at (46.04,16.22) {$0$};
  \node[lab,wingInk] at (82.38,17.93) {$p_+$};
  \node[lab,wingInk] at (4.21,5.38) {$t=U(s)^2/2,\quad u=\log t$};
  % Middle panel: the meridian plane and the orbit plane at $p$.
  \filldraw[fill=wingPath!17,draw=wingPath!83,line width=.6675bp]
    (144.48,51.18)--(200.87,42.84)--(200.87,68.1)--(144.48,76.43)--cycle;
  \filldraw[fill=wingOrbit!16,draw=wingOrbit!82,line width=.6675bp]
    (144.48,51.18)--(200.87,42.84)--(181.96,22.95)--(125.57,31.29)--cycle;
  \draw[wingPath,line width=.9066bp] (166.64,47.9)--(166.64,78.54);
  \fill[wingPath] (166.64,81.3)--(164.44,76.89)--(166.64,78.54)--(168.85,76.89)--cycle;
  \draw[wingOrbit,line width=.9066bp] (166.64,47.9)--(145.52,25.69);
  \fill[wingOrbit] (143.63,23.69)--(148.26,25.37)--(145.52,25.69)--(145.07,28.4)--cycle;
  \draw[wingInk,line width=.8568bp] (166.64,47.9)--(206.37,42.03);
  \fill[wingInk] (209.02,41.64)--(205.09,44.38)--(206.37,42.03)--(204.47,40.14)--cycle;
  \fill[wingInk] (166.64,47.9) circle[radius=1.545];
  \node[lab,wingPath] at (163.59,83.14) {$T$};
  \node[lab,wingOrbit] at (133.7,18.41) {$e_2$};
  \node[lab,wingInk] at (210.84,40.47) {$e_1$};
  \node[lab,wingInk] at (160.19,51.64) {$p$};
  \node[lab,wingPath] at (179.11,57.1) {$K_m$};
  \node[lab,wingOrbit] at (171.01,31.97) {$K_o$};
  \node[lab,wingInk] at (126.88,91.27) {\textbf{Two planes in} $T_p\Sigma_s$};
  \node[lab,wingInk] at (139.87,5.38) {$j(u)=K_m/K_o$};
  % Right panel: the band around the limit $\Theta(g)$.
  \draw[wingGrey,line width=.478bp] (247.15,22.39)--(247.15,64.5);
  \fill[wingGrey] (247.15,66.61)--(245.46,63.23)--(247.15,64.5)--(248.84,63.23)--cycle;
  \draw[wingGrey,line width=.478bp] (246.01,22.95)--(326.96,22.95);
  \fill[wingGrey] (329.07,22.95)--(325.69,24.64)--(326.96,22.95)--(325.69,21.26)--cycle;
  \fill[wingMirror!15]
    (251.12,64.91)--(252.67,64.08)--(254.23,63.33)--(255.78,62.64)--(257.34,62.01)--(258.9,61.43)--(260.45,60.89)--(262.01,60.39)--(263.56,59.93)--(265.12,59.5)--(266.67,59.09)--(268.23,58.72)--(269.79,58.37)--(271.34,58.03)--(272.9,57.72)--(274.45,57.43)--(276.01,57.15)--(277.57,56.88)--(279.12,56.63)--(280.68,56.4)--(282.23,56.17)--(283.79,55.96)--(285.34,55.75)--(286.9,55.56)--(288.46,55.38)--(290.01,55.2)--(291.57,55.03)--(293.12,54.86)--(294.68,54.71)--(296.24,54.56)--(297.79,54.41)--(299.35,54.28)--(300.9,54.14)--(302.46,54.02)--(304.01,53.89)--(305.57,53.77)--(307.13,53.66)--(308.68,53.55)--(310.24,53.44)--(311.79,53.34)--(313.35,53.24)--(314.91,53.14)--(316.46,53.05)--(318.02,52.96)--(319.57,52.87)--(321.13,52.79)--(322.68,52.71)--(324.24,52.62)--(324.24,42.04)--(322.68,41.96)--(321.13,41.88)--(319.57,41.8)--(318.02,41.71)--(316.46,41.62)--(314.91,41.52)--(313.35,41.43)--(311.79,41.33)--(310.24,41.22)--(308.68,41.12)--(307.13,41.01)--(305.57,40.89)--(304.01,40.77)--(302.46,40.65)--(300.9,40.52)--(299.35,40.39)--(297.79,40.25)--(296.24,40.11)--(294.68,39.96)--(293.12,39.8)--(291.57,39.64)--(290.01,39.47)--(288.46,39.29)--(286.9,39.11)--(285.34,38.91)--(283.79,38.71)--(282.23,38.49)--(280.68,38.27)--(279.12,38.03)--(277.57,37.78)--(276.01,37.52)--(274.45,37.24)--(272.9,36.95)--(271.34,36.63)--(269.79,36.3)--(268.23,35.95)--(266.67,35.57)--(265.12,35.17)--(263.56,34.74)--(262.01,34.28)--(260.45,33.78)--(258.9,33.24)--(257.34,32.66)--(255.78,32.03)--(254.23,31.34)--(252.67,30.59)--(251.12,29.76)--cycle;
  \draw[wingMirror!78,line width=.6575bp]
    (251.12,29.76)--(252.67,30.59)--(254.23,31.34)--(255.78,32.03)--(257.34,32.66)--(258.9,33.24)--(260.45,33.78)--(262.01,34.28)--(263.56,34.74)--(265.12,35.17)--(266.67,35.57)--(268.23,35.95)--(269.79,36.3)--(271.34,36.63)--(272.9,36.95)--(274.45,37.24)--(276.01,37.52)--(277.57,37.78)--(279.12,38.03)--(280.68,38.27)--(282.23,38.49)--(283.79,38.71)--(285.34,38.91)--(286.9,39.11)--(288.46,39.29)--(290.01,39.47)--(291.57,39.64)--(293.12,39.8)--(294.68,39.96)--(296.24,40.11)--(297.79,40.25)--(299.35,40.39)--(300.9,40.52)--(302.46,40.65)--(304.01,40.77)--(305.57,40.89)--(307.13,41.01)--(308.68,41.12)--(310.24,41.22)--(311.79,41.33)--(313.35,41.43)--(314.91,41.52)--(316.46,41.62)--(318.02,41.71)--(319.57,41.8)--(321.13,41.88)--(322.68,41.96)--(324.24,42.04);
  \draw[wingMirror!78,line width=.6575bp]
    (251.12,64.91)--(252.67,64.08)--(254.23,63.33)--(255.78,62.64)--(257.34,62.01)--(258.9,61.43)--(260.45,60.89)--(262.01,60.39)--(263.56,59.93)--(265.12,59.5)--(266.67,59.09)--(268.23,58.72)--(269.79,58.37)--(271.34,58.03)--(272.9,57.72)--(274.45,57.43)--(276.01,57.15)--(277.57,56.88)--(279.12,56.63)--(280.68,56.4)--(282.23,56.17)--(283.79,55.96)--(285.34,55.75)--(286.9,55.56)--(288.46,55.38)--(290.01,55.2)--(291.57,55.03)--(293.12,54.86)--(294.68,54.71)--(296.24,54.56)--(297.79,54.41)--(299.35,54.28)--(300.9,54.14)--(302.46,54.02)--(304.01,53.89)--(305.57,53.77)--(307.13,53.66)--(308.68,53.55)--(310.24,53.44)--(311.79,53.34)--(313.35,53.24)--(314.91,53.14)--(316.46,53.05)--(318.02,52.96)--(319.57,52.87)--(321.13,52.79)--(322.68,52.71)--(324.24,52.62);
  \draw[wingInk,line width=.6675bp,dash pattern=on 2.989bp off 2.092bp]
    (246.01,47.33)--(326.52,47.33);
  \node[lab,wingInk] at (236.8,91.27) {\textbf{The phase at infinity}};
  \node[lab,wingInk] at (243.85,75.68) {$j(u)^{-1}-2u-3\log u$};
  \node[lab,wingGrey] at (330.7,21.24) {$u$};
  \node[lab,wingInk] at (226.36,45.34) {$\Theta(g)$};
  \node[lab,wingInk] at (267,5.38) {$u\longrightarrow\infty$};
\end{tikzpicture}%
}%
}
\caption{Each interior point of the left
profile represents a round $S^2$ orbit. At a point $p$ of the maximal
orbit ($z=0$), the orthonormal frame consists of the meridian tangent
$T$ and orbit tangents $e_1,e_2$. The green and purple planes give
$K_m=\Sec_{\Sigma_s}(T,e_1)$ and $K_o=\Sec_{\Sigma_s}(e_1,e_2)$.
The blue band illustrates the $O((\log u)/u)$ error about $\Theta(g)$
in Proposition~\ref{cl:invariant}.}
\label{fig:equatorial-phase}
\end{figure}

To state the classification result, we consider a non-Bryant soliton with the
normalization (\ref{eq:normalization}). On a level $\Sigma_s=\{f=s\}$, let $U(s)$ be the
largest area radius of a spherical orbit, and let $t=U(s)^2/2$ and
$u=\log t$. The maximizing orbit is the equatorial sphere. At this
sphere let $K_m$ be the sectional curvature of a plane spanned by the
unit meridian tangent and an orbit tangent, and let $K_o$ be the
sectional curvature of the orbit tangent plane, both computed in the
induced metric on $\Sigma_s$ (Figure~\ref{fig:equatorial-phase}).
Put $j(u)=K_m/K_o$. By
Proposition~\ref{cl:invariant}, the normalized isometry class determines
the finite limit
\[
           \Theta(g)=\lim_{u\to\infty}
           \bigl(j(u)^{-1}-2u-3\log u\bigr).
\]
In fact, $\Theta(g)$ is the third term in the expansion of $j(u)^{-1}$ in Proposition~\ref{cl:invariant}. $\Theta(g)$ is the key quantity in our classification result.

\begin{theorem}\label{cl:main}
Every complete  four-dimensional $\kappa$-noncollapsed steady gradient Ricci soliton
with positive curvature operator is, up to scaling and isometry, the four-dimensional Bryant soliton
or a member of $\mathcal L$. With normalization (\ref{eq:normalization}), the map
\[
                \mathcal L\longrightarrow\mathbb R,\qquad
                [g]\longmapsto\Theta(g)
\]
is a bijection. 
\end{theorem}

\begin{rem}
Given two non-Bryant steady gradient Ricci solitons $(M_1,g_1,f_1)$ and $(M_2,g_2,f_2)$. Suppose they both belong to $\mathcal L$ and $\alpha(g_1)=\alpha(g_2)$. It remains open to us whether $\Theta(g_1)=\Theta(g_2)$.
\end{rem}

Theorem \ref{cl:main} consists of two parts, Proposition~\ref{cl:A} and Proposition~\ref{cl:B}. In  Proposition~\ref{cl:A},
we show that for every $\theta\in\mathbb R$ there is a steady construction limit
$h$ in Lai's family with $\Theta(h)=\theta$. In the proof, we take a compact
connected limit of the expander paths. The phase extends continuously
to its Bryant and product endpoints, where it takes the values $-\infty$
and $+\infty$ respectively, so its image contains every finite value. In Proposition~\ref{cl:B}, we show that if
 two normalized non-Bryant solitons in the class of
Theorem~\ref{th:ray} have the same phase, then there is a global isometry
between them which preserves their potentials.
Now let $g$ be any non-Bryant soliton in the theorem. By the realization result, there is
a particular Lai limit $h$ whose phase is equal to the phase of $g$, and by the rigidity result,
$g$ is isometric to $h$. Moreover, two normalized members of $\mathcal L$ are
isometric if and only if their phases agree. Therefore, the classification
holds no matter whether the original tip ratio $\alpha$ separates
the members of $\mathcal L$ or not.

The classification result in Theorem \ref{cl:main} can be generalized even if  the sectional curvature is nonnegative.

\begin{theorem}\label{th:classification-sec-nonnegative}
Let $(M^4,g,f)$ be a complete $\kappa$-noncollapsed steady gradient Ricci soliton
with nonnegative sectional curvature. Up to scaling, it is isometric to one of the following
\[
 \mathbb R^4,\;\mathrm{Bry}^3\times\mathbb R,\;\mathrm{Bry}^4,
 \;\mathcal L .
\]
Here $\mathrm{Bry}^n$ denotes the $n$-dimensional Bryant soliton.
\end{theorem}

 Under the hypotheses of Theorem \ref{th:classification-sec-nonnegative},
Hao and the  author~\cite{DH} show that $M$ is diffeomorphic to $\mathbb R^4$.
Hence, if the Ricci tensor has a null direction, then by the splitting theorem of Guan, Lu and Xu~\cite{GLX},
a Euclidean factor splits off from $M$ itself. In this case, we show in Proposition~\ref{secNN:models} that the product soliton is the only nonflat model. So we are left to consider soliton with nonnegative sectional curvature and positive Ricci curvature. We verify the $O(3)$ action under the weaker condition in the Appendix B. The $O(3)$-action with fixed axis together with the nonnegative sectional curvature and positive Ricci curvature will imply the positivity of curvature operator.

The paper is organized as follows. In Section~\ref{sec:geom}, we describe the geometry of the level sets. In
Section~\ref{sec:envelope}, we prove Theorem \ref{th:ray}.
In Section~\ref{sec:modes}, we derive the asymptotic estimates for the level sets.
These estimates are used in Section~\ref{sec:reflection} to compare the two sides of a level set and to
construct the reflection. In Sections~\ref{cl:phase} and~\ref{cl:selection}, we
define the phase and prove its realization. In Sections~\ref{cl:central} and~\ref{cl:caps},
we consider two ends with equal phase and compare their level
metrics from the equator to the poles. In Section~\ref{cl:ambient}, we obtain
a single identification of the four-dimensional ends. In
Section~\ref{cl:rigidity}, we prove the global rigidity and complete the proof
of Theorem~\ref{cl:main}. Finally, we prove
Theorem~\ref{th:classification-sec-nonnegative} in Section~\ref{sec:nonnegative}.

\medskip
\noindent\textbf{Acknowledgement.} The author would like to thank Professor Xiaohua Zhu and Professor Ziyi Zhao for their interests in this work. Professor Xiaohua Zhu and Professor Ziyi Zhao also give an independent proof  of Corollary \ref{cor:sphere-decay} for four-dimensional $\kappa$-noncollapsed steady gradient Ricci solitons with $O(3)$-symmetry in a preprint.

\medskip
\noindent\textbf{Declaration of AI use.} The main ideas of this paper were developed by the author. The idea of doing classification according to phase $\Theta(g)$ is suggested by  ChatGPT.  ChatGPT was used to check calculations, to draw pictures and to improve the language of the paper. The author independently verified the AI-assisted material and takes full responsibility for the mathematical arguments, results, and final content of this paper.

%\section{to be determined}

%On each level set of a four-dimensional Lai soliton, the angular
%$O(3)$ action has two poles on its fixed axis $\Gamma$ and round $S^2$
%orbits between them. The reflection preserves the level, exchanges its
%poles, and fixes its equatorial $S^2$ pointwise. These equatorial
%spheres, together with the tip, form the reflection hypersurface $H$
%shown in Figure~\ref{fig:lai-action}. 

%Both proposed models have an angular $O(3)$ action and a reflection
%interchanging the ends of its fixed axis. Fang~\cite{F} obtained an
%effective angular $O(3)$ action for every soliton in this class and
%showed that its metric blow-down is a closed Euclidean sector. Our
%ray blow-down theorem determines that sector.

\section{Geometry of the potential levels}\label{sec:geom}
We first recall Fang's theorem.

\begin{theorem}[Fang {\cite[Theorem~1.1]{F}}]\label{thm:fang}
Suppose $(M^4,g,f)$ is a complete $\kappa$-noncollapsed steady gradient Ricci soliton with positive curvature operator. Then there is an effective action of $O(3)$ by isometries
which leaves $f$ invariant and whose principal orbits are diffeomorphic to $S^2$.
Moreover, $f$ is proper and has exactly one critical point $o$,
at which it attains its global minimum.
\end{theorem}

Since $f$ has a critical point $o$ by Theorem \ref{thm:fang}, it follows from \cite[Proposition 2.3]{CaoChen12} that $f$ satisfies
\begin{align}\label{f linear}
    C_1 d(o,x)-C_2\le f(x)-f(o)\le d(o,x),
\end{align}
for some constants $C_1,C_2>0$, where the upper bound follows from $|\nabla f|\le1$ in the normalization \eqref{eq:normalization}.
After replacing $f$ by $f-f(o)$, we have $f(o)=\min_M f=0$. By (\ref{f linear}), every sublevel set
$\{f\le c\}=f^{-1}([0,c])$ is compact. Moreover, each $c>0$ is a regular
value. We will mainly work on the non-Bryant soliton and use the
potential-preserving angular $O(3)$ action in Theorem \ref{thm:fang}. Our purpose is to describe
the level set quantitatively by the radii of its sphere orbits, and to compare these
radii as the potential $f$ increases. The conclusions of this section which will be used in
Section~\ref{sec:envelope} are a scalar evolution equation, a maximum
orbit radius of order $\sqrt f$, and an upper bound of order $f$ for
the length of a meridian. The two sides estimate of a level set will not be treated
until the reflection is constructed in Section~\ref{sec:reflection}.

\subsection{The meridian of a single level set}Let $\Gamma$ be the $O(3)$-fixed axis as in  Figure \ref{fig:lai-action}. We have the following topological description of the level set.
\begin{lem}\label{lem:level-interval}
For every $c>0$, the level set $\Sigma_c=\{f=c\}$ is a smooth compact
three-sphere. It meets the fixed axis $\Gamma$ at two points $p_-(c)$ and
$p_+(c)$, one on each side of $o$. The quotient $\Sigma_c/O(3)$ is a
closed interval whose endpoints are these two points. Every interior point of
the interval represents one round sphere $O(3)/O(2)\cong S^2$.
\end{lem}
\begin{proof}
Note that the Hessian of $f$ at $o$ is positive definite. By the Morse lemma, $\Sigma_c$ is diffeomorphic to $S^3$ when $c$ is small.  On a slab
$0<a_0\le f\le b_0$, the smooth vector field
\[
                       Z=\frac{\nabla f}{|\nabla f|^2}
\]
satisfies $df(Z)=1$. The slab is compact, and it contains no critical point since $o$ is the unique critical point of $f$.
So the flow generated by $Z$ exists  on $\{f>0\}$, and it identifies any two level sets
in this region. Moreover, the flow commutes with the group action since $g$ and $f$ are invariant.
Thus every $\Sigma_c$ is a three-sphere.

In the equivariant coordinates of \cite[Proposition~7.6]{F}, the projection
$\pi:M\to Q=M/O(3)$ is given by $(x',x_4)\mapsto (|x'|,x_4)$. Since $\Gamma=\{x'=0\}$ is the fixed point set of  $O(3)$ and has dimension one, $\Gamma$
is a complete geodesic. The elements of $O(3)$ fix $o$ by the uniqueness of the minimum. Along each of the two
outgoing half-axes, $f$ increases strictly from zero to infinity, by the
same geodesic calculation. This gives exactly the two poles in the lemma.
At the other points of the level set, the stabilizer is $O(2)$. So the orbit has
dimension two, and its invariant metric is a multiple of the unit round
metric. Now we consider a pole. Note that the gradient of $f$ at a pole is nonzero and it is tangent to the axis.
Then it follows from the implicit function theorem that $x_4$ is locally a smooth
invariant function of $x'$ on the level set, and we see that the quotient of this level chart has the coordinate
$|x'|\ge0$. So the quotient of the connected compact level is a
one-dimensional manifold with exactly two boundary points, which means that it is a
closed interval.
\end{proof}

\begin{defi}\label{def:level-size}
If $p\notin\Gamma$, then we set
\[
           r(p)=\left(\frac{\operatorname{Area}_g(O(3)\cdot p)}
                                  {4\pi}\right)^{1/2},
 \qquad r\big|_{\Gamma}=0 .
\]
For the level set $\Sigma_c$, we also define
\[
 \begin{aligned}
 U(c)&=\max_{p\in\Sigma_c}r(p),\\
 D(c)&=d_{\,g|_{T\Sigma_c}}(p_-(c),p_+(c))
       =\inf_{\gamma}\int_0^1|\dot\gamma(t)|_g\,dt .
 \end{aligned}
\]
Here we take the infimum over all the piecewise smooth curves in $\Sigma_c$ which join its
two poles. If the value of the potential is
denoted by $f$, then we will also write $U(f)$ and $D(f)$.
\end{defi}

\begin{defi}\label{def:meridian plane}
Let $(M^4,g,f)$ be the $O(3)$-symmetric steady gradient Ricci soliton considered above. We fix a point $x\in M^4\setminus\Gamma$, where $\Gamma$ is the $O(3)$-fixed axis in $M^4$, and let $G_x=\{A\in O(3)|A(x)=x\}\cong O(2)$. The meridian plane through $x$ is defined by
\[P_x:=\{y\in M^4|A(y)=y,\text{ for any $A\in G_x$}\}.\]
In the equivariant coordinates of $\mathbb{R}^4=\mathbb{R}^3\oplus\mathbb{R}$, we denote $(x',x^4):=(x^1,x^2,x^3,x^4)$. Since the vectors in $\mathbb{R}^3$ fixed by $G_x\cong O(2)$ form the line $\mathbb{R}x'$ and the whole $O(3)$-action also fixes $\mathbb{R}$, we have
\[P_x=\mathbb{R}x'\times\mathbb{R}=\{(tx',s)|t,s\in\mathbb{R}\}.\]
Hence $P_x$ is diffeomorphic to $\mathbb{R}^2$. Moreover, $P_x$ is totally geodesic, since it is the fixed point set of the isometry group $G_x\cong O(2)$. For any $z=(z',z^4)\in M^4\setminus\Gamma$, an element of $O(3)$ which carries $\mathbb{R}x'$ onto $\mathbb{R}z'$ restricts to an isometry from $P_x$ onto $P_z$.
\end{defi}

\begin{defi}\label{meridian}
Let $\Sigma_c=\{f=c\}$ be a regular level set and let $x\in\Sigma_c\setminus\Gamma$. The meridian of $\Sigma_c$ through $x$ is the closure in $\Sigma_c$ of the connected component of $(P_x\cap\Sigma_c)\setminus\Gamma$ containing $x$. Since $G_x=\{A\in O(3)|A(x)=x\}\cong O(2)$ preserves $\Sigma_c$, we have
\[P_x\cap\Sigma_c=Fix(G_x)|_{\Sigma_c}.\]
Hence the meridian is a geodesic of the induced metric $g|_{T\Sigma_c}$ on $\Sigma_c$. When $\Sigma_c$ is compact with two poles $p_-(c),p_+(c)$, these poles divide $P_x\cap\Sigma_c$ into two meridian arcs. Each meridian arc minimizes the length between $p_-(c)$ and $p_+(c)$, so its length is $D(c)$ by Definition~\ref{def:level-size}. We orient each meridian arc from $p_-(c)$ to $p_+(c)$ and parametrize the meridian by arc-length,
\[\gamma:[0,D(c)]\to \Sigma_c,\]
with $\gamma(0)=p_-(c),\gamma(D(c))=p_+(c)$.

\end{defi}

\begin{figure}[htbp]
\centering
\IfFileExists{figures/level-measure.pdf}{%
  \includegraphics{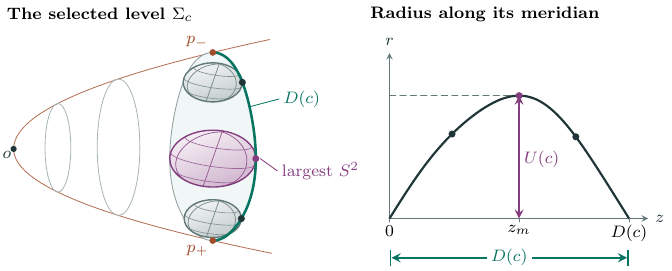}%
}{%
\begin{tikzpicture}[x=.95cm,y=.95cm,>=stealth,line cap=round,line join=round,
  font=\footnotesize,declare function={
    levelrad(\t)=sin(\t)-.14*sin(\t)^3*cos(\t);}]
  \node[anchor=west] at (0,2.40) {\textbf{The selected level $\Sigma_c$}};
  \node[anchor=west] at (6.48,2.40) {\textbf{Radius along its meridian}};
  % t = 180 z/D in degrees. Both panels use precisely the same profile.
  % The sample giving its maximum is the zero of the derivative of levelrad.
  \def\peak{97.387548}
  \def\xc{3.80}
  \def\a{.76}
  \pgfmathsetmacro{\topheight}{sqrt((\xc-.25)/1.20)}
  \pgfmathsetmacro{\botdepth}{sqrt((\xc-.25)/1.33)}
  \pgfmathsetmacro{\yc}{(\topheight-\botdepth)/2}
  \pgfmathsetmacro{\b}{(\topheight+\botdepth)/2}
  \draw[wingAxis!80,thin] plot[domain=0:1.95,samples=41] ({.25+1.20*\x^2},\x);
  \draw[wingAxis!80,thin] plot[domain=0:1.86,samples=41] ({.25+1.33*\x^2},{-\x});
  \fill[wingInk] (.25,0) circle(1.5pt) node[below left,inner sep=1pt] {$o$};
  \foreach \cx/\width in {1.04/.23,2.12/.38}{
    \pgfmathsetmacro{\hi}{sqrt((\cx-.25)/1.20)}
    \pgfmathsetmacro{\lo}{sqrt((\cx-.25)/1.33)}
    \draw[wingGrey!55] (\cx,{(\hi-\lo)/2})
                 ellipse[x radius=\width,y radius={(\hi+\lo)/2}];
  }
  \path[fill=wingPath!6]
    plot[domain=0:180,samples=91] ({\xc+\a*levelrad(\x)},{\yc+\b*cos(\x)})
    -- plot[domain=180:0,samples=91] ({\xc-\a*levelrad(\x)},{\yc+\b*cos(\x)}) -- cycle;
  \draw[wingGrey!72] plot[domain=0:180,samples=91]
       ({\xc-\a*levelrad(\x)},{\yc+\b*cos(\x)});
  \foreach \t/\col in {47/wingGrey,\peak/wingOrbit,140/wingGrey}{
    \pgfmathsetmacro{\rad}{\a*levelrad(\t)}
    \pgfmathsetmacro{\orbheight}{\yc+\b*cos(\t)}
    \begin{scope}[shift={(\xc,\orbheight)},yscale=.66]
       \WingOrbitSphere{(0,0)}{\rad}{\col}
    \end{scope}
  }
  % Split the plots at every marker, so that the vector path contains each point.
  \foreach \lo/\hi in {0/47,47/\peak,\peak/140,140/180}
    \draw[wingPath,line width=1.23pt] plot[domain=\lo:\hi,samples=39]
        ({\xc+\a*levelrad(\x)},{\yc+\b*cos(\x)});
  \foreach \t/\col in {47/wingInk,\peak/wingOrbit,140/wingInk}
    \fill[\col] ({\xc+\a*levelrad(\t)},{\yc+\b*cos(\t)}) circle(1.7pt);
  \fill[wingAxis] (\xc,\topheight) circle(1.6pt) node[above left,inner sep=2pt] {$p_-$};
  \fill[wingAxis] (\xc,{-\botdepth}) circle(1.6pt) node[below left,inner sep=2pt] {$p_+$};
  \draw[wingOrbit,thin] ({\xc+\a*levelrad(\peak)+.08},{\yc+\b*cos(\peak)})
         -- (4.95,-.39) node[right,inner sep=2pt] {largest $S^2$};
  \draw[wingPath,thin] ({\xc+\a*levelrad(65)},{\yc+\b*cos(65)})
         -- (4.98,.89) node[right,inner sep=2pt] {$D(c)$};
  \begin{scope}[shift={(6.95,-1.24)}]
    \def\Dwidth{4.27}
    \def\Rheight{2.17}
    \pgfmathsetmacro{\xm}{\Dwidth*\peak/180}
    \pgfmathsetmacro{\ym}{\Rheight*levelrad(\peak)}
    \draw[->,wingGrey] (0,0)--(4.61,0) node[right,text=wingInk] {$z$};
    \draw[->,wingGrey] (0,0)--(0,2.95) node[above,text=wingInk] {$r$};
    \draw[wingGrey,densely dashed] (0,\ym)--(\xm,\ym);
    \draw[wingOrbit,<->,line width=.9pt] (\xm,0)--(\xm,\ym)
          node[pos=.48,right,inner sep=2pt] {$U(c)$};
    \foreach \lo/\hi in {0/47,47/\peak,\peak/140,140/180}
       \draw[wingInk,line width=1.05pt] plot[domain=\lo:\hi,samples=39]
               ({\Dwidth*\x/180},{\Rheight*levelrad(\x)});
    \foreach \t/\col in {47/wingInk,\peak/wingOrbit,140/wingInk}
       \fill[\col] ({\Dwidth*\t/180},{\Rheight*levelrad(\t)}) circle(1.7pt);
    \foreach \xx/\lab in {0/0,\xm/z_m,\Dwidth/D(c)}
      {\draw[wingGrey] (\xx,0)--(\xx,-.06);
       \node[below,inner sep=2pt] at (\xx,-.04) {$\lab$};}
    \draw[wingPath,|<->|,line width=.95pt] (0,-.70)--(\Dwidth,-.70)
         node[midway,fill=white,inner sep=2pt] {$D(c)$};
  \end{scope}
\end{tikzpicture}%
}
\caption{The green curve on $\Sigma_c$ is a meridian. It joins $p_-$ to $p_+$
with fixed angular coordinate, and its $g$-length is $D(c)$. Each marked point on
the meridian lies on the sphere orbit drawn at that position. Note that the meridian is perpendicular to the orbit at
this point with respect to the level metric
$dz^2+r(z,c)^2g_{S^2}$, but this angle may not be preserved by the planar projection. In the
right panel, we plot the same illustrative radius profile as a function of the meridian
arclength $z$. The three marks in the two panels correspond to each other in the order from $p_-$ to $p_+$.
The violet orbit is the one which has radius $U(c)=r(z_m,c)$ and area $4\pi U(c)^2$,
and at the two endpoints the orbits are single points.}
\label{fig:level-measure}
\end{figure}

Recall that
$\pi:M\to Q=M/O(3)$ is the projection which was used in Lemma \ref{lem:level-interval}. Then we have $\pi(\Sigma_c)=\Sigma_c/O(3)$, and we will take the arclength $z$ along $\pi(\Sigma_c)$ as the parameter. We shall show in Proposition \ref{prop:level-profile} that $r$ is concave with respect to $z$. In Figure~\ref{fig:level-measure}, the minimizing meridian and the
maximal sphere are marked on the same level. Note that we can not tell from the strict concavity
where $z_m(c)$ lies relative to the arclength midpoint.

\begin{prop}\label{prop:level-profile}
Let $z$ be the arclength along $\pi(\Sigma_c)$ which is measured from $p_-(c)$
with respect to the quotient metric $h$. Then we have $0\le z\le D(c)$ and
\[
                g|_{T\Sigma_c}=dz^2+r(z,c)^2g_{S^2}.
\]
The function $r$ vanishes at the two endpoints, and there we have $r_z(0,c)=1$ and $r_z(D(c),c)=-1$.
In the interior, we have $r_{zz}<0$ and $|r_z|<1$. Moreover, there is exactly one
maximizing orbit. Its position $z_m(c)$ is in $(0,D(c))$ and it depends
smoothly on $c$.
\end{prop}
\begin{proof}
Let
$\pi:M\to Q=M/O(3)$ be the projection which is given by $(x',x_4)\mapsto (|x'|,x_4)$, and we choose a unit vector in the $x'$-space. Then the subgroup which fixes this vector
also fixes a meridian plane. We consider the half of this plane where $x'$ points along the chosen vector,
and it is easy to see that it meets each regular orbit once. At such a point, the stabilizer acts trivially
on the plane, and its action on the orbit tangent space is the standard $O(2)$ action.
We claim that the two spaces are orthogonal. In fact, if we pair a fixed vector with an orbit
vector, then we get an invariant linear functional on the standard representation,
and it is well known that such a functional must vanish. It follows that in these coordinates the metric away from the axis can be written as
$g=h+r^2g_{S^2}$, where $h$ is the metric on the chosen half-plane.
We remark that $|x'|$ is a label of the orbit, but $r$ is measured by $g$.
Now we restrict this metric to the inverse image of the quotient arc, and then we get the
level set metric in the proposition.

Suppose that we have a curve which joins the two poles and that its angular coordinate is $\theta(t)$. Then
on the regular part the speed of the curve is
$\sqrt{\dot z^2+r^2|\dot\theta|_{S^2}^2}\ge|\dot z|$.
So we see that the length is decreased by the projection to the orbit interval. This is also true for
the paths through a pole, which can be seen by continuity or by truncating the curve away from each endpoint.
On the other hand, if we take a lift with $\theta$ constant, then it reaches both poles and the equality is attained.
So we have proved that the length of the arc is $D(c)$. By the same argument,
we also see that near $p_-$ the parameter $z$ is the intrinsic distance from that pole. Note that the metric is smooth.
Then in the normal coordinates of the level set, the geodesic spheres around this pole have area
$4\pi z^2(1+O(z^2))$. But these geodesic spheres are just
the angular orbits, and so we get $r(z,c)/z\to1$. In the same way, we have
$r(z,c)/(D(c)-z)\to1$ for the corresponding ratio at the other end.

Next we compute the curvature of the level set. Let $N=\nabla f/|\nabla f|$. If $X,Y$ are vectors which are tangent to the level set,
then we can differentiate the unit normal and get
$\mathrm{II}(X,Y)=\langle\nabla_XN,Y\rangle=\Ric(X,Y)/|\nabla f|$.
We see that this form is positive definite. Now if $X,Y$ are orthonormal, then it follows from the Gauss equation
that the sectional curvature of the level set is obtained by adding $\mathrm{II}(X,X)\mathrm{II}(Y,Y)-\mathrm{II}(X,Y)^2>0$
to the positive ambient sectional curvature. Hence the level set
has positive sectional curvature. In the warped metric of the level set, we write $T=\partial_z$
and let $e_1,e_2$ be an orthonormal spherical pair. We vary a constant-angle
meridian by changing its angle on the unit sphere at unit speed. Then the
variation field $J$ commutes with $T$ and has length $r$. For an angular
field $V$ which is lifted independently of $z$, the Koszul formula gives
\[
  2\langle\nabla^{\Sigma_c}_T J,V\rangle=\partial_z\langle J,V\rangle
           =2(r_z/r)\langle J,V\rangle .
\]
Since $T$ is a unit geodesic field, the $T$-component of $\nabla^{\Sigma_c}_T J$ vanishes.
Hence $e_1=J/r$ is parallel, and the Jacobi
equation becomes $r_{zz}+r\Sec_{\Sigma_c}(T,e_1)=0$.
For the second formula below, we note that the orbit sphere has intrinsic
curvature $1/r^2$, and that its two principal curvatures in the level have
absolute value $|r_z|/r$, which can be seen by differentiating $r^2g_{S^2}$
in the $T$ direction. Thus we have
\[
 \Sec_{\Sigma_c}(T,e_1)=-r_{zz}/r,\qquad
 \Sec_{\Sigma_c}(e_1,e_2)=(1-r_z^2)/r^2 .
\]
Since these curvatures are positive, we get the strict concavity and the slope bound. Moreover, the one-sided
derivatives at the endpoints are equal to the limiting secant slopes computed above.
Thus $r_z$ decreases from $1$ to $-1$ and vanishes once, and at its zero we have
$r_{zz}<0$. The flow of $Z$ gives a smooth parameter for the family of
arcs. Integrating the smooth arclength element from the moving pole, we see that
$z$ is a smooth coordinate in the interiors of these arcs. Finally, applying the implicit function
theorem to $r_z(z,c)=0$, we obtain the smoothness of $z_m(c)$ and $U(c)$.
\end{proof}

\subsection{Quotient curvature and the level set error}Let $Q=M/O(3)$ and $h$ be the quotient metric, i.e. $(Q,h)=(M,g)/O(3)$. The following Lemma is due to \cite{F}
\begin{lem}\label{lem:level-quotient}
Let $Q^\circ$ be the regular region of $Q$. Then the soliton equations on $Q^\circ$ are
\begin{align}
 \nabla_h^2 f&=Kh-2r^{-1}\nabla_h^2r,&
 r\Delta_h r+|\nabla r|^2+r\langle\nabla f,\nabla r\rangle=1,\label{e:wp}\\
 Rr^2&\le C,&
 |\nabla^j\Rm|\le C_jR^{1+j/2}\quad(j\ge0).\label{e:der}
\end{align}
Here $K$ denotes the Gauss curvature of $h$. All contractions in
\eqref{e:wp} are taken with respect to $h$, while \eqref{e:der} is with respect to the ambient metric $g$.
Moreover, the radius of a regular orbit increases strictly along $\nabla f$.
\end{lem}
\begin{proof}
The equations (\ref{e:wp}) and (\ref{e:der}) follow from (7.12), (7.13) and (7.14) in \cite{F}.
Since $f$ is constant on the orbits,  the vertical Hessian of $f$ is
$\nabla^2 f(e,e)=\langle\nabla_h f,\nabla_h r\rangle/r$.
Since $\Ric(e,e)$ is positive, the soliton equation implies that $r$ is strictly increasing along $\nabla f$.
\end{proof}

In Proposition~\ref{prop:level-evolution}, we will study how the induced metric $g|_{T\Sigma_c}$ changes under the flow of $Z=\frac{\nabla f}{|\nabla f|^2}$. For this purpose, we compare the normalized tangential part of $|\nabla f|^{-2}\Ric_g|_{T\Sigma_c}$  with $\Ric_{\Sigma_c}:=\Ric_{g|_{T\Sigma_c}}$. Lemma~\ref{lem:level-error} shows that the difference of these two Ricci curvatures is a nonnegative symmetric 2-tensor $E$.
\begin{lem}\label{lem:level-error}
On every regular level set, the symmetric form
\[
 E=|\nabla f|^{-2}\Ric|_{T\Sigma_c}-\Ric_{\Sigma_c}
       =E^z\,dz^2+E^o\,g_{S^2}
\]
is nonnegative. For every integer $j\ge0$, we have
\begin{equation}
             |\nabla_h^jK|\le C_j R^{2+j/2} .\label{e:K}
\end{equation}
Moreover, on the part of a level where $r\ge1$, we have
\begin{equation}
 \begin{aligned}
  |\partial_z^jE^z|&\le C_jr^{-4-j},&
  |\partial_z^jE^o|&\le C_jr^{-2-j} \qquad(j\ge0),\\
  |\partial_z^jr|&\le C_jr^{1-j} &&(j\ge1).
 \end{aligned}\label{e:err}
\end{equation}
\end{lem}
\begin{proof}
Let $v=|\nabla f|$, $N=\nabla f/v$ and  $A=\Ric|_{T\Sigma_c}$. We denote by $A^2$ the square of $A$, where
one index is raised by the level set metric. If $X,Y$ are tangent vectors, then we put
$B(X,Y)=\langle R(X,\nabla f)\nabla f,Y\rangle$, and the same
formula will also be used for ambient vectors when we take derivatives. We know that $\mathrm{II}=A/v$. Then it follows from the contracted
Gauss equation that
\[
 v^2\Ric_{\Sigma_c}=v^2 A-B+(\operatorname{tr}_{\Sigma_c} A) A-A^2.
\]
Here we have used the fact that the normal direction is included in the ambient Ricci trace and it contributes
$B/v^2$, but the intrinsic trace of the level is taken only over
tangential directions. We subtract this identity from $A$ and use $1-v^2=R$. Then we get
\begin{align}\label{formula of E}
    v^2 E=B+\Ric(N,N) A+A^2.
\end{align}
Note that we have $B(X,X)=v^2\Sec(X/|X|,N)|X|^2$ for $X\in T\Sigma_c$ if $X\ne0$.
All three forms on the right-hand side of (\ref{formula of E}) are nonnegative. Let $T=\partial_z$. Since the metric $g$ is $O(3)$-symmetric and the $O(3)$-action preserves $f$, we have $E^z=E(T,T)$ and $E^o=r^2 E(e,e)$ for a $g$-unit orbit
vector $e$.

Next we estimate the normal curvature. Commuting
one derivative of $\nabla^2 f=\Ric$, we have
\begin{align}
     B(X,Y)=(\nabla_X\Ric)(\nabla f,Y)-(\nabla_{\nabla f}\Ric)(X,Y).\label{formula of B}
\end{align}
Differentiating $R+|\nabla f|^2=1$, we get
$\Ric(\nabla f,\cdot)=-\tfrac12\,\nabla R$. Hence,  the first term of (\ref{formula of B})  is
$$(\nabla_X\Ric)(\nabla f,Y)=-\tfrac12\nabla^2R(X,Y)-\Ric^2(X,Y).$$ On the other hand, by the evolution equation of the Ricci curvature under
$\partial_tg=-2\Ric$, we get
\[
 -\nabla_{\nabla f}\Ric-2\Ric^2
    =\mathcal L_{-\nabla f}\Ric
    =\Delta\Ric+2\Rm(\Ric)-2\Ric^2,
\]
where $\Rm(\Ric)$ is defined by
$$\Rm(\Ric)(X,Y)=\sum_{i,k}\langle R(X,e_i)e_k,Y\rangle\Ric(e_i,e_k)$$
in a $g$-orthonormal frame. Combining these identities, we obtain
\[
             B=\Delta\Ric+2\Rm(\Ric)-\tfrac12\nabla^2R-\Ric^2 .
\]
Thus, by \eqref{e:der}, we have
\begin{align}\label{estimate of B}
    |\nabla^j B|\le C_j R^{2+j/2}.
\end{align}

On $\{ f\ge1\}$, we have $Z(v^2)=2\Ric(\nabla f,\nabla f)/v^2>0$, and  $v$  attains a positive minimum on
$\Sigma_1$.
Note that both $\nabla f$ and the tangent directions of the meridians are
horizontal. Taking the trace of $B$ over the meridian plane, we get $v^2K$, and we write $v^2K=\tr_hB$. Moreover, the covariant differentiation of this horizontal trace is
the intrinsic differentiation on a totally geodesic plane. Since $v$ is
bounded away from zero on $\{f\ge1\}$, $v^{-2}$ is bounded there. For $k\ge1$,
the iterated chain rule for $(1-R)^{-1}$ gives
\begin{align}\label{estimate of v}
 |\nabla^k v^{-2}|
 &\le C_k\sum_{m=1}^k (1-R)^{-m-1}
   \sum_{\substack{a_1+\cdots+a_m=k\\ a_1,\ldots,a_m\ge1}}
   \prod_{i=1}^m |\nabla^{a_i}R|
 \le C'_k R^{1+k/2}.
\end{align}
Indeed, \eqref{e:der} bounds each product by $C_k R^{m+k/2}$, and
$0<R\le1$ by the steady normalization.
Then \eqref{e:K} follows from (\ref{estimate of B}), (\ref{estimate of v}) and $v^2K=\tr_hB$.
On the compact sublevel set $\{f\le1\}$,  the curvature derivatives are bounded and $R$ has a strictly
positive minimum, so the estimate is trivial there.

To estimate $E$, we extend tangential tensors by projecting each
argument with the orthogonal projection $I-N\otimes N$. For $k\ge1$,
$\nabla^k\nabla f$ is a derivative of $\Ric$ of order $k-1$. So we have
$$|\nabla^k N|\le C_k R^{(k+1)/2}$$ on $f\ge1$. Thus, if we differentiate a projection $k$ times, the result is bounded
by $C_kR^{k/2}$. Recall that an intrinsic derivative is the tangential
projection of the ambient derivative of a tangential tensor. When we iterate
this rule, each derivative falls either on the original tensor or on one of its
projections. By (\ref{formula of E}), (\ref{estimate of B}) and (\ref{estimate of v}), we get
\begin{align}
     |(\nabla^{\Sigma_c})^j E|\le C_j R^{2+j/2}.\label{level set estimate of E}
\end{align}
Applying the same reasoning to the uncontracted Gauss equation,
which has one ambient curvature term and two second fundamental form factors,
we also get
$$|(\nabla^{\Sigma_c})^j\Rm_{\Sigma_c}|\le C_j R^{1+j/2}.$$
Note that the part of $\{r\ge1\}$ below height one is compact and disjoint from
$o$. So both estimates hold there as well after we increase the constants.

Along a meridian, both $T$ and a parallel unit angular vector $e$ are
intrinsically parallel. So the ordinary $z$-derivatives of the tensor components
in this frame satisfy the same estimates. On $r\ge1$, we have
$R\le C r^{-2}$ by \eqref{e:der}. In particular, we get
$|\partial_z^k\Sec_{\Sigma_c}(T,e)|\le C_k r^{-2-k}$.
Now we start with $|r_z|\le1$ and differentiate
$r_{zz}=-r\Sec_{\Sigma_c}(T,e)$. At order $j\ge2$, each term is a
constant times $\partial_z^i r$ times a curvature derivative of order
$j-2-i$. By induction, this product is bounded by
$C r^{1-i}r^{-2-(j-2-i)}=C r^{1-j}$.  Hence, $$|\partial_z^i(r^2)|\le C_i r^{2-i}.$$ Finally, by (\ref{level set estimate of E}) and (\ref{e:der}), we get
$$|\partial_z^j E(T,T)|+|\partial_z^j E(e,e)|\le C_j r^{-4-j}.$$
Then the desired estimate follows by applying the product rule to $r^2E(e,e)$.
\end{proof}

\subsection{Variation with the potential}
\begin{figure}[htbp]
\centering
\IfFileExists{figures/mode-transport.pdf}{%
  \includegraphics{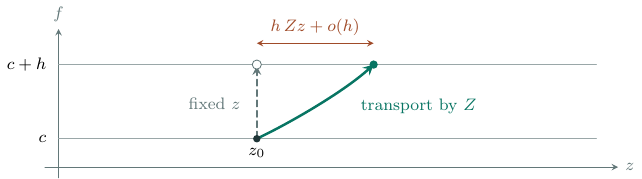}%
}{%
\begin{tikzpicture}[x=.92cm,y=.86cm,>=stealth,font=\footnotesize,
                    line cap=round,line join=round]
  \draw[->,wingGrey] (.3,.16)--(10.85,.16) node[right] {$z$};
  \draw[->,wingGrey] (.55,-.05)--(.55,2.88) node[above] {$f$};
  \draw[wingGrey!65] (.55,.72)--(10.45,.72);
  \draw[wingGrey!65] (.55,2.18)--(10.45,2.18);
  \node[anchor=east] at (.46,.72) {$c$};
  \node[anchor=east] at (.46,2.18) {$c+h$};
  \draw[wingGrey,densely dashed,->,line width=.85pt] (4.2,.74)--(4.2,2.15);
  \draw[wingPath,line width=1.2pt,->] (4.2,.72) .. controls
          (4.85,1.05) and (5.83,1.67) .. (6.35,2.18);
  \fill[wingInk] (4.2,.72) circle(1.8pt);
  \fill[white] (4.2,2.18) circle(2.2pt);
  \draw[wingGrey] (4.2,2.18) circle(2.2pt);
  \fill[wingPath] (6.35,2.18) circle(2pt);
  \draw[<->,wingAxis] (4.20,2.60)--(6.35,2.60)
     node[midway,above] {$h\,Zz+o(h)$};
  \node[anchor=east,wingGrey] at (4.03,1.40) {fixed $z$};
  \node[anchor=west,wingPath] at (5.98,1.36) {transport by $Z$};
  \node[anchor=north] at (4.20,.68) {$z_0$};
\end{tikzpicture}%
}
\caption{The two derivatives in the $(z,f)$ coordinate chart of
Proposition~\ref{prop:level-evolution}. Both paths start at the solid lower point.
The dashed path ends at the point of arclength $z_0$ on the next level, while the green
path follows the same moving orbit, whose new arclength is
$z_0+h\,Zz+o(h)$. The horizontal separation of the two endpoints gives the coefficient
of $\partial_z$ in $Z$. The axes only show coordinate values, and the angles in
this diagram are not the angles measured in the quotient metric.}
\label{fig:mode-transport}
\end{figure}

In Proposition~\ref{prop:level-profile}, a single level $\Sigma_c$ is described by the orbit function $r(z,c)$. In order to estimate the growth of its maximum $U(c)$, we need to compare the orbit functions on nearby levels. 
\begin{prop}\label{prop:level-evolution}
Let $z$ be the arclength along the quotient level set $\pi(\Sigma_c)$ measured from $p_-(c)$. Then the vector field $Z=\nabla f/|\nabla f|^2$
satisfies
\begin{align}
 Z(r^2)&=2-(r^2)_{zz}+2E^o
       =\frac{2r^2\Ric(e_o,e_o)}{|\nabla f|^2},\label{e:ZW}\\
 \partial_z(Zz)&=-2r_{zz}/r+E^z,\qquad (Zz)|_{z=0}=0,\label{e:beta}
\end{align}
where $e_o$ is an ambient unit vector tangent to the sphere orbit.
Let $y=r^2/2$. In the coordinates $(f,z)$, the equation becomes
\begin{equation}
 -\partial_f(y-f)-\partial_{zz}(y-f)=(Zz)y_z-E^o .\label{e:heat}
\end{equation}
 The same formulas also hold locally if $z$ is the
signed arclength from a reference point transported by $Z$, with $z=0$
at that point.
\end{prop}
\begin{proof}
We fix a regular value $c_0$ and denote by $\Phi_t$ the flow generated by $Z=\nabla f/|\nabla f|^2$. Since $df(Z)=1$, we have
$f(\Phi_t(p))=c_0+t$ for $p\in\Sigma_{c_0}$. Hence,
$\Phi_{c-c_0}$ identifies $\Sigma_{c_0}$ with $\Sigma_c$ when $c$ is
close to $c_0$. Note that
the $O(3)$-action preserves both $g$ and $f$, so its elements commute
with $\Phi_t$. Therefore, the flow maps sphere orbits to sphere
orbits and preserves the fixed-point set. Moreover, it maps
$p_-(c_0)$ to the chosen pole $p_-(c)$.

On $\Sigma_{c_0}$, we choose a unit-speed meridian starting from $p_-(c_0)$ and
let $u$ be its arclength parameter. By this meridian and the
$O(3)$-action, we identify the regular orbits with the same unit sphere
$S^2=O(3)/O(2)$, on which we use local polar coordinates $(\theta,\varphi)$.
Then we evolve $(u,\theta,\varphi)$ by $\Phi_{c-c_0}$. In the
resulting coordinates, we have $Z=\partial_c$, and the initial pole on each
level has  $u=0$. By equivariance, the stabilizer $O(2)$ fixes
$\partial_u$, and it has no nonzero invariant covector on the orbit
tangent plane. Thus $\partial_u$ is orthogonal to the orbits. Since
each orbit metric is $O(3)$-invariant, the induced metric can be written as
$$
    g|_{T\Sigma_c}=a(u,c)^2\,du^2+r(u,c)^2g_{S^2},
    \qquad a(u,c)=|\partial_u|_g,
$$
where $r$ is the orbit radius. Let $z$ be the arclength along the
meridian from $p_-(c)$, and let $T$ be the corresponding unit tangent. Then we have
$$
    z(u,c)=\int_0^u a(s,c)\,ds,
    \qquad T=\partial_z=a^{-1}\partial_u.
$$
In particular, $a(u,c_0)=1$ and $z(u,c_0)=u$.

Let $X,Y$ be tangent fields transported by $\Phi_t$. Then
$[Z,X]=[Z,Y]=0$. Since $df(X)=df(Y)=0$, the terms in $(\mathcal L_Zg)(X,Y)$
which come from differentiating $|\nabla f|^{-2}$ vanish.
Therefore, by the soliton equation $\nabla^2 f=\operatorname{Ric}_g$ and the
definition of $E$ in Lemma~\ref{lem:level-error}, we get
$$
\begin{aligned}
    (\mathcal L_Zg)(X,Y)=2|\nabla f|^{-2}\nabla^2 f(X,Y)&=2\bigl(\operatorname{Ric}_{\Sigma_c}+E\bigr)(X,Y).
\end{aligned}
$$

Next we fix a point on a regular orbit and choose the polar chart such that it
contains this point. Let $V=\partial_\theta$, which is the velocity at $s=0$ of
the curve $s\mapsto(u,\theta+s,\varphi,c)$. Since $\partial_\theta$
has unit length on the model sphere, we have $g(V,V)=r^2$, and $e_o=V/r$ is
a unit orbit vector. Since the angular labels are transported, we also have
$[Z,V]=[\partial_c,\partial_\theta]=0$ and
$(\mathcal L_Zg)(V,V)=Z(r^2)$. We choose $e_2$ tangent to the orbit such that
$(T,e_o,e_2)$ is an orthonormal frame  of $T\Sigma_c$. By the sectional curvatures of the
warped product $dz^2+r^2g_{S^2}$, we have
$$
\begin{aligned}
    \operatorname{Ric}_{\Sigma_c}(V,V)
    &=r^2\bigl(\sec_{\Sigma_c}(T,e_o)+\sec_{\Sigma_c}(e_2,e_o)\bigr)\\
    &=r^2\left(-\frac{r_{zz}}{r}+\frac{1-r_z^2}{r^2}\right)
      =1-\tfrac12(r^2)_{zz}.
\end{aligned}
$$
Note that $E=E^o g_{S^2}$ on the orbit, and hence $E(V,V)=E^o$. Now we evaluate the
metric variation identity on $(V,V)$ and obtain
$$
    Z(r^2)
    =2-(r^2)_{zz}+2E^o
    =\frac{2r^2}{|\nabla f|^2}\operatorname{Ric}_g(e_o,e_o).
$$
Here $Z(r^2)$ means the $c$-derivative with the orbit label $u$ fixed, while
the $z$-derivatives are taken with $c$ fixed.

The intrinsic Ricci curvature of $T=a^{-1}\partial_u$ is
$-2r_{zz}/r$, which is the sum over the two spherical directions. Taking the $du^2$
component of the metric variation, we get
\[
         \frac{Za}{a}=\Ric_{\Sigma_c}(T,T)+E^z=-2r_{zz}/r+E^z .
\]
Now we differentiate the integral defining $z$ at fixed $u$. Its upper limit is
constant because $Zu=0$. Then, by $\partial_z=a^{-1}\partial_u$, we get
$\partial_z(Zz)=(Za)/a$. Moreover, that integral is identically zero at $u=0$.
This proves \eqref{e:beta}. If we move the arclength origin to any transported
reference point, the lower limit of the integral is changed to the constant $u$-coordinate of this point.
So the same differentiation works in a signed chart there.

In Figure~\ref{fig:mode-transport}, we compare the endpoint obtained by transporting an orbit with the one obtained by keeping its arclength fixed. Since $Zf=1$, the chain
rule shows that $Z=\partial_f+(Zz)\partial_z$ on any orbit-invariant function of $(f,z)$. Dividing \eqref{e:ZW} by two
and substituting $Zy=y_f+(Zz)y_z$, we get
$y_f+(Zz)y_z=1-y_{zz}+E^o$. Moving the terms to the corresponding sides, we
obtain \eqref{e:heat}. Note that the principal part of \eqref{e:heat} is the forward
heat operator in the direction of decreasing $f$.
\end{proof}

\subsection{Maximum radius and length bounds}\ \\
We first give the growth estimate of the largest orbit.

\begin{prop}\label{prop:level-radius}
For $c>0$,  we have $$(U^2)'(c)>2.$$
Moreover, there is a constant $C$ such that
\begin{equation}
                        2f<U(f)^2\le Cf \mbox{~for~}f>0.\label{e:size}
\end{equation}
\end{prop}
\begin{proof}
By Proposition \ref{prop:level-profile}, $z_m(c)$ depends smoothly on $c$. We differentiate $U(c)^2=r(z_m(c),c)^2$. In the chain rule, the term containing
$z_m'(c)$ vanishes, since $(r^2)_z=0$ at the maximum $z=z_m$. Let $u_*$ be the transported label of this maximizing orbit, so that
$r(u_*,c)=U(c)$. When $s$ is sufficiently close to $c$, the label
$u_*$ determines an orbit on $\Sigma_s$. By the definition of $U(s)$, we have
$$
    U(s)^2-r(u_*,s)^2\ge 0,
    \qquad U(c)^2-r(u_*,c)^2=0.
$$
Hence, the function $s\mapsto U(s)^2-r(u_*,s)^2$ attains a local minimum
at $s=c$. Therefore, at any $c$ where $U^2$ is differentiable, the derivative
of this function vanishes at $s=c$, and we get
$$
    (U^2)'(c)
    =\left.\frac{d}{ds}\,r(u_*,s)^2\right|_{s=c}
    =Z(r^2)
$$
at the chosen orbit. Here the last equality holds because varying $s$
with $u_*$ fixed transports the orbit by the flow of $Z$.
By Proposition~\ref{prop:level-profile}, we have $r=U(c)>0$ and $r_{zz}<0$ at the maximum.
Hence, since $E^o\ge0$, the first expression in
\eqref{e:ZW} gives
\[
       (U^2)'(c)=2-2U(c)r_{zz}(z_m(c),c)+2E^o(z_m(c),c)>2.
\]
It follows that $U(c)^2-2c$ is strictly increasing on $(0,\infty)$.

Next we prove the upper bound. Note that $\Ric(e_o,e_o)\le R$, and $r^2 R\le C$ by
\eqref{e:der}. We also know that $|\nabla f|$ is bounded from below on $f\ge1$, see the
proof of Lemma~\ref{lem:level-error}. Then it follows from the second expression in
\eqref{e:ZW} that $(U^2)'(c)\le C$ there, and by integration we get
$U(c)^2\le C^{\prime}c$ for $c\ge1$.

When $c$ is small, we work in an equivariant coordinate neighborhood of $o$.
Note that $g$ is smooth and positive definite, and so it is uniformly comparable to
the Euclidean metric there. On the orbit of $(x',x_4)$, this implies that
$\operatorname{Area}_g(O(3)\cdot p)\le C |x'|^2$, and then we have
$r(p)^2\le C d_g(o,p)^2$. We also know that the Hessian at $o$ is positive definite.
So by Taylor's formula we have $f(p)\ge a d_g(o,p)^2$ in a smaller neighborhood.
By properness, all the sufficiently small levels are contained in this neighborhood,
and it follows that $U(c)^2\le Cc$ down to zero. On the remaining compact
range of positive values, we just use the continuity of $U$. In particular, we have
$U(c)^2-2c\to0$ as $c\to0$. Now we integrate the strict derivative
inequality on $[c/2,c]$ and get
$U(c)^2-2c>U(c/2)^2-c\ge0$, where the last inequality follows from the
monotonicity and the zero limit. So we have proved the strict lower bound
for every $c>0$.
\end{proof}

Finally, we have the following  estimate for $D(f)$.

\begin{prop}\label{prop:level-length}
For every $f>0$, the intrinsic distance between the poles satisfies
\begin{equation}
                         D(f)\le C(1+f).\label{e:length}
\end{equation}
\end{prop}
\begin{proof}
We evolve the whole meridian together with its two poles by the flow of $Z$. Then we see from equivariance
and horizontality that it is mapped to a constant-angle meridian of the new level set,
and it still realizes the distance by Proposition~\ref{prop:level-profile}.
Now we use the transported $u$ coordinate which was introduced in the proof  of Proposition~\ref{prop:level-evolution},
and we choose its two endpoint values to be $0,1$. Then, by the formula for $Za$ in the proof of Proposition \ref{prop:level-evolution}, we get
\[
 D'(c)=\frac{d}{dc}\int_0^1a(u,c)\,du
      =\int_0^{D(c)}\frac{\Ric(T,T)}{|\nabla f|^2}\,dz .
\]
In particular, $D$ is increasing. Note that $|\nabla f|^2$ is bounded
away from zero for $c\ge1$. In the orthonormal frame $T,N=\nabla f/|\nabla f|,e_1,e_2$, the plane spanned
by $T,N$ is horizontal, so its sectional  curvature is $K$. For each $i=1,2$,
since the second fundamental form is positive, the Gauss formula gives
$\Sec_M(T,e_i)\le\Sec_{\Sigma_c}(T,e_i)$. Summing up, we obtain
\[
        \Ric(T,T)\le K+\Ric_{\Sigma_c}(T,T).
\]

We first estimate the integral of the second term on a single
meridian. On every level set $\Sigma_c$ with $c\ge1$, the ambient curvature and
$\mathrm{II}=\Ric|_{T\Sigma_c}/|\nabla f|$ are uniformly bounded.
Then the Gauss equation implies that $0\le\Ric_{\Sigma_c}(T,T)\le C_0$ uniformly.
Let $L=D(c)$. If $L\le2$, then the desired integral is at most $2C_0$.
If $L>2$, we consider the piecewise linear function
$\phi(z)=\min\{z,1,L-z\}$ on $[0,L]$. We take a parallel orthonormal
pair perpendicular to $T$ within the three-dimensional level set and multiply
each field by $\phi$. These variation fields vanish at the endpoints. More precisely, let $\gamma\colon[0,L]\to\Sigma_c$ be the unit-speed meridian.
At some $z_0\in(0,L)$, we choose an orthonormal basis of
$\dot\gamma(z_0)^\perp\subset T_{\gamma(z_0)}\Sigma_c$ and extend
it by parallel transport along $\gamma$ to get $E_1,E_2$.
Since $\gamma$ is a geodesic, $\dot\gamma$ is also parallel. Hence,
$(\dot\gamma,E_1,E_2)$ is an orthonormal frame along $\gamma$.

For $i=1,2$, we set $V_i=\phi E_i$. Since $\phi(0)=\phi(L)=0$,
each $V_i$ vanishes at both endpoints, and it is admissible for the
fixed-endpoint index form. Since the meridian minimizes length in
$\Sigma_c$, we have $I(V_i,V_i)\ge 0$. Let
$D_z=\nabla^{\Sigma_c}_{\dot\gamma}$. Then $D_zE_i=0$, and hence
$D_zV_i=\phi' E_i$. Substituting this into the index form, we get
$$
\begin{aligned}
0\le I(V_i,V_i)
&=\int_0^L
  \left(
    |D_zV_i|^2
    -\big\langle R^{\Sigma_c}(V_i,\dot\gamma)\dot\gamma,V_i
     \big\rangle
  \right)\,dz\\
&=\int_0^L
  \left(
    |\phi'|^2-\phi^2\sec_{\Sigma_c}(\dot\gamma,E_i)
  \right)\,dz.
\end{aligned}
$$
Note that $E_1,E_2$ form an orthonormal basis of $\dot\gamma^\perp$
at each point. So we have
$$
\operatorname{Ric}_{\Sigma_c}(\dot\gamma,\dot\gamma)
=\sec_{\Sigma_c}(\dot\gamma,E_1)+\sec_{\Sigma_c}(\dot\gamma,E_2).
$$
Therefore,
$$
\int_0^L \phi^2
  \operatorname{Ric}_{\Sigma_c}(\dot\gamma,\dot\gamma)\,dz
\le 2\int_0^L |\phi'|^2\,dz=4.
$$
The cutoff $\phi$ is equal to one on $[1,L-1]$. The other two intervals have total
length two and contribute at most $2C_0$. Thus, we obtain
$$\int_0^{D(c)}\Ric_{\Sigma_c}(T,T)\,dz\le4+2C_0, ~\mbox{for~all }c\ge1.$$

For the term with $K$, we integrate in height. In the non-Bryant branch,
each of the two isometric halves of the smooth meridian double has
curvature integral $\int_QK\,dA_h$. Therefore, by Fang's angle formula
\cite[Theorem~7.10, (7.9)]{F}, we have $$\int_Q K\,dA_h
=\pi-\alpha\le\pi,$$ where $\alpha$ is the nonnegative blow-down
angle. We remark that this bound is used before that angle is determined. Since $f$
descends to $Q$, the length of its quotient gradient is
$|\nabla_h f|=|\nabla f|\le1$. Now we apply the coarea formula on $Q$, where $dz$ is the arclength
on each level arc, and get
\[
 \int_1^c\int_{\pi(\Sigma_s)} K\,dz\,ds
   =\int_{\{q\in Q:\,1\le f(q)\le c\}}K|\nabla_h f|\,dA_h
   \le\int_Q K\,dA_h \le\pi .
\]
Note that the two pole endpoints have zero arclength measure. Integrating the
first variation inequality from $1$ to $c$, we obtain
$D(c)-D(1)\le C(c-1)+C\pi$. For $0<c\le1$, we have
$D(c)\le D(1)$ by monotonicity. Combining these estimates, we get \eqref{e:length}.
\end{proof}

\section{The ray blow-down}\label{sec:envelope}
By Fang's sector theorem \cite[Theorem~1.3]{F}, in order to prove
Theorem~\ref{th:ray} it suffices to show that the scalar curvature decays uniformly.
In this section, we work on the non-Bryant branch. We first show that, for any $a>0$,
its maximal orbit cannot satisfy
$U(f)^2\le 2f+Cf^{1-a}$ on every sufficiently high level.
The reason is geometric. If this estimate were true, then $r^2$ would drop
too slowly along a meridian, whose length is at most $C(1+f)$. Then we show
that if the scalar curvature decay fails, such an upper bound must hold. In this step,
we only use the side of a level which is adjacent to a nondecaying pole.
We remark that the $O(3)$ action is sufficient for both steps.

\subsection{The poles and the large sphere orbits}
Let us first recall some results on dimension reduction. Let $x_i\to\infty$ be a sequence of points.
We rescale the soliton by the scalar curvature at these points such that the scalar curvature
at the new base point is one. On the non-Bryant soliton, it was proved in \cite[Theorems~1.2--1.3]{CMZ} that
a subsequence of the canonical flows
converges smoothly either to the three-dimensional Bryant soliton times a line
or to the shrinking round $S^2\times\mathbb R^2$ cylinder.
Moreover, the rotations converge on compact sets to the standard
rotations of the Bryant factor or of the spherical factor. We also know from \cite[Lemma~7.13]{F} that
at a limit fixed point there are actual axis points which converge to it
in the rescaled charts. So a pole can only converge to a point on the
Bryant tip line.

Using the dimension reduction, we now give asymptotic estimates for the curvature at the poles
and for the slopes on large orbits.

\begin{lemma}\label{lem:thick}
At either pole, we have $fR_{\rm tip}(f)\to\infty$. Moreover, for each fixed $d>0$, we have
$|\nabla_h r|\to0$ uniformly on $r\ge d\sqrt f$ as $f\to\infty$.
\end{lemma}
\begin{proof}
We parametrize one outgoing half of the axis by the arclength $l$. Since the gradient
of $f$ is fixed by the isotropy, it is tangent to the axis there.
Its component $(f\circ\gamma)'$ is strictly increasing from zero, because
$(f\circ\gamma)''=\Ric(\dot\gamma,\dot\gamma)>0$. Hence,
$R_{\rm tip}=1-((f\circ\gamma)')^2$ is decreasing and converges to a limit in $[0,1)$.
If this limit is positive, then the first assertion follows immediately.

Now we assume that the limit is zero. For any sequence of such poles, we denote by $R_i$
the scalar curvature at the base point. For the rescaled flow $g_i(t)=R_i g(t/R_i)$, the time
derivative of the scalar curvature at $(x_i,0)$ is
$-\langle\nabla f,\nabla R\rangle(x_i)/R_i^2$. Here one factor $R_i^{-1}$
comes from the scalar curvature and the other one comes from the rescaling of time.
We note that on the limiting Bryant product this derivative vanishes at the tip
line. In fact, the canonical flow moves the Bryant factor by the gradient of its potential,
which vanishes at the tip, and the scalar curvature is constant in
the line direction. Therefore, since the flows converge smoothly, the ratio above
tends to zero along every such sequence. Note that $|\nabla f|^2\to1$
on the chosen end. Hence, along the pole trajectory, we obtain
\[
 \frac{d}{df}\frac1{R_{\rm tip}}
       =\frac{-\langle\nabla f,\nabla R\rangle}
                    {|\nabla f|^2 R_{\rm tip}^2}\longrightarrow0 .
\]
It follows that
$1/(fR_{\rm tip})\to0$.

Next we turn from the poles to the large orbits.
For this purpose, we first need a lower bound of the curvature.
Recall that on steady Ricci solitons the scalar equation is
$$\Delta R+\langle\nabla f,\nabla R\rangle=-2|\Ric|^2.$$
Combining \eqref{e:der} with the scalar equation, we have $|\langle\nabla f,\nabla R\rangle|\le CR^2$.
Hence,
$$|Z(1/R)|\le C, \mbox{~on~}\{f\ge1\},$$  where $Z=\nabla f/|\nabla f|^2$ and $|\nabla f|$ has a
positive lower bound. By integrating the inequality above, we get
$R\ge c/(1+f)$. Together with the bound $Rr^2\le C$ in (\ref{e:der}), this implies
\[
              c_d\le fR\le C_d\qquad (r\ge d\sqrt f,\ f\ge1).
\]

Let $J=\{x\in M|r(x)\ge d\sqrt f(x),\ f(x)\ge1\}$. Suppose that $x_i\to\infty$ and $x_i \in J$.
We first assume that the sequence related to $x_i$ has a Bryant-product limit. Then we choose
a point on its tip line at finite rescaled distance. By the persistence of the axis and the smooth
convergence, there are actual poles $q_i$ at distance $O(R_i^{-1/2})$
from $x_i$, such that $R(q_i)/R_i$ tends to a positive finite number.
Let $f_i=f(x_i)$. Then we have $R_i\asymp1/f_i$. Since $|\nabla f|\le1$, we also have
$|f(q_i)-f_i|\le d_g(q_i,x_i)=O(\sqrt{f_i})$.
It follows that $f(q_i)R(q_i)$ stays bounded, which contradicts the pole estimate.

Therefore, we only need to consider the cylindrical subsequences. In the rescaled metrics, the orbit
radius is $\sqrt{R_i}\,r$. As in Lemma~\ref{lem:level-quotient}, its square is a fixed multiple
of the sum of the squared lengths of the three rotation generators.
Since the generators converge smoothly and the limiting
sphere has positive radius, the radii converge in $C^1$ near the base point.
On the cylinder, the radius of this sphere is spatially constant.
Note that $$|\nabla(\sqrt{R_i}\,r)|_{R_i g}=|\nabla r|_g
=|\nabla_h r|.$$ Since the limit radius is constant for any subsequence, we conclude that
$|\nabla_h r|\to0$ uniformly on J as $f\to\infty$.
\end{proof}

\subsection{A barrier for the meridian slope}
\begin{lemma}\label{lem:shift}
Suppose that $U(f)^2\le2f+Cf^{1-a}$ holds for all sufficiently large $f$, where
$0<a<1$. Let $m=\max\{3,2(1-a)/a\}$. Then, at every point where $r^2\ge f$ and $f$ is
sufficiently large, we have
\begin{equation}
             |(r^2)_z|\le Cf^{1/2-1/m} .\label{e:slope}
\end{equation}
\end{lemma}
\begin{proof}
Since $r$ is strictly concave, on each side of the maximal orbit we can use $r$ as a smooth
coordinate for $0<r<U(f)$. On that side, we regard $v=r_z^2$ as a function of $(f,r)$.
Then we have $0<v<1$, and by the chain rule $v_r=2r_{zz}<0$.
In these coordinates, the level metric is $dr^2/v+r^2g_{S^2}$.
In order to differentiate at fixed radius, we add to the normal transport the
tangential motion $j\partial_r$ which cancels $Zr$. By
\eqref{e:ZW} and $(r^2)_{zz}=2v+r v_r$, we get
\[
                j=-Zr=(v-1)/r+v_r/2-E^o/r .
\]
The level set Ricci curvature of the unit meridian tangent is $-2r_{zz}/r$.
Since $|\partial_r|^2=1/v$ and $2r_{zz}=v_r$, its $rr$ component
is $-v_r/(rv)$. For the same reason, we have $E_{rr}=E^z/v$. Moreover, when the Lie derivative
acts on the coefficient $1/v$, it contributes
$j(1/v)_r+2j_r/v=-jv_r/v^2+2j_r/v$.
Thus the first variation formula of Section~\ref{sec:geom} becomes
\[
       -v_f/v^2=-2v_r/(rv)+2E^z/v-jv_r/v^2+2j_r/v .
\]
Substituting the expression of $j$ and collecting terms, we obtain
\[
 \begin{split}
 -v_f&=\mathcal B[v]+2vE^z+r^{-1}v_r E^o-2v(E^o/r)_r,\\
 \mathcal B[v]&=vv_{rr}-\tfrac12 v_r^2+
                    r^{-2}(1-v)(rv_r+2v).
 \end{split}
\]
Note that $v_r E^o\le0$. Moreover, we have
$\sqrt v\,|(E^o)_r|=|(E^o)_z|\le Cr^{-3}$ for $r\ge1$, and also
$|E^o|\le Cr^{-2}$, $|E^z|\le Cr^{-4}$ by \eqref{e:err}.
Since $v\le\sqrt v$, it follows that
\[
                -v_f\le\mathcal B[v]+Cr^{-4}\sqrt v .
\]

Now we fix a sufficiently large barrier parameter $A$ and recall Brendle's barrier function
$\psi_A$ \cite[Proposition 2.7]{clBrendle}, whose properties are collected in
\cite[Proposition~2.19]{MMS}. We will use the positive profile
$\psi_A$ with the argument
$x=r/\sqrt{2f+kA^m}$, where $k$ will be chosen independently of $A$.
The profile is defined for $r_*/A\le x\le1+A^{-2}/100$, where $r_*\in(0,1)$ is a universal constant.
We first check the upper endpoint. For this, we maximize $Cf^{1-a}-\delta f/A^2$
over $f\ge0$. Since its critical point is of order $A^{2/a}$, we have
\[
 Cf^{1-a}\le\delta f/A^2+C_\delta A^{2(1-a)/a}.
\]
We choose $\delta$ small and $k>C_\delta$ large. Since $m\ge2(1-a)/a$,
it follows that $U^2$ is strictly less than
$(2f+kA^m)(1+A^{-2}/100)^2$ at all heights under consideration.
At the lower endpoint, we have $r\ge r_*\sqrt{k}\,A^{m/2-1}\ge1$ when $A$ is large.
Hence, the error estimate for $v$ holds on the whole comparison interval.

For the calculation, we temporarily write $H=2f+kA^m$. When $r$ is fixed, we have
$x_f=-x/H$ and $\partial_r=H^{-1/2}\partial_x$. So in the $r$ coordinate we find that the excess of
$-\partial_f\psi_A$ over $\mathcal B[\psi_A]$
is equal to $H^{-1}(x\psi_A'-\mathcal B_x[\psi_A])$, where the
subscript means that $r$ is replaced by $x$ everywhere in the operator.
It follows from the profile inequality cited above that this excess is at least
$cA^{-4}x^{-5}/H$ for $x\ge N/A$ and at least $cA/H$ for $r_*/A\le x\le N/A$,
where $N,r_*,c$ are independent of $A$. On the other hand, at a possible contact point, the error in
the equation for $v$ is at most $C H^{-2}x^{-4}\sqrt v$. In the outer
interval, we have the profile bound $\psi_A\le2A^{-2}x^{-2}$, and it implies that
$\sqrt v\le C/(Ax)$. Then we divide the error by the excess above and get
$CA^3/H$. In the inner interval, we use $\sqrt v\le1$ and $x^{-4}\le
A^4/r_*^4$, and we see that the ratio is again at most $CA^3/H$. Note that $m\ge3$. So this
is at most $C/k$, and we can make both ratios strictly less than one if we enlarge $k$.

Next we compare $v$ with $\psi_A$. At $x=r_*/A$, the profile is at least $3/2>v$, and at $r=U(f)$
we have $v=0<\psi_A$. When $A$ is fixed, the comparison region satisfies
$r\ge(\sqrt2\,r_*/A)\sqrt f$. Then we know from Lemma~\ref{lem:thick} that $v\to0$
uniformly there, and we also have $\psi_A\ge A^{-4}/32$. It follows that $v<\psi_A$ on every
sufficiently high terminal level. Now we let the comparison evolve
towards decreasing $f$. Suppose that there is a first contact. Then by the strict boundary
inequalities, we see that the contact point lies in the spatial interior. At this point, we have $v=\psi_A$,
$v_r=(\psi_A)_r$, $v_{rr}\le(\psi_A)_{rr}$ and
$-\partial_f(v-\psi_A)\ge0$. However, the highest derivative in $\mathcal B$
has the positive coefficient $v$. Then from the differential inequality and the
strict excess we get $-\partial_f(v-\psi_A)<0$, which is impossible. If the contact occurs at a join
of the profile, which is $C^1$ and
piecewise $C^2$, the same argument works with either one-sided Taylor
expansion. Note that this comparison takes place on compact height slabs. If a
spatial interval shrinks to an endpoint, the boundary inequality still
holds there. Taking the terminal height above a prescribed target, we conclude that
$v\le\psi_A$ on both sides at that target.

Finally, we choose $A=f^{1/m}$ at the target. On $r^2\ge f$, we have
$x\ge(2+k)^{-1/2}>N/A$ for all sufficiently high levels. Hence, by the outer
bound, we get $v\le C_k A^{-2}$. Since $r\le U(f)\le C\sqrt f$, we obtain
$|(r^2)_z|=2r\sqrt v\le C f^{1/2-1/m}$. This proves the lemma.
\end{proof}

\subsection{Heat estimates near the maximal orbit}
\begin{lemma}\label{lem:heat}
Let $0<b<b'<2$ and $0\le\gamma<1/2$. Assume that on all sufficiently
high levels we have $|(r^2)_z|\le Cf^\gamma$ wherever $r^2\ge bf$. Then we have
the following estimates on every sufficiently high level and wherever $r^2\ge b'f$,
\begin{align}
 |(r^2)_{zz}|&\le C f^{2\gamma-1}\log f,\label{e:hess}\\
 |(r^2)_z|&\le C f^{3\gamma-1/2}\qquad\text{if }\gamma<1/4 .
                                                    \label{e:improve}
\end{align}
Here the constants are allowed to depend on $b,b'$ and $\gamma$.
\end{lemma}
\begin{proof}
We first construct a chart on which the bounds persist.
Let $S$ be a large height, and we fix a target on this level with $r^2\ge b'S$.
We take this target as the origin of the arclength, and then we transport the origin by $Z$ for all
$f\ge S$.
On $r^2\ge bf$, we see that the assumption $|(r^2)_z|\le Cf^\gamma$ implies that
$|r_z|\le C f^{\gamma-1/2}$. Note that $r_{zz}<0$ and $E^o\ge0$. Then by
\eqref{e:ZW} we have
\[
                 Z(r^2)\ge 2-2r_z^2>b'
\]
for all sufficiently large $S$. Now if the moving origin crossed
$r^2-b'f=0$ downward for a first time, then the derivative along $Z$ would be nonpositive
at the crossing. But we know from the inequality above that this derivative is strictly positive. So
the moving origin stays in
$r^2\ge b'f$. On either side of the origin, the first drop of $r^2$ from its value there to $bf$
is at least $(b'-b)f$. If we integrate the assumed slope bound, then we see that such a drop
requires arclength at least $(b'-b) f^{1-\gamma}/C$.
Now we choose $\delta>0$ which is fixed and sufficiently small, and we set
$d(f)=\delta f^{1-\gamma}$. Then it is easy to see that the signed chart $|z|\le2d(f)$
is contained in $r^2\ge bf$ for every $f\ge S$. We also have
$r\asymp\sqrt f$ on this chart, where the upper bound follows from \eqref{e:size}.

Next we estimate the source term and introduce a cutoff.
In this chart, we have $Zz=0$ at $z=0$, and so we can integrate \eqref{e:beta} from zero.
Note that on the segment of integration $r_z$ is decreasing and the slopes at both endpoints
are $O(f^{\gamma-1/2})$. Hence, the integral of $|r_{zz}|$ is
at most $C f^{\gamma-1/2}$. Dividing by the lower bound for $r$ and
using $|E^z|\le Cf^{-2}$, we get
\[
 |Zz|\le C f^{\gamma-1}+C d(f)f^{-2}\le C f^{\gamma-1}.
\]
 Let
$F=(Zz)y_z-E^o$, where $y=r^2/2$ as in Proposition \ref{prop:level-evolution}.
By \eqref{e:err} and \eqref{e:beta}, we have $|y_z|\le C f^\gamma$, $|E^o|\le C/f$,
$|(Zz)_z|\le C/f$, $|y_{zz}|\le C$ and $|(E^o)_z|\le C f^{-3/2}$.
Differentiating $F$ once, we then obtain $|F|\le C f^{2\gamma-1}$ and $|F_z|\le C$.

Let $\eta$ be a smooth cutoff function which is equal to one on $[-1,1]$ and is supported in
$[-2,2]$. We define $\chi(z,f)=\eta(z/d(f))$ and
$w=\chi(y-f)$, and extend $w$ by zero to the whole $z$-line. Then we have
\[
 (-\partial_f-\partial_{zz})w
   =\chi F+
   \big(-\partial_f\chi-\partial_{zz}\chi\big)(y-f)
                         -2(\partial_z\chi)y_z .
\]
The last two terms are supported in $d(f)\le |z|\le2d(f)$, and they are bounded
by a uniform constant. In fact, we have $|y-f|\le Cf$ and
$|\partial_f\chi|\le C/f$, while the first and second spatial
derivatives of the cutoff are bounded by $C/d(f)$ and $C/d(f)^2$. Since $\gamma<1/2$, both
$f/d(f)^2$ and $f^\gamma/d(f)$ are bounded.
We also have $\|w_z(\cdot,f)\|_\infty\le C f^\gamma$, where the term in
which the derivative falls on the cutoff is included.

Now we integrate the heat equation.
Let $G_t(z)=(4\pi t)^{-1/2}\exp(-z^2/(4t))$. We apply Duhamel's formula from a terminal level
$S+T$ down to $S$. It expresses $w(\cdot,S)$ as the sum of
$G_T*w(\cdot,S+T)$ and the integral, for $0<t<T$, of $G_t$ convolved
with the right-hand side above at height $S+t$. Here the sign of the time
direction is determined by the operator $-\partial_f-\partial_{zz}$.
By scaling the Gaussian, we have
$\|\partial_z^jG_t\|_1\le C_j t^{-j/2}$ for $j=1,2$. If $|z|\le d(S)/2$, then the
terms which are supported on the cutoff annulus are at distance at least
$d(S+t)/2$. So after we take $j$ spatial derivatives, we see that their contribution is bounded by
\[
 C_j\int_0^\infty t^{-j/2}
   \exp\!\left(-c\frac{(S+t)^{2-2\gamma}}{t}\right)dt
                         =O_N(S^{-N})\quad\text{for every }N .
\]
Now we explain the last equality. When $t\le S$, we substitute $u=S^{2-2\gamma}/t$, whose lower endpoint is
$S^{1-2\gamma}$, and when $t\ge S$ we note that the exponent is at most $-ct^{1-2\gamma}$.
Then it is clear that both tails decay faster than any power of $S$.

We first estimate the second derivative. In the terminal convolution, we put one derivative
on $w$ and the other one on $G_T$. Then we see that its size is $O(T^{-1/2}(S+T)^\gamma)$, and it tends
to zero as $T\to\infty$. In the source integral, for $0<t<S^{-4}$, we put
one derivative on the cutoff times $F$, whose first derivative is uniformly
bounded. This part contributes at most $C\int_0^{S^{-4}}t^{-1/2}dt
=O(S^{-2})$. For later times, we put both derivatives on $G_t$. The two
remaining integrals satisfy
\[
 \int_{S^{-4}}^S t^{-1}(S+t)^{2\gamma-1}dt
       \le 5 S^{2\gamma-1}\log S,\qquad
 \int_S^\infty t^{2\gamma-2}dt=\frac{S^{2\gamma-1}}{1-2\gamma}.
\]
The cutoff is equal to one at the target. Since $(r^2)_{zz}=2w_{zz}$ there,
\eqref{e:hess} follows from these estimates.

Next we consider the first derivative. Here we estimate its \emph{difference} at two points
$z_1,z_2$ in $|z|\le d(S)/2$. When we move the derivative onto the terminal
data, what remains is the difference of two translates of $G_T$. By the fundamental theorem of
calculus, its $L^1$ norm is
at most $|z_1-z_2|\,\|\partial_zG_T\|_1$. Hence, for fixed
$S,z_1,z_2$, the terminal contribution tends to zero, and this is true even when $w_z$ itself
has a spatially constant part.
For the source term at either point, we put one derivative on $G_t$. By the scaling
$t=Su$, we have
\[
 \int_0^\infty t^{-1/2}(S+t)^{2\gamma-1}dt
 = S^{2\gamma-1/2}\int_0^\infty u^{-1/2}(1+u)^{2\gamma-1}du .
\]
The integrand is integrable at zero and behaves like $u^{2\gamma-3/2}$
at infinity. This is where the restriction $\gamma<1/4$ comes from. It follows that the
oscillation of $(r^2)_z$ on the interval $|z|\le d(S)/2$ is bounded by $C S^{2\gamma-1/2}$.
Finally, since $r$ is strictly concave, the set $\{r^2\ge b'S\}$ is a single interval
which contains the maximum. We join the target to the maximum by successive
steps of length at most $d(S)/2$, and at each step we recenter the construction above.
By \eqref{e:length}, the number of steps is at most $1+2D(S)/d(S)\le C S^\gamma$.
Since the slope is zero at the maximum, summing the local
oscillation bounds gives \eqref{e:improve}.
\end{proof}

\begin{figure}[htbp]
\centering
\IfFileExists{figures/mode-drop.pdf}{%
  \includegraphics{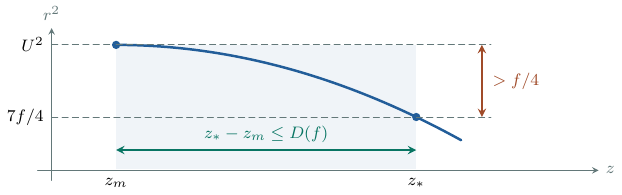}%
}{%
\begin{tikzpicture}[x=.95cm,y=.86cm,>=stealth,font=\footnotesize,
                    line cap=round,line join=round]
  \draw[->,wingGrey] (.8,.30)--(10.8,.30) node[right] {$z$};
  \draw[->,wingGrey] (1.05,.10)--(1.05,3.10) node[above] {$r^2$};
  \fill[wingMirror!7] (2.2,.32) rectangle (7.55,2.77);
  \draw[wingGrey,densely dashed] (1.05,2.77)--(2.2,2.77)
                 (1.05,1.35)--(7.55,1.35);
  \node[left] at (1.04,2.77) {$U^2$};
  \node[left] at (1.04,1.35) {$7f/4$};
  \draw[wingMirror,line width=1.2pt] plot[domain=2.2:8.35,samples=50]
      (\x,{2.77-1.42*((\x-2.2)/5.35)^2});
  \fill[wingMirror] (2.2,2.77) circle(2pt) (7.55,1.35) circle(2pt);
  \draw[wingGrey,densely dashed] (7.55,1.35)--(8.88,1.35)
                         (2.2,2.77)--(8.88,2.77);
  \draw[<->,wingAxis,thick] (8.72,1.35)--(8.72,2.77);
  \node[right,wingAxis] at (8.78,2.06) {$\displaystyle>f/4$};
  \node[below] at (2.2,.28) {$z_m$};
  \node[below] at (7.55,.28) {$z_*$};
  \draw[<->,wingPath,thick] (2.2,.70)--(7.55,.70)
          node[midway,above] {$z_*-z_m\le D(f)$};
\end{tikzpicture}%
}
\caption{The last step in the proof of Proposition~\ref{prop:noenv}, drawn on the side
where $z$ increases. Here $z_*$ denotes the first point where $r^2=7f/4$ is crossed.
The improved slope bound can be integrated on the shaded interval.
The vertical arrow shows the drop $U^2-7f/4>f/4$.
By integrating the slope bound, we see that this drop is at most $Cf^{-1/2}(z_*-z_m)$.}
\label{fig:mode-drop}
\end{figure}

In the following proposition, we show that $U(f)^2$ cannot have an upper bound of the form
$U(f)^2\le2f+Cf^{1-a}$ for any positive $a\in(0,1)$.
\begin{proposition}\label{prop:noenv}
On a non-Bryant soliton under the assumption of Theorem \ref{th:ray}, there is no $a>0$ such that
$U(f)^2\le2f+Cf^{1-a}$ holds on all sufficiently high level sets.
\end{proposition}
\begin{proof}
We argue by contradiction. Suppose that the inequality holds for some $a>0$. We first show that
the inequality then holds for $a=2/5$ and some constant $C$. For this, it suffices to deal with
the case $0<a<2/5$. When $a<2/5$, the exponent in
Lemma~\ref{lem:shift} is $m=2(1-a)/a>3$. Hence, on $r^2\ge f$ we can use
$\gamma=1/2-1/m$ in Lemma~\ref{lem:heat}. By \eqref{e:size}, the maximum lies in
$r^2\ge3f/2$ when $f$ is large. Since $(r^2)_z=0$ at that orbit, we have
$(U^2)'=Z(r^2)$ there. Then, by \eqref{e:hess},
\eqref{e:ZW} and $|E^o|\le C/f$, we get
\[
        (U^2)'\le2+Cf^{-2/m}\log f+C/f,\qquad 2/m=a/(1-a).
\]
For any $a'<a/(1-a)$, the two error terms are $O(f^{-a'})$.
Since $a'<1$, by integration we get $U^2\le2f+C_{a'}f^{1-a'}$.
In order to reach $2/5$ in finitely many steps, we choose
$1/a'=\max\{5/2,\,1/a-1/2\}$. This number lies strictly between $1/a-1$
and $1/a$, and so $a<a'<a/(1-a)$. With this choice, $1/a$ decreases by $1/2$
at each step until it reaches $5/2$.

Now we consider the case $a=2/5$. In this case, we have $m=3$ and $\gamma=1/6$.
Applying \eqref{e:improve} with $(b,b')=(1,3/2)$, we see that the slope exponent becomes
$3(1/6)-1/2=0$. Applying \eqref{e:improve} again with $(b,b')=(3/2,7/4)$ and
$\gamma=0$, we get $|(r^2)_z|\le C f^{-1/2}$ on $r^2\ge7f/4$.
We follow either side of a meridian from the maximum to the first point
with $r^2=7f/4$. Such a point exists since the value at the pole is zero. By \eqref{e:size}, the
drop satisfies $U^2-7f/4>f/4$. On the other hand, integrating the
last slope bound over this segment (see Figure~\ref{fig:mode-drop}), we find that the same drop is
bounded by $C f^{-1/2}D(f)\le C'\sqrt f$. These two estimates contradict each other
when $f$ is large. The proposition is proved.
\end{proof}

\subsection{Excluding a nondecaying end}\label{subsec:end}

\begin{lemma}\label{lem:ray-end}
Suppose that $R$ does not tend to zero uniformly at infinity. Then along one half of the fixed
axis the scalar curvature has a limit $L$ with $0<L<1$.
\end{lemma}
\begin{proof}
By assumption, there is a sequence of points $x_i\to\infty$ with $R_i=R(x_i)\ge\epsilon>0$.
Since $R\le1$, by passing to a subsequence, we may assume that $R_i\to R_*>0$.
We consider the curvature-normalized metrics $g_i=R_i g$. With respect to these metrics, the gradient
of the function $f-f(x_i)$ has norm at most $\epsilon^{-1/2}$, and we have
$\nabla^2_{g_i}f=\Ric_{g_i}$. By taking a diagonal subsequence, we get a potential
$f_\infty$ on the limit $(M_{\infty},g_{\infty})$ which satisfies the same Hessian equation. Hence, the limit is also a steady Ricci soliton, and it splits off a line. Recall that $R_*>0$. Then, by the uniqueness of $\kappa$-noncollapsed three-dimensional steady Ricci solitons \cite{Bre13}, the limit must be a Bryant product. We take a point on the
tip line of the limit which has finite distance from the limiting base point. As in
Lemma~\ref{lem:thick}, we can find points $q_i$ on the fixed axis which converge
to this point in the pointed charts. For all sufficiently large $i$, we have $R_i\ge R_*/2$ and hence
$d_g(x_i,q_i)\le C/\sqrt{R_i}\le C\sqrt{2/R_*}$. Thus $q_i$ tends to infinity.
Moreover, $R(q_i)/R_i\to c>0$, so $R(q_i)\to cR_*$ and $R(q_i)\ge cR_*/2>0$ for all sufficiently large $i$.
Note that infinitely many of the points $q_i$ lie on the same half-axis. Since the scalar
curvature has a positive lower bound along this subsequence, the limiting scalar curvature
on this half-axis is strictly positive.
\end{proof}

\begin{lemma}\label{lem:ray-side}
Suppose that the scalar curvature along one half of the fixed
axis tends to $L>0$. Then, on
the part of every sufficiently high level set between that pole and the
maximal orbit, we have
\[
                         |(r^2)_z|\le C_L .
\]
\end{lemma}
\begin{proof}
 On the regular quotient, we define
\[
                b=\frac{1-|\nabla_h r|^2}{r^2}+K.
\]
We first note that
\begin{align}
 \nabla_h^2(f+2\log r)&=Kh-2r^{-2}dr\otimes dr, \label{hess}\\
 |\nabla_h(f+2\log r)|^2&=1+4/r^2-2b. \label{e:mass}
\end{align}
 In fact, we get (\ref{hess}) by differentiating $\log r$ twice and then using the first equation in \eqref{e:wp}. To prove (\ref{e:mass}), we observe that the first formula in \eqref{e:wp} implies
\begin{align}\label{1}
    R=2K-4\Delta_h r/r+2(1-|\nabla_h r|^2)/r^2.
\end{align}
Moreover, by the second equation in \eqref{e:wp}, we have
\begin{align}\label{2}
    \langle\nabla_h f,\nabla_h r\rangle/r=(1-|\nabla_h r|^2)/r^2
-\Delta_h r/r.
\end{align}
The two Laplacian terms in (\ref{1}) and (\ref{2}) cancel, and we obtain \eqref{e:mass}.

By the proof of Proposition \ref{prop:level-profile}, $r_z$ changes its sign only once, namely at the maximum.  Hence, $(f,r)$ are
local coordinates on the part before $r$ reaches the maximum. Moreover, for each pair $0<r<U(f)$,
these coordinates determine a unique orbit on this side. In the following, we differentiate at fixed $r$.
Note that the coordinate vector $\partial_f$ satisfies $dr(\partial_f)=0$
and $d(f+2\log r)(\partial_f)=1$. If we evaluate
\eqref{hess} on $\partial_f,\nabla_h(f+2\log r)$, we get $K$.
Therefore, by differentiating the second identity \eqref{e:mass}, we obtain
\[
          -2\partial_f b=\partial_f|\nabla_h(f+2\log r)|^2=2K>0 .
\]
We also note that, by Proposition~\ref{prop:level-radius}, $U$ increases strictly to infinity.
So, for each fixed radius, the orbit on this side exists at every greater
height.

 Next we compare $b$ with its limiting value as $f\to\infty$. Let us first explain  how to take the limit. Suppose that $x_i$ are points on one half of the fixed
axis such that $x_i\to\infty$ and the scalar curvature $R_i$ at these points tends to $L$. Then $(M,g,p_i,f-f(p_i))$ converge to a Bryant product whose three-dimensional Bryant factor has tip curvature $L$.
As in Lemma~\ref{lem:ray-end}, the centered potentials converge on compact sets.
We denote by $t$ the line coordinate with $t=0$
at the base point, and by $f_B$ the Bryant potential which is normalized to vanish
at the tip. The limiting potential and the pullback of $f_B$ to the
product both have Hessian equal to the Ricci tensor of the product. So the gradient of
their difference is parallel, and hence $f_\infty=f_B+ct$.
Since $\nabla f_B=0$ and $|\nabla f_\infty|^2=1-L$ at the tip, we get
$|c|=\sqrt{1-L}>0$.

Let $s$ be the radial arclength of the Bryant factor. Then the metric of this factor has the form
$ds^2+\omega_L(s)^2g_{S^2}$, where the orbit radius $\omega_L$ increases
from zero to infinity and satisfies $\omega_L'>0$. The quotient of the product has the metric $ds^2+dt^2$,
and we have $K=0$ there. The limiting level set through the base point is the graph
$t=-f_B(s)/c$. This graph reaches each prescribed radius $r>0$ only once, namely at
$s=\omega_L^{-1}(r)$, and the derivative of the radius along its meridian is nonzero there.
In order to find the corresponding points on the original level sets, we follow
this graph from the pole to a point slightly beyond that crossing. Since $\partial_t f_\infty
=c\ne0$, the implicit function theorem can be applied on this compact arc. Then, by the smooth
convergence in meridian charts, we get a nearby level set arc
starting from the chosen pole, and the transversal crossing persists on it. Moreover, the radius
of this arc is increasing at the crossing. Since the radius is strictly concave along
the original meridian, the crossing must lie before the unique
maximum of the radius. Hence, the crossing is the orbit selected by the $(f,r)$ coordinates. Note that the first term of $b$ is the sectional curvature of the orbit tangent plane, and that $K$ is the horizontal sectional curvature.
By equivariance, the radius and its first derivatives converge at that regular orbit,
and the horizontal sectional curvature converges to zero.
Writing $\Phi_L(r)=(\omega_L'(\omega_L^{-1}(r)))^2$, we obtain
\[
              \lim_{f\to\infty} r^2 b(f,r)=1-\Phi_L(r)
                          \qquad (r>0\text{ fixed}).
\]
We remark that every sequence of heights has a subsequence which gives this same value.
This is because the tip curvature determines the Bryant scale, and the right-hand side
does not depend on the sign of $c$.
Together with the monotonicity obtained above, this implies $b(f,r)\ge (1-\Phi_L(r))/r^2$. By the definition of $b$, it follows that
$$|\nabla_h r|^2\le\Phi_L(r)+r^2K.$$ On the Bryant soliton, we have $$\Phi_L(r)\le C_L/(1+r^2).$$
Moreover, by \eqref{e:K} and \eqref{e:der}, we have $r^4K\le C(r^2 R)^2\le C$.
Multiplying by $4r^2$, we now get
$$|\nabla_h(r^2)|^2\le 4r^2\Phi_L(r)+4r^4K\le C_L$$ on the whole of
this side. The derivative in the direction of the unit meridian tangent satisfies the same norm
bound, and by continuity this bound extends to the maximum and to the pole.
\end{proof}

We can now prove Theorem~\ref{th:ray}.

\begin{proof}[Proof of Theorem~\ref{th:ray}]
The theorem is obvious in the Bryant case. So we assume that the soliton is non-Bryant, and we suppose that the scalar curvature does not decay uniformly. By Lemmas~\ref{lem:ray-end} and~\ref{lem:ray-side}, there is one side
of every high level on which the slope of $r^2$ is bounded. Let $s$ be the
arclength measured from the maximal orbit towards that pole, and let
$\ell=f^{1/4}$. Since $|r_s|\le1$, the length of a meridian which goes from radius $U$
to zero is at least $U$. By \eqref{e:size}, we have $U\asymp\sqrt f$.
Hence, the segment $0\le s\le\ell$ exists, and
$U-\ell\le r\le U$ holds on it. In particular, we have $r\asymp\sqrt f$ there.
By the third derivative estimate in \eqref{e:err}, on this segment we get
\[
 |(r^2)_{sss}|=|6r_s r_{ss}+2r r_{sss}|\le C/r\le C f^{-1/2}.
\]
Now we apply Taylor's formula with integral remainder to the function $w=(r^2)_s$.
Note that $w(0)=0$ and $|w(\ell)|\le C_L$. Then the identity
$w(\ell)=\ell w'(0)+\int_0^\ell(\ell-t)w''(t)\,dt$ implies
\[
 |(r^2)_{ss}(0)|\le C_L/\ell+\tfrac12\ell C f^{-1/2}
                   \le C_L f^{-1/4}.
\]
Next, at the
maximum orbit we have
$(U^2)'=Z(r^2)$. By \eqref{e:ZW} and \eqref{e:err}, we get
$$(U^2)'\le2+C_L f^{-1/4}+C/f.$$ Integrating this inequality,
we obtain $$U^2\le2f+C_L f^{3/4}.$$ However, this is impossible by
Proposition~\ref{prop:noenv} with $a=1/4$. Therefore, $R$ tends to zero uniformly. Finally, by Fang's
Theorem~1.3, the limiting sector angle is zero, and we have the
pointed Gromov--Hausdorff convergence to the ray for the full family
of scales and for every fixed base point.

\end{proof}

\section{Asymptotics of the potential level sets}\label{sec:modes}
In this section, we study the level sets of the non-Bryant soliton as $f(x)\to\infty$. There are three different length scales on such a level set. The radius of the largest
sphere is of order $\sqrt f$, the distance between the two poles is
of order $\sqrt{f\log f}$, and near either pole the level set
looks like a Bryant cap at scale $\sqrt{f/\log f}$. We will determine the
constants in these statements, and we will also estimate how the orbit radius changes
between these regions. The latter estimates will be used in
Section~\ref{sec:reflection}, where we have to compare two
orbits of equal radius before we know
that they are exchanged by an isometry.

\subsection{Selecting the cylindrical mode}
We assume that the soliton is non-Bryant.
We take the normal trajectory given by \cite[Proposition~3.6]{MMS}
as the origin of the signed meridian arclength $z$, which increases from the
minus pole to the plus pole. So $z=0$ corresponds to a transported orbit.
The maximizing orbit has its own coordinate $z_m(f)$. We also set $\xi=z/\sqrt f$.
Then a fixed $\xi$-interval describes the middle of the level set at the scale of the
largest sphere.

\begin{proposition}\label{prop:mode-selection}
The cylindrical profile belongs to the neutral alternative of
\cite[Proposition~3.41]{MMS}.
\end{proposition}
\begin{proof}
We first note that all the results in Sections~2.2 and~3
of \cite{MMS} hold under our assumption. By assumption, the perturbation is $r/\sqrt f-\sqrt2$, which is the difference from the
radius of the limiting cylinder. In terms of the logarithmic time $\tau=-\log f$,
the linearized operator in the Gaussian space is
$\partial_{\xi\xi}-(\xi/2)\partial_\xi+1$.
If we apply this operator to $1,\xi,\xi^2-2$, we get $1,\xi/2,0$, respectively.
Hence, the constant component and the linear component have positive eigenvalues, while
the eigenvalue of the quadratic component is zero. According to the dichotomy in \cite[Proposition~3.41]{MMS} for the
localized perturbation, either the positive component or the neutral
component dominates as $\tau\to-\infty$.

Suppose that the positive component dominates. We will derive a contradiction, and for this we only need two properties of the quantity
$\Gamma_+$ in \cite[Lemma~3.39]{MMS}. By its definition, $\Gamma_+$ is a supremum
over earlier times of the positive-mode energy plus a power of the
maximum radius excess, with a positive error allowance. As a consequence, we have
\[
 \left(U/\sqrt f-\sqrt2\right)_+^{\,8-1/200}\le\Gamma_+(\tau).
\]
Here $(a)_+=\max\{a,0\}$ denotes the positive part of a function $a$. Moreover, by the lemma and the assumption that the positive component dominates, we also have
$\Gamma_+(\tau-1)\le(e^{-1}+o(1))\Gamma_+(\tau)$.
We choose $T<0$ such that the factor is at most $e^{-1/2}$ for $\tau\le T$.
For $\tau<T-1$, let $n$ be the integer with $T-1<\tau+n\le T$.
By iteration, $\Gamma_+(\tau)$ is bounded by $e^{-n/2}\Gamma_+(\tau+n)$.
Note that the terminal factor is bounded by $\Gamma_+(T)$ and that $n=-\tau+O(1)$.
It follows that $\Gamma_+(\tau)\le C e^{\tau/2}=C f^{-1/2}$.
Taking the corresponding root in the inequality above, we get
$U/\sqrt f\le\sqrt2+C f^{-a}$ with $a=1/[2(8-1/200)]>0$.
Squaring this inequality and multiplying by $f$, we obtain $U^2\le2f+C'f^{1-a}$ on
every high level. This is impossible by Proposition~\ref{prop:noenv}.
Hence, the neutral alternative holds.
\end{proof}

\subsection{The center, the caps, and scalar curvature}

We begin with an estimate for the size of the maximal orbit.
\begin{proposition}\label{prop:mode-center}
As $f\to\infty$, on each fixed compact $\xi$-interval we have
\begin{equation}
 \frac{r}{\sqrt f}=\sqrt2-\frac{\xi^2-2}{4\sqrt2\,\log f}
       +o((\log f)^{-1})\quad\text{in every spatial }C^j .
 \label{e:neutral}
\end{equation}
Moreover, we have $z_m/\sqrt f\to0$ and
\[
          \frac{U(f)^2}{2}=f+\frac{f}{2\log f}(1+o(1)).
\]
\end{proposition}
\begin{proof}
By Proposition~\ref{prop:mode-selection}, condition (5.1) of
\cite{MMS} is satisfied. So Proposition~5.7 of
\cite{MMS} can be applied in our case. It follows that $(\log f)(r/\sqrt f-\sqrt2)$ converges in  the $C^\infty_{\rm loc}$ sense to the limit
$-(\xi^2-2)/(4\sqrt2)$, and this gives \eqref{e:neutral}.
The $C^1$ convergence is essential in what follows. For any fixed $\delta>0$, the
derivative in $\xi$ is positive at $-\delta$ and negative at $\delta$
when $f$ is large. Since $\partial_\xi(r/\sqrt f)=r_z$ and $r_{zz}<0$,
the unique maximum lies between these two points. Letting $\delta\to0$,
we get $z_m/\sqrt f=o(1)$. Now we can evaluate the uniform
expansion, for example on $[-1,1]$, at this moving maximum. In this way we obtain
$U/\sqrt f=\sqrt2+1/(2\sqrt2\log f)+o(1/\log f)$. The cross term
in the square of this expansion is $1/\log f$, which gives the last formula of the proposition.
\end{proof}

We recall that all the results in Section 5 of \cite{MMS} hold if (A1)-(A3) and the uniform scalar curvature decay are satisfied and the neutral mode dominates. If the soliton has positive curvature operator, all these conditions are satisfied by Theorem \ref{thm:fang}, Theorem \ref{th:ray} and Proposition \ref{prop:mode-selection}. Hence, we have the following estimates.

\begin{proposition}\label{prop:mode-regions}
Let $0<\theta<1/2$ be fixed and sufficiently small. For any given $\epsilon>0$, we
can choose $M$ and then a height such that on all higher levels we have
\begin{align}\label{3}
 (1-\epsilon)\frac{z^2}{2\log f}\le 2f-r^2
 \le(1+\epsilon)\frac{z^2}{2\log f}
 \quad (|z|\ge M\sqrt f,\ r\ge\theta\sqrt{2f}).
\end{align}
On either side, the distance from the reference orbit to the crossing
$r=\theta\sqrt{2f}$, divided by $\sqrt{f\log f}$, tends to
$2\sqrt{1-\theta^2}$, and for the distance to the pole the limit is $2$.
In particular, we have $D(f)=(4+o(1))\sqrt{f\log f}$, and at each pole we have
$R_{{\rm tip},\pm}=(1+o(1))(\log f)/f$.

Let $x=r\sqrt{(\log f)/f}$. On either cap, $r_z^2$, regarded as a function of $x$,
converges smoothly on compact $x$-intervals to $\Phi(x)$,
which is the squared radial slope of the three-dimensional Bryant soliton
with tip scalar curvature one. In the same coordinates, the ambient scalar
curvature divided by $(\log f)/f$ converges smoothly on each compact
cap to the scalar curvature of this Bryant factor.
\end{proposition}
\begin{proof}
Recall that all the results in Section 5 of \cite{MMS} hold in our case. By Propositions~5.12 and~5.13 of \cite{MMS}, we have \[
 \frac{M^2-C(\theta)}{M^2}\frac{z^2}{2\log f}\le 2f-r^2
 \le\frac{M^2+C(\theta)}{M^2-2}\frac{z^2}{2\log f}
 \quad (|z|\ge M\sqrt f,\ r\ge\theta\sqrt{2f}).
\]
If we take $M$ large enough, we get (\ref{3}). Since $r(z,f)$ is strictly concave and
$r(0,f)/\sqrt{2f}\to1$, there is one crossing of the threshold on each side. So we may assume
\begin{align}
    \{z: r(z,f)\ge \theta \sqrt{2f}\}=[z_1(\theta,f),z_2(\theta,f)].\notag
\end{align}
The points $z_1,z_2$ both satisfy (\ref{3}). Then, by $r(z_1,f)=\theta\sqrt{2f}=r(z_2,f)$ and (\ref{3}), we get
\begin{align*}
    \frac{4(1-\theta^2)}{1+\epsilon}\le \frac{z_i^2}{f\log f}\le \frac{4(1-\theta^2)}{1-\epsilon},~i=1,2.
\end{align*}
Letting $\epsilon\to 0$, we obtain the limit.

Next we observe that the
geometric definition of the two tips in \cite[Corollary~2.16]{MMS}
is invariant under the isometries of the level which are induced by the angular action.
Hence, each tip is fixed by $SO(3)$. On a level set, the tips in \cite[Corollary~2.16]{MMS} are the
two poles of Lemma~\ref{lem:level-interval}. We choose one side, and we denote by $d_\theta$ the distance to the radius
$\theta\sqrt{2f}$ on this side and by $d_{\rm tip}$ the distance to
its pole, where both distances are measured from the reference orbit. By \cite[Proposition~5.16]{MMS}, we have $d_{\rm tip}=(2+o(1))\sqrt{f\log f}$, and hence $D(f)=(4+o(1))\sqrt{f\log f}$. Moreover, $R_{{\rm tip},\pm}=(1+o(1))(\log f)/f$ follows from  \cite[Lemma 5.17,Proposition~5.18]{MMS}.

 It remains to consider the caps. By \cite[Proposition~2.13]{MMS}, the pointed level sets rescaled by $R_{\rm tip}$
 converge smoothly to the Bryant soliton at each tip.
As in \cite[Lemma~7.13]{F}, the rotations also converge equivariantly.
On the Bryant soliton $ds^2+\omega(s)^2g_{S^2}$, the radius is strictly
increasing on each bounded $s$-interval. Hence,
$\Phi(\omega(s))=\omega'(s)^2$ defines a smooth function of the radius.
The same change of coordinates can be used on each compact cap of the
approximating level sets. Finally, since
$\sqrt{R_{\rm tip}}/\sqrt{(\log f)/f}\to1$, the convergence still holds on every fixed
compact interval after we replace the scaled radius by $x$.

To compare the ambient and intrinsic scalar curvatures, let $k_c=(\log c)/c$
on $\Sigma_c$, write $\gamma=g|_{T\Sigma_c}$ and $\mathcal S=\operatorname{Scal}_\gamma$,
and take $N=\nabla f/|\nabla f|$. The traced Gauss equation gives
\[
 \mathcal S-R|_{\Sigma_c}
 =-2\Ric(N,N)+(\operatorname{tr}_\gamma\mathrm{II})^2-|\mathrm{II}|_\gamma^2 .
\]
With $B$ as in the proof of Lemma~\ref{lem:level-error}, we have
$\Ric(N,N)=|\nabla f|^{-2}\operatorname{tr}_g B$ and
$\mathrm{II}=\Ric|_{T\Sigma_c}/|\nabla f|$. By \eqref{estimate of B},
\eqref{estimate of v}, \eqref{e:der} and the tangential projection estimates
in that proof, differentiation of the Gauss identity gives
\[
 |\nabla_\gamma^j(R|_{\Sigma_c}-\mathcal S)|_\gamma
 \le C_j R^{2+j/2}\qquad (j\ge0)
\]
on all sufficiently high levels. On each fixed cap the smooth convergence of
$k_c\gamma$ gives $|\mathcal S|\le C k_c$. Since $R\to0$ uniformly by
Theorem~\ref{th:ray}, the estimate for $j=0$ gives
$R\le Ck_c+C_0R^2\le Ck_c+R/2$ for large $c$. Thus $R=O(k_c)$ there.
In the metric $k_c\gamma$, the norm of the $j$th covariant derivative of
$(R|_{\Sigma_c}-\mathcal S)/k_c$ is therefore bounded by
$C_j k_c^{-1-j/2}R^{2+j/2}=O_j(k_c)$. As $\mathcal S/k_c$ is the scalar
curvature of $k_c\gamma$, the ambient scalar curvature $R/k_c$ has the
same smooth Bryant limit on the fixed cap, including at the pole.
\end{proof}

Now we can give a uniform curvature decay estimate.
\begin{lemma}\label{lem:mode-levelR}
On the non-Bryant soliton with the normalization $R+|\nabla f|^2=1$, it holds that
\[
              \lim_{c\to\infty}\frac{c}{\log c}\sup_{\Sigma_c}R=1 .
\]
\end{lemma}
\begin{proof}
We set $k_c=(\log c)/c$ and fix a constant $A>0$. By Proposition~\ref{prop:mode-center}, we have
$U(c)\sqrt{k_c}\to\infty$. Hence, when the level set $\Sigma_c$ is sufficiently high, the set
$\{r\sqrt{k_c}\le A\}$ is the union of two caps, one at each pole.
It follows from Proposition~\ref{prop:mode-regions} that, on each component of $\{r\sqrt{k_c}\le A\}$,
$R/k_c$ converges uniformly to the scalar curvature of the unit-tip Bryant soliton
on $0\le x\le A$.  Since the scalar curvature of the unit-tip Bryant soliton is no greater than one, the supremum of $R/k_c$ on $\{r\sqrt{k_c}\le A\}$
tends to one. On the rest of the level set, by \eqref{e:der} we have
\[
                  \frac{R}{k_c}\le\frac{C}{k_c r^2}\le\frac{C}{A^2}.
\]
Now we fix $A$ such that $C/A^2<1$. Letting $c\to\infty$, we see that
the upper limit is at most one. The lower limit is obtained at either pole,
since $R_{{\rm tip},\pm}(c)/k_c\to1$.
\end{proof}

As a direct consequence of Lemma \ref{lem:mode-levelR}, we obtain Corollary~\ref{cor:sphere-decay}.

\begin{proof}[Proof of Corollary~\ref{cor:sphere-decay}]

Using (\ref{f linear}), the uniform curvature decay and \cite[Lemma 2.2(1)]{DZ20b}, we get
\begin{align}
    \frac{f(x)}{d(x,o)}\to1,~\mbox{as}~f(x)\to\infty,
\end{align}
where we have used the normalization  $|\nabla f|^2+R=1$.  It follows that
\[
    \frac{(\log f(x))/f(x)}{(\log d(x,o))/d(x,o)}
      =\frac{d(x,o)}{f(x)}
          \left(1+\frac{\log(f(x)/d(x,o))}{\log d(x,o)}\right)\longrightarrow1
\]
uniformly on distance spheres.

Then, on the non-Bryant soliton, we get from Lemma~\ref{lem:mode-levelR} that
\[
        \limsup_{d\to\infty}\frac{d}{\log d}\sup_{\partial B_g(o,d)} R\le1 .
\]

Fix one outgoing half of the symmetry axis $\Gamma$ from $o$. Along this
half-axis the distance to $o$ is continuous, starts at zero, and is unbounded,
since $f$ tends to infinity there and is bounded on each closed metric ball.
Thus for every sufficiently large $d$ there is an axis point $x_d$ with
$d_g(o,x_d)=d$. Properness of $f$ gives $f(x_d)\to\infty$, so the distance
asymptotic above gives $f(x_d)/d\to1$. Since $x_d$ is a pole of its
potential level, Proposition~\ref{prop:mode-regions} gives
\[
 \frac{d}{\log d}R(x_d)
 =\frac{d}{f(x_d)}\frac{\log f(x_d)}{\log d}
   \frac{f(x_d)R(x_d)}{\log f(x_d)}\longrightarrow1 .
\]
Consequently,
\[
 \liminf_{d\to\infty}\frac{d}{\log d}
      \sup_{\partial B_g(o,d)}R\ge1 .
\]
Together with the upper bound, this proves the asserted limit.

On the Bryant soliton, the scalar curvature decays exactly linearly.
\end{proof}

\subsection{Recovering slopes between the center and the caps}
In this subsection and in Subsection~\ref{subsec:differentiated-sheets},
we assume again that the soliton is on the non-Bryant soliton.
As in \eqref{e:heat}, we let $y=r^2/2$, and we set $q=U^2/2-y$.
All the limits below are taken along sequences of points on level sets $f\to\infty$.
By a sheet we mean the part of the orbit interval from a given pole to the
maximum. We will first estimate the slopes in the region where the radius is
still a positive fraction of $\sqrt f$. After that, we consider the region where
that fraction tends to zero, but the point has already left every bounded
Bryant cap.

\begin{figure}[htbp]
\centering
\IfFileExists{figures/mode-secants.pdf}{%
  \includegraphics{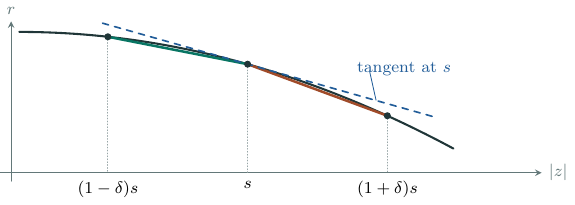}%
}{%
\begin{tikzpicture}[x=.91cm,y=.89cm,>=stealth,font=\footnotesize,
                    line cap=round,line join=round]
  \draw[->,wingGrey] (.65,.22)--(10.72,.22) node[right] {$|z|$};
  \draw[->,wingGrey] (.85,.06)--(.85,3.10) node[above] {$r$};
  \draw[wingInk,line width=1.0pt] plot[domain=1:9.08,samples=65]
                      (\x,{2.90-.034*(\x-1)^2});
  \coordinate (A) at (2.65,{2.90-.034*1.65^2});
  \coordinate (B) at (5.25,{2.90-.034*4.25^2});
  \coordinate (C) at (7.85,{2.90-.034*6.85^2});
  \foreach \p in {A,B,C}
     \draw[wingGrey!66,densely dotted] (\p)--(\p |- 0,.22);
  \draw[wingPath,line width=1.2pt] (A)--(B);
  \draw[wingAxis,line width=1.2pt] (B)--(C);
  \draw[wingMirror,thick,dashed] (2.55,{2.90-.034*4.25^2-.289*(2.55-5.25)})
                              --(8.78,{2.90-.034*4.25^2-.289*(8.78-5.25)});
  \foreach \p in {A,B,C} \fill[wingInk] (\p) circle(1.7pt);
  \node[below] at (2.65,.20) {$(1-\delta)s$};
  \node[below] at (5.25,.20) {$s$};
  \node[below] at (7.85,.20) {$(1+\delta)s$};
  \node[wingMirror,anchor=west] at (7.15,2.20) {tangent at $s$};
  \draw[wingMirror,thin] (7.52,2.15)--(7.64,1.60);
\end{tikzpicture}%
}
\caption{How a derivative is recovered from three values of the radius on the
same level. The plot should be read from the reference orbit toward either pole,
where $s=|z|$. By concavity, the slope of the dashed tangent
lies between the slope of the green left chord and the slope of the rust right
chord. Both chords remain in the region where the radius estimate is valid.}
\label{fig:mode-secants}
\end{figure}

We first give the asymptotics of $r_z$ in the intermediate region.
\begin{lemma}\label{lem:mode-secant}
Suppose that $\liminf r/\sqrt{2f}>0$ and $|\xi|\to\infty$. Then we have
\begin{equation}
        r_z=-(1+o(1))\frac{z}{2r\log f}. \label{e:pslope}
\end{equation}
Moreover, $q\sim z^2/(4\log f)$ lies in this region.
\end{lemma}
\begin{proof}
By passing to a subsequence, we may assume that $z$ has one sign. Let $s=|z|$.
In this proof, we measure distance from the reference orbit along the chosen
direction. Note that the function $r$ is also concave in $s$ on this side.
We fix small numbers $0<\theta_1<\theta_2<\liminf r/\sqrt{2f}$, and
we let $s_i$ be the distance to the outer crossing $r=\theta_i\sqrt{2f}$
on this side. By Proposition~\ref{prop:mode-regions}, we have
$s_i/\sqrt{f\log f}\to2\sqrt{1-\theta_i^2}$. Moreover, the base point lies
at $s\le s_2$. Since the two limiting crossing distances satisfy a strict
inequality, for any sufficiently small fixed $\delta>0$ we have
$(1+\delta)s_2<s_1$ on all late levels. Hence, the points at
$(1-\delta)s,s,(1+\delta)s$ all satisfy $r\ge\theta_1\sqrt{2f}$.
In addition, their distances divided by $\sqrt f$ all tend to infinity. We need
this last property, because the intermediate estimates require both
conditions at the same time.

In the following formula, the height is fixed and $r$ is regarded as a function
of the outward coordinate $s$. For $h=\pm\delta s$, by
Proposition~\ref{prop:mode-regions}, we have
\[
 r(s+h)^2-r(s)^2=-\frac{2sh+h^2+o_\delta(s^2)}{2\log f}.
\]
Here $o_\delta(1)$ denotes a quantity which tends to zero as
$f\to\infty$ for each fixed $\delta>0$, and $o_\delta(s^2)$ denotes
$s^2$ times such a quantity. In fact, once $\delta$ is fixed, the proposition
can be applied at all three points with arbitrarily small relative error.
We also have $s^2\le C f\log f$, and all three
squared radii are at least $2\theta_1^2 f$. Hence, factoring the difference
of squares, we get $$\frac{r(s+h)}{r(s)}=1+O_{\theta_1}(\delta)+o_\delta(1).$$
Next we divide the first identity above by $h(r(s+h)+r(s))$. For both choices
of $h$, the chord slope obtained in this way, divided by $-s/(2r(s)\log f)$, is
equal to $1+O_{\theta_1}(\delta)+o_\delta(1)$. On the other hand, since $r$ is
concave, we have the two-sided bound
\[
 \frac{r((1+\delta)s)-r(s)}{\delta s}
 \le r_s(s)\le
 \frac{r(s)-r((1-\delta)s)}{\delta s}.
\]
We first take the limit in height with $\delta$ fixed, and then we let
$\delta\downarrow0$. Since $r_z=\operatorname{sgn}(z)r_s$, taking the signs
into account we obtain \eqref{e:pslope} on either side. Figure~\ref{fig:mode-secants}
explains why estimates of the radius values are enough for this step.

Finally, we write $2q=(U^2-2f)+(2f-r^2)$. By Proposition~\ref{prop:mode-center},
the first term is $O(f/\log f)$, while the second term is asymptotic to
$z^2/(2\log f)$. Thus the absolute value of their ratio is $O(f/z^2)=o(1)$.
This proves the last assertion of the lemma.
\end{proof}

Then, we get the One-sided control of the second derivative of $y$.
\begin{lemma}\label{lem:mode-second}
There exists a fixed $L_0$ with the following property. For each fixed
$L\ge L_0$ and for all sufficiently high levels, we have
\[
               y_{zz}\le C_L\frac{(\log f)^2}{f}
                    \qquad\text{\quad where }r\sqrt{\log f/f}\ge L .
\]
\end{lemma}
\begin{proof}
\emph{Boundary values on the two caps.}
We begin with the sign of the corresponding quantity on the Bryant soliton.
For the unit-tip Bryant metric, we write $x=\omega(s)$ and denote
its potential by $f_B$. Then the spherical soliton equation reads
$f_B'\omega\omega'=1-\omega'^2-\omega\omega''$.
By \cite[(2.6)--(2.8)]{F}, the right-hand side of this equation tends to one and
$f_B'\to1$. Hence, we have $x^2\Phi(x)=(\omega\omega')^2\to1$.
In this normalization, by the function-value expansion in
\cite[Appendix~B, Proposition~1(4)]{AIK}, we have
$\Phi=x^{-2}+2x^{-4}+O(x^{-6})$. Note that the sphere factor has dimension two
here. Let $A(x)=x^2\Phi(x)$. Then $A=1+2x^{-2}+O(x^{-4})$. It follows that
$A>1$ for large $x$ and $A\to1$. Thus, for any sufficiently large $a$, there
exists $b>a$ such that $A(b)<A(a)$. By the mean value theorem, there is a point
of $(a,b)$ at which $A'<0$.

Next, at any critical point of $A$ we have $\Phi'=-2\Phi/x$. Recall the radial
equation \cite[Appendix~B, Proposition~1(1)]{AIK},
\[
 \Phi\Phi''-\tfrac12(\Phi')^2
             +(1-\Phi)(x\Phi'+2\Phi)/x^2=0 .
\]
At such a point, the last term of this equation vanishes. Since $\Phi>0$, we
get $\Phi''=2\Phi/x^2$. In particular, we have
\[
      A''=2\Phi+4x\Phi'+x^2\Phi''=-4\Phi<0
                      \qquad\text{whenever }A'=0 .
\]
Now we choose $L_0$ such that $A'(L_0)<0$. Suppose that $A'$ has a first zero
to the right of $L_0$. Since $A''<0$ at that zero, we would have $A'>0$
immediately to its left, which contradicts the choice of the first zero. Hence,
we conclude that $\Phi+(x/2)\Phi'=A'/(2x)<0$ for every $x\ge L_0$.

On either cap, differentiating along a level set, we have
\[
 y_{zz}=r_z^2+rr_{zz}=r_z^2+\tfrac{x}{2}\partial_x(r_z^2).
\]
For each fixed $L\ge L_0$, by the smooth cap convergence in
Proposition~\ref{prop:mode-regions}, the value of this quantity at $x=L$ tends to
$\Phi(L)+(L/2)\Phi'(L)<0$. Therefore, the boundary value is strictly negative
on both caps for all sufficiently high level sets.

For the parabolic comparison, we also need the terminal values. Namely, we
claim that $\sup_{x\ge L}(y_{zz})_+\to0$ as $f\to\infty$. To prove this claim, we take
an arbitrary sequence in that region. If $x$ stays bounded, then we can pass to a
convergent subsequence and use the negativity on the Bryant limit. If $x\to\infty$,
then we fix $A>L$, and we see that the base point, viewed from its own pole, eventually lies
beyond $x=A$ on the same sheet. Note that by concavity the nonnegative
outward-from-pole slope decreases to zero at the maximum. So we know that the limsup of
its value at the base point is at most $\sqrt{\Phi(A)}$.
Now we let $A\to\infty$. It follows from $rr_{zz}\le0$ that the positive part
of $y_{zz}$ is bounded by $r_z^2\to0$. Then the uniformity follows from this
sequential argument.

\emph{Evolution and comparison.}
We apply $\partial_{zz}$ to \eqref{e:heat} at fixed signed arclength.
For the moment, we denote the unknown of the resulting scalar equation by
$H=y_{zz}$. It follows from the product rule that the drift term becomes
$(Zz)H_z+2\partial_z(Zz)H+\partial_{zz}(Zz)y_z$.
Then by using \eqref{e:beta}, $y_z=rr_z$ and $H_z=3r_zr_{zz}+rr_{zzz}$,
we get
\[
 \begin{aligned}
 -H_f-H_{zz}
   ={}&(Zz-2r_z/r)H_z-4r_{zz}^2
                    +4r_z^2 r_{zz}/r\\
      &+2E^zH+(E^z)_z y_z-(E^o)_{zz}.
 \end{aligned}
\]
Note that both explicit terms in $r_{zz}$ are nonpositive. We also see from \eqref{e:err} that
the last three terms have absolute value $\le Cr^{-4}$, where we have used
$|H|\le C$, $|y_z|\le r$, $|(E^z)_z|\le Cr^{-5}$ and
$|(E^o)_{zz}|\le Cr^{-4}$. So we know that this source is
at most $C_L(\log f/f)^2$ throughout $x\ge L$.

Now we fix a terminal height $T>f$. If we regard decreasing height as the time
direction, then the spatial operator is forward parabolic. At a height $v$
between $f$ and $T$, we subtract from $H$ the spatially constant function
\[
   \sup_{x\ge L} H_+(\cdot,T)+
             C_L\int_v^T(\log t/t)^2\,dt .
\]
Then the new function satisfies the homogeneous subsolution inequality.
Moreover, it is nonpositive at $v=T$ and on both moving boundaries $x=L$.
Note that these boundaries are smooth on the slab, since the radius crosses
the threshold transversely on each sheet. We also know that the interior contains the
maximal orbit, where $z$ is a regular coordinate. Then by the scalar maximum
principle, we see that the difference is bounded by zero throughout. Finally, we let
$T\to\infty$ and use the terminal estimate. Note that the remaining integral is
equal to $((\log f)^2+2\log f+2)/f$, and so the claim follows.
\end{proof}

\begin{figure}[htbp]
\centering
\IfFileExists{figures/mode-match.pdf}{%
  \includegraphics{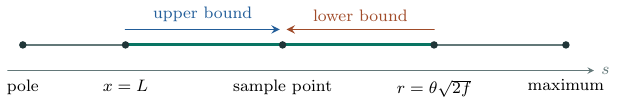}%
}{%
\begin{tikzpicture}[x=.95cm,y=.82cm,>=stealth,font=\footnotesize,
                    line cap=round,line join=round]
  \draw[->,wingGrey] (.45,.34)--(10.90,.34) node[right] {$s$};
  \draw[wingGrey,line width=.8pt] (.72,.87)--(10.40,.87);
  \draw[wingPath,line width=1.5pt] (2.55,.87)--(8.05,.87);
  \foreach \x in {.72,2.55,5.35,8.05,10.40}
     \fill[wingInk] (\x,.87) circle(1.8pt);
  \node[below] at (.72,.30) {pole};
  \node[below] at (2.55,.30) {$x=L$};
  \node[below] at (5.35,.30) {sample point};
  \node[below] at (8.05,.30) {$r=\theta\sqrt{2f}$};
  \node[below] at (10.40,.30) {maximum};
  \draw[->,wingMirror,thick] (2.55,1.19)--(5.30,1.19)
             node[midway,above] {upper bound};
  \draw[->,wingAxis,thick] (8.05,1.19)--(5.42,1.19)
             node[midway,above] {lower bound};
\end{tikzpicture}%
}
\caption{The matching argument on a single sheet of a non-Bryant soliton. The
locations are ordered by the arclength $s$ from the chosen pole. At the left
comparison orbit, we have $x|r_z|\to L\sqrt{\Phi(L)}$ by the Bryant limit. At the
right one, we have $x|r_z|\to\sqrt{1-\theta^2}$ by the intermediate profile.
Using a bound on the increase of $r|r_z|$ along the green segment, we carry
an upper bound from the left and a lower bound from the right.}
\label{fig:mode-match}
\end{figure}

Now we give the estimate in the region matching the Bryant and intermediate slopes.
\begin{proposition}\label{prop:mode-match}
Assume that the soliton is on the non-Bryant soliton. Then on either side, for every
sequence with $x\to\infty$ and $r/\sqrt{2f}\to0$, we have
\begin{equation}
                           x|r_z|\longrightarrow1 .\label{e:capslope}
\end{equation}
\end{proposition}
\begin{proof}
We denote by $s$ the arclength which is measured inward from the chosen pole toward the maximal
orbit. On this sheet we have $r_s=|r_z|$, and hence $(r|r_z|)_s=y_{ss}=y_{zz}$.
We fix $L\ge L_0$ and a small $\theta>0$. Along the given sequence, $x$ tends
to infinity and the ratio $r/\sqrt{2f}$ tends to zero. Hence, by the
monotonicity of $r$, for large $f$ a point of the sequence lies after $x=L$ and
before $r=\theta\sqrt{2f}$. These two comparison orbits are shown
in Figure~\ref{fig:mode-match}. By Lemma~\ref{lem:mode-second}, the total
possible increase of $r|r_z|$ on the segment between them is at most
$C_L D(f)(\log f)^2/f$. Using the length estimate of
Proposition~\ref{prop:mode-regions}, we see that after multiplication by
$\sqrt{(\log f)/f}$ this error is $O_L((\log f)^3/f)=o(1)$.

We first integrate the upper derivative bound from $x=L$ to the point
under consideration. This shows that the value of $x|r_z|$ at this point is at
most the cap-end value plus $o(1)$. Then we integrate again from that point to
the outer threshold, and we see that the value at this point is at least the
threshold value minus $o(1)$. At $x=L$, the smooth cap convergence gives
$L\sqrt{\Phi(L)}$. At $r=\theta\sqrt{2f}$, by \eqref{e:pslope} and the crossing
location in Proposition~\ref{prop:mode-regions}, we get
\[
 x|r_z|=(1+o(1))\frac{|z|}{2\sqrt{f\log f}}
                          \longrightarrow\sqrt{1-\theta^2}.
\]
Therefore, for the original sequence we have
\[
 \sqrt{1-\theta^2}\le\liminf x|r_z|
        \le\limsup x|r_z|\le L\sqrt{\Phi(L)} .
\]
Here the upper endpoint tends to one as $L\to\infty$, and the lower endpoint
tends to one as $\theta\downarrow0$. Note that both auxiliary parameters were
fixed before we took the limit in height, so we do not use any uniform
convergence on expanding Bryant caps.
\end{proof}

We remark that 
$x|r_z|\to\infty$ on the Bryant soliton. So the non-Bryant assumption is necessary.

\subsection{Higher derivatives  estimates for the comparison of the two sheets}\label{subsec:differentiated-sheets}
\begin{lemma}\label{lem:outer}
On either sheet, in each of the two regimes of \eqref{e:pslope}
and \eqref{e:capslope}, we have
\begin{equation}
                 r|r_{zz}|+r^2|r_{zzz}|=o(|r_z|). \label{e:pder}
\end{equation}
Moreover, if two points of a common level have the same radius in either regime,
then the ratio of their absolute slopes tends to one.
\end{lemma}
\begin{proof}
\emph{The moving chart.}
We take an arbitrary sequence of base points in one of the two regimes, and we
denote a base height and radius by $f_0,r_0$. Let $p_0=|r_z|$ at that point.
By the slope estimates above, we have $p_0>0$ eventually, $r_0\to\infty$,
and $r_0p_0\to\infty$. In fact, we have
$r_0p_0\sim |z|/(2\log f_0)$ in the intermediate region, and
$r_0p_0\sim\sqrt{f_0/\log f_0}$ in the cap overlap. We fix a small $c>0$. We
follow the base point by $Z$ for $f_0\le f\le f_0+c r_0^2$, and we let
$z_{\rm rel}$ be the signed arclength from this transported point, with the same
orientation as $z$. We will work in the chart $|z_{\rm rel}|\le c r_0$.
We claim that $r$ is comparable to $r_0$ throughout the chart. To see this,
note first that the center radius is increasing by \eqref{e:ZW}. Moreover, by
\eqref{e:ZW} and \eqref{e:err}, we have $|Z(r^2)|\le C$ at radii at least one.
Hence, the square of the center radius stays between $r_0^2$ and
$(1+Cc)r_0^2$. Since $|r_z|\le1$ and the pole radius is zero, each pole is at
arclength at least the center radius. Therefore, the chart exists, and the
spatial variation of $r$ from its center is at most $c r_0$.

\emph{Control of the normalized slope.}
We claim that $r_z/p_0$ is uniformly bounded in this chart, and that
its oscillation on each spatial slice tends uniformly to zero.
We first consider the cap overlap. In this case we have $r_0^2=o(f_0)$, and
hence $f/f_0\to1$ throughout the chart. Since the radii are comparable, we
still have $x\to\infty$ and $r/\sqrt{2f}\to0$ uniformly there. Moreover, the
entire chart remains on the same sheet, because the maximum has radius
comparable to $\sqrt f$. It follows that the formula \eqref{e:capslope} holds
throughout the chart with uniform relative error. Otherwise, there would be a
sequence of points violating it which still satisfies its hypotheses.
In particular, we have $\sup|r_z|\to0$.
On one time slice, the radii vary by at most $2cr_0\sup|r_z|=o(r_0)$,
and hence their ratios tend to one. Then both claims follow from the formula
$|r_z|=(1+o(1))\sqrt{f/\log f}/r$, and the slope has a constant sign on the
chart.

Next we consider the intermediate region. In this case we have
$r_0\asymp\sqrt{f_0}$ and $f\asymp f_0$,
and $r\ge d\sqrt f$ holds on the whole chart for a fixed $d>0$. We also
have to check that the chart remains far from the original reference at $z=0$.
Since the radius of that reference is asymptotic to $\sqrt{2f}$, by concavity
the segment joining it to the moving center lies in the same thick
region, where $d$ is possibly smaller. On this segment we have
$|r_{zz}|/r+|E^z|\le C/f$. Integrating \eqref{e:beta} on this segment, we get
$|Zz|\le C|z|/f$ at the center. Since $f/f_0$ is bounded, by Gronwall's
inequality the value of $|z|$ at the center stays between fixed multiples of
its initial value. In particular, we have $|z|/r_0\to\infty$ there. If we move
to any other point on its spatial slice, then $z$ changes by at most $c r_0$,
so both its sign and this divergence are preserved. By Lemma~\ref{lem:thick},
we have $\sup|r_z|\to0$ there. Hence, the radii on any one slice again have
ratios tending to one. Now, by the same sequential argument, \eqref{e:pslope}
holds uniformly. Note that in the expression $-z/(2r\log f)$, the ratios
of $z$ and of $r$ at two points on that slice both tend to one,
and their ratios to the original values stay bounded above and
away from zero. So we have proved the claim on the normalized slope.

\emph{The control of values gives control of derivatives here.}
The limits with relative error which we obtained above are only estimates of values.
In order to get estimates of derivatives, we will use the equation on the chart.
In the following formula, the spatial derivatives and the height derivative are
taken in the coordinates $(z_{\rm rel},f)$. Note that $r_z=\partial_{z_{\rm rel}}r$,
because the two arclength origins differ by a translation on each
level. Then by \eqref{e:ZW} and the chain rule in
Proposition~\ref{prop:level-evolution}, we have
\[
 -\partial_f r=r_{zz}+(r_z^2-1)/r+(Zz_{\rm rel})r_z-E^o/r .
\]
We differentiate this equation once in the spatial variable, and we get
\[
 \begin{aligned}
 -\partial_f r_z={}&r_{zzz}+(Zz_{\rm rel}+2r_z/r)r_{zz}\\
 &+\big(\partial_z(Zz_{\rm rel})+(1-r_z^2)/r^2\big)r_z
 -(E^o/r)_z .
 \end{aligned}
\]
We set $a=z_{\rm rel}/r_0$, $t=(f_0-f)/r_0^2$ and $v=r_z/p_0$.
Then the time variable satisfies $-c\le t\le0$ and it increases in the heat
direction, and the estimate that we want is at the latest time $t=0$. In these
units, the principal part is $v_t-v_{aa}$. The drift coefficient and the
potential coefficient are the two coefficients in parentheses above, multiplied
by $r_0$ and by $r_0^2$, respectively. For the drift, we use $Zz_{\rm rel}=0$
at the moving origin and integrate \eqref{e:beta} over a segment of length at
most $c r_0$, and we get $|Zz_{\rm rel}|\le C/r_0$. For its spatial derivatives,
we can just use \eqref{e:beta}. Then we see from the bounds in \eqref{e:err} that
the $a$-derivatives of both coefficients of every fixed order are bounded.
The forcing term is $-(r_0^2/p_0)(E^o/r)_z$. Note that each original spatial
derivative contributes a factor $r_0^{-1}$, and so its $j$th $a$-derivative
is $O_j(r_0^{-2}/p_0)$, which tends to zero because $r_0p_0\to\infty$.

For completeness, we include the interior estimate which is used at $t=0$.
We choose nested space-time rectangles ending at $0$, which have positive
separation from the lateral and initial boundaries of the larger
rectangle. After we differentiate $j$ times in space, we see that the leading terms
in the equation for $\partial_a^j v$ are its second derivative
and a bounded coefficient times its first derivative. It follows from the product
rule that all the other terms are bounded by
$C_j(1+\sum_{i=0}^j|\partial_a^i v|)$. Now we test this equation against
$\partial_a^j v$ times the square of a smooth cutoff function, which is equal
to one on the smaller rectangle and zero near the initial and lateral
boundaries. Then we integrate by parts in $a$ and get the positive integral of
$|\partial_a^{j+1}v|^2$. By the inequality
$2|AB|\le\epsilon A^2+\epsilon^{-1}B^2$, we know that the drift and cutoff terms with one top
derivative are bounded by half of this integral plus lower-order
squares. Then we integrate in time up to any slice and take the supremum over
the slices, and we get
\[
 \sup_t\int_{I'}|\partial_a^j v|^2\,da
   +\iint_{Q'} |\partial_a^{j+1}v|^2\,da\,dt
 \le C_j\left(1+\sum_{i=0}^j\iint_Q|\partial_a^i v|^2\,da\,dt\right),
\]
where $I'$ is the spatial interval of the smaller rectangle $Q'\subset Q$.
We start the induction with the uniform bound on $v$. At each step of the
induction, we see that the space-time bound on the next derivative gives the right-hand
side on a smaller rectangle. We iterate this up to $j=5$, and then we have uniform
$L^2$ bounds on the final slice for the derivatives up to order five near
$a=0$. Then it follows from the one-dimensional Sobolev embedding that $\partial_a^3v$ and
$\partial_a^4v$ are bounded pointwise there.

At $t=0$, we denote by $\operatorname{osc}v$ the oscillation on that smaller
interval. Then for any sufficiently small fixed $h>0$, by using centered first and
second differences and Taylor's formula, we have
\[
 |v_a(0,0)|\le \frac{\operatorname{osc}v}{2h}+C h^2,\qquad
 |v_{aa}(0,0)|\le \frac{2\operatorname{osc}v}{h^2}+C h^2 .
\]
We know from the chart estimates that the oscillation tends to zero. So we can first take the limit along
the sequence, and then let $h\downarrow0$. In the original units,
these two derivatives are $r_0r_{zz}/p_0$ and $r_0^2r_{zzz}/p_0$. This proves
\eqref{e:pder}.

Finally, we fix the same $(f,r)$ on opposite sheets. In the cap overlap,
$x$ is the same at the two points, so the assertion follows directly from
\eqref{e:capslope}. In the intermediate region, the radius fraction is the same.
Suppose that the opposite point has bounded $|\xi|$ along a subsequence. Then
\eqref{e:neutral} would give $q=O(f/\log f)$ there. On the other hand, on the
first side we have
$\dfrac{q}{f/\log f}\sim\dfrac{z^2}{4f}=\dfrac{\xi^2}{4}\to\infty$. This is
impossible, since both points have exactly the same $q$. Thus we have
$|\xi|\to\infty$ on both sides, and the ratio of their absolute
$z$-coordinates tends to one. Applying \eqref{e:pslope}
with the common $r$, we finish the proof.
\end{proof}

\section{Global involution on the soliton}\label{sec:reflection}
\begin{figure}[htbp]
\centering
\IfFileExists{figures/ax-fold.pdf}{%
  \includegraphics{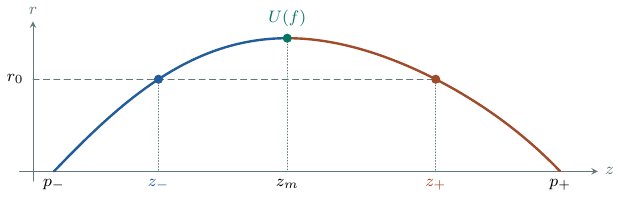}%
}{%
\begin{tikzpicture}[x=.91cm,y=.83cm,>=stealth,line join=round,line cap=round,font=\footnotesize]
  \draw[->,wingGrey] (.42,.22)--(11.18,.22) node[right] {$z$};
  \draw[->,wingGrey] (.67,.02)--(.67,3.28) node[above] {$r$};
  \foreach \lo/\hi/\col in {0/.46086/wingMirror,.46086/1/wingAxis}{
    \draw[\col,line width=1.15pt] plot[domain=\lo:\hi,samples=54]
      ({1.06+9.42*\x},{.22+10.8*(\x*(1-\x)+.65*\x^2*(1-\x)^2*(1-2*\x))});
  }
  \draw[wingGrey,densely dashed] (.67,2.10)--(8.164,2.10);
  \draw[wingGrey,densely dotted] (3.005,2.10)--(3.005,.22)
                                   (8.164,2.10)--(8.164,.22)
                                   (5.401,2.937)--(5.401,.22);
  \foreach \x/\col in {3.005/wingMirror,8.164/wingAxis}
      \fill[\col] (\x,2.10) circle(2.2pt);
  \fill[wingPath] (5.401,2.937) circle(2.2pt);
  \node[anchor=south,wingPath] at (5.401,3.02) {$U(f)$};
  \node[anchor=east] at (.63,2.10) {$r_0$};
  \node[below] at (1.06,.21) {$p_-$};
  \node[below,wingMirror] at (3.005,.21) {$z_-$};
  \node[below] at (5.401,.21) {$z_m$};
  \node[below,wingAxis] at (8.164,.21) {$z_+$};
  \node[below] at (10.48,.21) {$p_+$};
\end{tikzpicture}%
}
\caption{The pairing at equal radius on one potential level. The horizontal
line singles out two sphere orbits, one on each sheet, and at the maximum the
two coordinates coincide. In the graph we only see the strict concavity, which is used
to select the pair. The fact that the meridian distances from $z_-$ and
$z_+$ to $z_m$ are equal will follow from the isometry.}
\label{fig:ax-fold}
\end{figure}

The involution we are looking for must preserve $f$, map each round orbit to an orbit
of the same radius, and exchange the two poles. On a positive level set,
these conditions determine a map of the orbit interval. Namely, we match the two
points where each value $0<r<U(f)$ is attained, and we fix the maximum. This pairing is shown in
Figure~\ref{fig:ax-fold}. But at this stage we do not know yet that it preserves lengths.

It suffices to work on the non-Bryant soliton. We will prove that the scalar
$b=(1-|\nabla_h r|^2)/r^2+K$ which is introduced in Section~\ref{subsec:end}
determines all three coefficients of the quotient metric. There are two
values $b_+(f,r)$ and $b_-(f,r)$, which are indexed by the adjacent pole, and they
satisfy the same elliptic equation. We first show by a comparison on high levels
that their difference is bounded by $C_\nu e^{-\nu f}$. After that we change
variables by using $b_f=-K<0$. More precisely, we fix a radius and a value of $b$, and we
compare the two heights at which this value is attained. Note that the equation for a height
is invariant under adding a constant. So we can use the maximum principle for this equation, and we will
deduce equality at every height from the decay estimate. Only after this step do
we lift the exchange of orbits to an isometry of $M$.

\subsection{An elliptic equation at equal radius}\label{sec:eq}
A \emph{sheet} means the part of the positive-level quotient between
one designated pole and the maximal orbit. By Proposition~\ref{prop:level-profile},
we have $r_z\ne0$ in the interior of a sheet, so $(f,r)$ are independent
coordinates there. In this subsection, the derivative $b_f$ is taken with $r$ fixed, and
$b_r$ is taken with $f$ fixed. When the scalar curvature $R(f,r)$ appears in a formula
on a sheet, it is evaluated on that same sheet.

\begin{lemma}\label{lem:ax-metric}
On each sheet, we have
\begin{equation}
 (r^3 b)_r=\tfrac12 r^2R,\qquad
 b(f,r)=\frac1{2r^3}\int_0^r v^2R(f,v)\dif v,\qquad
 b_f=-K.\label{e:b}
\end{equation}
The value at the pole is $R_{\rm tip}/6$, and the two profiles have the same
value at the maximal orbit. Moreover, the inverse metric in $(f,r)$ is given by
\begin{equation}
 \begin{aligned}
 |\nabla f|^2 &=1-6b-2r b_r,\\
 h^{fr} &=r(2b+b_f)+\tfrac12 r^2 b_r,\\
 |\nabla r|^2 &=1-r^2(b+b_f).
 \end{aligned}\label{e:metric}
\end{equation}
Here the gradients of functions on $Q$ are taken with respect to $h$, and
$h^{fr}=\langle\nabla f,\nabla r\rangle$.
\end{lemma}
\begin{proof}
We differentiate the second identity in \eqref{e:mass}, and we replace
$d|\nabla(f+2\log r)|^2$ by twice the first identity there, evaluated
on $\nabla(f+2\log r)$ and the free argument. Multiplying by $-r^2/2$
and adding $2rb\,dr$, we obtain
\[
 d(r^2b)=-r^2K\,df+
 \left(2\langle\nabla f,\nabla r\rangle
               +2(|\nabla r|^2-1)/r\right)dr
         =-r^2K\,df-2\Delta_h r\,dr .
\]
Here the last equality follows from the spherical equation in \eqref{e:wp}. Thus, we have
$b_f=-K$ at fixed radius, and
$(r^3b)_r=r^2b-2r\Delta_h r$ at fixed height. Substituting the definition of $b$ and
the warped scalar formula $R=2K-4\Delta_h r/r+
2(1-|\nabla r|^2)/r^2$, we get the first identity in \eqref{e:b}.
Since both the angular and the horizontal sectional curvatures stay bounded at a
smooth pole, we have $r^3b\to0$ there. Then the second
identity follows by integration in $r$. By the substitution $v=sr$ in the integral,
we get the continuous limit
$$R(f,0)\int_0^1s^2\,ds/2=R_{\rm tip}/6.$$
At the maximum, the formula defining $b$ only involves a single regular point
of $Q$ with $r=U>0$. So its two limiting values coincide.

Expanding $(r^3b)_r$, we get $R=6b+2rb_r$. Together with
$|\nabla f|^2=1-R$, this gives the first metric entry. The last entry follows by using
$K=-b_f$ in the definition of $b$. For the mixed entry, we use
$r h^{fr}=1-|\nabla r|^2-r\Delta_h r$ and
$-2\Delta_h r=(r^2b)_r$. After substitution we obtain the middle line.
\end{proof}

Now, we have the following  profile equation which can extend smoothly to the pole.

\begin{proposition}\label{prop:ax-b}
On each open sheet, we have the equation
\begin{equation}
 \begin{aligned}
 \mathcal Q_b[b]={}&|\nabla f|^2 b_{ff}+2h^{fr} b_{fr}
       +|\nabla r|^2(b_{rr}+4b_r/r)\\
 &+(|\nabla f|^2+4h^{fr}/r)b_f+h^{fr}(b_r+2b/r)=0 .
 \end{aligned}\label{e:bq}
\end{equation}
The coefficients of this equation are given by \eqref{e:metric}. If we regard it as a radial equation
in five auxiliary Euclidean variables, then it extends smoothly to the pole and is
elliptic there at every positive height.
\end{proposition}
\begin{proof}
We evaluate $\nabla_h^2 f$ on $\nabla f,\nabla r$ in the horizontal equation of \eqref{e:wp}.
By the metric compatibility, we get
\[
 \langle\nabla f,\nabla|\nabla r|^2\rangle
 +\tfrac r2\langle\nabla r,\nabla|\nabla f|^2\rangle
 -rK h^{fr}=0 .
\]
In the coordinates $(f,r)$, the two differential operators on the left hand side are
$|\nabla f|^2\partial_f+h^{fr}\partial_r$ and
$h^{fr}\partial_f+|\nabla r|^2\partial_r$. For example, we have
$\partial_f|\nabla r|^2=-r^2(b_f+b_{ff})$ and
$\partial_r|\nabla f|^2=-8b_r-2rb_{rr}$.
Inserting these two derivatives and the other two derivatives from \eqref{e:metric},
we get the following second-order terms $$-r^2(|\nabla f|^2 b_{ff}+2h^{fr}b_{fr}
+|\nabla r|^2b_{rr}).$$ Moreover, the terms involving $b_f,b_r,b$, after dividing by
$-r^2$, are $$(|\nabla f|^2+4h^{fr}/r)b_f+
(4|\nabla r|^2/r+h^{fr})b_r+2h^{fr}b/r.$$
This proves \eqref{e:bq}.

Next we consider a pole. The smooth orbit-invariant scalar $R$ on a level is even in
the signed meridian arclength at a pole. This is because opposite radial directions are
mapped to each other by the isotropy. Moreover, the radius, when extended with a sign, is smooth
and odd, and its derivative is one with respect to the arclength oriented from the pole.
To see this, we note that in geodesic polar coordinates its square is $s^2$ times a positive
smooth even function which is equal to one at zero. Hence, $R$ and the integral
$b=\tfrac12\int_0^1 t^2R(f,tr)\,dt$ are smooth even functions of
$r$, and they depend smoothly on $f$ on positive slabs. Now we write $r=|v|$ with $v\in\mathbb R^5$.
Then the radial Laplacian in \eqref{e:bq} is $\Delta_v$, and we have
$r b_r=v\cdot D_v b$. The coefficient vector of the mixed derivative is
$(h^{fr}/r)v$, and the function $h^{fr}/r=2b+b_f+(v\cdot D_vb)/2$ is smooth.
All the lower-order terms which look singular can be extended by the same formulas.

It remains to check the ellipticity. On $r>0$, the radial principal block is the Gram matrix of $df,dr$.
Note that the component of $\nabla r$ perpendicular to
$\nabla f$ is $r_z\partial_z$. So the determinant of this matrix in the orthonormal meridian frame is
$|\nabla f|^2r_z^2$, which is positive away from the
maximum. Moreover, each of the four spherical directions in $v$ has coefficient
$|\nabla r|^2>0$. At $v=0$, the mixed vector vanishes and the spatial
block tends to the identity, since $r_z^2\to1$ and $Zr\to0$, where the latter follows from
$\langle\nabla f,\nabla r\rangle=r\Ric(e_o,e_o)$.
We also have $|\nabla f|^2>0$ since the height is positive. So after the lift, the
pole becomes an interior point for the maximum principle. We note that the dimension five
corresponds to the coefficient $4/r$ in \eqref{e:bq}.
\end{proof}

Now, we study the equation for the difference.
\begin{lemma}\label{lem:ax-linear}
As before, let $y=r^2/2$, and let $u=1-2yb$. By subtracting the equations
on paired sheets, we get linear homogeneous equations for $b_+-b_-$
and for $u_+-u_-$, which we denote by $\mathcal L_b$ and $\mathcal L_u$
respectively. For any test function $h$ on $r>0$, they satisfy
$\mathcal L_u(2yh)=2y\mathcal L_b h$. Moreover, their principal matrices are
positive definite on the open sheets.
\end{lemma}
\begin{proof}
Using $\partial_r=r\partial_y$, we can rewrite \eqref{e:metric} as
\[
 |\nabla f|^2=1+(u-1)/y+2u_y,\qquad |\nabla r|^2=u+u_f,
 \qquad h^{fy}=1-|\nabla r|^2-yu_y .
\]
Substituting these into the same horizontal identity, or multiplying
\eqref{e:bq} by $-r^2$, we get
\[
 \begin{aligned}
 \mathcal Q_u[u]={}&|\nabla f|^2(u_{ff}+u_f)+2h^{fy}u_{fy}
             +2y|\nabla r|^2u_{yy}\\
 &+(1-yu_y)u_y+|\nabla r|^2(1-u)/y=0 .
 \end{aligned}
\]
Here we have used that $b=(1-u)/(2y)$ gives $rb_r=-2b-u_y$ and
$b_{rr}+4b_r/r=-u_{yy}-u_y/(2y)-(1-u)/(2y^2)$
when written out after the substitution, where
the radial Laplacian acts on $b$.
Note that each of these expressions is at most quadratic in the unknown and its first and
second derivatives. For any scalar quantity $A$ which is defined on the two paired sheets, we define
\[
\overline{A}(f,r):=
\frac{A_+(f,r)+A_-(f,r)}{2}.
\]
Since $y=r^2/2$, we use the same convention in $(f,y)$ coordinates.
The bar is always applied to the whole quantity under it. For example, we have
\[
\overline{|\nabla q|^2}
=\frac{|\nabla q_+|^2+|\nabla q_-|^2}{2},
\qquad
\overline{\langle\nabla f,\nabla q\rangle}
=
\frac{
\langle\nabla f,\nabla q_+\rangle+
\langle\nabla f,\nabla q_-\rangle
}{2}.
\]
For a product, we have the exact identity
$a_+d_+-a_-d_-=\bar a(d_+-d_-)+\bar d(a_+-a_-)$,
where $\bar{a}=\frac{a_++a_-}{2}$, and similarly for $\bar{b}$. We apply this identity to each product.
Then the difference is given by
the derivative of the polynomial at the mean of the two sets of
values and derivatives, applied to their difference.  In particular,
the principal matrix is the arithmetic mean of two positive definite
matrices. So it is also positive definite. Finally, differentiating the identity
$\mathcal Q_u[1-r^2b]=-r^2\mathcal Q_b[b]$ in a test direction,
we get the identity for the linear operators stated in the lemma.

\end{proof}

For later use, we write down the $u$ operator in these coordinates. It is
\[
 \mathcal L_u=\overline{|\nabla f|^2}\,\partial_{ff}
 +2\overline{h^{fy}}\,\partial_{fy}
 +2y\overline{|\nabla r|^2}\,\partial_{yy}
 +C_f\partial_f+C_y\partial_y+C_0 ,
\]
where
\begin{align}
 C_f&=1+2\bar u_y-2\bar u_{fy}+2y\bar u_{yy},\nonumber\\
 C_y&=1-2y\bar u_y-2y\bar u_{fy}+2\bar u_f+2\bar u_{ff},\nonumber\\
 C_0&=(1-2\bar u+\bar u_{ff})/y-2\bar u_{fy}+2y\bar u_{yy}.\label{e:lin}
\end{align}
We remark that at the fold these coefficients will be estimated after averaging the
sheets, because the individual radius derivatives can be singular.

\subsection{A strict supersolution on both sheets}\label{sec:bar}
The difference $b_+-b_-$ satisfies a homogeneous elliptic equation, but the
zeroth-order coefficient of this equation may be positive or negative. In order to get
a comparison principle for both signs of $b_+-b_-$, we divide $b_+-b_-$ by a positive
strict supersolution, so that the zeroth-order coefficient becomes negative.
More precisely, suppose that $P>0$ and $\mathcal L_bP<0$. Applying the product rule to $P$ times
the quotient, we see that the zeroth-order coefficient in the equation for
$(b_+-b_-)/P$ is $(\mathcal L_bP)/P<0$.
We need to find one such $P$ which is adapted to the caps and to the fold at the same time.
We will use the exponent $5/4$, which lies between $1$ and $3/2$. With this exponent, the leading
singular term in the central calculation is negative, and the barrier
dominates the order-$3/2$ vanishing of $b_+-b_-$ at the place where the sheets meet.
We will verify these requirements in the next three local estimates
and in the sign calculation after them. Recall that $x=r\sqrt{(\log f)/f}$
and $q=(U^2-r^2)/2$. Here $x$ measures the radius in tip units, and $q$
vanishes at the fold. We look for a function
$J_0(x)>0$, which is even and smooth at zero and is equal to $x^{-6}$ for large $x$, such that
\begin{equation}
 H_b=f^{-7/8}J_0(x)(q/f)^{5/4},\qquad H_u=2yH_b
 \label{e:bar}
\end{equation}
are strict supersolutions for their respective operators on all sufficiently high levels, away from
the fold. We also check the additional exponential factor which is needed at the lower comparison boundary.

\subsubsection{Bounded Bryant radii}
\begin{lemma}\label{lem:ax-cap}
There exists a positive even function $J_0$, which is smooth at zero and equal to $x^{-6}$
outside a compact interval, such that
$\mathcal L_b H_b<0$ holds on each fixed cap $0\le x\le A$ for all sufficiently high levels.
\end{lemma}
\begin{proof}
Let $b_B(x)=(1-\Phi(x))/x^2$ be the Bryant profile in \eqref{e:b}.
We first show that $fb/\log f\to b_B$ on fixed $x$-intervals together with all spatial
derivatives in $x$, and that this convergence also holds at the pole.
For this cap calculation, we set $k=(\log f)/f$. Substituting $v=rs$ in the integral
in \eqref{e:b} and using $r=x/\sqrt k$, we obtain
\[
 \frac{b}{k}=\frac12\int_0^1 s^2
       \frac{R(f,sx/\sqrt k)}{k}\,ds .
\]
By Proposition~\ref{prop:mode-regions}, the
ambient scalar curvature $R/k$ converges smoothly on both caps to the unit-tip Bryant
scalar. Note that in the product limit the horizontal plane is spanned by the
line and the Bryant radial direction, so its sectional curvature is
zero. Moreover, the squared radius gradient converges to $\Phi$. Hence, by the definition
of $b$, the limit of the integral is equal to $b_B$ for $x>0$.
Using the integral over $0\le s\le1$ and the smooth evenness at the pole,
we get the convergence with all spatial derivatives down to $x=0$.
Moreover, we have $b_f=-K=O((\log f)^2/f^2)$
together with its scaled spatial derivatives, and $b_{ff}=O((\log f)^{5/2}/f^{5/2})$. To see this,
we note that at fixed $r$ the
differentiation vector is $Z-(Zr/r_z)\partial_z$. This vector is bounded on a fixed Bryant cap,
since $|r_z|$ has a
positive lower bound there and $Zr=O(r(\log f)/f)$. Then we use \eqref{e:K}, which also
applies at $r=0$ by evenness. Finally, the
factors $(q/f)^{5/4}$ tend to one on these caps together with their scaled spatial derivatives.

Now we consider the exact Bryant product. We choose a nonzero affine line component of the
potential, and we rescale the
Bryant factor so that the total normalization holds. Since the quotient has $K=0$, we have $b_f=0$.
Substituting
$b=k b_B(\sqrt k\,r)$ into \eqref{e:bq} and dividing by $k^2$, we get the radial equation
\[
 (1-x^2 b_B)(b_B''+4b_B'/x)
       +(2xb_B+x^2 b_B'/2)(b_B'+2b_B/x)=0 .
\]
Here $k$ is the tip scalar curvature of the product. In the next calculation, we take the
radial limit of
the linearized operator on the original soliton. For this we only use $K=O((\log f)^2/f^2)$ there.

Let $L_B$ be the linearization of this radial equation at $b_B$. We claim that on any bounded cap,
\begin{equation}
 \frac{f}{\log f}\,\frac{\mathcal L_b H_b}{H_b}\longrightarrow \frac{L_B J_0}{J_0}.
 \label{e:captest}
\end{equation}
In fact, we have $|\nabla f|^2=1+O((\log f)/f)$, $|\nabla r|^2\to\Phi$ and
$h^{fr}/\sqrt{(\log f)/f}\to2xb_B+x^2b_B'/2$. For a radial test $J_0(x)$,
by linearizing \eqref{e:metric}, we get
\[
 \begin{aligned}
 \delta|\nabla f|^2&=-6J_0-2xJ_0',&
 \delta|\nabla r|^2&=-r^2(J_0+\partial_f J_0),\\
 \delta h^{fr}&=r(2J_0+\partial_f J_0)+r^2\partial_r J_0/2 .
 \end{aligned}
\]
After dividing by $(\log f)/f$, the leading terms in \eqref{e:bq} are exactly
the variations of the two radial products above. On compact caps, the coefficients of
$b_f=O((\log f)^2/f^2)$ are $O(1)$, and we have $(\delta h^{fr})b_{fr}=O((\log f)^2/f^2)$.
The first $f$-derivative of the test with all factors in
\eqref{e:bar} is of order $O(1/f)=o((\log f)/f)$ relative to its size, its second derivative is
$O(f^{-2}+|(U^2/2)''|/f)$, and the mixed term is $O((\log f)/f^2)$. The bound $(U^2/2)''=O(1/f)$
will be verified below, independently,
from the central evolution. By the five-dimensional
lift, these estimates also hold at zero.

We now explain how to follow the limit in \eqref{e:captest} so that no
powers of $k$ are lost in the terms which involve a variation of the coefficients.
We write $b=k B(x)$ on one fixed level set, and for the moment we omit all $f$
derivatives. Then the two radial products in \eqref{e:bq} become $k^2$
times the radial expression displayed above, with $B$ in place of $b_B$.
If we vary $b$ by a radial test $J$, then $B$ varies by $J/k$. Therefore, their
linearization is $k$ times the radial linearization on $J$.
Note that the scaled spatial profile $H_b/f^{-7/8}$ tends to $J_0$ in every
fixed spatial norm by (\ref{e:bar}). Next, the terms omitted above are $o(k)$ relative
to $H_b$. For example, we have $b_{fr}=O(k^{5/2})$, and a coefficient variation
$\delta h^{fr}/H_b=O(k^{-1/2})$ gives $O(k^2)$. The factors
$1/x$ which occur at the pole multiply smooth odd derivatives, and the quotients of these
derivatives by $x$
are bounded by their first derivative. Thus, the convergence is uniform
on closed caps. The estimate on $(U^2/2)''$ which we used here follows directly
from the scalar evolution in the next lemma, and it does not depend on $J_0$.

From the Bryant scaling family $k b_B(\sqrt k\,x)$, we get the positive Jacobi function $W_B=b_B+x
b_B'/2=-\Phi'/(2x)$, which satisfies $W_B(0)>0$. In fact, differentiating the exact family in $k$ at
one, we get $L_B W_B=0$. The equality with $-\Phi'/(2x)$ follows by
differentiating $x^2 b_B=1-\Phi$. Since
$\Phi'(\omega(s))=2\omega''(s)<0$ and
$\Phi(x)=1-x^2/6+O(x^4)$ at the scalar-one tip
\cite[(2.7)]{F}, this function is strictly positive up to zero.
For $J_0=W_B v$, we compute
\[
 W_B^{-1}L_B(W_Bv)=\Phi v''+\gamma_B v',\quad
 \gamma_B=4\Phi/x+3xb_B+x^2b_B'+2\Phi W_B'/W_B .
\]
At the pole, we have $\gamma_B=4/x+O(x)$. At infinity, we apply
Lemma~\ref{lem:bryant-derivatives} through order three to the expansion
$\Phi=x^{-2}+2x^{-4}+O(x^{-6})$ which was used in Lemma~\ref{lem:mode-second}.
Then we get $W_B=x^{-4}+4x^{-6}+O(x^{-8})$ together with the corresponding first and
second derivative bounds. Inserting $b_B=(1-\Phi)/x^2$ and
$W_B'/W_B=-4/x+O(x^{-3})$ into the expression above, we obtain
$\gamma_B=x^{-1}+O(x^{-3})$. Therefore, for $v=x^{-6}/W_B$, we have
\[
 v=x^{-2}+O(x^{-4}),\qquad
 v'=-2x^{-3}+O(x^{-5}),\qquad v''=6x^{-4}+O(x^{-6}).
\]
It follows that $\Phi v''+\gamma_B v'=-2x^{-4}+O(x^{-6})<0$ for all
sufficiently large $x$, and we also have $v'<0$ there.
We choose $I>0$ such that $I'/I=\gamma_B/\Phi$ and $I\sim x^4$ at zero. On the outer interval,
$-Iv'$ is positive and strictly increasing. We extend it smoothly to a
strictly increasing
positive function on $x>0$ which is equal to $cIx$ near zero. This can be done by choosing $c$
small and interpolating its positive
derivative between two fixed join points by a smooth positive function with the prescribed integral.
Then we integrate $v'$ inward, starting from its prescribed value at the outer join. For the
function obtained in this way, we have $v>0$ and $(Iv')'<0$, and it is even at zero.
Moreover, the limit of $\Phi v''+\gamma_B v'$ at zero is $-5c$. Dividing by the positive integrating
factor, we get $\Phi v''+\gamma_Bv'=\Phi(Iv')'/I$. Therefore, we have constructed $J_0$ such that
$L_BJ_0<0$ on every fixed cap. By continuity, the maximum of
$L_BJ_0/J_0$ on each such compact
cap is strictly negative. The lemma now follows from this strict sign and \eqref{e:captest}.
\end{proof}

\subsubsection{The coefficients at the maximum}
At the maximum, the radius coordinate degenerates. We first work in a signed
coordinate which is regular on both sides. Then we take the even part in order to
recover the coefficients of the linearized equation at equal $y$.

\begin{lemma}\label{lem:ax-center}
On each bounded central interval, the two sides have a smooth signed
coordinate $\eta=\operatorname{sgn}(z-z_m)\sqrt{q\log f/f}$.
The means at equal radius extend across the fold as smooth functions
of $\eta^2$. Moreover, the estimates  below hold uniformly,
\begin{equation}
 \begin{aligned}
 \frac{\overline{|\nabla q|^2}}q&=\frac{1+o(1)}{\log f},&
 (\mathcal L_u-C_0)q&=-\frac{1+o(1)}{2\log f},\\
 |\overline{\langle\nabla f,\nabla q\rangle}|&\le C\eta^2,&
 |\overline{\langle\nabla y,\nabla q\rangle}|&\le C\eta^2+Cq/\log f .
 \end{aligned}\label{e:centralND}
\end{equation}
and we have $(U^2/2)''=O(f^{-1})$.
\end{lemma}
\begin{proof}
In this calculation, all estimates are considered on $|\xi|\le L$, where $L$ is arbitrary and fixed.
We recenter at $\xi_m$
and set $\eta=\operatorname{sgn}(z-z_m)\sqrt{q\log f/f}$. By Taylor's formula with integral remainder and
\eqref{e:neutral}, we have
\[
 \eta=(\xi-\xi_m)(1/2+o(1))
\]
in every fixed spatial $C^j$ norm, which means that the coefficient converges in this sense together
with its derivatives. To see this, we note that $q\log f/f$ has a double zero and that its second
derivative is $1/2+o(1)$. By Taylor's formula, a smooth function with this
double zero can be written as $(\xi-\xi_m)^2$ times the integral of its second
derivative against $1-t$, for $0\le t\le1$. This integral tends to
$1/4$ together with its spatial derivatives, so its positive square root tends
to $1/2$. Note also that the two points at the same $q$ have exactly opposite values of $\eta$.
Thus, $\eta$ and its inverse are
uniformly regular on the compact intervals which we consider. For an even function written as
$v(f,\eta)=V(f,s)$ with $s=\eta^2$,
the notation $C^j(\eta^2)$ will mean $C^j$ regularity in $s$ at fixed $f$.
We also write $\partial_{\eta^2}v$ for $\partial_s V$ evaluated at
$s=\eta^2$. We remark that the even part of a $C^{2j}$
remainder of size $o(1)$ in $\eta$ is $o(1)$ in $C^j(\eta^2)$ up to $\eta=0$. This follows by
applying $v'(\eta)/\eta=\int_0^1
v''(t\eta)\dif t$ for even $v$ and iterating.

Throughout the level, in a meridional frame, we have
\begin{equation}
 \begin{aligned}
 u-r_z^2&=|\nabla f|^2(Zr)^2-r^2K,\\
 |\partial_z^j(u-r_z^2)|&\le C_j r^{-2-j}\quad(r\ge1).\end{aligned}\label{e:p2}
\end{equation}
In fact, we have $u=|\nabla r|^2-r^2K$, and the squared norms of the tangential and normal components of
$\nabla r$ are $r_z^2$ and $|\nabla f|^2(Zr)^2$. The derivative
bounds follow from $\langle\nabla f,\nabla r\rangle=r\Ric(e_o,e_o)$,
\eqref{e:K} and \eqref{e:der}. Thus, by \eqref{e:neutral}, we get $(\log f)^2\bar u\to \eta^2/2$ in
every fixed $C^j(\eta^2)$. More explicitly, we have
$$(\log f)r_z=-\xi/(2\sqrt2)+o(1),\qquad\xi=\xi_m+(2+o(1))\eta$$ with $\xi_m=o(1)$.
Moreover, the error in \eqref{e:p2} is $O(f^{-1})$ together with its scaled spatial
derivatives. Squaring and taking the even part, we obtain the limit above.
At fixed $f$, we have $\partial_y=-(\log f/f)\partial_{\eta^2}$ on even
functions, and the second derivative of the limiting affine function
of $\eta^2$ is zero. In particular, we have
\[
 \bar u_y=-\frac{1+o(1)}{2f\log f},\qquad
 \bar u_{yy}=o(f^{-2}).
\]
Next we deal with the $f$-derivatives. Here we differentiate first and take the average afterwards.
On each open sheet, we have $u_f=2yK$.
Recall that the bar denotes the arithmetic mean at the two paired
points with the same $(f,y)$ coordinates, which was defined in Lemma~\ref{lem:ax-linear}.
Let $K_+$ and $K_-$
be the values of $K$ at these points. Then we obtain
\[
\bar u_f=\overline{2yK}
=\frac{2yK_+ + 2yK_-}{2}
=y(K_+ + K_-).
\]
We apply the even-part argument to this mean $\bar u_f$. Since $r\asymp\sqrt f$, the derivatives
$\partial_z^j K$ carry an extra factor $O(f^{-j/2})$ in addition to the $O(f^{-2})$ bound on
$K$. This follows from \eqref{e:K} in a meridian frame and the tangent projection
estimates of Lemma~\ref{lem:level-error}. At first order, a derivative in $\eta$
is at most $O(\sqrt f)$ times a $z$ derivative,
and the higher rescaled coordinate changes are bounded. Multiplying by
$2y=O(f)$, we see that any fixed number of $\eta$
derivatives of $2yK$ is also $O(f^{-1})$. Hence, after taking the mean, we get $O(f^{-1})$
with bounded derivatives of the same order in $\eta^2$. It follows that $\bar u_{fy}=O(\log f/f^2)$.
Now we give the details for the remaining time
derivatives. By \eqref{e:ZW}, we have $Zy=1-y_{zz}+E^o$. At the maximum, we have
$$(U^2/2)'=Zy=1+(1+o(1))/(2\log f).$$ In the central region, by \eqref{e:neutral} and \eqref{e:err}, we have
$|\partial_z(Zy)|=O(1/(\sqrt f\log f))$ and $|\partial_{zz}(Zy)|=O(1/(f\log f))$.
Since
$[Z,\partial_z]=-\partial_z(Zz)\partial_z$, the maximum moves relative
to $Z$ with speed $-\partial_z(Zy)/y_{zz}=O(f^{-1/2})$. We also have
\[
 Z y_{zz}=\partial_{zz}(Zy)-2\partial_z(Zz)y_{zz}-\partial_{zz}(Zz)y_z.
\]
The right-hand side is $O(f^{-1})$. To check the remaining error derivative, we use
$E^o=r^2 E(e,e)$ and the expression for $|\nabla f|^2E$ in the proof
of Lemma~\ref{lem:level-error}. Here we have $R=O(f^{-1})$ by $Rr^2\le C$, and
$|Z|=O(1)$. Moreover, the ambient differentiated tensor is $O(R^{5/2})$.
This is obtained by using the bounds in that lemma when we differentiate $B$, the two Ricci
products, and their projections. Note that a unit angular field can be chosen
parallel along a horizontal curve, namely the fixed spherical lift divided
by its radius. We also have $Zr=O(f^{-1/2})$ by \eqref{e:ZW}. Hence, we get
$$ZE^o=O(f R^{5/2})+O(R^2)=O(f^{-3/2}).$$
Now we differentiate $1-y_{zz}+E^o$ along the moving maximum. Its derivative
along $Z$ is $O(f^{-1})$ by the commutator calculation, and the extra
meridian term is its $z$-derivative times a speed $O(f^{-1/2})$.
Thus, we have $(U^2/2)''=O(f^{-1})$.
At fixed signed $\eta$, the differentiation vector has the form $Z+v\partial_z$ with
$v=O(f^{-1/2})$. To see this, away from zero we solve
\[
 (U^2/2)'-Zy-vy_z=q(1/f-1/(f\log f)).
\]
By the bound on $\partial_z(Zy)$, the difference satisfies $(U^2/2)'-Zy=O(|\eta|/\log f)$,
while $|y_z|\asymp\sqrt{f}|\eta|/\log f$. The
quotient extends through zero, and the same estimate holds there. Hence, by \eqref{e:K} and the
bounded rate of change of $y$, we get
$\left.\partial_f(2yK)\right|_{\eta}=O(f^{-3/2})$. If we differentiate $\eta^2=q\log f/f$ at fixed
$y$, we get
$(U^2/2)'(\log f)/f+q(1-\log f)/f^2=O((\log f)/f)$ on this region. We take the average first and
then use the $O(1/f)$ bound for the
derivative of $\overline{2yK}$ with respect to $\eta^2$. This gives $\bar u_{ff}=O(f^{-3/2}+\log f/f^2)$.

When we apply the product rule to $q^{5/4}$, we need the gradient of $q=U^2/2-y$
and its pairings with $\nabla f$ and $\nabla y$. In the central region, we have \eqref{e:centralND}. 

We now prove estimates \eqref{e:centralND}. On each sheet, we have $$|\nabla q|^2=|\nabla f|^2((U^2/2)'-Zy)^2+y_z^2.$$
As we approach the double zero, \eqref{e:neutral} gives
$y_z^2/q=(1+o(1))/\log f$ uniformly. The first term is
$O(\eta^2/(\log f)^2)=o(q/\log f)$.
For the second estimate, we use
$$(\mathcal L_u-C_0)q=\overline{|\nabla f|^2}(U^2/2)''+C_f(U^2/2)'-C_y.$$
By \eqref{e:lin}, we have $C_f=1+o(1/\log f)$,
$C_y=1+(1+o(1))/\log f$, and $yC_0=1+o(1)$.
The pairing $\langle\nabla f,\nabla q\rangle=|\nabla f|^2((U^2/2)'-Zy)$
vanishes at the maximum, and its first two $\eta$-derivatives are bounded.
Therefore, its even part is $O(\eta^2)$. Finally, we have
$$\langle\nabla y,\nabla q\rangle=(Zy)\langle\nabla f,\nabla q\rangle-y_z^2,$$
and the same even-part estimate can be applied to the first term.
This proves \eqref{e:centralND}, including its uniform interpretation
as $q$ tends to zero on an open sheet.
\end{proof}

\subsubsection{The coefficients outside a bounded central interval}
\begin{lemma}\label{lem:ax-outer}
In the two regions of Lemma~\ref{lem:outer}, both the following sheetwise estimates
\eqref{e:outer2} and the averaged estimates \eqref{e:outerC} hold:
\begin{equation}
 u=o(1),\quad yu_y=o(|r_z|),\qquad y^2u_{yy}=o(1).\label{e:outer2},
\end{equation}
and
\begin{equation}
 \begin{gathered}C_f=1+O(1/y),\qquad C_y=1+o(|r_z|),\qquad yC_0=1+o(1),\\
 \overline{|\nabla f|^2}=1+o(1),\qquad
 \overline{h^{fy}}=O(1),\qquad
 \overline{|\nabla r|^2}=O(r_z^2+r^{-2}).\end{gathered}\label{e:outerC}
\end{equation}
Moreover,
 \begin{align}
 \overline{|\nabla q|^2}&=O(1+yr_z^2),&
 \overline{\langle\nabla f,\nabla q\rangle}&=O(1),\label{5}\\
 \overline{\langle\nabla y,\nabla q\rangle}&=O(1+yr_z^2),&
 (\mathcal L_u-C_0)q&=O((\log f)^{-1}+y^{-1})+o(|r_z|).\label{6}
 \end{align}

\end{lemma}
\begin{proof}
We consider either of the two regimes of Lemma~\ref{lem:outer}, always with $x\to\infty$. We first
note that all estimates which we make for the two sheets can be averaged. This is because, for a
common $y$, they belong to the same regime. Moreover, in
\eqref{e:pslope} their $|z|$ and $|r_z|$ are asymptotically equal by $q\sim z^2/(4\log f)$, and in
\eqref{e:capslope} both values of $|r_z|$ are asymptotic to $1/x$. On these sequences, we have
$r_z\to0$ and
$r|r_z|\to\infty$. By \eqref{e:p2} and $\partial_y=(rr_z)^{-1}\partial_z$, we get \eqref{e:outer2}.

For the main term $r_z^2$, the first derivative in $y$ is $2r_{zz}/r$, and the second
derivative is
$2r_{zzz}/(r^2r_z)-2r_{zz}/r^3$. For the remainder, by the chain rule we have
\[
 (u-r_z^2)_{yy}
 =\frac{\partial_{zz}(u-r_z^2)}{r^2r_z^2}
 -\frac{(r_z^2+rr_{zz})\,\partial_z(u-r_z^2)}{r^3r_z^3}.
\]
After multiplying by $y^2$, the two remainder terms are bounded by
\[
 O\!\left(\frac1{r^2 r_z^2}\right),\qquad
 O\!\left(\frac{r_z^2+r|r_{zz}|}{r^2|r_z|^3}\right).
\]
Both of them tend to zero by \eqref{e:pder}, since
$r|r_{zz}|=o(|r_z|)$ and $r|r_z|\to\infty$.
For the first derivative, the remainder multiplied by $y$ is
$O(r^{-2}/|r_z|)=o(|r_z|)$. At fixed $y$, the differentiation vector is
$Z-(Zr/r_z)\partial_z$. Its norm is bounded, because $Zr=O(1/r)$ by
\eqref{e:ZW} and $r|r_z|\to\infty$. Hence, by \eqref{e:K}, we have
$u_f=2yK=O(r^{-2})$, $u_{ff}=O(r^{-3})$, and $y u_{fy}=O(r^{-2}+r^{-2}/|r_z|)=o(|r_z|)$.
Inserting these bounds into \eqref{e:lin}, we obtain \eqref{e:outerC}. Here we used that $h^{fy}=r\langle\nabla f,\nabla r\rangle=r^2\Ric(e_o,e_o)$ is bounded.
Since $q=U^2/2-y$, (\ref{5}) and (\ref{6}) follows from the estimates of the principal coefficients.

\end{proof}

\subsubsection{Barriers for $\mathcal L_b$ and $\mathcal L_u$}
\begin{proposition}\label{prop:ax-super}
For each $0<\nu<1$, there exists $S_0$ such that both $H_b$ and
$e^{-\nu(f-S_0)}H_b$ are positive strict supersolutions for
$\mathcal L_b$ on $f\ge S_0$, $0\le r<U(f)$, where we use the lift at
the pole. On $r>0$, multiplying by $2y$ gives the corresponding
supersolutions for $\mathcal L_u$.
\end{proposition}
\begin{proof}
When $x$ is large, we separate the explicit weight
$V=f^{7/8}(\log f)^{-3}y^{-2}$, so that $H_u$ is a fixed positive
multiple of $Vq^{5/4}$. At fixed $y$, we have
$V_f/V=(7/8-3/\log f)/f$ and $V_{ff}/V=O(f^{-2})$.
The exponents are obtained by substituting $x^{-6}=f^3 r^{-6}(\log f)^{-3}$
in \eqref{e:bar} and multiplying by $2y$. Note that for a second order
operator with principal matrix $a^{ij}$, the cross term on a product
$vw$ is $2a^{ij}v_i w_j$. In particular, we have
$\partial_{ij}q^p=pq^{p-1}q_{ij}+p(p-1)q^{p-2}q_iq_j$.
Taking $p=5/4$, we get

 \begin{align}
 \frac{\mathcal L_u V}{V}
 =\overline{|\nabla f|^2}\frac{V_{ff}}{V}
    -\frac{4\overline{h^{fy}}}{y}\frac{V_f}{V}
    +\frac{12\overline{|\nabla r|^2}}{y}
    +C_f\frac{V_f}{V}-\frac{2C_y}{y}+C_0,\label{e:sign}\\
 \frac{\mathcal L_u H_u}{H_u}
 =\frac{\mathcal L_u V}{V}+\frac{5}{4q}\bigg(
    (\mathcal L_u-C_0)q+2\frac{V_f}{V}\overline{\langle\nabla f,\nabla q\rangle}-\frac4y\overline{\langle\nabla y,\nabla q\rangle}\bigg)
   +\frac{5}{16q^2}\overline{|\nabla q|^2}.\notag
 \end{align}

We first consider the case $r/\sqrt{2f}\to0$. In this case, we have $q\sim f$ and $y=o(f)$.
Multiplying the first identity of \eqref{e:sign} by $y$ and using
\eqref{e:outerC}, we get $-1+o(1)$, and the other terms times $y$ are $o(1)$. Next we consider the
case that $r/\sqrt{2f}\to\rho_0>0$ and
$|\xi|\to\infty$. Then $f\mathcal L_u V/V$ tends to $7/8-1/\rho_0^2\le-1/8$. If $q\asymp f$, then all
correction terms times $f$ tend to zero, since $r_z=o(1)$, $\overline{|\nabla r|^2}=o(1)$ and
$\overline{|\nabla q|^2}=O(1+yr_z^2)$. In the remaining
case $q=o(f)$, we have $q\sim z^2/(4\log f)$ and $|r_z|\asymp |z|/(\sqrt f\log f)$. In particular,
we have
\[
 \frac{f}{q\log f}=O(|\xi|^{-2}),\qquad
 \frac{f|r_z|}{q}=O(|\xi|^{-1}),\qquad
 \frac{f(1+yr_z^2)}{q^2}=O\left(\frac{(\log f)^2}{f\xi^4}+|\xi|^{-2}\right).
\]
By the bounds on the two gradient pairings, we get the same limit for the mixed
terms in \eqref{e:sign}.

It remains to determine the sign uniformly up to the maximal orbit. We fix
$L<\infty$ and consider $|\xi|\le L$ with $\eta\ne0$. For a pair at equal
$y$, by \eqref{e:neutral} and the strict concavity, the other point stays in a fixed central
interval, which is possibly larger. The constants in the following estimates
may depend on that interval. Since $y/f\to1$ uniformly and
$\overline{|\nabla r|^2}=o(1)$, we get from the first identity of \eqref{e:sign} that
\[
 f\,\frac{\mathcal L_u V}{V}=-\tfrac18+o(1).
\]
We claim that the errors in \eqref{e:centralND} are uniform even as $\eta\to0$.
In fact, we have $q=f\eta^2/\log f$ exactly, and
$y_z=-2f\eta\eta_z/\log f$. Since
$\sqrt{f}\,\eta_z=1/2+o(1)$ uniformly through $\eta=0$, we get
$(\log f)y_z^2/q=4f\eta_z^2=1+o(1)$.
The other term $|\nabla f|^2((U^2/2)'-Zy)^2$ in $|\nabla q|^2$ is
$O(\eta^2/(\log f)^2)$, and it contributes $O(f^{-1})$ to
$(\log f)\overline{|\nabla q|^2}/q$. By \eqref{e:centralND} and $V_f/V=O(f^{-1})$, we have
\[
 \frac{f}{q}\left|2\frac{V_f}{V}\overline{\langle\nabla f,\nabla q\rangle}
        -\frac4y\overline{\langle\nabla y,\nabla q\rangle}\right|
 \le C_L\left(\frac{\log f}{f}+\frac1{\log f}\right).
\]
Inserting these estimates into \eqref{e:sign} and using $f/q=(\log f)/\eta^2$, we obtain
\[
 \left|f\frac{\mathcal L_u H_u}{H_u}+\frac18+\frac{5}{16\eta^2}\right|
 \le o_L(1)(1+\eta^{-2}).
\]
Note that the coefficient of $\eta^{-2}$ is $-(5/4)/2+5/16=-5/16$.
If we take $f$ so large that the error coefficient is at most $1/16$, then we get
\[
 f\frac{\mathcal L_u H_u}{H_u}\le-\frac1{16}-\frac1{4\eta^2}
 \qquad(|\xi|\le L,\ \eta\ne0).
\]
The lower bound on $f$ depends only on $L$. We use the operator on the
open sheets $q>0$, and the limiting boundary value at $q=0$ will be verified below.

Now let us take any sequence of test points escaping to $f=\infty$. We first pass to a
subsequence on which $x$ is bounded or tends to infinity. The case of bounded $x$ is
covered by \eqref{e:captest}. In the second case, we take $r/\sqrt{2f}$ to be convergent.
If the limit is zero, we are in the first outer regime. If the limit is positive, we take a
further subsequence with $|\xi|\to\infty$ or with $|\xi|\le L$. The former is
the second outer regime, and the latter is covered by the central estimate. Hence, there
exists $S_0$ such that $\mathcal L_b H_b<0$ wherever $f\ge S_0$ and
$0\le r<U(f)$.

Next we fix $0<\nu<1$. We show that, after increasing $S_0$, the function $e^{-\nu(f-S_0)}H_b$ is
also a strict supersolution. For
$x\to\infty$, when we multiply by this factor, the quantity
\[
 \nu^2\overline{|\nabla f|^2}
 -\nu\left(C_f+2\overline{|\nabla f|^2}\frac{V_f}{V}
 -\frac{4\overline{h^{fy}}}{y}
 +\frac{5}{2q}\overline{\langle\nabla f,\nabla q\rangle}\right)
\]
is added to \eqref{e:sign}. The expression in the parentheses tends to one. Here $q\to\infty$ in
the outer region,
and we have $\overline{\langle\nabla f,\nabla q\rangle}/q=O((\log f)/f)$ near the fold.
On bounded caps, we differentiate \eqref{e:bq}. The
$\partial_f$ coefficient of its linearization is
$\overline{|\nabla f|^2}+4\overline{h^{fr}}/r+O((\log f)/f)$,
since in particular $r^2(b_{rr}+4b_r/r)=O((\log f)/f)$.
We also have $\partial_f H_b/H_b=O(f^{-1})$ and
$\overline{h^{fr}}\,\partial_r H_b/H_b=O((\log f)/f)$. Since $|\nabla f|^2\to1$,
the added term tends to $\nu^2-\nu<0$. Suppose that one of the strict
inequalities failed arbitrarily far out. Then, after passing to a subsequence, a sequence of
failing points
would belong to one of the regimes
treated above. This proves the uniform height assertion for both tests.
\end{proof}

\subsubsection{Exponential decay of the equal-height profile difference}
\begin{proposition}\label{prop:ax-exp}
For each $0<\nu<1$, the equal-height profile difference satisfies
\begin{equation}
 \sup_{0\le r\le U(f)}|b_+(f,r)-b_-(f,r)|\le C_\nu e^{-\nu f}.
 \label{e:exp}
\end{equation}
on all sufficiently high level sets.
\end{proposition}
\begin{proof}
We first consider the boundary at infinity. By Lemma~\ref{lem:mode-levelR}, we have
$R\le C(\log f)/f$ on every sufficiently high level set. Thus, by
\eqref{e:b}, we get $|b_+-b_-|\le C(\log f)/f$. For $|\eta|\ge1$, we have $q/f\ge1/\log f$ and
$x\le C\sqrt{\log f}$. Hence, we have $H_b\ge c f^{-7/8}(\log f)^{-17/4}$
and $|b_+-b_-|/H_b\to0$. Here we used $J_0(x)\ge c(1+x)^{-6}$, which
follows from its positivity on compact intervals and its prescribed tail.
For $0<|\eta|\le1$, by
\eqref{e:neutral} and the strict concavity, both paired points lie in a fixed central interval.
Their separation is $O(\sqrt f\,|\eta|)$. Note that the values of $b$ agree
at $r=U$. Integrating the difference of \eqref{e:b} from $r$ up to $U$, we have the exact identity
\[
 r^3(b_+-b_-)(f,r)
  =-\tfrac12\int_r^{U(f)}v^2\big(R_+(f,v)-R_-(f,v)\big)\,dv .
\]
For $v\ge r$, the two points are closer to the maximum than the
original pair. By \eqref{e:der}, the scalar difference at any integration radius is then at most
$C f^{-3/2}\sqrt f\,|\eta|$. We also have $U-r\le Cq/\sqrt{f}$. Multiplying by $r^2$ in the
preceding integral, we are left with a factor of
at most $CU^2/r=O(\sqrt f)$ in front of the integration length.
Using $q=f\eta^2/\log f$ in this estimate and $J_0=x^{-6}$ here,
we obtain
\[
 |u_+-u_-|\le C|\eta|^3/\log f,\qquad H_u\asymp f^{1/8}(\log f)^{-17/4}|\eta|^{5/2}.
\]
In particular, we have
$$|u_+-u_-|/H_u\le Cf^{-1/8}(\log f)^{13/4}|\eta|^{1/2}.$$ 
Hence,
the ratio tends to zero uniformly for $0<|\eta|\le1$. We now fix $S_0>1$
so large that the supersolution estimates hold. On each compact slab
$S_0\le f\le S$, the same fold argument, where we use the nondegenerate maximum,
gives $|b_+-b_-|=O(q^{3/2})$. By \eqref{e:bar}, the factor
$H_b/q^{5/4}=f^{-17/8}J_0(x)$ is bounded from above and bounded away from zero
on that slab. Hence, we have $|b_+-b_-|/H_b=O(q^{1/4})$, and its extension to
the fold has value zero.

We choose a constant $C$ which bounds $|b_+-b_-|/H_b$ at $f=S_0$. Then we compare $\pm(b_+-b_-)$
with $(Ce^{-\nu(f-S_0)}+\epsilon)H_b$ on a
truncated slab. By the preceding uniform limit, the top boundary of this slab is controlled for all
sufficiently large truncations. At the fold, we may truncate at $q=\delta$ and then
let $\delta\to0$ on
the fixed slab. In fact, by the signed square-root coordinate below,
the separation of the paired points is bounded by $C\sqrt q$, and we have $U-r=2q/(U+r)=O(q)$
there. Moreover, all constants can be chosen uniformly on that compact positive slab. At the pole,
we apply the five-dimensional lift of Section~\ref{sec:eq}, which has no
boundary there. In a regular chart, we divide $\pm(b_+-b_-)$ by the positive
comparison function. Since the comparison function
is a strict supersolution, we get a quotient with the same principal matrix and
a strictly negative zeroth-order coefficient. At a positive interior maximum, the gradient
vanishes and the second-order term is nonpositive, which is a contradiction.
Therefore, the quotient is at most one everywhere. Finally, we let the upper truncation tend to
infinity and then let $\epsilon$ tend to zero. Since $H_b$ is bounded on the high level sets, this proves (\ref{e:exp}).
\end{proof}

\subsection{A lower bound on a sequence of slabs}\label{sec:lower}
Suppose that we fix the radius and raise the height from $f_1$ to $f_2$.
Then $b$ changes by $-\int_{f_1}^{f_2}K\,df$. Therefore, if $K$ is very small,
a substantial height difference could be nearly invisible
in \eqref{e:exp}. For our purpose we only need a lower bound on selected
unit-height slabs, and at height $S$ we can afford any loss $e^{-o(S)}$.

We first give a curvature lower bound at one pole.
\begin{lemma}\label{lem:ax-seed}
We fix one of the two halves of the axis. Then there is a sequence $S_j\to\infty$
such that the smallest eigenvalue $k$ of the ambient curvature operator satisfies
$k\ge cS_j^{-2}$ on the ball of radius $S_j^{-3}$ centered at the pole
of height $S_j$ on that half-axis.
\end{lemma}
\begin{proof}
We denote by $l$ the outgoing unit arclength on the fixed axis.
Note that at a pole the isotropy representation on $T_pM$ is the direct sum
of the axial line and the standard $\mathbb R^3$. We claim that every two-plane
spanned by the axial tangent and a unit vector in $\mathbb R^3$
has curvature $K$. To see this, we rotate its transverse vector into the totally
geodesic meridian. On two-forms, the two summands are the mixed
three-space and $\Lambda^2\mathbb R^3$. Note that the central element $-I$
of the angular group acts on these summands with opposite signs, and so
the curvature operator has no cross term. Moreover, on each summand all
unit directions are equivalent under rotations. Then by polarization we see that the quadratic
form on each summand is scalar. We denote by $K_2$ the eigenvalue
on $\Lambda^2\mathbb R^3$. Now we count the six planes and get
$R_{\rm tip}=6(K+K_2)$ and $\Ric_{ll}=3K$.

Along this half-axis we have $df/dl>0$ and
$R_{\rm tip}=1-(df/dl)^2$. We differentiate the second identity and divide by $df/dl$, and then
we get $dR_{\rm tip}/df=-2\,d^2f/dl^2=-6K$.
On the other hand, by Proposition~\ref{prop:mode-regions}, we have
$R_{\rm tip}(S)-R_{\rm tip}(2S)=(1/2+o(1))\log S/S$.
So if we integrate the derivative identity from $S$ to $2S$ and divide by
the length of the interval, then we can find some $S_j\in[S,2S]$ such that
$K\ge c\log S/S^2$ there. We let $S$ run through the powers of four.
At this pole, we have $K\le CR_{\rm tip}^2=o(R_{\rm tip})$ by \eqref{e:K},
and hence $K_2=R_{\rm tip}/6-K\ge cR_{\rm tip}$. It follows that both eigenvalues
are at least $cS_j^{-2}$. We also know from \eqref{e:der} and $R\le1$ that
the covariant derivative of the ambient curvature operator is bounded
globally. So under parallel transport from the pole, its quadratic
form on a unit two-form changes by at most $CS_j^{-3}$ everywhere in the
ball of the statement. Then the lemma follows after we decrease $c$.
\end{proof}

Then, we have the lower bound for $K$ on the unit-height slab.
\begin{lemma}\label{lem:slabs}
There exist sequences $S_j\to\infty$ and $a_j=o(S_j)$ such that
$K\ge e^{-a_j}$ holds everywhere on $S_j-1\le f\le S_j$.
\end{lemma}
\begin{proof}
We begin with the differential inequality which will be used to propagate the lower
bound for $k$. We apply Hamilton's evolution equation for the curvature operator \cite{Ham}
to the steady flow in a moving orthonormal frame. This gives
\[
 (\Delta+\nabla_{\nabla f})\mathcal R
              =-2(\mathcal R^2+\mathcal R^\#).
\]
Note that the right-hand side is
nonpositive as a quadratic form. Under the identification of two-forms
with skew endomorphisms, in an orthonormal eigenbasis with
eigenvalues $\lambda_i>0$ the sharp term evaluated on a form $w$
is a sum of nonnegative multiples of
$\lambda_i\lambda_j\langle[e_i,e_j],w\rangle^2$.

The function $k$ can have repeated eigenvalues, but it still satisfies
$(\Delta+\nabla_{\nabla f})k\le0$ when tested from below.
To see this, we extend a minimizing unit eigenform at the test point by radial
parallel transport in normal coordinates. Then we see that the first covariant
derivative and the rough Laplacian of this extension vanish at the center,
because both covariant derivatives vanish along every
radial geodesic, and this gives each
diagonal term in the trace. Note that the Rayleigh quotient of this
extension is smooth, it is at least $k$ everywhere near the center, and it is equal to $k$
at the center. So any smooth lower test for $k$ also touches this
Rayleigh quotient from below at the center. It follows that the drift Laplacian of the test
is no greater than that of the quotient, and the latter is equal to the curvature
equation evaluated on the chosen form at the center.

Next we choose $0<\delta<1/8$ such that every ball of radius $4\delta$ is
within the injectivity radius of its center. We can make such a uniform choice
by Cheeger--Gromov--Taylor \cite{CGT}, since we have the curvature bound and we also have,
by noncollapsing, uniform lower volume bounds
for balls of one sufficiently small fixed radius. Let $d$ be the distance to a
center. By Ricci comparison, we have $\Delta d\le3/d$. Since
$|\nabla f|\le1$ and the derivative of $d^{-8}$ is negative, we get
\[
 (\Delta+\nabla_{\nabla f})d^{-8}
 \ge 72d^{-10}-8d^{-9}(3/d+1)
 =8 d^{-10}(6-d)>0 \qquad(0<d\le4\delta).
\]
Now suppose that $k\ge h>0$ on $B(p,\epsilon)$ for some $0<\epsilon\le\delta$.
On $\epsilon\le d\le4\delta$, we compare $k$ with the function
\[
 h\,\frac{d^{-8}-(4\delta)^{-8}}{\epsilon^{-8}-(4\delta)^{-8}} .
\]
This function is equal to $h$ on the inner boundary and to zero on the outer boundary.
Suppose that $k$ fell below it somewhere. Then we translate the function downward by
the size of the most negative difference. The translated function would touch $k$
from below at an interior point, and its drift Laplacian is strictly positive there.
This contradicts the test inequality. At $d=2\delta$, the numerator of the comparison
function is a fixed positive number and its denominator is at most $\epsilon^{-8}$.
Thus we have $k\ge c_\delta h\epsilon^8$ on $B(p,2\delta)$. After decreasing
$c_\delta$, this also holds on the inner ball.

We first apply this with $\epsilon=S_j^{-3}$ and $h=cS_j^{-2}$ which are given by
Lemma~\ref{lem:ax-seed}. Then the enlarged ball has the lower bound
$c'S_j^{-26}$. Next we move its center along a path by at most $\delta$, and we see that the
new $\delta$ ball is contained in the old $2\delta$ ball. So in every subsequent
application, which is done with $\epsilon=\delta$, we lose just one fixed
factor in $(0,1)$. Note that any point on $f=S_j$ can be joined to the chosen pole
by a constant-angle meridian whose length is at most
$D(S_j)=O(\sqrt{S_j\log S_j})$. Moreover, any point on the slab can first be moved
to that level set by the flow of $Z$, and the distance travelled is at most a fixed constant,
since $Zf=1$ and $|Z|=1/|\nabla f|$ is uniformly bounded there.
So we need at most $C(1+\sqrt{S_j\log S_j})$ moves for each point.
It follows that
\[
 K\ge k\ge c S_j^{-26}\exp(-C\sqrt{S_j\log S_j}) .
\]
Here we get the first inequality when we evaluate the curvature operator on the unit
meridian two-form. Note that both $\log S_j$ and $\sqrt{S_j\log S_j}$ are $o(S_j)$.
Then we take negative logarithms and the lemma follows.
\end{proof}

\subsection{Inverse heights and the global involution}\label{sec:global}
In this subsection, we use the fact that \eqref{e:bq} holds on a sheet. Recall that
at equal $(f,r)$ the zeroth-order coefficient of the linearized equation did not have
a fixed sign. On the other hand, at equal $(r,b)$, the inverse equation is preserved
when we add a constant to the unknown height. So we can compare a height difference
by the strong maximum principle.

\begin{lemma}\label{lem:ax-inverse}
On a sheet, we denote by $t_\pm(r,\beta)$ the height at which $b_\pm(f,r)=\beta$.
We consider the common domain where both heights are at most $S$. If the maximum of
$|t_+-t_-|$ on this domain is strictly positive, then it can only be attained at a point
where one of the heights is equal to $S$.
\end{lemma}
\begin{proof}
\emph{The domain and its boundary.}
We know from Proposition~\ref{prop:level-radius} that the function $U$ is strictly increasing from
zero to infinity. For fixed $r>0$, the two sheets first meet at height
$U^{-1}(r)$, and after that they extend to every larger height. The starting values
of $b$ on the two sheets agree there. Note that $b_f=-K<0$, and so each profile is strictly
decreasing afterwards. Then we see that the range which is attained on both sheets up to height
$S$ is
\[
 \begin{gathered}
 0\le r\le U(S),\qquad
 \max\{b_+(S,r),b_-(S,r)\}\le\beta\le
                       b_+(U^{-1}(r),r).
 \end{gathered}
\]
At $r=0$ we take $U^{-1}(0)=0$. Then the right endpoint of the $\beta$
interval has the value $1/6$. This is because radii and heights approaching the
vertex stay in arbitrarily small sublevel sets of the proper function $f$,
and the integral in \eqref{e:b} tends to $R(o)/6=1/6$.
We also see that for fixed $0\le r<U(S)$ the interval has positive length, since the profiles are
strictly decreasing on both sheets. Note that the upper boundary is the fold,
where both inverse heights are equal to $U^{-1}(r)$, and the lower boundary
is the graph in the left inequality, on which at least one of the heights
is equal to $S$. These two boundaries meet at $r=U(S)$.

We note that the direct profiles are continuous up to the boundaries.
Now we consider a convergent sequence of inverse arguments. The corresponding
heights lie in the compact interval $[0,S]$. Then by continuity, we know that any subsequential
limit of the heights must solve the limiting direct equation,
and this solution is unique because of the strict monotonicity on the closed height
interval. At $(0,1/6)$, if the limiting height were positive, then we would
get $R_{\rm tip}=1$ instead, which contradicts $R=1-|\nabla f|^2<1$
on positive levels. It follows that both inverses are continuous on the whole of this
compact set.

\emph{Ellipticity of the inverse equation.}
We use the radial lift $r=|v|$ of Proposition~\ref{prop:ax-b}. At interior
points, including $v=0$ at regular positive height, we have $b_f<0$.
Then it follows from the implicit function theorem that the inverse is smooth there.
On either sheet, the inverse coordinate change is
$T(\beta,v)=(t(\beta,v),v)$, its Jacobian $J=DT$ is invertible,
and we have $b(T(\beta,v))=\beta$. We take first derivatives of this relation and get
$b_f=1/t_\beta$ and $D_v b|_f=-D_v t/t_\beta$.
If we differentiate the composition once more, then we have the matrix identity
\[
       J^{\mathsf T}D^2b\,J=-b_f D^2t .
\]
Here we have used that the second derivative of the right-hand coordinate $\beta$
is zero and that $t$ is the only nonlinear component of $T$.
So if we denote by $A$ the lifted principal matrix for $b$, then the principal matrix
in the inverse equation is
$-t_\beta^{-1}J^{-1} A J^{-\mathsf T}$. This matrix is positive
definite since $t_\beta<0$. Note that the substitution in all the lower-order terms involves
only $v,\beta,Dt$, since \eqref{e:bq} contains no explicit $f$.
In particular, the inverse equation does not depend on $t$ itself.

\emph{The maximum principle at a contact point.}
Suppose that a signed difference $t_+-t_-$ or its negative attains a
positive maximum in the lifted interior. Then the two gradients agree at that
point, and the inverse principal matrix is positive definite at their common gradient.
Note that this matrix depends on the independent variables and
on the gradient, but not on the Hessian. So on a sufficiently
small ball, we see that it remains uniformly positive for all linear interpolations
between the two gradients. We now use the mean value formula, that is, we differentiate
the inverse equation along the segment between the derivatives of the two solutions.
In this way we get a linear homogeneous equation for their difference
on this ball. The principal matrix of this equation is the average just described,
the first-order coefficients are locally bounded, and the
zeroth-order coefficient is zero. By the strong maximum principle,
the difference is therefore constant on the ball.

It follows that the set where this maximum is attained is open in the interior,
and by continuity it is also relatively closed. We claim that the interior is connected.
In fact, in the lifted ball $|v|<U(S)$, both boundary graphs are continuous
and there is a positive vertical gap between them. So we can move vertically to their
midpoint, follow the continuous midpoint graph along a path in the ball, and then move
vertically to the destination. Hence, if there were an interior maximum, the difference
would be equal to the same nonzero constant throughout. This contradicts the continuity
at the upper boundary, where the inverse heights agree. Since the lift
also takes care of $r=0$, a positive maximum can only be attained on the lower boundary.
\end{proof}

\begin{figure}[htbp]
\centering
\IfFileExists{figures/ax-heights.pdf}{%
  \includegraphics{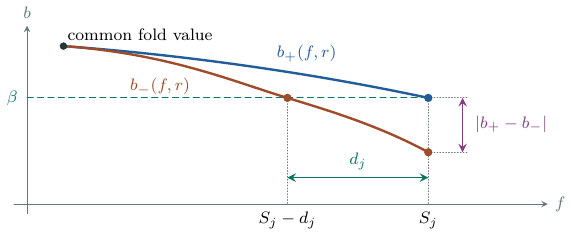}%
}{%
\begin{tikzpicture}[x=.88cm,y=.86cm,>=stealth,line cap=round,line join=round,font=\footnotesize]
  \draw[->,wingGrey] (.55,.22)--(10.82,.22) node[right] {$f$};
  \draw[->,wingGrey] (.80,.03)--(.80,3.73) node[above] {$b$};
  \draw[wingMirror,line width=1.15pt] (1.50,3.33) .. controls
       (3.78,3.16) and (6.23,2.81) .. (8.52,2.31);
  \draw[wingAxis,line width=1.15pt] (1.50,3.33) .. controls
       (3.53,3.1786) and (4.77,2.63) .. (5.81,2.31)
       .. controls (6.64,2.05) and (7.58,1.76) .. (8.52,1.24);
  \draw[wingPath,densely dashed] (.80,2.31)--(8.51,2.31);
  \draw[wingGrey,densely dotted] (5.81,2.31)--(5.81,.22)
                                  (8.52,2.31)--(8.52,.22);
  \draw[<->,wingPath,thick] (5.81,.74)--(8.52,.74)
      node[midway,above] {$d_j$};
  \draw[wingGrey,densely dotted] (8.52,1.24)--(9.30,1.24)
                                  (8.52,2.31)--(9.30,2.31);
  \draw[<->,wingOrbit,thick] (9.18,1.24)--(9.18,2.31);
  \node[anchor=west,wingOrbit] at (9.28,1.79) {$|b_+-b_-|$};
  \node[anchor=east,wingPath] at (.76,2.31) {$\beta$};
  \node[anchor=south,wingMirror] at (6.18,2.87) {$b_+(f,r)$};
  \node[anchor=north,wingAxis] at (3.36,2.84) {$b_-(f,r)$};
  \node[below] at (5.81,.20) {$S_j-d_j$};
  \node[below] at (8.52,.20) {$S_j$};
  \node[anchor=south west,inner sep=1pt] at (1.53,3.40) {common fold value};
  \fill[wingInk] (1.50,3.33) circle(1.8pt);
  \fill[wingMirror] (8.52,2.31) circle(2pt);
  \fill[wingAxis] (5.81,2.31) circle(2pt) (8.52,1.24) circle(2pt);
\end{tikzpicture}%
}
\caption{A hypothetical difference at a fixed $r$. The horizontal arrow
is the difference $d_j$ of the two heights at the same $b=\beta$, and
the vertical arrow is the difference of the two values at the same
$f=S_j$. Since $-b_f=K$ along the lower curve, the vertical difference is
equal to the integral of $K$ over the horizontal interval. This
configuration is excluded by Proposition~\ref{prop:ax-equal}. One can
also reverse the roles of the two sheets.}
\label{fig:ax-heights}
\end{figure}

As a consequence, we get the exact equality at the same radius and height on two sheets.
\begin{proposition}\label{prop:ax-equal}
For every $f>0$ and every $0\le r\le U(f)$, we have
$b_+(f,r)=b_-(f,r)$.
\end{proposition}
\begin{proof}
We note that the common truncated inverse domains increase with $S$, since for
every fixed radius both lower boundary values of $b$ decrease.
Suppose that their maximum difference were $d_0>0$ at some height cutoff $S_1$.
Then, at any larger $S_j$ of Lemma~\ref{lem:slabs}, the maximum $d_j$
would be at least $d_0$. By Lemma~\ref{lem:ax-inverse}, we can find
a pair at some common $(r,\beta)$ whose heights are $S_j$ and $S_j-d_j$.
On the sheet with the smaller height, the fixed-radius profile continues
from $S_j-d_j$ up to $S_j$, and its initial value is the value of the other
sheet at $S_j$. Hence, by the fundamental theorem of calculus, we get
the vertical difference shown in Figure~\ref{fig:ax-heights}, namely
\[
 |b_+(S_j,r)-b_-(S_j,r)|
    =\int_{S_j-d_j}^{S_j}K(v,r)\,dv
    \ge \min\{d_0,1\}\,e^{-a_j}.
\]
Here the integral is taken on the sheet with the smaller height. Note that its last
$\min\{d_0,1\}$ units lie in $[S_j-1,S_j]$, where the lower bound
holds for every radius. On the other hand, at the same point $(S_j,r)$,
we get $C e^{-S_j/2}$ instead from Proposition~\ref{prop:ax-exp}.
So we have $\min\{d_0,1\}\exp(S_j/2-a_j)\le C$, which is impossible since
$a_j=o(S_j)$. Hence, we conclude that all inverse discrepancies vanish on every finite domain.

It remains to recover the direct values. We fix $S>0$ and $0\le r\le U(S)$.
At the lower boundary $\beta=\max\{b_+(S,r),b_-(S,r)\}$, one of the inverse
heights is $S$ by construction. Then by what we have just shown, both of them are equal
to $S$, and it follows that both direct values at that height are equal to this same $\beta$.
This proves the proposition, including its two endpoints.
\end{proof}

\begin{figure}[htbp]
\centering
\IfFileExists{figures/ax-reflection.pdf}{%
  \includegraphics{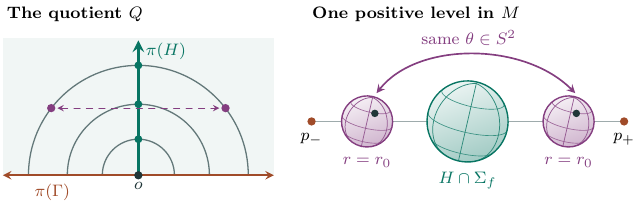}%
}{%
\begin{tikzpicture}[x=.94cm,y=.94cm,>=stealth,line cap=round,line join=round,font=\footnotesize]
  \node[anchor=west] at (.0,2.78) {\textbf{The quotient $Q$}};
  \node[anchor=west] at (5.50,2.78) {\textbf{One positive level in $M$}};
  \begin{scope}[shift={(2.50,-.12)}]
    \fill[wingPale] (-2.44,0) rectangle (2.44,2.48);
    \foreach \R in {.65,1.28,1.98}{
       \draw[wingGrey,line width=.62pt] (\R,0)
          arc[start angle=0,end angle=180,radius=\R];}
    \draw[wingAxis,thick,<->] (-2.44,0)--(2.44,0);
    \draw[wingPath,line width=1.45pt,->] (0,0)--(0,2.43);
    \foreach \y in {.65,1.28,1.98}{\fill[wingPath] (0,\y) circle(1.9pt);}
    \draw[wingOrbit,dashed,<->] (-1.46,1.21)--(1.46,1.21);
    \fill[wingOrbit] (-1.57,1.21) circle(2pt) (1.57,1.21) circle(2pt);
    \node[wingPath,anchor=west,inner sep=1pt] at (.09,2.23) {$\pi(H)$};
    \node[wingAxis,below] at (-1.55,-.03) {$\pi(\Gamma)$};
    \fill[wingInk] (0,0) circle(2pt) node[below] {$o$};
  \end{scope}
  \begin{scope}[shift={(5.62,.85)}]
    \draw[wingGrey!62] (0,0)--(5.63,0);
    \foreach \p in {0,5.63}{\fill[wingAxis] (\p,0) circle(2pt);}
    \WingOrbitSphere{(1.00,0)}{.46}{wingOrbit}
    \WingOrbitSphere{(2.81,0)}{.73}{wingPath}
    \WingOrbitSphere{(4.63,0)}{.46}{wingOrbit}
    \fill[wingInk] (1.14,.145) circle(1.8pt) (4.77,.145) circle(1.8pt);
    \draw[wingOrbit,<->,thick] (1.17,.51) .. controls (1.91,1.44) and
                        (3.82,1.44) .. (4.75,.51);
    \node[wingOrbit,above] at (2.83,1.24) {same $\theta\in S^2$};
    \node[below] at (0,-.09) {$p_-$};
    \node[below] at (5.63,-.09) {$p_+$};
    \node[wingPath,below] at (2.81,-.79) {$H\cap\Sigma_f$};
    \node[wingOrbit,below] at (1.00,-.52) {$r=r_0$};
    \node[wingOrbit,below] at (4.63,-.52) {$r=r_0$};
  \end{scope}
\end{tikzpicture}%
}
\caption{The reflection, after the equality of the profiles has been proved.
On the left, the arcs are levels in $Q$ and every interior dot
is a whole $S^2$. The purple pair is exchanged by the reflection. The
green curve consists of the maximal orbits and $o$, and the
rust boundary is the $O(3)$-fixed axis. On the right, we show a single
height, where the black dots on the purple spheres are matched angular
points. We see that the reflection fixes every point of the green sphere and exchanges
the poles. The three-dimensional reflection hypersurface $H$ consists of
the green spheres over all positive heights and $o$.}
\label{fig:ax-reflection}
\end{figure}

Now, we are able to prove the existence of the global involution.
\begin{proposition}\label{prop:ax-isom}
If we exchange the two sheets at fixed $(f,r)$ and match the same
angular point, we obtain an isometric involution $\sigma$ on $M$. On every
positive level set, $\sigma$ fixes the maximal sphere pointwise and exchanges
the two poles. Moreover, $\sigma$ commutes with the angular $O(3)$.
\end{proposition}
\begin{proof}
\emph{The open sheets.} We choose a unit vector in the transverse
$\mathbb R^3$ in Fang's global equivariant model. On its open half
meridian, each orbit has one representative, with the same
stabilizer $O(2)$. Given a direction $\theta\in S^2$, we apply any
$A\in O(3)$ which takes the chosen vector to $\theta$. Note that two different
choices differ by that stabilizer, and so the angular coordinate which we get in
this way is well-defined. As in Proposition~\ref{prop:level-profile}, it
gives the same identification on both sheets and at every height. We define
\[
 \sigma(f,r,\theta,-)=(f,r,\theta,+),\qquad
 \sigma(f,r,\theta,+)=(f,r,\theta,-).
\]
These maps between the open sheets are smooth and they are inverse to each
other. We know from Proposition~\ref{prop:ax-equal} that the functions $b$ on the two
sheets are equal on the open $(f,r)$ domain, and hence the same holds for
$b_f,b_r$. Then all entries of $h^{-1}$ in \eqref{e:metric} agree, and so the
inverse matrices $h$ agree as well. Moreover, both angular blocks are equal
to $r^2g_{S^2}$, and the horizontal directions are orthogonal to the
angular directions. It follows that $\sigma$ preserves the full metric
$h+r^2 g_{S^2}$ there. In particular, we see that the equality of radii by itself was
only the angular part of the required isometry.

\emph{Crossing a maximal orbit.} At a fixed positive height, we use a
smooth local arclength $z$ on the quotient. Note that the defining function
$U(f)^2-r(z,f)^2$ and its $z$-derivative both vanish at
$z=z_m(f)$, and its second derivative there is $-2U r_{zz}>0$.
By Taylor's formula along each level set, we can write this function as
$(z-z_m(f))^2$ times a smooth positive function of $(z,f)$.
Hence, the signed square root
\[
 s=\operatorname{sgn}(z-z_m(f))\sqrt{U(f)^2-r^2}
\]
is smooth and has nonzero $z$-derivative through the maximum. Since the
two paired points have the same $s^2$ and opposite $s$, the extension is
given by $(f,s,\theta)\mapsto(f,-s,\theta)$. This map is smooth and
fixes the whole sphere at $s=0$ pointwise. By continuity, it also preserves
the metric there. Note that this signed coordinate is local to the fold.
In Figure~\ref{fig:ax-reflection}, we show the resulting map on a whole
level set and on its sphere orbits.

\emph{The axis and the vertex.} So far we have obtained a bijective
Riemannian isometry of the regular-orbit set $M\setminus\Gamma$ with its
intrinsic length metric. We claim that for two points there, this distance
agrees with $d_g$. To see this, we choose a smooth path between the two
points whose length is arbitrarily close to the ambient distance, or we
take a minimizing geodesic. By relative transversality, this path admits
arbitrarily small $C^1$ perturbations which are fixed near the endpoints
and transverse to the closed embedded axis. Since the domain of such a path
has dimension one and the axis has codimension three, such a path does not
meet the axis. Moreover, by $C^1$ continuity, its length changes arbitrarily
little. The reverse distance inequality holds because every path in the
complement is also a path in $M$. Note that the complement is dense and
$M$ is complete. Therefore, by applying the isometry to any approximating
Cauchy sequence, we get a unique distance-preserving extension to $M$.
The inverse extends in the same way, so the extension is onto and
involutive. By the Myers--Steenrod theorem, a distance-preserving bijection
of smooth Riemannian manifolds is a smooth metric isometry. 

When $r\to0$ at fixed positive $f$, a sheet converges to its
own pole, and hence the two poles are exchanged. The vertex $o$ is fixed,
since it is the unique minimum of $f$. On regular orbits, we see that the
formulas preserve $f$ and commute with $A\theta$ for every $A\in O(3)$.
Then by continuity, these identities also hold on the axis. 
\end{proof}

Let $l$ be the
arclength coordinate from the minus pole along a whole level. Note that the induced
isometry of the interval exchanges the two endpoints and preserves the
distance to each of them, and so it sends $l$ to $D(f)-l$. It follows that its
unique fixed orbit is at $l=D(f)/2$. Note that this is a consequence of
the proposition above, which is not used
for the estimates in Section~\ref{sec:modes}.

To verify the $O(3)\times\mathbb{Z}_{2}$ action on the whole steady gradient Ricci soliton, we recall the definition  of such an action from \cite{Lai}. 
\begin{defi}[The Euclidean $O(3)\times\mathbb{Z}_{2}$ action]
We write $(x',x_{4})\in\mathbb{R}^{3}\times\mathbb{R}$, where $x'=(x_{1},x_{2},x_{3})$. For $A\in O(3)$ and $\varepsilon\in\mathbb{Z}_{2}=\{+1,-1\}$, the action is defined by
\[
(A,\varepsilon)\cdot(x',x_{4})=(Ax',\varepsilon x_{4}).
\]
Note that the $O(3)$ factor fixes the axis $\{x'=0\}$ pointwise, and for $x'\neq0$ its orbit through $(x',x_{4})$ is $\lvert x'\rvert S^{2}\times\{x_{4}\}$. We also see that the reflection $(I,-1)$ fixes $\{x_{4}=0\}$ pointwise and exchanges the two ends of the axis.
\end{defi}

\begin{defi}[The isometric $O(3)\times\mathbb{Z}_{2}$ action on manifold]
We say that a Riemannian manifold $(M^{4},g)$ is $O(3)\times\mathbb{Z}_{2}$-symmetric if it admits a smooth isometric action of this group and a diffeomorphism $\Phi:M\to\mathbb{R}^{3}\times\mathbb{R}$ such that for every $p$ and $(A,\varepsilon)\in O(3)\times\mathbb{Z}_{2}$, we have
\[
\Phi((A,\varepsilon)\cdot p)=(Ax',\varepsilon x_{4})
\qquad\text{when }\Phi(p)=(x',x_{4}).
\]
For a soliton $(M,g,f)$, we say that the action preserves the potential if $f((A,\varepsilon)\cdot p)=f(p)$ holds for every group element and every $p$.
\end{defi}

Now, we are ready to prove Theorem~\ref{th:ref}.

\begin{proof}[Proof of Theorem~\ref{th:ref}]
We first consider the non-Bryant soliton and let $H=\operatorname{Fix}(\sigma)$.
This set consists of $o$ and
the maximal sphere on each $f>0$. Away from $o$, by the fold chart above, we see that
there are two angular directions and one height direction which are tangent
to this set, and the sign is reversed in the remaining direction. At $o$,
the derivative of $\sigma$ is $-1$ along $T_o\Gamma$, because $\sigma$
exchanges the two outgoing rays of a unit-speed geodesic. On the orthogonal
standard $\mathbb R^3$, the derivative commutes with $O(3)$, and hence it
is a scalar. In fact, the stabilizer of each unit vector fixes just the
span of this vector, and all the spans are related by rotations. Note that the
derivative is an orthogonal involution, and so its scalar is $+1$
or $-1$. Suppose that the latter holds. Then we would have $d\sigma_o=-I$,
and by the equivariance of the exponential map, $o$ would be isolated in
the fixed set. But this contradicts the fact that the fixed maximal spheres
approach $o$ as $f\to0$ by properness. So we know that $d\sigma_o$ fixes the
transverse three-plane. We also know that an isometry preserves the initial data of
each geodesic. Then in normal coordinates based at any fixed point, the fixed set
is the exponential image of the $+1$ eigenspace. This also shows that
the three-dimensional hypersurface which we get in this way is totally geodesic
and meets $\Gamma$ orthogonally at $o$.

Now we define the $O(3)\times\mathbb Z_2$ action by
$(A,+1)\mapsto A$ and $(A,-1)\mapsto A\circ\sigma$.
It is easy to see that this is a homomorphism, since $\sigma$ commutes with $O(3)$ and
$\sigma^2=\mathrm{id}$. We also note that $A$ fixes the two
poles of each positive level individually, and so every $(A,-1)$ exchanges them.
Moreover, the angular action is effective, and so among the
elements $(A,+1)$ only $(I,+1)$ acts trivially. Hence, we see that the product action
is effective. On the Bryant soliton, we use the radial $O(4)$ action and the
standard angular $O(3)$ which fixes one coordinate line. The reflection in
the three-plane orthogonal to this line belongs to $O(4)$, commutes with
that $O(3)$, and exchanges the poles. Then the same kernel argument applies.

It remains to verify the global equivariant diffeomorphism. In both
cases, the common fixed set of $G=O(3)\times\mathbb Z_2$ consists of
$o$ alone. This is because an $O(3)$-fixed point lies on $\Gamma$, and on
each positive level the reflection exchanges the two axis points. Note that $G$
is compact and acts by isometries. Then by Fang~\cite[Proposition~7.5]{F}
we get a $G$-equivariant diffeomorphism $T_oM\to M$ for the
differential action. Let us check the hypotheses of that proposition. The
manifold $M$ is complete and connected, and positive curvature operator
implies strictly positive sectional curvature. Moreover, compactness would
contradict $\Delta f=R>0$ after integration. So the fixed point
given by that proposition must be $o$.
On $(T_o\Gamma)^\perp\oplus T_o\Gamma\cong\mathbb R^3\oplus\mathbb R$,
the angular group acts on the first summand by its standard representation
and fixes the second summand. As we have shown above, the reflection acts
by $+1$ and $-1$ on these two summands, respectively. On the Bryant branch,
we can see this in the radial model. Under this identification,
the tangent action is $(A,\varepsilon)(v,t)=(Av,\varepsilon t)$. Hence, we know that
the inverse of Fang's diffeomorphism gives the required global map
$M\to\mathbb R^3\times\mathbb R$.
\end{proof}

\section{The phase of a reflected end}\label{cl:phase}

By
Theorem~\ref{th:ref}, we know that the end admits a reflection. On the non-Bryant
soliton, we will use the neutral-mode and slope estimates of
Section~\ref{sec:modes} to define the phase which appears in
Theorem~\ref{cl:main}. Throughout this section, all orbit radii are area
radii.

\subsection{Reformulation of the level set equation}
Let $s=f$ and $t=U(s)^2/2$. By Theorem~\ref{th:ref}, the sphere where
the maximum is attained lies in the reflection hypersurface. We choose the
signed level arclength $z$ measured from that sphere. Since the normal gradient
field $Z=\nabla f/|\nabla f|^2$ is tangent to the hypersurface, the
integration constant in \eqref{e:beta} vanishes there. We set
\[
 y=r^2/2,\quad u=\log t,\quad \tau=-u,\quad
 \xi=z/\sqrt t,\quad w=y/t,\quad
 \nu=dt/ds,\quad \lambda=ds/dt .
\]
Thus $w$ is even and $w(0)=1$. In the following, primes always denote
$\xi$-derivatives. Let $j=-w''(0)$ and $k=w^{(4)}(0)$. At the equator, the two
intrinsic sectional curvatures are $-r_{zz}/r$ and $1/r^2$, and their ratio is
$j$. In particular, we have $j=1$ on round level sets and $j=0$ on product level sets.

We denote by $\gamma$ the induced level metric. When tensors on nearby
levels are pulled back along $Z$, we set
\[
 G=\tfrac12\gamma^{-1}\partial_s\gamma,\qquad
 T=\Ric_\gamma^\sharp,\qquad E=G-T .
\]
Thus, the tensor $E$ here is the endomorphism which we get from the
error tensor of Lemma~\ref{lem:level-error}. We use a
subscript $m$ or $o$ to denote its meridional or spherical eigenvalue. In the
formulas below, we set $e=r^2E_o$, $b=tE_m$, and $V=\sqrt t\,Zz$. By
Proposition~\ref{prop:level-evolution}, we have $Zy=1-y_{zz}+e$ and
$\partial_z(Zz)=-2r_{zz}/r+E_m$.
Applying the chain rule, first at fixed $z$ and then at fixed $\xi$, we get
\begin{align}
 \nu w_\tau&=w''+\nu(w-\xi w'/2)-1+Vw'-e,\label{cl:w}\\
 V'&=-w''/w+(w')^2/(2w^2)+b,\qquad V(0)=0,\label{cl:V}\\
 \nu&=1+j+e(0),\qquad
 \nu j_u=k-2j^2-2jb(0)-e''(0).\label{cl:clock}
\end{align}
Here the last two identities are obtained by evaluating \eqref{cl:w} and its
second derivative at zero. In particular, $t$ is an increasing parameter.
By Proposition~\ref{prop:mode-center}, after shifting the old arclength
origin by $o(\sqrt s)$ to its maximum and using $t/s\to1$, we get
$j\sim(2u)^{-1}$. Moreover, by Proposition~\ref{prop:mode-regions}, we have
$w\ge c_a>0$ on $|\xi|\le a\sqrt u$ for every $a<2$. On a band where
$w\ge c$, the estimate \eqref{e:err} gives a bound for every fixed spatial
derivative of $w$, and it bounds every such derivative of $e,b$ by
$C_{k,c}e^{-u}$. For the moment, these constants may depend on the soliton.

\subsection{Definition of phase.}
Later we apply the estimates after translating the argument \(u\) to define the phase. We now make
precise the strip and margin terminology used in the next lemma.
\begin{defi}[Strips and margins]\label{def:phase-strip}
We fix \(L>0\), \(d>0\), \(\rho>0\), and
\(\alpha\in\{1/3,1/2\}\). We put
\[
 R=\rho L^\alpha
\]
and we define the full strip by
\[
 \mathscr S(L,d,\rho,\alpha)
 =
 \left\{
 (u,\xi):
 L\le u\le(1+d)L,\quad |\xi|\le R
 \right\}.
\]
For \(0<\vartheta<1/2\), we also define the \(\vartheta\)-interior strip by
\[
 \mathscr S_\vartheta(L,d,\rho,\alpha)
 =
 \left\{
 (u,\xi):
 L+\vartheta dL\le u\le(1+d-\vartheta d)L,\quad
 |\xi|\le(1-\vartheta)R
 \right\}.
\]
Then we see that the distances from the two \(u\)-boundaries of
\(\mathscr S_\vartheta\) to the corresponding \(u\)-boundaries of
\(\mathscr S\) are both \(\vartheta dL\), and the distance from each
lateral boundary is \(\vartheta R\).

More generally, if we have
\[
 \mathscr S'=[u_-',u_+']\times[-R',R']
 \subset
 \mathscr S=[u_-,u_+]\times[-R,R],
\]
then we define
\[
 \operatorname{mar}(\mathscr S';\mathscr S)
 =
 \min\left\{
 u_-'-u_-,
 u_+-u_+',
 R-R'
 \right\}.
\]
This is the smallest Euclidean distance from a boundary component of
\(\mathscr S'\) to the corresponding boundary component of
\(\mathscr S\). We say that the margin is fixed when the ratios
\[
 \frac{u_-'-u_-}{L},
 \qquad
 \frac{u_+-u_+'}{L},
 \qquad
 \frac{R-R'}{R}
\]
are bounded below by fixed positive numbers. All constants below may
depend on these fixed ratios, but are independent of \(L\).

\end{defi}

\begin{lemma}\label{cl:strip}
Assume that on
\(\mathscr S(L,d,\rho,\alpha)\) one has
\[
 w\ge c>0,
 \qquad
 0\le j(u)\le \frac{C}{L},
\]
and, for every fixed integer \(m\ge0\),
\[
 \|\partial_\xi^m w\|_{L^\infty(\mathscr S)}\le C_m,
 \qquad
 \|\partial_\xi^m e\|_{L^\infty(\mathscr S)}
 +
 \|\partial_\xi^m b\|_{L^\infty(\mathscr S)}
 \le C_m e^{-u}.
\]
Set
\[
 P=w_\xi,\qquad
 Q=w_{\xi\xi},\qquad
 A=w_{\xi\xi\xi},\qquad
 H=w_{\xi\xi\xi\xi}.
\]
Then, on every fixed interior strip
\(\mathscr S_\vartheta(L,d,\rho,\alpha)\),
\[
 |H|\le CL^{-3},
 \qquad
 |A|\le C|\xi|L^{-3},
 \qquad
 |Q+j|\le C\xi^2L^{-3},
 \qquad
 |P|\le \frac{C|\xi|}{L}.
\]
Moreover, at the central observation point,
\[
 |k+8j^3|\le CL^{-4}.
\]
\end{lemma}

\begin{proof}
We divide the proof into eight steps.

\emph{Step 1. Differentiated equations.}

For \(m=3,4\), let \(D_m=\partial_\xi^m w\). Since
\[
 \partial_\xi^m(\xi w_\xi)
 =
 \xi\,\partial_\xi D_m+mD_m,
\]
differentiation of \eqref{cl:w} gives
\[
 \partial_\tau D_m
 =
 \nu^{-1}D_{m,\xi\xi}
 -\frac{\xi}{2}D_{m,\xi}
 +\left(1-\frac m2\right)D_m
 +\nu^{-1}\partial_\xi^m(VP)
 -\nu^{-1}\partial_\xi^m e,
 \qquad m=3,4.
 \tag{6.4}\label{eq:rawDm}
\]
Here $V$ is the full coefficient in \eqref{cl:w}, and the product derivative
contains the term $V D_{m,\xi}$. When we expand the differentiated factors
of $V$ using \eqref{cl:V}, we collect the terms containing $b$ or its
spatial derivatives together with $-\nu^{-1}\partial_\xi^m e$ in
$\mathcal E_m$. The strip hypotheses give
\[
 |\mathcal E_m|\le C e^{-u}.
\]
For the principal nonlinear term, equation \eqref{cl:V} gives
\[
 V_\xi^{(0)}
 =
 -\frac Qw+\frac{P^2}{2w^2}.
\]
Further differentiation gives
\[
 V_{\xi\xi}^{(0)}
 =
 -\frac Aw+\frac{2PQ}{w^2}-\frac{P^3}{w^3},
\]
and
\[
 V_{\xi\xi\xi}^{(0)}
 =
 -\frac Hw+\frac{3PA}{w^2}
 +\frac{2Q^2}{w^2}
 -\frac{7P^2Q}{w^3}
 +\frac{3P^4}{w^4}.
\]

For \(m=3\), Leibniz' rule yields
\[
 \partial_\xi^3(VP)
 =
 VH+3V_\xi A+3V_{\xi\xi}Q+V_{\xi\xi\xi}P.
\]
The coefficient of \(H\) is \(V-P/w\). The coefficient of \(A\) is
\[
 3\left(-\frac Qw+\frac{P^2}{2w^2}\right)
 +3\left(-\frac Qw+\frac{P^2}{w^2}\right)
 =
 \frac{9P^2}{2w^2}-\frac{6Q}{w}.
\]
Define
\[
 c_3=\frac{9P^2}{2w^2}-\frac{6Q}{w}.
\]
The remaining terms are
\[
\begin{aligned}
 F_3
 &=
 3\left(
 \frac{2PQ}{w^2}-\frac{P^3}{w^3}
 \right)Q
 \\
 &\quad+
 \left(
 \frac{2Q^2}{w^2}
 -\frac{7P^2Q}{w^3}
 +\frac{3P^4}{w^4}
 \right)P
 \\
 &=
 \frac{3P^5}{w^4}
 -\frac{10P^3Q}{w^3}
 +\frac{8PQ^2}{w^2}.
\end{aligned}
\]
Thus \(A=D_3\) satisfies
\[
 \begin{aligned}
 \partial_\tau A
 &=
 \nu^{-1}A_{\xi\xi}
 +\left(
 -\frac{\xi}{2}+\frac{V-P/w}{\nu}
 \right)A_\xi
 \\
 &\quad+
 \left(-\frac12+\frac{c_3}{\nu}\right)A
 +\frac{F_3}{\nu}
 +\mathcal E_3.
 \end{aligned}
 \tag{6.5}\label{eq:A}
\]

For \(m=4\), one more differentiation gives
\[
 V_{\xi\xi\xi\xi}^{(0)}
 =
 -\frac{w_{\xi\xi\xi\xi\xi}}{w}
 +\frac{4PH}{w^2}
 +\frac{7QA}{w^2}
 -\frac{13P^2A}{w^3}
 -\frac{18PQ^2}{w^3}
 +\frac{33P^3Q}{w^4}
 -\frac{12P^5}{w^5}.
\]
Using
\[
 \partial_\xi^4(VP)
 =
 Vw_{\xi\xi\xi\xi\xi}
 +4V_\xi H
 +6V_{\xi\xi}A
 +4V_{\xi\xi\xi}Q
 +V_{\xi\xi\xi\xi}P,
\]
the coefficient of \(w_{\xi\xi\xi\xi\xi}\) is \(V-P/w\), and the coefficient
of \(H\) is
\[
 4\left(-\frac Qw+\frac{P^2}{2w^2}\right)
 -\frac{4Q}{w}+\frac{4P^2}{w^2}
 =
 \frac{6P^2}{w^2}-\frac{8Q}{w}.
\]
Define
\[
 c_4=\frac{6P^2}{w^2}-\frac{8Q}{w}.
\]
The remaining terms are
\[
\begin{aligned}
 F_4={}&
 \frac{8Q^3}{w^2}
 -\frac{6A^2}{w}
 +\frac{31APQ}{w^2}
 -\frac{19AP^3}{w^3}
 \\
 &-\frac{46P^2Q^2}{w^3}
 +\frac{45P^4Q}{w^4}
 -\frac{12P^6}{w^5}.
\end{aligned}
\]
Consequently,
\[
 \begin{aligned}
 \partial_\tau H
 &=
 \nu^{-1}H_{\xi\xi}
 +\left(
 -\frac{\xi}{2}+\frac{V-P/w}{\nu}
 \right)H_\xi
 \\
 &\quad+
 \left(-1+\frac{c_4}{\nu}\right)H
 +\frac{F_4}{\nu}
 +\mathcal E_4.
 \end{aligned}
 \tag{6.6}\label{eq:H}
\]

Because
\[
 \nu=1+j+e(u,0),
 \qquad
 0\le j\le CL^{-1},
 \qquad
 |e(u,0)|\le Ce^{-u},
\]
we have, for \(L\) sufficiently large,
\[
 \frac12\le\nu\le2.
 \tag{6.7}\label{eq:nu-bounds}
\]

\emph{Step 2. The initial estimate for \(P\).}

Set
\[
 \rho_0(\xi)=\sqrt{2w(u,\xi)}.
\]
We know from the level-set geometry that \(\rho_0\) is concave in \(\xi\). Then for
\(|\xi|\le(1-\vartheta)R\), it follows from the secant inequalities that
\[
 \frac{\rho_0(\xi)-\rho_0(-R)}{\xi+R}
 \ge \rho_0'(\xi)
 \ge
 \frac{\rho_0(R)-\rho_0(\xi)}{R-\xi}.
\]
Note that
\[
 \xi+R\ge\vartheta R,
 \qquad
 R-\xi\ge\vartheta R,
\]
and \(\rho_0\) is uniformly bounded. So we get
\[
 |\rho_0'(\xi)|\le\frac{C}{R}.
\]
We differentiate \(\rho_0^2=2w\) and then we have
\[
 2\rho_0\rho_0'=2P,
 \qquad
 P=\rho_0\rho_0'.
\]
Hence, we get
\[
 |P|\le \frac{C}{R}.
 \tag{6.8}\label{eq:Pinitial}
\]

Now we see that the equation \eqref{cl:V} can be integrated explicitly. Note that
\[
 \left(\frac Pw\right)'
 =
 \frac Qw-\frac{P^2}{w^2},
\]
and so we have
\[
 -\frac Qw+\frac{P^2}{2w^2}
 =
 -\left(\frac Pw\right)'
 -\frac12\left(\frac Pw\right)^2.
\]
Then by using \(V(u,0)=0\), we get
\[
 V(u,\xi)
 =
 -\frac{P(u,\xi)}{w(u,\xi)}
 -\frac12\int_0^\xi
 \left(\frac{P(u,\eta)}{w(u,\eta)}\right)^2\,d\eta
 +\int_0^\xi b(u,\eta)\,d\eta.
 \tag{6.9}\label{eq:V-integral}
\]
Since \(w\ge c\), \eqref{eq:Pinitial} implies
\[
 \left|\frac Pw\right|\le\frac{C}{R},
\]
\[
 \left|
 \int_0^\xi
 \left(\frac Pw\right)^2\,d\eta
 \right|
 \le\frac{C|\xi|}{R^2}
 \le\frac{C}{R},
\]
and
\[
 \left|\int_0^\xi b\,d\eta\right|
 \le CRe^{-u}.
\]
Hence
\[
 |V|\le\frac{C}{R}+CRe^{-u}.
 \tag{6.10}\label{eq:V-bound}
\]

\emph{Step 3. Interpolation and the parabolic comparison.}

For \(m=3,4\), the reaction coefficients in \eqref{eq:A} and
\eqref{eq:H} are negative once \(P\) and \(Q\) are sufficiently small.
The required smallness of \(Q\) follows from interpolation.

On a smaller strip with a fixed spatial margin, the bound for
$P^{(4)}=\partial_\xi^5w$ gives
\[
 |P'|
 \le
 C\left(h^{-1}\|P\|_\infty+h^3\right),
 \qquad 0<h\le1.
 \tag{6.11}\label{eq:interpolation}
\]
For $L$ large this margin keeps $\xi\pm h$ and $\xi\pm2h$ in the larger
strip for $0<h\le1$.
Expanding at $\xi$ through degree three gives
\[
 P'(\xi)=
 \frac{P(\xi-2h)-8P(\xi-h)+8P(\xi+h)-P(\xi+2h)}{12h}
 +O\!\left(h^3\|P^{(4)}\|_\infty\right).
\]
The constant, quadratic and cubic terms cancel, and the quotient is
bounded by $Ch^{-1}\|P\|_\infty$ on the larger strip.

If \(R=\rho L^{1/3}\), choose \(h=L^{-1/12}\). Then
\[
 h^{-1}\|P\|_\infty
 \le CL^{1/12}L^{-1/3}
 =CL^{-1/4},
 \qquad
 h^3=L^{-1/4},
\]
and hence
\[
 |Q|\le CL^{-1/4}.
 \tag{6.12}\label{eq:Q-narrow-initial}
\]
If \(R=\rho L^{1/2}\), choose \(h=L^{-1/6}\). Then
\[
 h^{-1}\|P\|_\infty
 \le CL^{1/6}L^{-1/2}
 =CL^{-1/3},
 \qquad
 h^3=L^{-1/2},
\]
and therefore
\[
 |Q|\le CL^{-1/3}.
 \tag{6.13}\label{eq:Q-wide-initial}
\]

For later use, we let
\[
 \mathscr P_m
 =
 \partial_\tau
 -\nu^{-1}\partial_{\xi\xi}
 -\left(
 -\frac{\xi}{2}+\frac{V-P/w}{\nu}
 \right)\partial_\xi
 -\lambda_m,
\]
where
\[
 \lambda_3=-\frac12+\frac{c_3}{\nu},
 \qquad
 \lambda_4=-1+\frac{c_4}{\nu}.
\]
Once \(|Q|\) is bounded by \eqref{eq:Q-narrow-initial} or
\eqref{eq:Q-wide-initial}, we have
\[
 \lambda_3\le-\frac14,
 \qquad
 \lambda_4\le-\frac12.
 \tag{6.14}\label{eq:negative-reaction}
\]

The lateral boundary is controlled by
\[
 \Psi(\xi)=\exp\bigl(\gamma(\xi^2-(R')^2)\bigr).
\]
Since
\[
 \Psi_\xi=2\gamma\xi\Psi,
 \qquad
 \Psi_{\xi\xi}
 =(2\gamma+4\gamma^2\xi^2)\Psi,
\]
we obtain
\[
 \begin{aligned}
 \frac{\mathscr P_m\Psi}{\Psi}
 &=
 -\frac{2\gamma+4\gamma^2\xi^2}{\nu}
 -2\gamma\xi
 \left(
 -\frac{\xi}{2}
 +\frac{V-P/w}{\nu}
 \right)
 -\lambda_m.
 \end{aligned}
\]
The coefficient of \(\xi^2\) is
\[
 \gamma-\frac{4\gamma^2}{\nu},
\]
which is positive for sufficiently small \(\gamma\), by \eqref{eq:nu-bounds}.
The constant contribution is positive by \eqref{eq:negative-reaction}.
So we see that a sufficiently large multiple of \(\Psi\) is a supersolution near
the lateral boundary.

If a source term is bounded by \(B\), then a constant multiple \(CB\) is also a
supersolution, because we have
\[
 \mathscr P_m(CB)
 =
 -\lambda_mCB
 \ge\frac14CB.
\]
Note that the upper \(u\)-boundary is separated from the interior by the fixed
\(u\)-margin \(\vartheta dL\), and this is a positive
parabolic time separation since \(\tau=-u\). It follows that the negative reaction contributes
a factor \(e^{-cL}\) to the influence of that boundary. Then by the parabolic
maximum principle, we get interior estimates from the source size and
the fixed margins.

\emph{Step 4. Bootstrap on the narrow strip.}

We assume that \(R=\rho L^{1/3}\). It follows from \eqref{eq:Pinitial} and
\eqref{eq:Q-narrow-initial} that
\[
 |P|\le CL^{-1/3},
 \qquad
 |Q|\le CL^{-1/4}.
\]
We see that the three terms in \(F_3\) satisfy
\[
 \frac{|P|^5}{w^4}\le CL^{-5/3},\qquad
 \frac{|P|^3|Q|}{w^3}\le CL^{-5/4},
\]
and
\[
 \frac{|P||Q|^2}{w^2}\le CL^{-5/6}.
\]
So we have
\[
 |F_3|\le CL^{-5/6}.
\]
Then we apply the parabolic estimate to \eqref{eq:A} and get
\[
 |A|\le CL^{-5/6}.
\]
Note that \(w\) is even, and so
\[
 A(u,0)=0.
\]
It follows that
\[
 |Q(u,\xi)+j(u)|
 =
 \left|\int_0^\xi A(u,\eta)\,d\eta\right|
 \le CRL^{-5/6}
 \le CL^{-1/2}.
\]
Then by \(j\le C/L\), we get
\[
 |Q|\le CL^{-1/2}.
 \tag{6.15}\label{eq:Q-narrow-first}
\]

By using this improved estimate, we have
\[
 \frac{|P||Q|^2}{w^2}\le CL^{-4/3},
 \qquad
 \frac{|P|^3|Q|}{w^3}\le CL^{-3/2},
 \qquad
 \frac{|P|^5}{w^4}\le CL^{-5/3}.
\]
Hence, we get
\[
 |F_3|\le CL^{-4/3}.
\]
The equation for \(A\) gives
\[
 |A|\le CL^{-4/3},
\]
and integration over a spatial interval of length \(O(L^{1/3})\) yields
\[
 |Q+j|\le CL^{-1}.
\]
Consequently,
\[
 |Q|\le CL^{-1}.
 \tag{6.16}\label{eq:Q-narrow}
\]

Since \(P(u,0)=0\),
\[
 |P(u,\xi)|
 \le
 \int_0^{|\xi|}|Q(u,\eta)|\,d\eta
 \le\frac{C|\xi|}{L}.
 \tag{6.17}\label{eq:P-final}
\]
Substituting \eqref{eq:Q-narrow} and \eqref{eq:P-final} into \(F_3\),
\[
 \frac{|P||Q|^2}{w^2}
 \le
 C\frac{|\xi|}{L^3}
 \le CL^{-8/3},
\]
\[
 \frac{|P|^3|Q|}{w^3}
 \le
 C\frac{|\xi|^3}{L^4}
 \le CL^{-3},
\]
and
\[
 \frac{|P|^5}{w^4}
 \le
 C\frac{|\xi|^5}{L^5}
 \le CL^{-10/3}.
\]
Therefore
\[
 |F_3|\le CL^{-8/3}.
\]
The equation for \(A\) and then the equation for \(H\) give
\[
 |A|\le CL^{-8/3},
 \qquad
 |F_4|\le CL^{-3},
 \qquad
 |H|\le CL^{-3}.
 \tag{6.18}\label{eq:H-final}
\]

\emph{Step 5. Bootstrap on the wide strip.}

Assume now \(R=\rho L^{1/2}\). Initially,
\[
 |P|\le CL^{-1/2},
 \qquad
 |Q|\le CL^{-1/3}.
\]
The dominant term in \(F_3\) is \(PQ^2/w^2\), so
\[
 |F_3|
 \le
 CL^{-1/2}L^{-2/3}
 =
 CL^{-7/6}.
\]
Consequently,
\[
 |A|\le CL^{-7/6},
\qquad
 |Q+j|\le CRL^{-7/6}\le CL^{-2/3}.
\]
Thus
\[
 |Q|\le CL^{-2/3}.
\]
Using this bound once more,
\[
 |F_3|
 \le
 CL^{-1/2}L^{-4/3}
 =
 CL^{-11/6},
\]
so
\[
 |A|\le CL^{-11/6},
 \qquad
 |Q+j|\le CRL^{-11/6}\le CL^{-1}.
\]
Therefore
\[
 |Q|\le CL^{-1}.
\]
The same integration as in \eqref{eq:P-final} gives
\[
 |P|\le\frac{C|\xi|}{L}.
\]
A third application of the equation for \(A\) gives
\[
 |A|\le CL^{-5/2}.
\]
Substitution into \(F_4\) gives
\[
 |F_4|\le CL^{-3},
\]
and hence
\[
 |H|\le CL^{-3}.
\]
Thus the bound $|H|\le CL^{-3}$ in \eqref{eq:H-final} holds for both choices
\(\alpha=1/3\) and \(\alpha=1/2\).

\emph{Step 6. Integration from the central point.}

Because \(w\) is even,
\[
 P(u,0)=0,
 \qquad
 A(u,0)=0,
 \qquad
 Q(u,0)=-j(u).
\]
Using this common bound for $H$,
\[
 |A(u,\xi)|
 =
 \left|\int_0^\xi H(u,\eta)\,d\eta\right|
 \le C|\xi|L^{-3}.
\]
Integrating again,
\[
 |Q(u,\xi)+j(u)|
 =
 \left|\int_0^\xi A(u,\eta)\,d\eta\right|
 \le C\xi^2L^{-3}.
\]
Finally,
\[
 \begin{aligned}
 |P(u,\xi)|
 &=
 \left|\int_0^\xi Q(u,\eta)\,d\eta\right|\\
 &\le
 \int_0^{|\xi|}
 \left(
 \frac{C}{L}+\frac{C\eta^2}{L^3}
 \right)\,d\eta\\
 &\le
 \frac{C|\xi|}{L}.
 \end{aligned}
\]
This proves the first four estimates in the lemma.

\emph{Step 7. The estimate for \(k+8j^3\).}

From \eqref{cl:clock} and the estimates already proved,
\[
 |j_u|\le CL^{-2}.
 \tag{6.19}\label{eq:ju}
\]
At \(\xi=0\),
\[
 P=0,\qquad A=0,\qquad Q=-j,\qquad w=1.
\]
Therefore
\[
 F_4(u,0)=8Q(u,0)^3=-8j^3.
\]
We next estimate the error away from the center. Since
\[
 |Q+j|\le C\xi^2L^{-3},
 \qquad
 |P|\le C|\xi|L^{-1},
 \qquad
 |A|\le C|\xi|L^{-3},
\]
and
\[
 |w-1|
 =
 \left|\int_0^\xi P(u,\eta)\,d\eta\right|
 \le C\frac{\xi^2}{L},
\]
we have
\[
 \left|\frac1{w^m}-1\right|
 \le C\frac{\xi^2}{L}
\]
for every fixed \(m\). Since \(j\le C/L\), it follows that
\[
 \left|\frac{Q^3}{w^2}+j^3\right|
 \le CL^{-4}(1+\xi^2).
\]
The remaining terms satisfy
\[
 \frac{A^2}{w}\le C\xi^2L^{-6},\qquad
 \frac{|APQ|}{w^2}\le C\xi^2L^{-5},\qquad
 \frac{|AP^3|}{w^3}\le C\xi^4L^{-6},
\]
\[
 \frac{|P^2Q^2|}{w^3}\le C\xi^2L^{-4},\qquad
 \frac{|P^4Q|}{w^4}\le C\xi^4L^{-5},\qquad
 \frac{|P|^6}{w^5}\le C\xi^6L^{-6}.
\]
Because \(\xi^2\le CR^2\le CL\), all these terms are bounded by
\[
 CL^{-4}(1+\xi^2).
 \tag{6.20}\label{eq:F4-center-error}
\]

Let
\[
 Y=H+8j^3.
\]
The derivative of the added term \(8j^3\) is bounded by \(CL^{-4}\)
because of \eqref{eq:ju}. Combining this fact with
\eqref{eq:F4-center-error}, we find that \(Y\) satisfies a parabolic
equation whose source is bounded by
\[
 CL^{-4}(1+\xi^2).
\]
Choose
\[
 B(\xi)=CL^{-4}(K+\xi^2).
\]
Then
\[
 B_\xi=2CL^{-4}\xi,
 \qquad
 B_{\xi\xi}=2CL^{-4}.
\]
Using the negative reaction term and the inward linear drift in the
operator for \(H\), one obtains, after choosing \(K\) sufficiently large,
\[
 \mathscr P_4B
 \ge C_1L^{-4}(1+\xi^2).
\]
The lateral and upper-time boundary values are controlled by the same
fixed-margin barriers used above. Applying the parabolic maximum principle
to \(Y-B\) and to \(-Y-B\), we get
\[
 |Y(u,\xi)|\le CL^{-4}
\]
at the central observation point. Since
\[
 Y(u,0)=H(u,0)+8j(u)^3=k(u)+8j(u)^3,
\]
we conclude that
\[
 |k+8j^3|\le CL^{-4}.
\]
This completes the proof.
\end{proof}

\noindent\emph{One-sided form of Lemma~\ref{cl:strip}.}
Under the same hypotheses, the four spatial estimates hold on
\[
 L\le v\le (1+d-\vartheta d)L,\qquad
 |\xi|\le (1-\vartheta)R,\qquad R=\rho L^\alpha ,
\]
where $v$ denotes the height variable. At $\xi=0$ one also has
$|k(v)+8j(v)^3|\le C L^{-4}$ throughout this interval, including $v=L$. The constants depend only on the fixed
strip parameters and the bounds in the hypotheses.

\begin{proof}
In forward parabolic time $\tau=-v$, the initial face is $v=(1+d)L$
and the terminal face is $v=L$. Choose a fixed integer $N$ large enough
to allow one shrinkage for each interpolation and comparison in the
proof above, including the last comparison for $Y=H+8j^3$. Use the strips
\[
 [L,(1+d)L-q\vartheta dL/N]\times
 [-(1-q\vartheta/N)R,(1-q\vartheta/N)R],\qquad 0\le q\le N ,
\]
with integer $q$. All have the same terminal face. At each comparison,
write $\tau_0$ for the outer initial time and $R'$ for the outer half-width.
By \eqref{eq:negative-reaction}, a fixed multiple of
$\exp[-(\tau-\tau_0)/4]$ controls the initial values. On the next strip
it is $O(e^{-cL})$. The lateral barrier
$\exp[\gamma(\xi^2-(R')^2)]$ is $O(e^{-cR^2})$ there, since the
half-width decreases by $\vartheta R/N$. The boundary multipliers are
uniform by the assumed derivative bounds. For both values of $\alpha$
these errors are $O(L^{-p})$ for every fixed $p$.
The comparisons for $A$ and $H$, the integrations from $\xi=0$, and the
comparison for $Y$ with $CL^{-4}(K+\xi^2)$ therefore give the same
estimates up to $\tau=-L$. Each comparison uses only the initial and
lateral data of its outer strip.
\end{proof}

By the estimates of the fourth derivatives, we are able to define the phase $\Theta$. 

\begin{proposition}\label{cl:invariant}
Each normalized non-Bryant soliton has a finite equatorial phase
$\Theta$, which is determined by
\[
 j^{-1}=2u+3\log u+\Theta+O((\log u)/u).
\]
This number depends only on the isometry class of the soliton.
\end{proposition}
\begin{proof}
We apply Lemma~\ref{cl:strip} on a thick strip around each observation
height. At the equator, we get $k=-8j^3+O(u^{-4})$. Moreover, the geometric
errors in \eqref{cl:clock} are $O(e^{-u})$. Since $j=O(u^{-1})$, we see from
the clock that
\[
 j_u=\bigl(-2j^2-8j^3+O(u^{-4})\bigr)
          \bigl(1-j+O(u^{-2})\bigr)
     =-2j^2-6j^3+O(u^{-4}).
\]
Since $j\sim1/(2u)$, we have $(1/j)_u=2+6j+O(u^{-2})$.
Integrating once, we get $1/j=2u+O(\log u)$. Substituting this into the
right-hand side, we obtain $2+3/u+O((\log u)/u^2)$. The last error term is
integrable, and its tail is of size $O((\log u)/u)$. This proves the expansion
in the proposition.

Next we check the invariance under an isometry of the normalized metrics.
If we pull back the second potential by such an isometry, we get a function
which has the same Hessian as the first potential. Hence the difference of
these two functions has parallel gradient. Since the Ricci curvature is
positive, there is no nonzero parallel vector, and the additive constant is
fixed by the minimum values. Thus, the isometry identifies levels with the
same $s$. On a sufficiently high level, we have
$j\sim(2u)^{-1}$, so $K_m\ne K_o$ at the equator. There, the intrinsic Ricci
eigenvalues are $2K_m$ in the meridian direction and $K_m+K_o$
on the orbit plane. By Proposition~\ref{prop:level-profile}, at any regular
point of the target level set where there are two distinct eigenvalues, the double
eigenspace is also tangent to the sphere orbit. At a pole, the
isotropy acts transitively on unit tangent directions, so the intrinsic
Ricci tensor is scalar there. It follows that the image of the
equatorial sphere is everywhere tangent to the target sphere orbits
and is contained in one of them. Since it is compact and the dimensions are
equal, it is the whole orbit. Note that the equator is totally geodesic in the
level, and that the principal curvatures of a target orbit in its level are
$r_z/r$ up to sign. We also know from the strict concavity that there is exactly one orbit with
$r_z=0$. Hence, we see that the isometry maps the equator to the equator. So it preserves
the area radius $U$ and both sectional curvatures there, and then it also
preserves $u$, $j$ and $\Theta$.
\end{proof}

\section{Realization of the phase on Lai's expander paths}\label{cl:selection}

\begin{defi}[The compact continuum of marked construction limits]
We denote by \(\mathscr X_{\mathrm{mk}}\) the space of normalized marked pointed
smooth isometry classes
\[
 [M,g,f,o,\mathfrak m,\lambda],\qquad \lambda\ge0 .
\]
Here \(o\) is the marked tip, and the soliton equation and normalizations are
\[
 \nabla^2 f=\Ric+\lambda g,\qquad R(o)=1,\qquad
 f(o)=\min_M f=0,\qquad \nabla f(o)=0.
\]
The contracted Bianchi identity gives $dR=-2\Ric(\nabla f,\cdot)$,
so $d(R+|\nabla f|^2)=2\lambda\,df$. Evaluation at $o$ yields
\[
 R+|\nabla f|^2-2\lambda f=1 .
\]
The steady objects form the slice $\lambda=0$.
The marking \(\mathfrak m\) is an adapted orthonormal frame at \(o\), which is taken modulo
the compact group of frame changes preserving the isotropy representation.
Two marked tuples
\[
 (M,g,f,o,\mathfrak m,\lambda)
 \quad\text{and}\quad
 (M',g',f',o',\mathfrak m',\lambda')
\]
represent the same point of \(\mathscr X_{\mathrm{mk}}\) if $\lambda=\lambda'$ and there is a
potential-preserving isometry
\[
 F:(M,g)\longrightarrow(M',g')
\]
such that
\[
 F(o)=o'
\]
and \(dF_o\) maps the frame class \(\mathfrak m\) to the frame class
\(\mathfrak m'\).

The topology on \(\mathscr X_{\mathrm{mk}}\) is the marked pointed smooth
Cheeger--Gromov topology together with convergence of the soliton parameter. Thus a sequence
\[
 [M_j,g_j,f_j,o_j,\mathfrak m_j,\lambda_j]
 \longrightarrow
 [M_\infty,g_\infty,f_\infty,o_\infty,\mathfrak m_\infty,\lambda_\infty]
\]
if $\lambda_j\to\lambda_\infty$ and, for every \(A<\infty\), there are pointed embeddings
\[
 \Phi_j:B_{g_\infty}(o_\infty,A)\longrightarrow M_j
\]
such that
\[
 \Phi_j(o_\infty)=o_j,
\]
the pulled-back metrics, potentials, and marked frames converge in
\(C^\infty\) on \(B_{g_\infty}(o_\infty,A)\), and the convergence is
compatible as \(A\) increases. Fix a metric
\(d_{\mathrm{mk}}\) inducing this topology.

For every \(i\), let
\[
 \{E_{i,\mu}:0\le\mu\le1\}
\]
be the \(i\)-th normalized expander path. In the notation $E_{i,\mu}$ we
include the normalized potential and its soliton parameter $\lambda_{i,\mu}>0$.
Choose the adapted tip frames continuously along each path, which is possible
since their frame bundle over $[0,1]$ is trivial. Its marked image is
\[
 C_i
 :=
 \left\{
 [E_{i,\mu},o_{i,\mu},\mathfrak m_{i,\mu}]
 :
 0\le\mu\le1
 \right\}
 \subset\mathscr X_{\mathrm{mk}}.
\]
Each $C_i$ is compact and connected as the continuous image of $[0,1]$.
To see that $K=\overline{\bigcup_i C_i}$ is compact, take a sequence of
path points with indices $i_j$ and parameters $\mu_j$. A subsequence with
bounded indices lies in a finite union of compact paths. Otherwise pass to
$i_j\to\infty$. In the proof of \cite[Theorem~1.1]{Lai}, the asymptotic
volume ratios satisfy
\[
 \sup_{0\le\mu\le1}\operatorname{AVR}(E_{i,\mu})\longrightarrow0 .
\]
Thus \cite[Lemma~2.3 and equation~(2.5)]{Lai} gives
$\lambda_{i_j,\mu_j}\to0$ and a subsequence converging smoothly at the
tips, with the normalized potentials, to a steady soliton. Pull the tip
frames back by the convergence maps. The metrics at the tip converge, so
compactness of the orthonormal frame fiber gives a further subsequence
of frames converging to a frame adapted to the limiting action.
This is convergence in $\mathscr X_{\mathrm{mk}}$. Every sequence in the
union therefore has a subsequence converging in the ambient metric space,
which proves compactness of its closure $K$.

The nonempty compact subsets of $K$ form a compact space for the Hausdorff
metric. Choose a subsequence \(i_k\to\infty\) for which the compact connected sets
\(C_{i_k}\) converge in the Hausdorff metric induced by
\(d_{\mathrm{mk}}\). We define
\[
 \boxed{
 \mathcal C
 :=
 \lim_{k\to\infty}^{\,H} C_{i_k}
 }
\]
by the equivalent formula
\[
 \mathcal C
 =
 \bigcap_{N=1}^{\infty}
 \overline{\bigcup_{k\ge N}C_{i_k}}
 \subset\mathscr X_{\mathrm{mk}}.
\]
This means that
\[
 x\in\mathcal C
\]
if and only if there are integers \(k_j\to\infty\) and parameters
\(\mu_j\in[0,1]\) such that
\[
 [E_{i_{k_j},\mu_j},o_{i_{k_j},\mu_j},
   \mathfrak m_{i_{k_j},\mu_j}]
 \longrightarrow x
\]
in \(\mathscr X_{\mathrm{mk}}\). Along such a sequence the soliton parameter
tends to zero. Hence every point of $\mathcal C$ lies in the steady slice
and satisfies $R+|\nabla f|^2=1$.

\end{defi}

\begin{lemma}
By the definition above, we see that \(\mathcal C\) is the continuum of
all marked pointed smooth limits which we get by taking points on the selected
expander paths and letting the path index go to infinity. In particular, a point
of \(\mathcal C\) is not an arbitrary normalized steady soliton, but it is a limit
of a sequence
\[
 E_{i_{k_j},\mu_j},
 \qquad
 i_{k_j}\longrightarrow\infty,
 \qquad
 \mu_j\in[0,1].
\]
\end{lemma}

\begin{proof}
The equivalence follows from the definition of the Hausdorff limit.
We first suppose that \(x\in\mathcal C\). Then by
\[
 \mathcal C
 =
 \bigcap_{N=1}^{\infty}
 \overline{\bigcup_{k\ge N}C_{i_k}},
\]
we see that for every \(j\) there are \(k_j\ge j\) and a point
\(x_j\in C_{i_{k_j}}\) such that
\[
 d_{\mathrm{mk}}(x_j,x)<\frac1j.
\]
We also know from the definition of \(C_{i_{k_j}}\) that there is
\(\mu_j\in[0,1]\) such that
\[
 x_j=
 [E_{i_{k_j},\mu_j},o_{i_{k_j},\mu_j},
   \mathfrak m_{i_{k_j},\mu_j}].
\]
It follows that
\[
 [E_{i_{k_j},\mu_j},o_{i_{k_j},\mu_j},
   \mathfrak m_{i_{k_j},\mu_j}]
 \longrightarrow x.
\]

Conversely, we suppose that
\[
 [E_{i_{k_j},\mu_j},o_{i_{k_j},\mu_j},
   \mathfrak m_{i_{k_j},\mu_j}]
 \longrightarrow x
\]
with \(i_{k_j}\to\infty\). Then for every \(N\), we have \(i_{k_j}\ge N\)
for all sufficiently large \(j\). So we get
\[
 x\in
 \overline{\bigcup_{k\ge N}C_{i_k}}
\]
for every \(N\). Then we take the intersection over \(N\) and we see that
\(x\in\mathcal C\).

The endpoint statements follow from the construction of the expander paths.
In fact, the parameter values which correspond to the two endpoints converge to the
four-dimensional Bryant soliton and to
\(\mathbb R\times\operatorname{Bry}^3\), respectively. Note that these endpoint
sequences are included in the sets \(C_{i_k}\), and so both limits belong to the
Hausdorff limit \(\mathcal C\).
\end{proof}

The following lemma gives the relation of $\mathcal C$ and $\mathcal L$.

\begin{lemma}\label{cl:continuum}
 The set $\mathcal C$ is compact and connected, and it contains the four-dimensional Bryant soliton
and $\mathbb R\times\mathrm{Bryant}^3$. Moreover, every other member of this set is a
positive non-Bryant soliton in the class of Theorem~\ref{th:ray}.
\end{lemma}
\begin{proof}
The set \(\mathcal C\) is compact because it is a Hausdorff limit of compact
sets. It is connected because each \(C_{i_k}\) is connected and a Hausdorff
limit of nonempty compact connected subsets of a metric space is connected.
The endpoint limits of the expander paths give distinguished points
\[
 [\operatorname{Bry}^4]\in\mathcal C,
 \qquad
 [\mathbb R\times\operatorname{Bry}^3]\in\mathcal C.
\]

By \cite[Lemma~2.3, Definition~2.1 and Corollary~3.5]{Lai}, the limits
have the standard marked action on $\mathbb R^4$ and nonnegative curvature
operator, and they are noncollapsed at all scales. For clarity, we remark that
the degeneracy at the marked point can also be read off from the curvature
operator. Its two blocks on $\mathbb R\wedge\mathbb R^3$ and
$\Lambda^2\mathbb R^3$ are scalar and do not couple, by the central inversion
of $O(3)$. If both eigenvalues are positive, then we get positivity throughout by
Hamilton's strong maximum principle. Suppose that only the first one is
positive. Then its image could not be a Lie subalgebra, since brackets of mixed
rotations generate the second summand. This contradicts the same principle
applied to the ancient canonical flow \cite{Ham}. If the first one is
zero, then the axial Ricci null direction splits off a parallel line.
The remaining three-dimensional factor is rotationally symmetric,
steady and nonflat, and hence it is Bryant. The cover splits globally since the
marked manifold is diffeomorphic to $\mathbb R^4$. Finally, the flat case is
ruled out, since the scalar curvature at the mark is one.
\end{proof}

\subsection{Uniform estimates away from the poles}
We give the uniformity argument in detail, because compactness at the tip
alone says nothing about a strip whose height moves to infinity. For the
positive members, the tangent shrinker is the standard $S^3\times\mathbb R$ or
$S^2\times\mathbb R^2$ by \cite[Corollary~5.4 and Theorem~1.5]{CMZ}.
By \cite[Theorems~1.1 and~1.4]{clEntropy}, the asymptotic Nash entropy
agrees with the shrinker entropy. Then the smaller one of the entropies of
these two standard shrinkers gives a single lower bound for the whole family.
Since the curvature operator is nonnegative, Hamilton's trace Harnack
inequality holds and we have $|\Rm|\le C_4 R$. With these common constants,
we get from
\cite[Appendix A, Theorems~A.2 and~A.10 and Corollary~A.9]{CMZ}
the pointed smooth compactness at the base scalar-curvature scale,
the long-range curvature bound and the curvature derivative bounds.

\begin{lemma}\label{cl:uniform-orbit-radius}
There is a constant \(C<\infty\), independent of the positive members of
\(\mathcal C\), such that
\[
 R(p)\,r(p)^2\le C
\]
at every point \(p\) of every positive member of \(\mathcal C\).
\end{lemma}

\begin{proof}
Suppose, for contradiction, that no uniform constant exists. Then there are
positive members
\[
 (M_i,g_i,f_i)\in\mathcal C
\]
and points \(p_i\in M_i\) such that
\[
 Q_i:=R_{g_i}(p_i)>0,
 \qquad
 \widetilde r_i(p_i):=\sqrt{Q_i}\,r_i(p_i)\longrightarrow\infty.
 \tag{7.1}\label{cl:uniform-radius-failure}
\]
We rescale by
 $\widetilde g_i:=Q_i g_i$.
The scalar curvature rescales according to $R_{\widetilde g_i}=Q_i^{-1}R_{g_i}$,
so $R_{\widetilde g_i}(p_i)=1.$
The orbit radius rescales according to
$\widetilde r_i=\sqrt{Q_i}\,r_i.$ The uniform curvature estimates and noncollapsing estimates established
before this lemma imply pointed smooth compactness at the scale
\(Q_i^{-1/2}\). Passing to a subsequence, the pointed canonical flows
$(M_i,\widetilde g_i(\tau),p_i)$
 converge smoothly on compact subsets to a complete pointed ancient flow
$(M_\infty,\widetilde g_\infty(\tau),p_\infty).$
At time zero, the limit satisfies
$R_{\widetilde g_\infty}(p_\infty)=1.$

\emph{Step 1. The orbit distribution becomes parallel.}

Fix \(A<\infty\). Since \(|\widetilde\nabla\widetilde r_i|\le1\), every
point \(x\in B_{\widetilde g_i}(p_i,A)\) satisfies
\[
 \widetilde r_i(x)
 \ge
 \widetilde r_i(p_i)-d_{\widetilde g_i}(p_i,x)
 \ge
 \widetilde r_i(p_i)-A.
\]
By \eqref{cl:uniform-radius-failure},
\[
 \inf_{B_{\widetilde g_i}(p_i,A)}\widetilde r_i
 \longrightarrow\infty.
 \tag{7.4}\label{cl:radius-on-balls}
\]
Thus, for sufficiently large \(i\), the ball
\(B_{\widetilde g_i}(p_i,A)\) contains only regular \(O(3)\)-orbits.

We denote by \(P_i\) the orthogonal projection onto the tangent space of the
two-dimensional orbit through a regular point. In the warped-product
description $ \widetilde g_i
 =
 \widetilde h_i+\widetilde r_i^{\,2}g_{S^2},$
we see that the mixed components of the Levi-Civita connection contain the factor
$ \frac{\widetilde\nabla\widetilde r_i}{\widetilde r_i}.$
It follows that the covariant derivative of the projection satisfies
\[
 |\widetilde\nabla P_i|_{\widetilde g_i}
 \le
 C_4\frac{|\widetilde\nabla\widetilde r_i|}
          {\widetilde r_i}.
\]
Then by using \(|\widetilde\nabla\widetilde r_i|\le1\) and
\eqref{cl:radius-on-balls}, we get
\[
 \sup_{B_{\widetilde g_i}(p_i,A)}
 |\widetilde\nabla P_i|_{\widetilde g_i}
 \le
 \frac{C_4}{\widetilde r_i(p_i)-A}
 \longrightarrow0.
 \tag{7.5}\label{cl:projection-derivative}
\]

We choose a subsequence such that \(P_i(p_i)\) converges. Let
\(\gamma_i:[0,1]\to B_{\widetilde g_i}(p_i,A)\) be a curve of length at
most \(L\). Then by the fundamental theorem of calculus for tensor fields along
\(\gamma_i\), we have
\[
 \left|
 P_i(\gamma_i(1))
 -
 \mathcal P_{\gamma_i}P_i(\gamma_i(0))
 \right|
 \le
 L\sup_{B_{\widetilde g_i}(p_i,A)}
 |\widetilde\nabla P_i|_{\widetilde g_i},
\]
where \(\mathcal P_{\gamma_i}\) denotes parallel transport. By
\eqref{cl:projection-derivative}, the right-hand side tends to zero.
The connections converge smoothly on compact subsets, so the parallel
transports converge as well. Therefore \(P_i\) converges locally uniformly
to a tensor \(P_\infty\) satisfying
\[
 \nabla^{\widetilde g_\infty}P_\infty=0.
 \tag{7.6}\label{cl:parallel-projection}
\]
Because each \(P_i\) is self-adjoint and idempotent,
\[
 P_i^*=P_i,
 \qquad
 P_i^2=P_i,
\]
the same identities hold for \(P_\infty\). Its rank is two. Hence,
 $E_\infty:=\operatorname{Im}P_\infty$
is a parallel two-dimensional subbundle of \(TM_\infty\). The sectional curvature of an orbit plane is
$\operatorname{Sec}_{\mathrm{orb},i}
 =
 \frac{1-|\widetilde\nabla\widetilde r_i|^2}
      {\widetilde r_i^{\,2}}.$
Since the sectional curvature is nonnegative,
\[
 0\le \operatorname{Sec}_{\mathrm{orb},i}
 \le \frac1{\widetilde r_i^{\,2}}.
\]
We see from \eqref{cl:radius-on-balls} that this curvature tends uniformly to zero on
every fixed compact ball. So every two-plane which is contained in
\(E_\infty\) has zero sectional curvature. Now we pass to the complete simply connected cover and apply the de Rham
decomposition theorem, and then we get a product splitting
\[
 \widetilde M_\infty
 =
 \mathbb R^2\times N^2.
 \tag{7.7}\label{cl:two-flat-factor}
\]

\emph{Step 2. The case \(Q_i\to0\).} We first assume that $ Q_i\longrightarrow0.$
We define
$ F_i:=\sqrt{Q_i}\bigl(f_i-f_i(p_i)\bigr).$
Note that \(\widetilde g_i=Q_i g_i\), and so the gradient norm satisfies
 $|\widetilde\nabla F_i|_{\widetilde g_i}^2
 =
 |\nabla f_i|_{g_i}^2.$
Then by the normalization \(R_{g_i}+|\nabla f_i|_{g_i}^2=1\), we have
\[
 |\widetilde\nabla F_i|_{\widetilde g_i}^2
 =
 1-R_{g_i}
 =
 1-Q_iR_{\widetilde g_i}.
 \tag{7.8}\label{cl:F-gradient}
\]
On every fixed compact ball, we know that the rescaled curvature is uniformly bounded.
Then it follows from \(Q_i\to0\) and the equation \eqref{cl:F-gradient} that
\[
 |\widetilde\nabla F_i|_{\widetilde g_i}^2
 \longrightarrow1.
 \tag{7.9}\label{cl:F-unit-gradient}
\]

By the soliton equation, we have
$\nabla^2_{g_i}f_i=\operatorname{Ric}_{g_i}.$
Note that the connection is unchanged by a constant rescaling, and so we get $ \widetilde\nabla^2 F_i
 =
 \sqrt{Q_i}\,\operatorname{Ric}_{g_i}.$
Then by using \(|\operatorname{Rm}_{g_i}|\le C_4R_{g_i}\), we have on every fixed
rescaled ball $ |\widetilde\nabla^2F_i|_{\widetilde g_i}
 \le C\sqrt{Q_i}
 \longrightarrow0.$ After we pass to a subsequence, we have
 $F_i\longrightarrow F_\infty$
smoothly on compact subsets, and the limit satisfies
\[
 |\widetilde\nabla F_\infty|=1,
 \qquad
 \widetilde\nabla^2F_\infty=0.
 \tag{7.10}\label{cl:F-parallel}
\]
So we see that \(\widetilde\nabla F_\infty\) is a unit parallel vector field. We also know that each \(f_i\) is constant on \(O(3)\)-orbits, and so we have $ dF_i|_{\operatorname{Im}P_i}=0.$
Then we pass to the limit and get
$ dF_\infty|_{E_\infty}=0.$
It follows that
$ \widetilde\nabla F_\infty\in E_\infty^\perp.$
Note that the two-dimensional bundle \(E_\infty\) supplies two independent parallel
vector fields, and \(\widetilde\nabla F_\infty\) supplies a third one.
Hence, we see that the universal cover has at least three independent parallel
directions. In dimension four, the remaining one-dimensional factor is
automatically flat, and so the entire universal cover is flat. But this contradicts
the condition that the scalar curvature at
the base point is one. So the case \(Q_i\to0\) is impossible.

\emph{Step 3. The case \(Q_i\ge q>0\).}

Now we may assume that, after passing to a subsequence, $Q_i\ge q>0.$
By the normalization, we have
$ R_{g_i}\le1. $
It follows that
\[
 R_{\widetilde g_i}
 =
 Q_i^{-1}R_{g_i}
 \le q^{-1}.
 \tag{7.12}\label{cl:rescaled-scalar-bound}
\]
Note that every time slice of the rescaled canonical flow is isometric to the
time-zero slice. So \eqref{cl:rescaled-scalar-bound} holds on every
time slice. Then by the nonnegativity of the curvature operator and
$ |\operatorname{Rm}|\le C_4R,$
we get a uniform curvature bound on the entire rescaled ancient limit. We also see that the functions
$f_i-f_i(p_i)$
have uniformly bounded rescaled gradients, because
\[
 |\widetilde\nabla f_i|_{\widetilde g_i}
 =
 Q_i^{-1/2}|\nabla f_i|_{g_i}
 \le q^{-1/2}.
\]
Their Hessians are equal to the rescaled Ricci tensors and they are uniformly bounded
by the curvature estimate. So after we add constants such that the
functions vanish at \(p_i\), a subsequence converges smoothly on compact
sets to a steady potential \(F_\infty\) on the limit.
By the splitting \eqref{cl:two-flat-factor}, we can split one line from
the simply connected cover, that is,
\[
 \widetilde M_\infty
 =
 \mathbb R\times N^3.
 \tag{7.13}\label{cl:one-line-splitting}
\]
Note that the mixed and line components of
$\widetilde\nabla^2F_\infty=\operatorname{Ric}_{\widetilde g_\infty}$
vanish. Hence, if we denote by \(x\) the coordinate on the displayed
\(\mathbb R\)-factor, then we have
\[
 F_\infty(x,y)=ax+f_N(y)+b,
 \qquad
 \nabla_N^2 f_N=\operatorname{Ric}_N.
 \tag{7.14}\label{cl:split-potential}
\]
So we see that \(N^3\) is a complete three-dimensional steady gradient Ricci soliton.
The second parallel direction which comes from \(E_\infty\) lies in \(TN\), and so
\(N\) itself splits a line.

Now we show that the noncollapsing estimate descends from the four-dimensional product to
\(N\). In fact, suppose that on \(B_N(y,r)\) the curvature is bounded by
\(r^{-2}\). Then we have
\[
 B_{\mathbb R\times N}((0,y),r)
 \subset [-r,r]\times B_N(y,r),
\]
and hence
$ \operatorname{Vol}_{\mathbb R\times N}
 B_{\mathbb R\times N}((0,y),r)
 \le
 2r\,\operatorname{Vol}_N B_N(y,r).$
We also know from the four-dimensional \(\kappa\)-noncollapsing estimate that
$ \operatorname{Vol}_{\mathbb R\times N}
 B_{\mathbb R\times N}((0,y),r)
 \ge \kappa r^4.$
Then we combine the two inequalities and get
$\operatorname{Vol}_N B_N(y,r)
 \ge \frac{\kappa}{2}r^3.$ So \(N\) is a three-dimensional noncollapsed ancient solution with
nonnegative curvature and with a line splitting. By the three-dimensional classification, we know that the only noncompact
nonflat \(\kappa\)-solution which splits a line is the round cylinder.
Hence, we have $N=\mathbb R\times S^2$ with a round metric on \(S^2\). Let \(S^2\) be a compact round sphere slice in \(N\). Note that the slice is
totally geodesic. Then we restrict \eqref{cl:split-potential} to \(S^2\) and get
\[
 \Delta_{S^2}(f_N|_{S^2})
 =
 \operatorname{tr}_{TS^2}\nabla_N^2f_N
 =
 \operatorname{tr}_{TS^2}\operatorname{Ric}_N
 =
 R_{S^2}.
\]
But the round sphere has $R_{S^2}>0.$
Then we integrate over \(S^2\) and get
\[
 0
 =
 \int_{S^2}\Delta_{S^2}(f_N|_{S^2})\,d\mu_{S^2}
 =
 \int_{S^2}R_{S^2}\,d\mu_{S^2}
 >
 0,
\]
which is a contradiction. Both alternatives lead to contradictions. Therefore the assumption
\eqref{cl:uniform-radius-failure} is false, and there exists a uniform
constant \(C\) such that $ Rr^2\le C.$
on every positive member of \(\mathcal C\).
\end{proof}

As a consequence, on $w\ge c$ we have $R\le C_c/t$ uniformly.
By Lemma~\ref{lem:level-error}, with the common curvature-derivative
bounds in its proof, we get uniform scaled spatial bounds for $w$ and
the $C^k$ bounds $C_{c,k}/t$ for $e,b$. Also, the lower bound
$|\nabla f|^2=1-R\ge1/2$ there is uniform for large $t$.
These are the uniform bounds which are needed in Lemma~\ref{cl:strip}.

\subsection{Matching translated ends}
We fix one neutral reference $g$ and let $h$ vary through the neutral members
of $\mathcal C$. We choose $a,b_0\ge0$ such that
\[
 a-b_0=(\Theta_h-\Theta_g)/2,\qquad \min(a,b_0)=0.
\]
We also set $J_1(u)=j_h(u+b_0)$ and $J_2(u)=j_g(u+a)$.
\begin{lemma}\label{cl:translate}
There exist constants $C_g,u_g$, which are independent of $h,a,b_0$, such that
$|J_1(u)-J_2(u)|\le C_g u^{-3}$ for $u\ge u_g$.
\end{lemma}
\begin{proof}
For each fixed $h$, we have $u^2(J_1-J_2)\to0$ by Proposition~\ref{cl:invariant},
since the constant terms in the two inverse expansions
agree. We choose $U_h$ so large that for $u\ge U_h$ we have
$|J_1-J_2|\le u^{-2}$ and the translated profile of $h$ satisfies $w\ge1/2$
on $|\xi|\le4u^{1/3}$. By the fixed-end intermediate estimate, such a height exists,
but it may depend on $h$ and on the shifts. Starting
from this tail, we continue the difference bound and the narrower
condition $w\ge1/2$ on $|\xi|\le u^{1/3}$ towards smaller $u$.
On the continuation interval, we have $J_1,J_2\le C_g/u$. Moreover, the incoming strip
$u\le v\le(1+d)u$, $|\xi|\le u^{1/3}$ is available at each observation
height, including a candidate lower endpoint. Apply the one-sided form of
Lemma~\ref{cl:strip} with $L=u$, $\alpha=1/3$ and $\rho=1$. For $h$ the
continuation gives $w\ge1/2$ on this strip. For the fixed reference $g$ the
same bound holds at heights $v+a$ for $u\ge u_g$, by its fixed-end
intermediate estimate. The spatial derivative bounds preceding this
subsection are common to both profiles. At a translated height $v+s$,
where $s=b_0$ or $a$, the error bounds are $C_m e^{-(v+s)}\le C_m e^{-v}$.
Writing $k_1(u)=k_h(u+b_0)$ and $k_2(u)=k_g(u+a)$, we therefore have
$k_i(u)=-8J_i(u)^3+O_g(u^{-4})$ at the terminal face. At the equator,
\eqref{cl:clock} now gives
\[
 J_i'(u)=\frac{-2J_i(u)^2-8J_i(u)^3+O_g(u^{-4})}
                   {1+J_i(u)+O(e^{-u})}
        =-2J_i(u)^2-6J_i(u)^3+O_g(u^{-4}).
\]
The constants in the remainders are independent of $h,a,b_0$. Then the difference $D=J_1-J_2$
satisfies $D'+AD=O_g(u^{-4})$, where
\[
 A=2(J_1+J_2)+6(J_1^2+J_1J_2+J_2^2)
 \le 2/u+C_g(1+\log u)/u^2 .
\]
Here we used that $J_2\le1/(2u)+O_g((1+\log u)/u^2)$ uniformly for $a\ge0$.
It follows that $\exp(\int_u^v A(s)\,ds)\le C_g(v/u)^2$ for every
$v\ge u$, since $(1+\log s)/s^2$ is integrable at infinity.
Now we integrate from infinity. Since $v^2D(v)\to0$, the terminal term vanishes,
and we see that $|D(u)|$ is bounded by $C_g\int_u^\infty(v/u)^2v^{-4}\dif v
=C_g/u^3$.

We also need to improve the width condition at a possible first contact.
Let $H_0=1-w$, $R_*=u^{1/3}$, and, for the translated solution, let
$\mathcal D=\partial_\tau-\nu^{-1}\partial_{\xi\xi}
-(-\xi/2+V/\nu)\partial_\xi-1$.
Then we have $\mathcal D H_0=(e-e(0)-j)/\nu\le Ce^{-u}$.
It follows from the formula for $V$ after \eqref{cl:V}, the monotonicity of $w$ and
$|w'|\le2$ that $|V|\le C$ on the positive half strip. In fact, we see that
the quadratic integral there is bounded by using $\int_0^\xi |w'|\le1$. On
$A_0\le\xi\le R_*$, we use the function
\[
 B_0=C\xi^p/u+\exp(\gamma(\xi^2-R_*^2)),\qquad 2<p<3 .
\]
For the first term, the leading factor of $\mathcal D$ is $p/2-1>0$,
and the errors are $O(1/A_0+1/u)$. For the exponential term, the quadratic
coefficient is $\gamma-4\gamma^2/\nu>0$ when $\gamma$ is small, and we can absorb the
constant and linear errors if we increase $A_0$.
Moreover, we note that shrinking $R_*$ in forward $\tau$ contributes with a favorable sign.
At $A_0$, we know from the one-sided form of Lemma~\ref{cl:strip} that $H_0$ is bounded by $C A_0^2/u$, and at
$R_*$ the exponential is equal to one. On an initial slice $U\ge U_h$, the
wider band persists on the incoming strip. So we can apply
the one-sided form of Lemma~\ref{cl:strip} with that spatial margin and integrate from $\xi=0$, and then
we get $H_0(U,\xi)\le C\xi^2/U$ up to $U^{1/3}$. Hence, if we take the constant
in $B_0$ sufficiently large, then $B_0$ dominates these initial values for
$\xi\ge A_0$. Then by the finite-strip comparison on the moving strip, we have
$H_0\le B_0$ on the whole continuation interval below $U$.

Now we fix $\ell>0$ with $2\ell<1/4$. At height $u$, we see that the inner edge of the
last interval satisfies
\[
 H_0(R_*-\ell)\le B_0(R_*-\ell)
 \le C u^{p/3-1}
       +\exp\!\bigl(-2\gamma\ell u^{1/3}+\gamma\ell^2\bigr).
\]
Note that the right-hand side tends to zero since $2<p<3$. So we can increase $u_g$, by using only
the reference end and the common strip bounds, such that this expression
is less than $1/8$, $R_*-\ell>A_0$, and $C_g/u<1/2$ for $u\ge u_g$.
Then, at a candidate first contact, by the monotonicity of $H_0$ towards the pole,
the reflection, and $|w'|\le2$ on the last interval, we have
\[
 |J_1-J_2|\le C_g u^{-3}<\tfrac12u^{-2},\qquad
 \inf_{|\xi|\le R_*} w
 \ge 1-H_0(R_*-\ell)-2\ell>\tfrac58 .
\]
Note that the strip edge is still in the regular region since $w>0$ there.
Therefore, by smoothness, both conditions extend past a proposed lower
endpoint which is greater than $u_g$. Thus the assertion holds down to $u_g$.
Only the initial slice $U\ge\max\{U_h,u_g\}$ depends on the individual
comparison, while the threshold $u_g$ and the bound $C_g$ are independent
of $h,a,b_0$.
\end{proof}

\subsection{Detection of the round endpoint}
\phantomsection\label{cl:finite-radius}
\noindent\emph{Continuity at finite orbit radius.}
Suppose $h_i\to h_\infty$ in the marked topology on $\mathcal C$ and
$t_i\to t_*>0$. Then
\[
 j_{h_i}(\log t_i)\longrightarrow j_{h_\infty}(\log t_*),
\]
where the reading is $1$ on $\operatorname{Bry}^4$ and $0$ on
$\mathbb R\times\operatorname{Bry}^3$.

\begin{proof}
Choose representatives of the converging tip frames. For each member let
$e_1,e_2,e_3$ span the marked transverse three-plane and let
$\gamma(\ell)=\exp_o(\ell e_1)$. The reflection fixes $e_1$, hence
fixes this geodesic pointwise. On a positive member its orbit at
$\ell>0$ is the equator of its level by Theorem~\ref{th:ref}.
Indeed, $f\circ\gamma$ starts with zero derivative and
$(f\circ\gamma)''=\Ric(\dot\gamma,\dot\gamma)>0$, so it strictly
increases to infinity. The radius there is $U(f\circ\gamma)$, and
Proposition~\ref{prop:level-radius} shows that its square increases
from zero to infinity with positive derivative for $\ell>0$.

On the product the marked axis is the line factor, so $\gamma$ is
a radial geodesic in the Bryant factor at line coordinate zero.
Again $f\circ\gamma$ increases to infinity by positive radial Ricci
curvature. Each positive level is the product of the line with a round
sphere. Its orbit radius $r$ is independent of the line coordinate $z$.
The local first variation \eqref{e:ZW} thus gives
$d(r^2)/df=2+2E^o\ge2$ along $\gamma$. Here $E^o\ge0$ follows from
the Gauss identity \eqref{formula of E} on this nonnegatively curved
product. Thus the limiting ray crosses $r=\sqrt{2t_*}$ at a unique
finite parameter $\ell_*>0$ with positive radial derivative in every case.

To locate this crossing in the approximating spaces, vary $e_1$ through
$\cos(a)e_1+\sin(a)e_b$ for $b=2,3$. The resulting Jacobi fields $J_b$
along $\gamma$ have $J_b(0)=0$ and $D_\ell J_b(0)=e_b$. By exponential
equivariance the variations are unit angular rotations, so at a regular
orbit the fields are orthogonal, span its tangent plane and have length
$r$. On a fixed interval $[0,S]$ with $S>\ell_*$, smooth convergence of
the marked metrics and the geodesic and Jacobi equations give smooth
convergence of the rays and these fields in the pointed charts.
In particular $r_i=|J_{2,i}|$ converges in $C^1$ near $\ell_*$.
By transversality there are $\ell_i\to\ell_*$ with
$r_i(\ell_i)=\sqrt{2t_i}$. All the selected points therefore lie
within distance $S$ of their tips for large $i$.

At the limiting point $df_\infty\ne0$, so the tangent spaces to the levels
converge by smooth convergence of the potentials. The Jacobi fields
identify the orbit planes, and their orthogonal complements within the
level tangent spaces identify the meridian lines. The Gauss equation,
with second fundamental form $\nabla^2 f/|\nabla f|$ restricted to the
level, now gives convergence of both intrinsic curvatures $K_m,K_o$.
Since $K_o=1/(2t_i)$ there, their ratio converges as asserted.
The round levels of $\operatorname{Bry}^4$ give $K_m=K_o$, while
the line direction on each product level gives $K_m=0$.
\end{proof}

\begin{lemma}\label{cl:round}
Suppose that $h_i$ are neutral construction limits, $t_i\to\infty$, and
$j_{h_i}(\log t_i)\ge\eta>0$. Then every marked limit of the $h_i$ is the
four-dimensional Bryant soliton.
\end{lemma}
\begin{proof}
At the equator, we have $r=\sqrt{2t_i}$ and $-rr_{zz}\ge\eta$. For a distance
$c_\eta\sqrt{t_i}$ on each side, $r/\sqrt{t_i}$ stays bounded below
and the scaled spatial derivatives just proved are uniformly bounded.
Hence, by Taylor's theorem, the outward slope is at most $-c'_\eta$
at the end of a shorter such interval. By concavity, this bound still holds
up to the pole. The level has diameter at most $C_\eta\sqrt{t_i}$, because
one may go to a pole along meridians. At the equator, we have $R\asymp_\eta
1/t_i$ by Gauss, the orbital curvature $1/(2t_i)$ and the thick
bound. Here, for two positive quantities $A,B$, we write $A\asymp_\eta B$ if there is
a constant $C_\eta\geq1$ depending on $\eta$ such that $C_\eta^{-1}B\leq A\leq C_\eta B$.
Then the ambient two-point estimate gives
$c_\eta/t_i\le R\le C_\eta/t_i$ on the whole level, where we apply it also
with the two points reversed. By Lemma~\ref{lem:level-error}, the intrinsic scalar
curvature satisfies the same bounds, since the two-normal curvature
and the squared second fundamental form are $O_\eta(t_i^{-2})$.
The same conclusions, including the scaled derivative bounds, hold at
every high level with this lower bound for $j$.

We consider an ancient pointed smooth limit, where the spatial metrics are scaled
by $1/t_i$ and the base points are at the equators. The compactness follows from the common
curvature-scale bounds above and \cite[Appendix A]{CMZ}. At time zero,
$(f-s_i)/\sqrt{t_i}$ tends to a function with parallel unit gradient on
bounded balls. This is because its rescaled Hessian is $O(t_i^{-1/2})$ there by the
two-point estimate. Note that the smooth limit of the levels is contained in the transverse
slice through the base point. More explicitly, all the points of the levels lie
in a fixed ball, where the normalized functions converge as submersions,
and every subsequential limit of points lies on the same zero level.
The product splitting given by the parallel gradient has connected zero
level. The intrinsic curvature in its tangent planes is equal to the ambient
curvature there. At the equator, the three transverse curvature
eigenvalues are bounded below by a positive number in rescaled units.
The ambient curvature operator has nonnegative spectrum, it has rank at most
three on the split slice, and it has rank three at the equator. Then by Hamilton's
strong maximum principle on the ancient limit, we know that the rank is constant
on the time slice \cite{Ham}. Hence, we see that the transverse curvature
operator is positive everywhere on the bounded part which contains the
level limits. Now we take a contradiction sequence based at the least-pinched points of
these same levels, and we get a uniform positive lower bound for
the smallest intrinsic curvature eigenvalue divided by the trace.

Next we carry this information to fixed levels. Let
$S$ be the trace of the intrinsic curvature operator and let $x\le y\le z$
be its eigenvalues divided by $S$, so that $x+y+z=1$. We choose a small
$\delta>0$ which is less than half of the lower bound above. On a level where
$x\ge\delta$, the equatorial curvatures are $j/(2t),j/(2t),1/(2t)$,
and so $j$ is bounded below in terms of $\delta$. It follows that the bounds of the first
paragraph hold there. In the inward time $v=-s$, the metric variation is
$-2\Ric_\gamma-2\gamma E$. Then the corresponding error in the evolution of the curvature
endomorphism is $O_\delta(t^{-3})$, since the variation formula uses
$\nabla^2 E$ and $\Rm_\gamma E$. When we divide by $S\asymp_\delta t^{-1}$,
we lose one power of $t$. The three-dimensional Ricci flow reaction for
the smallest normalized eigenvalue, divided by $2S$, is
\begin{align*}
 x^2+yz-x(x^2+y^2+z^2+xy+xz+yz)
   &=x^2(1-3x)+(1+x)(y-x)(z-x).
\end{align*}
We can see this from the reaction $2(\alpha^2+\beta\gamma)$ for the
unnormalized eigenvalue $\alpha$. At a spatial minimum, the normalized
diffusion contributes a Laplacian and a drift $2\nabla\log S$,
with the favorable sign. In other words, at a contact we test the
curvature form minus $\delta S$ times the identity on its null vector,
which is parallel to first order. This also handles the case of a multiple eigenvalue.
At $x=\delta<1/3$, we see that the positive reaction is at least $c_\delta/t$,
and the absolute error is at most $C_\delta/t^2$.
It follows that a first return to $\delta$ in inward time is impossible for
$t\ge T_\delta$. So we know that $j$ is bounded below by a fixed positive constant on every
level $T_\delta\le t\le t_i$. By Lemma \ref{cl:finite-radius}, the same lower bound holds on the marked limit at
each fixed $t\ge T_\delta$. The product has $j=0$, and on a neutral end
$j(\log t)\to0$ by Proposition~\ref{cl:invariant}. Then, by Lemma~\ref{cl:continuum}, the only
possibility left is the Bryant soliton.
\end{proof}

The following proposition realizes any prescribed phase.
\begin{proposition}\label{cl:A}
The phase extends to a continuous map
$\mathcal C\to[-\infty,+\infty]$ which takes the value $-\infty$ at Bryant and
$+\infty$ at the product. For each $\theta\in\mathbb R$, there exists a
positive non-Bryant member $h\in\mathcal C$ such that $\Theta(h)=\theta$.
In particular, the phase of any normalized non-Bryant soliton satisfying
the hypotheses of Theorem~\ref{th:ray} is realized by a fixed such $h$.
\end{proposition}
\begin{proof}
Let $g$ be any fixed neutral reference, and consider a convergent sequence
$h_i$ of neutral members. We first assume that $\Theta_{h_i}\to\theta\in\mathbb R$.
Then the shifts $a_i,b_i$ of Lemma~\ref{cl:translate} converge to finite
$a,b$. For fixed $u\ge u_g$, apply \hyperref[cl:finite-radius]{continuity
at finite orbit radius} with $t_i=e^{u+b_i}\to e^{u+b}$. Passing to the
limit in Lemma~\ref{cl:translate} gives
\[
 |j_{h_\infty}(u+b)-j_g(u+a)|\le C_g u^{-3},
 \qquad h_\infty=\lim_i h_i .
\]
The limiting
reading is asymptotic to $1/(2u)$ as $u\to\infty$, so the limit is neither of the two
endpoints. Taking inverses, the $O(u^{-3})$ error becomes $O(u^{-1})$.
Then Proposition~\ref{cl:invariant} gives
$\Theta_{\rm lim}=\Theta_g+2(a-b)=\theta$. Next we assume that
$\Theta_{h_i}\to+\infty$. Then $a_i\to\infty$ and $b_i=0$. At each fixed
$u\ge u_g$ we apply the same continuity statement with $t_i=e^u$.
Since $j_g(u+a_i)\to0$, Lemma~\ref{cl:translate} gives
$|j_{h_\infty}(u)|\le C_g u^{-3}$, and only the product is possible. Finally, we assume that
$\Theta_{h_i}\to-\infty$. Then $a_i=0$ and $b_i\to\infty$. We fix $u_1$
so large that $j_g(u_1)-C_gu_1^{-3}>0$. Applying the lemma at $u_1$, we get
a uniform lower bound for $j_{h_i}(u_1+b_i)$, and then by Lemma~\ref{cl:round}
the limit is Bryant. Constant endpoint sequences cause no
difficulty. By the compactness of the extended real line, we now get sequential
continuity, and hence continuity. The image is connected and contains both
infinite endpoints, so it contains every finite value. We note that a neutral
reference exists, for example by \cite[Theorem~1.1]{Lai} with a tip
ratio strictly between zero and one. For a specified target, we can
equally well use the target itself as the reference. Then a point in the inverse
image of its finite phase is the required fixed member.
\end{proof}
\section{Two ends with the same phase}\label{cl:central}

We fix normalized neutral ends $g_1,g_2$ with $\Theta_1=\Theta_2$.
In the rest of the proof, all the constants are allowed to depend on this pair. By a
difference $\delta$ we mean the first quantity minus the second one, where both are evaluated
at the same $t$ and at the given spatial coordinate. In particular,
the two potential heights $s_i(t)$ need not be equal. We first make the comparison
in $\xi$. We will use a definite gain in the exponent twice, and then we do not need
any infinite-order asymptotic expansion.

\begin{proposition}\label{cl:gain}
Let $p=0$ or $3/5$. Assume that on every strip $|\xi|\le a\sqrt u$, $a<2$,
the differences of the profiles and of their fixed spatial derivatives
are $O((1+u)^N e^{-pu})$, and the differences of $e,b$ and of their
first four derivatives are $O((1+u)^N e^{-(1+p)u})$.
Here the powers $N$ are allowed to depend on $a$ and on the derivative order. Then for
$\sigma=p+2/3$, we have
\[
 |\delta j|+|\delta\lambda|
  +\sup_{|\xi|\le\sqrt u/4}\sum_{i=0}^4
       |\partial_\xi^i\delta w| \le C e^{-\sigma u}.
\]
Moreover, on any slightly smaller strip, the assertion holds through any fixed
order if we replace $\sigma$ by any slightly smaller exponent.
\end{proposition}
\begin{proof}
We freeze $u=L$ and use the independent height variable $v$ on $[L,5L]$.
We set $R_0=\sqrt{7L/2}$, $\alpha=1/5$, and
\[
 M(v)=\sup_{|\xi|\le R_0}e^{-\alpha\xi^2}
                    |\partial_\xi^4(w_1-w_2)(v,\xi)|.
\]
Note that $\sqrt{7/2}<2$, and so there is room to apply Lemma~\ref{cl:strip} on a larger strip even
at $v=L$. At other heights, we can apply it at
height $v$ on its own incoming strip, and $v$ and $L$ are comparable.
We write $H_i=\partial_\xi^4 w_i$. On the enlarged spatial strip, we have
$\|H_i\|_\infty\le CL^{-3}$ by Lemma~\ref{cl:strip}, and by the
fixed-order spatial bounds we also have
$\|\partial_\xi^6 H_i\|_\infty
 =\|\partial_\xi^{10}w_i\|_\infty\le C$.
For $k=1,2$, we differentiate the degree-five interpolation polynomial on
a stencil of spacing $h$ in this strip. The derivative weights are
$O(h^{-k})$, and the remainder in Taylor's formula is
$O(h^{6-k}\|\partial_\xi^6 H_i\|_\infty)$. Thus, on $|\xi|\le R_0$, we have
\[
 \|\partial_\xi^k H_i\|_\infty
 \le C\bigl(h^{-k}\|H_i\|_\infty
            +h^{6-k}\|\partial_\xi^6H_i\|_\infty\bigr).
\]
Taking $h=L^{-1/2}$, we get $\|H_i'\|_\infty=O(L^{-5/2})$ and
$\|H_i''\|_\infty=O(L^{-2})$. Note that the spatial margin contains the
stencil for all sufficiently large $L$.

The Taylor data of $\delta w$ through order three at zero are
$(0,0,-\delta j,0)$. Hence, by integration, the weighted
suprema of its derivatives through order three are bounded by $C(|\delta j|+M)$.
For example, the integration kernels are bounded in integral by a
constant depending on $\alpha$, since
$\xi^2-y^2\ge|\xi|(|\xi|-|y|)$ for $0\le |y|\le|\xi|$.
By the integrated formula for $V$ in the proof of Lemma~\ref{cl:strip},
the same bound holds for $\delta V$ after adding
$C(1+L)^N e^{-(1+p)v}$. These estimates hold on both halves of the strip.

Now we subtract the fourth derivative equations, keeping the diffusion, drift
and reaction of the first solution. Then the differences of the diffusion
and of the extra drift multiply $H_2''$ and $H_2'$. The difference of the reaction
multiplies $H_2=O(L^{-3})$, and we have $\delta\nu=\delta j+\delta e(0)$.
In $F_4$, the partial derivatives in $w,P,Q$ are $O(L^{-2})$ along
the segment between the two lists of arguments. This follows term by
term from $w\ge c$, $|P|\le C/\sqrt L$, $|Q|\le C/L$ and
$|A|\le C L^{-5/2}$. For the remaining partial derivative, we need the
vanishing of $A$ at the reflection sphere. In fact, we have
\[
 (F_4)_A=-12A/w+31PQ/w^2-19P^3/w^3,\qquad
 |(F_4)_A|\le C|\xi|L^{-2}.
\]
Here we used $|P|\le C|\xi|/L$, $|A|\le C|\xi|/L^3$ and
$\xi^2\le 7L/2$. Since $\delta A(0)=0$, we get
\[
 |\xi|e^{-\alpha\xi^2}|\delta A(\xi)|
 \le M(v)|\xi|\int_0^{|\xi|}e^{-\alpha(\xi^2-y^2)}\dif y
 \le M(v)/\alpha .
\]
In this way we have estimated the product, including its behavior at lateral distances
of order $\sqrt L$. Each term which is produced by differentiating $b$ in
\eqref{cl:V} is a derivative of $b$ times an individual profile
derivative, or an individual $b$-derivative times a profile difference.
The latter terms can be absorbed in $CL^{-2}(|\delta j|+M)$. The former difference
terms, as well as $\delta e^{(4)}$, have the assumed geometric error
bound. Therefore, the weighted source is at most
\[
 CL^{-2}(|\delta j(v)|+M(v))+C(1+L)^N e^{-(1+p)v}.
\]

Conjugating by $e^{-\alpha\xi^2}$, the reaction becomes at most
\[
 -1+2\alpha/\nu_1+O(L^{-1})
 +(4\alpha^2/\nu_1-\alpha)\xi^2+O(L^{-1/2}|\xi|).
\]
The negative quadratic term absorbs the last term at the cost of $O(L^{-1})$.
Hence the damping is larger than $1/3$ even after $CL^{-2}M$ is absorbed.
To be precise, we compare both signs with a scalar supersolution
which depends only on the forward time $-v$. At an interior first contact, its
value bounds $M$, and we solve the resulting scalar differential inequality
with the smaller damping $1/3$. A constant multiple of the supremum
of the lateral data and of the geometric source is a supersolution
for these contributions. The initial data at $v=5L$ contribute
$C(1+L)^N e^{-4L/3}$. At $|\xi|=R_0$, the weight contributes
$e^{-7L/10}$ in addition to the known $e^{-pL}$. It follows that
for every $\epsilon>0$,
\begin{equation}
 |\delta k(L)|\le M(L)\le
 C\int_L^{5L}e^{-(v-L)/3}v^{-2}|\delta j(v)|\dif v
 +C_\epsilon e^{-(p+7/10-\epsilon)L}.                 \label{cl:duhamel}
\end{equation}
Since $p+7/10\le13/10<4/3$, the upper-time contribution can be included
in the last term for both choices of $p$.

Let $J=\delta(j+j^2/2)$. Subtracting \eqref{cl:clock}, we get
\[
 J_u+aJ=\delta k+O((1+u)^N e^{-(1+p)u}),\qquad
 a=\frac{2(j_1+j_2)}{1+(j_1+j_2)/2}.
\]
For the product $e_1(0)\delta j_u$, we use \eqref{cl:clock} first.
Namely, the assumed bounds for $\delta k,\delta j$ and the errors imply
$\delta j_u=O((1+u)^N e^{-pu})$. All the other product differences are
directly of the order appearing in the equation above. By Proposition~\ref{cl:invariant},
we have $a=2/u+O((1+\log u)/u^2)$, so the integrating factor $I'=aI$
can be normalized by $I/u^2\to1$. Also by Proposition~\ref{cl:invariant},
equal phases give $j_1^{-1}-j_2^{-1}=o(1)$. Since
$j_i\sim(2u)^{-1}$, it follows that $u^2\delta j\to0$ and
$I J=I(1+(j_1+j_2)/2)\delta j\to0$. This terminal condition determines the
integration constant. In particular,
$S(u)=\sup_{v\ge u}|I(v)J(v)|$ is finite. We integrate the
equation from infinity, insert \eqref{cl:duhamel}, and use
$I(v)|\delta j(y)|\le C S(u)$ for $u\le v\le y\le5v$.
Then the part involving $S(u)$ is bounded by
\[
 C S(u)\int_u^\infty v^{-2}
       \int_v^{5v}e^{-(y-v)/3}\dif y\dif v
 \le 3C S(u)/u.
\]
The same bound holds if we start at any height above $u$. Taking the
supremum of the integrated estimate, we obtain
\[
 S(u)\le C S(u)/u+C_\epsilon u^2e^{-(p+7/10-\epsilon)u}.
\]
By absorption, we get a bound for $\delta j$, and then \eqref{cl:duhamel} gives
a bound for $M$ with the same exponent. When we remove the weight on
$|\xi|\le\sqrt u/4$, we lose only $e^{u/80}$. Then the integrated Taylor bounds
and the equatorial formula for $\lambda$ give every gain strictly below
$7/10-1/80=11/16$. In particular,
\[
 \frac{11}{16}-\frac23=\frac1{48}>0 .
\]
Choosing $\epsilon$ and the losses used to absorb fixed polynomial
factors such that their total is less than $1/48$, we obtain the exponent
$\sigma=p+2/3$ in the proposition. For the last assertion, we interpolate spatially
on a slightly larger interval. For a function with small supremum $d$ and bounded
derivative of order $n$, the finite-difference bound
$C(h^{-k}d+h^{n-k}\|D^n w\|_\infty)$ loses arbitrarily little in the
exponent if we first choose $h=e^{-\eta u}$ and then choose $n$ large.
\end{proof}
\section{Comparison up to the two poles}\label{cl:caps}

Let $\ell$ be the arclength measured from the pole toward the equator. On a closed
cap $0\le\rho=r/\sqrt t\le\rho_b<\sqrt2$, we set
$V_c=dr/d\ell$ and $l=1/V_c$. The function $l$ extends smoothly and evenly
to $\rho=0$, where it takes the value one. In these coordinates, the scaled
metric is given by $t^{-1}\gamma=l^2 d\rho^2+\rho^2g_{S^2}$,
while in the central chart it is given by $d\xi^2+2wg_{S^2}$.

In what follows, we will need bounds which are polynomial in $u$ for every fixed
spatial derivative of the individual cap data. We give the details at the pole. Let
$L_c=\ell/\sqrt t$ and $v(L_c)=r/\sqrt t$. Then we have $v(0)=0$, $v'(0)=1$,
$|v'|\le1$, and $v''=-(tK_m)v$, where $K_m$ denotes the meridional sectional
curvature of the level. By Lemmas~\ref{lem:level-error}
and~\ref{lem:mode-levelR} and \eqref{e:der}, the meridional
derivatives of $tK_m$ in the $L_c$ variable are bounded by powers of $1+u$.
Here we use a parallel unit orbital vector along the meridian and extend it along
the geodesic through the pole. Note that the curvature contractions are smooth
there. Differentiating the equation, we get bounds for all derivatives of $v$.
By Lemma~\ref{lem:mode-secant}, the slope satisfies $v'\asymp_{\rho_b}u^{-1/2}$
at $\rho_b$, and by concavity the same lower bound holds up to
the pole. Then the required bounds on $l$ follow by inverting the odd extension of $v$.
Using the same parallel frame, the argument above also shows that scaled covariant
derivative bounds on an invariant tensor imply polynomial coordinate bounds on its
radial eigenvalues. On a central thick
band, these conversions follow directly from $w\ge c$.

\subsection{An inequality for the inverse slope}
We first compute the first variation of the radial arclength at a fixed true radius.
It is given by
\[
 (\partial_s l)_r=lG_m-\partial_r(l rG_o).
\]
To see this, we note that a material radial line element evolves by $G_m$, while its
radius moves with speed $rG_o$. Then the formula follows by subtracting the flux at the
two ends of an interval of radii, or alternatively by differentiating the arclength
from the pole at fixed $r$.
Note that the intrinsic Ricci eigenvalues are $-2V_c(V_c)_r/r$ and
$(1-V_c^2)/r^2-V_c(V_c)_r/r$. Inserting them into the formula above and changing
the variable to $\tau=-\log t$ at fixed $\rho$, we obtain
\begin{align}
 l_\tau&=-\rho l_\rho/2+\lambda\{\mathcal P(l)+F_c\},\label{cl:l}\\
 \mathcal P(l)&=\partial_\rho\{(l-V_c)/\rho-(V_c)_\rho\}
                   +2(V_c)_\rho/\rho,\notag\\
 F_c&=\partial_\rho(\rho l\,tE_o)-l\,tE_m.\notag
\end{align}
By parity and $l(0)=1$, the terms which contain quotients are smooth at zero.

Next we choose $0<\rho_a<\rho_b<\sqrt2$ so close to $\sqrt2$ that
the entire annulus, together with some room at its two ends, is contained in
$0<\xi<\sqrt u/4$ on both ends. This is possible, since by inverting
Proposition~\ref{prop:mode-regions} on compact subintervals we have
$\xi(\rho)/\sqrt u\to2\sqrt{1-\rho^2/2}$. We choose these constants once and
fix them from now on.
Suppose that a central profile difference, including its first derivative, has size
$O(e^{-\sigma u})$ there. Then the inverse coordinates and the inverse
slopes at equal $\rho$ differ by $O((1+u)^N e^{-\sigma u})$.
For the inverse coordinate, we apply the mean value theorem to
$\sqrt{2w(\xi)}$, whose slope has magnitude at least $c/\sqrt u$ on the
enlarged annulus. For the inverse slope, we use the same estimate,
differentiate the profiles at each point, and then take reciprocals.

\begin{lemma}\label{cl:kato}
Let $D=l_1-l_2$. For any $K>0$, there exist a smooth positive weight
$\psi=\rho e^{B\rho^2}$ on the open half cap and a cutoff $\chi$, which is equal
to one on $[0,\rho_a]$ and vanishes near $\rho_b$, such that for
$Y=\int_0^{\rho_b}\chi\psi|D|\dif\rho$ and for all sufficiently large $u$, we have
\[
 (\partial_\tau+K)Y\le
 C\left(\int_{\rho_a}^{\rho_b}|D|\dif\rho+
 \|\delta F_c\|_\infty+(1+u)^N|\delta\lambda|\right).
\]
\end{lemma}
\begin{proof}
Let $A_c=1/(l_1l_2)$. Then we have $\delta V_c=-A_cD$. Expanding
\eqref{cl:l}, we get the identity
\[
 \delta\mathcal P=(A_cD)''+\big((1-A_c)D/\rho\big)'
                         -2A_cD/\rho^2 .
\]
The remaining source terms in the difference equation are
$\lambda_1\delta F_c+\delta\lambda\,(\mathcal P(l_2)+F_{c,2})$,
and the second coefficient here is polynomially bounded. On each annulus
$\rho>0$, we multiply the difference equation by the sign of $D$ and obtain the
absolute-value inequality in the distributional sense. Here we have used $A_c>0$ and
$\operatorname{sgn}(A_cD)(A_cD)''\le |A_cD|''$, and for the drift and
divergence terms we have used the ordinary chain rule for the absolute value.
Before the derivatives of the cutoff are taken into account, the
adjoint acting on $\psi$ is
\[
 \lambda_1\{A_c\psi''-(1-A_c)\psi'/\rho-2A_c\psi/\rho^2\}
                    +(\rho\psi)'/2 .
\]
Dividing it by $\psi$, we get
\[
 -\lambda_1(1+A_c)/\rho^2+1+
 B(8\lambda_1A_c-2\lambda_1+\rho^2)+4\lambda_1 A_c B^2\rho^2 .
\]
Since $\lambda_1\to1$ and $\rho_b^2<2$, we can take $B$ so large
that $1+B(-2+\rho_b^2)<-K-2$. On a sufficiently small fixed pole
interval, by $0<A_c\le1$, all the positive terms are dominated by the singular
negative term. On the complementary closed annulus, we have $A_c\to0$
uniformly by the slope asymptotic. This proves the adjoint bound.
The terms which contain derivatives of the cutoff are bounded by
$C\int_{\rho_a}^{\rho_b}|D|$. Finally, the integration down
to zero introduces no boundary mass, since $D=O(\rho^2)$, $D'=O(\rho)$,
$\psi=O(\rho)$ and $\psi'=O(1)$. Equivalently, one can first insert a cutoff
at $\rho=\delta$ and then let $\delta$ decrease to zero.
\end{proof}

\begin{lemma}\label{cl:cap-transfer}
Suppose that the central estimate of Proposition~\ref{cl:gain} holds with
$\sigma>m$, and that $\|\delta F_c\|_\infty=O(e^{-\sigma' u})$ for some
$\sigma'>m$. Then $\delta l$ on every slightly smaller cap and
$\delta w$ on each thick band $|\xi|\le a\sqrt u$, $a<2$, are bounded by
$C_k e^{-mu}$ through each prescribed finite spatial order $k$.
Moreover, the scaled half-lengths from the equator to the pole differ by $O(e^{-mu})$.
\end{lemma}
\begin{proof}
We choose an intermediate exponent which lies strictly between $m$ and both input
exponents. By the inverse comparison on the annulus, the right-hand side of the
inequality in Lemma~\ref{cl:kato} has that exponent.
Note that $Y$ grows at most polynomially as $\tau\to-\infty$.
Integrating the scalar inequality, we get the same
decay for $Y$. Here we have used that the integral of $e^{-K(\tau-v)}e^{\sigma_0 v}$ is
$e^{\sigma_0\tau}/(K+\sigma_0)$. On a closed subcap where $\chi$ is bounded
below, we have for small $h>0$
\[
 |\partial_\rho^k D(\rho)|\le C_{n,k}h^{-k-2}Y+
                  C_{n,k}h^{n+1-k}(1+u)^{N_n}.
\]
To prove this formula, we integrate the Taylor polynomial up to degree $n$ against
a smooth kernel on $[\rho+h,\rho+2h]$ whose moments reproduce the
$k$th derivative at $\rho$. Such a kernel can be obtained by inverting the
moment matrix on polynomials of degree at most $n$ against a positive
smooth interior weight. After scaling, the supremum of this kernel is $C h^{-k-1}$,
and the remaining factor $h^{-1}$ comes from $\psi\ge c h$ on its support. The
remainder term is given by Taylor's formula. We note that this argument also works
at $\rho=0$.
Now we first choose $h=e^{-\eta u}$ with $\eta$ small compared with the exponent
margin, and then choose $n$ large. This still leaves a strict margin above $m$.
The cutoffs can be moved within the overlap, and on its remainder we use the
central bounds.

Next we choose an anchor $\rho_*$ in the overlap where both comparisons apply.
For $0\le\rho\le\rho_*$, we use the formula
\[
 \xi_i(\rho)=\xi_i(\rho_*)+\int_\rho^{\rho_*}l_i(y)\dif y .
\]
At zero, this is the scaled half-length. Since the slope of the inverse $\rho_i(\xi)$
is bounded by one, we get the same absolute bound on $w_1-w_2$ at equal $\xi$
throughout the common part of the sheets, where on the equatorial side of the anchor
we use the central estimate. In particular,
every fixed thick band is covered when $u$ is large. The estimate for the
higher derivatives follows by interpolation on a slightly larger band,
where we use the remaining exponent margin. We note that this argument can be
carried out with any fixed finite number of derivatives.
\end{proof}

\subsection{The first comparison and its effect on the normal error}
By \eqref{e:err}, the hypotheses of Proposition~\ref{cl:gain} are satisfied
for $p=0$. On the whole level, by Lemmas~\ref{lem:level-error}
and~\ref{lem:mode-levelR}, we have
$|\nabla^k E|\le C_k\Lambda^{2+k/2}$ with
$\Lambda=C(1+\log t)/t$. These estimates also hold at comparable
lower heights. Together with the pole coordinate bounds above, they imply
$|F_{c,i}|\le C(1+u)^N/t$ in \eqref{cl:l}. Then Proposition~\ref{cl:gain} gives
the central exponent $2/3$. Moreover, the cap source admits every exponent below
$1$. Since $2/3-3/5=1/15$, we can apply Lemma~\ref{cl:cap-transfer} with $m=3/5$
and get
\begin{equation}
 \delta l,\quad \delta w,\quad \delta\lambda=O(t^{-3/5})
                 \label{cl:first}
\end{equation}
in their respective charts, with any fixed number of spatial derivatives
on the profiles. We remark that the cap bound holds on any prescribed closed cap
below $\sqrt2$, because we may first place the overlap closer to the equator.

In order to improve \eqref{cl:first}, we need to consider the difference of the normal
errors. We denote $q_0=|\nabla f|^2$, $N=q_0^{-1/2}$, and
$\mathcal S=\Scal_\gamma$. Here we use the subscript in $q_0$ to distinguish it from
the radius drop in Section~\ref{sec:modes}. We define
\[
 \mathcal P_2(\gamma)=-D_\gamma T[2\gamma T]
                       -\tfrac12(\nabla_\gamma^2\mathcal S)^\sharp.
\]
Here $D_\gamma T[h]$ denotes the first metric variation of the Ricci
endomorphism in the direction $h$. Note that $\mathcal P_2$ is an intrinsic tensor,
which consists of quadratic curvature terms and second covariant derivatives of
the curvature.

\begin{lemma}\label{cl:normal}
For every $k\ge0$, we have uniformly on a fixed neutral end
\[
 |\nabla_\gamma^k(E-\mathcal P_2(\gamma))|
                       \le C_k\Lambda^{3+k/2}.
\]
\end{lemma}
\begin{proof}
In the normal gradient coordinates, the ambient metric is
$N^2 ds^2+\gamma$, and the mixed shape operator is $\sqrt{q_0}G$. It follows from the
normal soliton equation and the trace identity that
$(q_0)_s/2+q_0\tr G=1-q_0$. Then we combine this with the traced Gauss equation
$\mathcal S=R-2\Ric_{nn}+q_0((\tr G)^2-\tr(G^2))$ and
$\Ric_{nn}=(q_0)_s/2$, and we get the first one of the identities
\begin{equation}
 q_0=\frac{1+\mathcal S}{(1+\tr G)^2-\tr(G^2)},\qquad
 q_0\partial_sG+G=T-N^{-1}(\nabla_\gamma^2 N)^\sharp .
 \label{cl:normal-eq}
\end{equation}
To prove the second one, we use the tangential Ricci formula for the mixed shape
operator $B=\sqrt{q_0}G$, which reads
$\Ric^\sharp|_{\rm tan}=T-N^{-1}\partial_sB-(\tr B)B
-N^{-1}(\nabla_\gamma^2N)^\sharp$. Here the left-hand side means its
tangential compression, and by the soliton equation we know that it is equal to $q_0G$.
After transposition, we see that the coefficient of $G$ is
$q_0+(q_0)_s/2+q_0\tr G=1$, and then the second identity follows.

We note that all the estimates in this paragraph remain valid after we take $k$ further spatial
derivatives, with an additional factor $\Lambda^{k/2}$. By Gauss and the level error
bounds, we have $T,G=O(\Lambda)$ and $G-T=O(\Lambda^2)$. We also have
$N=1+\mathcal S/2+O(\Lambda^2)$ by the algebraic formula, and the differentiated remainders can be
estimated by the product rule. Next we consider $\partial_sT=D_\gamma T[2\gamma G]$.
In this formula, the error which is caused by the substitution
$G=T$ is $O(\Lambda^3)$. In fact, the variation of the connection is
$\gamma^{-1}*\nabla(2\gamma G)$, and the variation of the Ricci endomorphism
contains $\nabla^2G$ and curvature times $G$. Here a star means a
finite sum of metric contractions. We use the same formulas to differentiate the
curvature polynomial $\mathcal P_2$, and we get
$\partial_s\mathcal P_2=O(\Lambda^3)$. This is because in the quadratic terms
$\partial_s\Rm$ has order two in curvature units, and in the differentiated
terms the commutation of $\partial_s$ with the connection brings in derivatives
of $G$. Now we substitute these estimates into \eqref{cl:normal-eq}, and we see that
$W=G-T-\mathcal P_2$ satisfies $q_0\partial_s W+W=O(\Lambda^3)$,
and the same bound holds for derivatives of the right-hand side.

Now we integrate this equation along each $Z$-trajectory on $[s/2,s]$. On this
interval, the curvature scale $\Lambda$ is comparable to its value at the upper end.
Moreover, $q_0$ tends to one uniformly there, and by $\partial_s\gamma=2\gamma G$
we see that the tensor norms change at a rate which is at most $C\Lambda$ times the norm.
So for large $s$, the scalar comparison for $|W|$ has damping at least $1/2$,
and the lower data decay exponentially in $s$. For the
derivatives, we commute $\nabla_\gamma^k$ through this equation. A connection
commutator contains a derivative of $G$ and a lower derivative of $W$.
When derivatives fall on $q_0$, they are multiplied by lower derivatives of
$\partial_s W$, which are controlled by the equation and by the lower orders
that have already been improved. Note that $|\nabla^j q_0|\le C_j\Lambda^{1+j/2}$ for $j>0$, and so
these source terms are $O(\Lambda^{3+k/2})$. Then the lemma follows by induction and
the same scalar comparison.
\end{proof}

\begin{proposition}\label{cl:second}
The difference $\delta l$ on each closed cap and the difference $\delta w$ on each
central thick band, with any prescribed number
of spatial derivatives, as well as the scalar difference $\delta\lambda$,
are $O(t^{-6/5})$. Moreover, the scaled half-lengths differ by $O(t^{-6/5})$.
\end{proposition}
\begin{proof}
In this proof, we denote by $\bar\gamma_i$ the rescaled level metric
$t^{-1}\gamma_i$. Under a constant scaling of the metric, we have
$t\mathcal P_2(\gamma_i)=t^{-1}\mathcal P_2(\bar\gamma_i)$.
In common spatial charts, the curvature polynomials of the scaled metrics are
controlled by the profiles in \eqref{cl:first}. At a pole, we use the coordinates
$X=\rho\omega$, in which the coefficients are
\[
 (\bar\gamma_i)_{ab}=\delta_{ab}
       +\frac{l_i^2-1}{\rho^2}X_aX_b .
\]
Note that for an even smooth function $H$ with $H(0)=0$, we have
$H(\rho)/\rho^2=\int_0^1(1-v)H''(v\rho)\dif v$.
So if we take $H=\delta(l^2)$, then we can control all the required radial derivatives.
To estimate Cartesian derivatives of a smooth even radial function, we subtract
its even Taylor polynomial at zero before differentiating. Then the
remainder and its derivatives vanish to the corresponding orders.
In this way, we get Cartesian $C^k$ bounds, uniformly at the origin, from finitely
many one-dimensional bounds. Since $l\ge1$, we lose at most powers of $1+u$ when we
invert the metric coefficients. In the central chart, we use
$w\ge c_a$ instead. By Lemma~\ref{cl:normal}, we now have in either chart
\[
 \|\delta(tE)\|_{C^k}\le C_k(1+u)^{N_k}t^{-1-3/5}.
\]
In fact, the contribution of its intrinsic polynomial is $t^{-1}$ times
\eqref{cl:first}, and after multiplying by
$t$ each individual remainder is $O((1+u)^N/t^2)$ in scaled coordinates. The radial
and orbital eigenvalues can be read off along a fixed radial ray in the smooth chart.
The product differences in $e=2twE_o$, $b=tE_m$, and in \eqref{cl:l}
satisfy the same bound, because each individual $tE$ and its derivatives
are $O((1+u)^N/t)$. Now we apply
Proposition~\ref{cl:gain} for the second time, with $p=3/5$, and get
$\sigma=3/5+2/3=19/15$. Moreover, the cap source has every exponent strictly
below $1+3/5=8/5$. Since $19/15-6/5=1/15$, the transfer in
Lemma~\ref{cl:cap-transfer} can be applied with $m=6/5$, together with the
fixed number of spatial derivatives in the statement.
\end{proof}
\section{One identification of the four-dimensional ends}\label{cl:ambient}

Let $d_i(t)$ be the distance in the level from its equator to its pole.
Then $2d_i$ is the full length which was denoted by $D$ in Section~\ref{sec:geom}.
We use the coordinate $x=z_i/d_i(t)$, with $-1\le x\le1$. As the angular variable,
we take the common $O(3)/O(2)$ coordinate from Section~\ref{sec:global}.
We identify the two ends by preserving $(t,x,\omega)$, and we denote by $\Phi$ the map
from the first end to the second end. At the equator, this is a
smooth signed arclength coordinate. Near a pole, the map $\Phi$ sends the radial
geodesic in each level to the corresponding one on the other level
by $\ell_2=(d_2/d_1)\ell_1$. In the level exponential charts, it is a
linear similarity on the tangent spaces. Since these charts, their centers
and this similarity depend smoothly on $t$, the map $\Phi$ and its
inverse extend smoothly across both polar curves.

\begin{proposition}\label{cl:end}
For this identification, there is a constant $c\in\mathbb R$ such that
\[
 |\Phi^*g_2-g_1|_{g_1}=O(f_1^{-11/10}),\qquad
 f_2\circ\Phi-f_1=c+O(f_1^{-1/5})
\]
uniformly on the complete end.
\end{proposition}
\begin{proof}
By Proposition~\ref{cl:second}, we have
$|\delta(d/\sqrt t)|\le Ct^{-6/5}$, while
$d_i/\sqrt t\asymp\sqrt{\log t}$. Note that the angular coefficient has to be
compared in the relative sense, even at points which are arbitrarily close to the
poles. For this purpose, we define $L_i(\rho)=\int_0^\rho l_i(y)\dif y$ at a common
cap radius. Then we have
\[
 |L_1(\rho)-L_2(\rho)|\le C\rho\,t^{-6/5},
 \qquad L_i(\rho)\ge\rho .
\]
Thus, at the same positive radius, the two scaled distances to the pole
have relative error $O(t^{-6/5})$. Note that at common $x$ their ratio is
$d_2/d_1$. In addition, by concavity and $r(0)=0$, each metric satisfies
$0\le d\log r/d\log\ell=\ell r_\ell/r\le1$.
We first take the same radius in both metrics, and then move in the second
metric to its prescribed distance from the pole. Integrating the last
inequality in $\log\ell$, we get $|\log(r_2/r_1)|\le Ct^{-6/5}$.
Here we use overlapping subcaps so that both points stay in the chart.
On the central side, $r_i/\sqrt t$ is bounded below. In this case, we compare $w_i$
and use $\xi_i=xd_i/\sqrt t$, and we obtain the same assertion.

It remains to reconstruct the coefficient of $dt$ and the cross term.
Note that the normal rate of change of a meridian line element is $G_m$.
Hence, at a material point, we have $Zz=\int_0^zG_m\,dz'$ and
$\partial_s d_i=\int_0^{d_i}G_m\,dz'$. It follows that in these coordinates
the metric of the end is given by
\begin{align}
 g_i&=\frac{\lambda_i^2}{q_{0,i}}dt^2+
      d_i^2(dx-\lambda_i v_i\,dt)^2+r_i^2g_{S^2},\label{cl:metric}\\
 v_i(t,x)&=\int_0^xG_{m,i}(t,y)\dif y
                 -x\int_0^1 G_{m,i}(t,y)\dif y .\notag
\end{align}
Here we have $Zx=v_i$ and $Zt=\nu_i$, and this determines the sign of the cross term.
Note also that the angular direction is horizontal under the normal motion in the
fixed meridian chart.

Now we work in the two level charts of Section~\ref{cl:caps}. By
Proposition~\ref{cl:second} with two spatial derivatives, the intrinsic curvature
comparison for $tT$ has order $O((1+u)^N t^{-6/5})$. Moreover, the normal
error estimate which was used in the proof of Proposition~\ref{cl:second} already
has the better exponent $1+3/5$. Thus, the same $6/5$ bound holds for $\delta(tG)$.
In order to evaluate at common $x$, we use the difference of $\xi_i$ in the
central chart, and in the cap chart we solve
$L_i(\rho_i)=(1-|x|)d_i/\sqrt t$. Since $L_i'\ge1$, the absolute
difference of the evaluation radii is $O((1+u)^N t^{-6/5})$. Note that each scaled
tensor component and its first spatial derivative contribute at most a
power of $1+u$ in either chart. This proves the same bound on
$\delta(tG_m)$ at common $x$. By the algebraic identity
\eqref{cl:normal-eq}, we get the corresponding bound on $\delta q_0$.
In fact, this bound holds with an additional factor of $t^{-1}$.

In the orthonormal coframe of the first metric, the shift difference is
measured by $d_2|\delta(\lambda v)|$. By \eqref{cl:metric}, the boundedness of
$\lambda_i$, \eqref{cl:clock}, and Proposition~\ref{cl:second}, this quantity is
at most $C\sqrt{tu}\,(1+u)^Nt^{-1-6/5}$. Moreover, the lapse
$\lambda_1/\sqrt{q_{0,1}}$ is bounded from above and from below. Now all the relative
coefficient differences in \eqref{cl:metric} are at most
$C(1+u)^Nt^{-6/5}$. Since $(1+\log t)^N=O(t^{1/10})$ and $s_1\sim t$,
we get the metric estimate in the proposition. Finally, we have
$\partial_t(s_2-s_1)=\lambda_2-\lambda_1=O(t^{-6/5})$, which is integrable
with tail $O(t^{-1/5})$. This proves the estimate for the potential for
the same map $\Phi$.
\end{proof}

\begin{lemma}\label{cl:two-ric}
Let $\mu_1\le\cdots\le\mu_4$ be the eigenvalues of the ambient Ricci
tensor on a fixed normalized neutral end. Then there is a positive constant $b$
depending on the end such that
\[
          \mu_1+\mu_2\ge b/(s\log s)
\]
uniformly at high levels.
\end{lemma}
\begin{proof}
The meridional sectional curvature of the level is given by
$K_m=-r_{zz}/r=(r_z^2-y_{zz})/r^2$. On
$|\xi|\le\sqrt u/4$, by Lemma~\ref{cl:strip}, we have
$y_{zz}=w''=-j+O(\xi^2/u^3)\le-c/u$. At the edge of this band, we have
$|r_z|\ge c/\sqrt u$ by the intermediate
profile and \eqref{e:pslope}. By concavity,
this lower bound still holds towards each pole. Outside the fixed Bryant
core of Lemma~\ref{lem:mode-second}, the same lemma gives
$y_{zz}\le C(\log s)^2/s=o(1/u)$. Inside the fixed core, by the smooth
tip-scale convergence of Proposition~\ref{prop:mode-regions} and the
strict positivity on the compact Bryant core, we have $K_m\ge c u/t$,
and this also holds at the smooth pole. Since $r^2\le2t$, we have proved that
$K_m\ge c/(tu)$ on the whole level. The three intrinsic Ricci
eigenvalues of the level are $2K_m,K_m+K_o,K_m+K_o$, where $K_o>0$ denotes the
orbital sectional curvature. By Lemma~\ref{lem:level-error}, the
tangential compression of the ambient Ricci tensor is at least
$q_0\Ric_\gamma$. Hence, by eigenvalue interlacing, we get
$\mu_2\ge q_0 K_m$. The lemma now follows from $q_0\to1$, $\mu_1>0$ and $s\sim t$.
\end{proof}

\section{Global rigidity and the classification}\label{cl:rigidity}

Throughout this section, we denote the fixed pair by $(M,g,f)$ and
$(\widehat M,h,\hat f)$. After interchanging the two ends if necessary, we may
assume that the constant $c$ in Proposition~\ref{cl:end} is nonnegative.
Note that the inverse identification is asymptotically isometric and the two
heights differ by a bounded amount. Hence, the inverse identification has the
same decay.

\begin{lemma}\label{cl:initialmap}
There is a smooth proper map $F_0:M\to\widehat M$ of degree
$\pm1$ such that $\hat f\circ F_0-f\to c$ uniformly at infinity.
Moreover, for every orthonormal pair $e_1,e_2$ in the source, we have
\[
        |dF_0(e_1)|_h^2+|dF_0(e_2)|_h^2\le2 .
\]
\end{lemma}
\begin{proof}
We first extend $\Phi$ from the end to a smooth map defined over the compact core.
For example, using the global smooth $\mathbb R^4$ coordinates of the target
from Section~\ref{sec:global}, we can cut the map off to a constant
on an interior annulus. The extension is proper and has degree $\pm1$, since
a regular target point which is sufficiently far away has a unique preimage on the
end. Let $\varphi_v$ be the flow of $-\nabla\hat f$. This flow is complete, and
the vector field has norm at most one. We also set $s=f$.

We first consider a frozen time $T\ge0$. For an orthonormal pair at height $s$,
the sum of the squared lengths of its images under $\Phi$ is initially
$\le2+Cs^{-11/10}$. Along the flow $\varphi_v$, the derivative of this sum is
$-2\sum_{j=1}^2\Ric_h(d(\varphi_v\Phi)e_j,d(\varphi_v\Phi)e_j)$.
We claim that for $0\le v\le\min(T,1)$ the Gram matrix of the two image vectors
is bounded below by a fixed positive multiple of the identity. In fact, this holds
at $v=0$ by Proposition~\ref{cl:end}, and the Ricci tensor of the target is
bounded. Since $|d\hat f/dv|\le1$, their heights differ from $s$ by $O(1)$.
By orthogonal diagonalization, the trace of the Ricci tensor on the plane
spanned by them is at least the sum of the two smallest eigenvalues of the
Ricci tensor. Hence, by Lemma~\ref{cl:two-ric} and the positivity at all later
times, we get
\begin{equation}
 \sum_{j=1}^2|d(\varphi_T\Phi)e_j|^2
 \le 2+Cs^{-11/10}-b_1\frac{\min(T,1)}{s\log s}.
 \label{cl:frozen}
\end{equation}

Next we let the time depend on $s$ and define $F_0=\varphi_{T(s)}\Phi$. By the
chain rule, besides the frozen differential we get the additional term
$-\nabla\hat f(F_0)T'(s)\,df$.
Since the frozen differential is bounded and $|df|,|d\hat f|\le1$, when
$|T'|\le1$ the contribution of this term to the left side of
\eqref{cl:frozen} is at most $C_1|T'|$, and this bound is independent of $T$.
Let $\chi$ be a smooth cutoff function which is equal to zero for $s\le S$,
equal to one for $s\ge2S$, and takes values between zero and one. We start with
$T=H\ge1$ on $s\le S$ and solve
\[
 T'=-\epsilon\frac{\chi(s)}{s\log s}\frac{T}{1+T}
 \qquad (s>S).
\]
Integrating $(T+\log T)'$, we see that $T>0$ and $T\to0$.
Moreover, for $s\ge S$ we have
\[
 |T'|\le\frac{\epsilon\min(1,T)}{s\log s},\qquad
 \min(1,T)\ge\left(\frac{\log S}{\log s}\right)^\epsilon .
\]
Here the second bound is obtained by integrating
$(\log T)'\ge-\epsilon/(s\log s)$ and using $H\ge1$.
We first choose $C_1\epsilon<b_1/3$, and then we choose
$S$ so large that the positive error in \eqref{cl:frozen} is at most a third of
its negative term. Since $s^{-1/10}(\log s)^{1+\epsilon}\to0$, this can be done
uniformly in $H$. Finally, the inward flow of
$\Phi(\{s\le S\})$ stays in a compact sublevel set of the target, on which
$\Ric_h\ge m h$ for some $m>0$. Thus, we have $|d\varphi_H|\le e^{-mH}$
there. By choosing $H$ large, we get the desired inequality on the core as well.
The maps $\varphi_{\theta T(s)}\Phi$, $0\le\theta\le1$, form a proper
homotopy, because their displacements are bounded by $H$. Moreover, the
potential changes at infinity by at most $T(s)\to0$, which is what we need.
\end{proof}

\begin{proposition}\label{cl:B}
Suppose that the normalized reflected neutral ends have equal phase. Then there
exists a global isometry $F:(M,g)\to(\widehat M,h)$ such that $\hat f\circ F=f$.
\end{proposition}
\begin{proof}
We start from the map given by Lemma~\ref{cl:initialmap} and consider the equation
\begin{equation}
 \partial_v F=\tr_g\nabla dF+dF(\nabla f)-\nabla\hat f(F).
 \label{cl:hmhf}
\end{equation}
The short-time existence of a solution with bounded differential follows from
\cite[Appendix A, Lemmas~A.23--A.24]{clBK}. To see this, let $\psi_g,\psi_h$
be the flows of the negative potential gradients of the two solitons. Then
$\psi_h(v)^{-1}\circ F_v\circ\psi_g(v)$ satisfies the harmonic map heat
flow between the canonical Ricci flows $\psi_g(v)^*g,\psi_h(v)^*h$.
Both of these flows are complete, and by \eqref{e:der} their curvature and its
fixed derivatives are uniformly bounded. Moreover, the time derivatives of
their metrics are equal to minus twice the Ricci tensor. The initial map is
smooth and globally Lipschitz. Note that these bounds remain uniform when the
starting time is shifted.

Next we prove an estimate for the differential, which allows us to continue the
solution for all time.
Let $B=F^*h$, and let $a_1\ge\cdots\ge a_4\ge0$ be its ordered eigenvalues
relative to $g$. At a given point, we take a singular frame, so that
$dF(e_i)=\sqrt{a_i}\,\hat e_i$. Here we complete a target orthonormal frame
also in the case when the rank drops. We write $K_{ik}=\Sec_h(\hat e_i,\hat e_k)$.
With respect to the tensor operator $\partial_v-\Delta-\nabla_{\nabla f}$,
the diagonal Bochner reaction for $B_{ii}$ is
$2\sum_{k\ne i}a_i(a_k-1)K_{ik}$, and there is an additional term
$-2|\nabla dF(\cdot,e_i)|^2$. In fact, since $\nabla^2 f=\Ric_g$, the domain
curvature term $-2B(\Ric_g^\sharp e_i,e_i)$ cancels with the derivative of the
source drift. The target potential contributes
$-2a_i\Ric_h(\hat e_i,\hat e_i)$, and by commuting the two derivatives of
$dF$ we get $2\sum_k a_i a_k K_{ik}$. When $a_1+a_2=2$, the sum of the
first two reactions divided by two is
\[
 -\tfrac12(a_1-a_2)^2K_{12}
  +\sum_{i=1,2;\ k=3,4}a_i(a_k-1)K_{ik}\le0 ,
\]
since $a_3,a_4\le1$. Note that the set of symmetric forms whose trace on every
orthonormal two-frame is at most two is closed and convex. This set is preserved
by the tensor maximum principle. To check this, at a supporting plane we extend
its frame by parallel transport to first order and test the trace.
Alternatively, we can use the equivalent null-vector test on $\Lambda^2$ for the
additive action of $B^\sharp$ on two-forms. On a short time interval on which the
differential is bounded, the reaction in this test is uniformly Lipschitz in the
squares of the singular values. More explicitly, if $a_1+a_2=2+d$ with $d>0$, we
replace the four $a_i$ by $2a_i/(2+d)$ in the reaction. Since the curvature and
the $a_i$ are bounded, the error is at most $C d$. This also gives the
barrier which we need on the noncompact manifold. Namely, we subtract
$\zeta e^{Cv}(1+f)$ from the two-frame trace bound, where we increase $C$ so
that it dominates that Lipschitz constant and
$(\Delta+\nabla f\cdot\nabla)(1+f)=1$. Then the first contact can only occur in
a compact set, and afterwards we let $\zeta\downarrow0$.
It follows that $a_1+a_2\le2$ and $|dF|_{\rm op}\le\sqrt2$. Applying the
short-time existence theorem cited above again and again, we obtain a solution
for every $v\ge0$ which satisfies the same estimate. Moreover, the local
derivative estimates of this solution at positive time are uniform on unit time
windows away from their initial slices.

Let $L=\Delta_g+\nabla f\cdot\nabla$ and
$P=\hat f\circ F-f$. Since $L f=1$, applying the scalar chain rule to
\eqref{cl:hmhf}, we get the exact equation
\begin{equation}
 (\partial_v-L)P
   =R_h(F)-\tr_g F^*\Ric_h
   =\sum_{i<k}(2-a_i-a_k)K_{ik}\ge0 .                 \label{cl:scalar}
\end{equation}
On the other hand, for $0<\beta<1$ we have
$L e^{-\beta f}
=(-\beta+\beta^2|\nabla f|^2)e^{-\beta f}
\le-\beta(1-\beta)e^{-\beta f}$.
By the uniform potential limit, for any $\eta>0$ the initial lower bound
$P(0,\cdot)\ge c-\eta-C_\eta e^{-\beta f}$ holds. Then the scalar comparison
principle gives
\begin{equation}
 P(v,x)\ge c-\eta-
 C_\eta e^{-\beta f(x)-\beta(1-\beta)v}.              \label{cl:lowerP}
\end{equation}
On the noncompact domain, in this comparison we may add
$\zeta e^{5v}(1+f)^2$ to the supersolution. Here we use
$L(1+f)^2=2(1+f)+2|\nabla f|^2$ and $P\ge-f$.
We first let the compact truncation exhaust $M$, and then we let $\zeta\downarrow0$.

By the inequality \eqref{cl:lowerP}, the homotopy is proper on each finite
time interval. This is because a uniform upper bound for $\hat f(F_v(x))$
then implies a uniform bound for $f(x)$. So we have $\deg F_v=\deg F_0\ne0$,
and every $F_v$ takes the value $\hat o$, which is the minimum point of the
target. We choose $x_v\in F_v^{-1}(\hat o)$. At this point we have $P=-f(x_v)$,
and hence $\limsup_{v\to\infty}f(x_v)\le\eta-c$ for all $\eta>0$. Note that $c\ge0$. Then
we must have $c=0$, and $x_v\to o$ by the properness and the uniqueness of the
minimum of the source. It follows from the Lipschitz bound that $F_v(o)\to\hat o$.

Now we translate the time by an arbitrary sequence $v_i\to\infty$. By the last
conclusion and the Lipschitz estimate, we know that on every compact space-time
set the images stay in a common compact set. Then we use the interior parabolic
estimates for \eqref{cl:hmhf}, or alternatively \cite[Appendix A]{clBK} in the
moving metrics, and we take a diagonal subsequence to get a smooth eternal limit.
It follows from \eqref{cl:lowerP} that the potential difference of this limit is nonnegative, and
it vanishes at $(o,0)$. Then by the strong parabolic minimum principle, it is
identically zero on every earlier slice. Note that all $K_{ik}>0$ in
\eqref{cl:scalar}, and so we have $a_i+a_k=2$ for every pair. In dimension four,
we subtract these identities and get $a_i=1$ for all $i$. This means that a fixed
earlier slice is a local isometry. Recall that a local isometry from a complete
connected source is a covering of the connected target, which follows by
continuation of lifts of geodesics. We also know from
Section~\ref{sec:global} that the target is simply connected, and so this covering is a global isometry. Finally, since the
potential difference is zero, we get the identity of potentials stated in the
proposition.
\end{proof}

\begin{proof}[Proof of Theorem~\ref{cl:main}]
We fix the scale and the additive constant of the potential as in
Section~\ref{sec:intro}. If the soliton is rotationally symmetric, then it is
the Bryant soliton. Let $g$ be any other soliton in the theorem. Then by
Theorem~\ref{th:ref} and
Propositions~\ref{prop:mode-selection}--\ref{prop:mode-regions}, we know that it has a
reflected neutral end, and its finite phase $\theta$ is defined by
Proposition~\ref{cl:invariant}. By Proposition~\ref{cl:A}, there is a fixed
construction limit $h$ with $[h]\in\mathcal L$ and
$\Theta(h)=\theta$. Then by Proposition~\ref{cl:B}, we see that $g$ is identified with this
$h$ by an isometry which preserves the normalized potentials. Hence, every
soliton in the class is one of the models in the statement.

It remains to identify the parameter on $\mathcal L$. We know from
Proposition~\ref{cl:invariant} that isometric representatives have the same phase.
We also know from Proposition~\ref{cl:A} that every real value is attained, and from
Proposition~\ref{cl:B} that any two representatives with the same value are
isometric. So the map $[h]\mapsto\Theta(h)$ is well defined,
surjective and injective. The proof of the theorem is complete.
\end{proof}

\section{Classification under nonnegative sectional curvature}\label{sec:nonnegative}

 If the Ricci tensor has a null direction, then the soliton splits globally.
On the non-Bryant soliton with $\Ric>0$, we get a fixed axis from the angular
action of Proposition~\ref{app:angular-action}.
Then by using the isotropy at the axis, we will show that the curvature operator is
positive, and so we can apply Theorem~\ref{cl:main}.

We first state a dimension reduction result.

\begin{proposition}\label{secNN:models}
Let $(M^4,g,f)$ be a complete $\kappa$-noncollapsed steady gradient Ricci soliton
with nonnegative sectional curvature.
Suppose that the Ricci tensor has a null direction at some point. Then, up to
scaling, the metric is either Euclidean or the product of a three-dimensional Bryant
soliton with a line. Suppose instead that the Ricci curvature is positive. Then
either the soliton is the four-dimensional Bryant soliton, or, after passing to
subsequences, the curvature-normalized backward
limits at infinity are shrinking cylinders $S^2\times\mathbb R^2$ or
Bryant products $\mathrm{Bry}^3\times\mathbb R$.
\end{proposition}

\begin{proof}
If the soliton is nonflat, we scale the metric so that $R+|\nabla f|^2=1$. Since
the sectional curvature is nonnegative, we get $|\Rm|\le C_4 R$  and hence the curvature
is bounded. Note that all time slices of the canonical flow are isometric.
Then the spatial noncollapsing on the final time slice gives the parabolic
noncollapsing under a parabolic curvature bound. In the case of a scalar curvature
bound, we use $|\Rm|\le C_4 R$ and reduce the radius of the ball by a
dimensional factor.

By the result of Hao and the author~\cite[Theorem~1.3]{DH}, a complete
$\kappa$-noncollapsed steady gradient Ricci soliton with nonnegative
sectional curvature is diffeomorphic to Euclidean space. In particular it is simply connected.
By the theorem of Guan, Lu and Xu~\cite[Theorem~1.1(I)]{GLX}, the Ricci tensor
has constant rank when the sectional curvature is nonnegative. If the Ricci
tensor has a null vector at one point, then the universal cover splits
isometrically with a Euclidean factor of dimension
$\dim\ker\operatorname{Ric}$. Since $M$
is simply connected, we get a splitting of $M$ itself. Splitting off a
line, we have $M=N^3\times\mathbb{R}$. From the line and mixed equations in
$\operatorname{Hess} f=\operatorname{Ric}$, we get $f(x,z)=f_N(x)+az+b$ and
$\operatorname{Hess}_N f_N=\operatorname{Ric}_N$. Suppose that
$|\operatorname{Rm}_N|\le r^{-2}$ on $B_N(x,r)$. Then the product ball centered
at $(x,0)$ satisfies the same bound, and it is contained in
$B_N(x,r)\times(-r,r)$. Hence, we have
$\operatorname{Vol}_N B_N(x,r)\ge(\kappa/2)r^3$. If $N$ is nonflat, then by
Brendle's classification \cite{Bre13} we have $N=\mathrm{Bry}^3$ up to scale.
In the flat case, by completeness and simple connectivity we get
$\mathbb{R}^4$. Finally, two potentials differ by a function with parallel
gradient. Thus, the potentials on the Bryant product are
$f_{\mathrm{Bry}^3}+az+b$, and the potentials in the flat case are affine.

Next we assume that $\operatorname{Ric}>0$. By the results of Chan--Ma--Zhang
\cite[Theorems 1.2, 1.3, 1.5 and Corollary 5.4]{CMZ}, we have the following
alternatives for the limits at infinity. If a tangent flow at infinity is
$S^3\times\mathbb{R}$, then the soliton is $\mathrm{Bry}^4$. Otherwise, the
tangent flow is $S^2\times\mathbb{R}^2$. For every escaping sequence, after
rescaling by the scalar curvature at the basepoint and passing to a subsequence,
the backward limit is a shrinking $S^2\times\mathbb{R}^2$ or
$\mathrm{Bry}^3\times\mathbb{R}$. Note that in the latter case the basepoint
need not be a tip. We remark that these results can be applied before a
critical point is established.
\end{proof}

The ray blow-down and the uniform scalar decay apply on the branch with positive curvature operator. Here we only have nonnegative sectional curvature and positive Ricci curvature. The next proposition deals with the case of positive Ricci curvature.

\begin{proposition}\label{secNN:axis}
Consider a complete connected four-dimensional steady soliton with bounded curvature, $\sec\ge0$ and $\operatorname{Ric}>0$. Suppose that it admits an isometric $O(3)$ action which is equivariantly diffeomorphic to the standard
action on $\mathbb R\oplus\mathbb R^3$, trivial on the first factor. Then its curvature operator is positive.
\end{proposition}
\begin{proof}
On the regular set, we write $g=h+u^2g_{S^2}$. We choose an orthonormal horizontal
eigenframe of $\operatorname{Hess}_h u$ with eigenvalues $\mu_1,\mu_2$, together
with an arbitrary vertical orthonormal frame. By the warped-product formulas, the
curvature operator is diagonal with respect to the six coordinate two-forms, and
its eigenvalues are
\[
 K_h,\quad \frac{1-|\nabla^h u|^2}{u^2},\quad
 -\frac{\mu_1}{u},\quad -\frac{\mu_1}{u},\quad
 -\frac{\mu_2}{u},\quad -\frac{\mu_2}{u}.
\]
Since these basis elements are simple, all the eigenvalues are sectional
curvatures, and hence they are nonnegative. By continuity, we also have
$\mathcal R\ge0$ on the axis, and by pullback the same holds at all times of the
canonical flow.

At a point on the axis, we write $T_s M=\mathbb{R} e_0\oplus V$, where
$V=\mathbb{R}^3$. By the isotropy, the curvature operator acts as a scalar on
each of $e_0\wedge V$ and $\Lambda^2V$. This is because the normal inversion has
opposite signs on the two summands and the $SO(3)$ representations are
irreducible. The first eigenvalue $a$ satisfies
$3a=\operatorname{Ric}(e_0,e_0)>0$. Therefore, the image of the curvature
operator contains $e_0\wedge e_i$, $1\le i\le3$. By Hamilton's tensor strong
maximum principle for a nonnegative curvature operator \cite[Section 8]{Ham},
the image has constant spatial rank and is closed under the Lie bracket at
positive elapsed time on a steady flow. For Hamilton's strong maximum principle
for the curvature operator,
we also refer to \cite[Theorem~6.60]{CLN}. Note that every time slice of the
canonical flow is isometric to the original soliton. In
$\Lambda^2\simeq\mathfrak{so}(4)$, we have
\[
 [e_0\wedge e_i,e_0\wedge e_j]=\pm e_i\wedge e_j\qquad(i\ne j).
\]
Thus, the image contains a basis at the axis point. Since the spatial rank is
constant, the rank is six everywhere, which proves the positivity.
\end{proof}

\Needspace{6\baselineskip}
\noindent\hypertarget{proof:sec-nonnegative}{}Now we are ready to prove
Theorem~\ref{th:classification-sec-nonnegative}.

\begin{proof}[Proof of Theorem~\ref{th:classification-sec-nonnegative}]
Since the sectional curvature is nonnegative, we have $\Ric\ge0$. If $\Ric$ has a null direction
at some point, then by Proposition~\ref{secNN:models}, $(M,g)$ is isometric up to scale
to $\mathbb R^4$ or $\mathrm{Bry}^3\times\mathbb R$.

Now we suppose that $\Ric>0$ everywhere. In the non-Bryant case, by
Proposition~\ref{app:angular-action}, we have the standard angular $O(3)$
action with a fixed axis. The curvature is bounded as in the proof of
Proposition~\ref{secNN:models}. Hence, by Proposition~\ref{secNN:axis}, the
curvature operator is positive everywhere. Then, by Theorem~\ref{cl:main},
$(M,g)$ is isometric up to scale to a member of $\mathcal L$, and we also get
the phase bijection $\Theta:\mathcal L\to\mathbb R$ stated in the theorem.
\end{proof}

\appendix
\section{Finite derivatives of the Bryant profile}\label{app:bryant-profile}

\begin{lemma}\label{lem:bryant-derivatives}
Consider the Bryant soliton with metric $ds^2+a(s)^2g_{S^m}$, $m=2$ or $3$,
and let $x=a(s)$ and $\Phi(x)=a'(s)^2$. Let $c_{2\ell}$ denote the
coefficients of its expansion at infinity given in
\cite[Appendix~B, Proposition~1(4)]{AIK}. Then for every fixed $N\ge1$
and every integer $j\ge0$, we have
\[
 \frac{d^j}{dx^j}\left(\Phi(x)-\sum_{\ell=1}^{N}c_{2\ell}x^{-2\ell}\right)
                  =O_{N,j}(x^{-2N-2-j}).
\]
\end{lemma}
\begin{proof}
For $j=0$, the lemma follows from the expansion cited above at every truncation
order, and we have $c_2>0$. To prove the derivative bounds, we use the radial
equation in part~(1) of the same proposition, namely
\[
 \Phi\Phi''-\tfrac12(\Phi')^2
        +(m-1-\Phi)\frac{\Phi'}{x}
        +2(m-1)\frac{(1-\Phi)\Phi}{x^2}=0 .
\]
We first show that every fixed derivative of $\Phi$ grows at most polynomially.
Note that $\Phi\asymp x^{-2}$. Moreover, the radial sectional curvature
$K_{\rm rad}=-a''/a$ is bounded, since the Bryant metric has positive
sectional curvature and bounded scalar curvature. Since $a'>0$, we have
$\Phi'=2a''=-2xK_{\rm rad}=O(x)$. Then, solving the radial equation above for
$\Phi''$ and using $\Phi^{-1}=O(x^2)$, we get a polynomial bound.
Differentiating the equation repeatedly, we obtain the same conclusion for the
higher derivatives by induction. In fact, each new term is a product of
derivatives which are already bounded, powers of $1/x$, and powers of $1/\Phi$.

Now we fix $j\ge1$ and set $D=2N+2+j$. For a smooth function $u$, we write
$\Delta_h u(t)=u(t+h)-u(t)$. Applying the fundamental theorem of calculus
$j$ times, we get
\[
 h^{-j}\Delta_h^j u(x)=
 \int_{[0,1]^j}u^{(j)}\bigl(x+h(t_1+\cdots+t_j)\bigr)\,dt_1\cdots dt_j .
\]
If $jh\le1$, this quantity differs from $u^{(j)}(x)$ by at most
$jh\sup_{[x,x+1]}|u^{(j+1)}|$. We choose $h=x^{-p}$, where $p$ is so large
that this error for $\Phi$ is $O(x^{-D})$. Here we use the polynomial
bound proved above. We also require that $p\ge D$. Then, for any fixed truncation
$P_K(x)=\sum_{\ell=1}^K c_{2\ell}x^{-2\ell}$, the same difference error
is $O_K(x^{-D})$, since the derivatives of the truncation are bounded on
$[1,\infty)$.

Next we choose $K\ge N$ so large that $2K+2-pj\ge D$. At each of the $j+1$
sampling points, the function-value remainder $\Phi-P_K$ is
$O_K(x^{-2K-2})$. Since the sum of the absolute values of the forward-difference
coefficients is $2^j$, we have
$h^{-j}|\Delta_h^j(\Phi-P_K)(x)|=O_K(x^{-2K-2+pj})=O_K(x^{-D})$.
Combining the three errors, we obtain $\Phi^{(j)}-P_K^{(j)}=O(x^{-D})$.
On the other hand, the $j$-th derivative of each term of $P_K-P_N$ is of order
at most $x^{-2N-2-j}=x^{-D}$. This proves the lemma for the original truncation.
\end{proof}

\section{Angular symmetry under nonnegative sectional curvature}\label{app:fang}

In this appendix, we construct the standard isometric $O(3)$ action on the
non-Bryant solitons of Theorem~\ref{th:classification-sec-nonnegative}
with $\sec\ge0$ and $\Ric>0$. When the curvature operator is positive, the
angular action was constructed by Fang~\cite[Theorem~1.1]{F}. However, in the proof
of our classification, the positivity of the curvature operator is
proved in Proposition~\ref{secNN:axis} by using the fixed axis of
this action. So we need to follow Fang's construction of angular Killing
fields under the weaker hypotheses. In the proof, we supply the convexity and
vanishing arguments, and we also keep track of the domains on which the local fields
are compared. The action and its fixed axis are given in
Proposition~\ref{app:angular-action}. When we refer to Fang's
results, we use the numbering of~\cite{F}.

\subsection{An auxiliary exhaustion and an integral vanishing lemma}

On the nonflat branch of Proposition~\ref{secNN:models}, we normalize the
soliton by $R+|\nabla f|^2=1$ and we write $X=\nabla f$. On the canonical flow
we have the curvature bound $|\Rm|\le C_4R\le C_4$, and then by Shi's estimates we get
uniform derivative bounds. Moreover, by the noncollapsing condition and the
Cheeger--Gromov--Taylor estimate, we have the injectivity radius bound which
is required in Fang's Lemma~2.1. We will use this bounded local geometry below.

In \cite[Section~6]{F}, Fang uses positive sectional curvature to get a
strictly convex exhaustion, and he glues the local fields on the levels of
this exhaustion. In our case, properness comes from the Cheeger--Gromoll
function, and we get strict convexity from
$\nabla^2 f=\Ric>0$ after we add
$\log(1+e^f)$ to that function. Note that we can carry out this construction
before the properness of $f$ has been established.

\begin{lemma}\label{secNN:exhaustion}
On a complete noncompact normalized steady soliton with $\sec\ge0$ and $\operatorname{Ric}>0$,
there exists a smooth proper Lipschitz function $\psi$ whose Hessian is positive definite.
\end{lemma}
\begin{proof}
We fix a point $o\in M$ and consider the Cheeger--Gromoll function
\[
 b(x)=\sup_{\gamma(0)=o}\lim_{t\to\infty}\big(t-d(x,\gamma(t))\big),
\]
where the supremum is taken over unit-speed rays. It is known that this function is
convex and $1$-Lipschitz, see \cite{CGAnn}. We now show that it is proper. Let $x_i$ be any
escaping sequence. After passing to a subsequence, we may assume that the initial
directions of minimizing segments from $o$ converge to the initial direction of a ray
$\gamma$. Let $r_i=d(o,x_i)$ and let $\theta_i$ be the angle with that ray. By hinge
comparison, we have $d(x_i,\gamma(t))^2\le r_i^2+t^2-2r_it\cos\theta_i$. For fixed $i$,
we let $t\to\infty$ and get $b(x_i)\ge r_i\cos\theta_i\to\infty$.

Next we let $\chi(s)=\log(1+e^s)$ and $v=b+\chi(f)$. Since $\chi\ge0$ and $0<\chi'<1$,
the function $v$ is proper and $2$-Lipschitz. Moreover, the Hessian
\[
 \operatorname{Hess}(\chi(f))=\chi'(f)\operatorname{Ric}+\chi''(f)\,df\otimes df
\]
has a positive lower bound relative to $g$ on every compact set. So we see that $v$ is strictly
convex in the local quantitative sense of Greene--Wu. Then by their approximation theorem \cite[Theorem~2(a)]{GW}, there is a smooth strictly convex function $\psi$ such that $|\psi-v|<1$
and such that the Lipschitz constant of $\psi$ is at most three. It follows that the function $\psi$ is
proper and it has a unique nondegenerate minimum. Finally, we subtract the minimum value
from $\psi$.
\end{proof}

In the proof of \cite[Proposition~7.3]{F}, Fang uses positive sectional curvature
to get $R-2\lambda_{\max}(\Ric)>0$. In our case the sectional curvatures are
allowed to vanish, and we keep the eigenvalues of $h$ in the curvature term.
Then the strictness which we need for the integral argument comes from the
equality case in \eqref{secNN:gap} and $\Ric>0$.

\begin{lemma}\label{secNN:vanishing}
Let $n\ge3$. Suppose that $\sec\ge0$ and $\operatorname{Ric}>0$ on a complete normalized
steady soliton with bounded local geometry as above, and let $X=\nabla f$. Let
$h\in L^1(d\mu_g)$ be a smooth symmetric two-tensor such that
$(\Delta+\nabla_X)h+2\operatorname{Rm}(h)=0$. Then $h=0$.
\end{lemma}
\begin{proof}
We choose an orthonormal frame which diagonalizes $h$. Let $h_i$ be the eigenvalues of
$h$ and let $K_{ij}$ be the sectional curvatures. With Fang's convention, we have
$\langle\operatorname{Rm}(h),h\rangle=2\sum_{i<j}K_{ij}h_ih_j$, and hence
\begin{equation}\label{secNN:gap}
 R|h|^2-2\langle\operatorname{Rm}(h),h\rangle
 =2\sum_{i<j}K_{ij}\left((h_i-h_j)^2+\sum_{k\ne i,j}h_k^2\right).
\end{equation}
Suppose that equality holds. Then for any positive $K_{ij}$ we must have $h_i=h_j$ and
$h_k=0$ for $k\notin\{i,j\}$. For an index $k$ outside this pair, since
$\operatorname{Ric}_{kk}>0$, there is a second positive pair containing $k$. By the
equality condition for this second pair, the common value must vanish. So we see that the
left side is strictly positive for every nonzero $h$.

Let $u=|h|$. We define $G$ to be the left side of \eqref{secNN:gap} divided by $u$ when
$u>0$, and we set it to be zero elsewhere. Then $G$ is nonnegative and continuous.
We first apply Kato's inequality to $(u^2+\varepsilon^2)^{1/2}$, and then we get
$(\Delta+X)u+Ru\ge G$ in the distributional sense. We take smooth compact cutoffs
$0\le\eta_A\le1$ which are equal to one on exhausting balls and whose gradient and
Laplacian are uniformly bounded. We can construct such cutoffs from smoothed distance as
in Fang's Proposition 7.3. Then by $\operatorname{div}X=R$, we have
\[
 0\le\int\eta_A G\,d\mu_g
 \le\int u(\Delta\eta_A-X\eta_A)\,d\mu_g\longrightarrow0
\]
by integrability and $|X|\le1$. So we have proved the claim, and we note that we do not need a uniform
spectral gap at infinity.
\end{proof}

\subsection{Necks and caps}
In this subsection and in the next two subsections, we work on the non-Bryant branch
with $\Ric>0$. By Proposition~\ref{secNN:models}, we have the same comparison models as in
\cite{F}. They are the shrinking cylinder $S^2\times\mathbb R^2$ and
the Bryant product $\mathrm{Bry}^3\times\mathbb R$. The estimates below on necks
and on cap collars are obtained by smooth comparison with these
models. At the tips, we use $\Ric>0$ to orient the tip arcs.

\paragraph{Necks and tilted caps.}
The scale-normalized potential satisfies Fang's identities (2.5). On a tip-normalized
Bryant product, the time-zero limit of this potential is
$\lambda f_{\mathrm{Bry}^3}+az+c$, $\lambda^2+a^2=1$. In fact, after subtracting
$\lambda f_{\mathrm{Bry}^3}$ the gradient is parallel, and the normalization is obtained
by evaluating at the tip. Changing the sign of $z$, we may take $a\ge0$. For cylindrical
charts, by integrating the spherical Hessian identity we get $\lambda=O(\eta)$ and an
affine limiting flat direction as in Fang's Lemma 3.5. In Proposition 4.3, the positive
quadratic estimate is computed on the Bryant product, and then it is transferred under a
quantitative bound on the difference of curvature tensors. By Corollary 4.6, every tilt
is allowed on collars away from the tip, including $a=0$. All of these are conditions on
the fixed comparison charts, and all of them are supplied by the preceding models.

\paragraph{Locating tips.}
At the normalized Bryant tip, the three-dimensional Ricci block is $I/3$ and the Bryant
block of $\operatorname{Hess} R$ is $-2I/9$. These are the spectral gap and the transverse
inverse which appear in Proposition 5.7. On the actual intrinsic tip set, its spectral
projection satisfies $0=P_E\nabla R=-2\operatorname{Ric}^\sharp P_E X$, and hence
$P_E X=0$. Also, in Corollary 5.8 the smallest Ricci eigenvalue is positive. Since $f$
is strictly convex, it has at most one critical point. If this critical point exists, we
may include it in the compact exclusion. This verifies the strictness which is needed to
orient the tip arcs. For the tube, projection and common CMC leaf constructions, we use
smooth model convergence and the local inverses. In particular, the local maxima of
scalar curvature on a regular level in this construction do not require a global maximum
of scalar curvature.

\subsection{Local symmetry estimates on necks and caps}

In the neck step of Proposition~\ref{secNN:fang521}, we apply Fang's
Corollary~3.15 with gain $\theta=1/2$. After the terminal domain is
fixed, we choose $2T=4^J$, where $J$ is so large that the bounds on
$T$ in that corollary are satisfied. Then we choose the cylindrical threshold.
On this mesh, the rotation fields are constructed by Fang from the local
soliton identities, smooth cylinder comparison, and the assumed
symmetry bounds on the input neighborhoods, see
\cite[Theorem~3.4, Lemma~3.8 and Corollary~3.15]{F}.

The limits in Proposition~\ref{secNN:models} give us comparison on
backward time intervals. On the other hand, the cylindrical neighborhoods in
\cite[Proposition~5.10]{F} contain positive normalized times as well.
We will continue the comparison on an exhaustion of the complete model
through the final time in Fang's definition.

\begin{proposition}[Fang, Proposition~5.10]\label{secNN:fang510}
Let \(0<\eta_{\mathrm{cyl}}<1/2\) be fixed, where cylindrical neighborhoods are understood
as in Definition~3.1. Then there are \(b_0<\infty\) and a compact set \(K\) with the
following property. Every \(x\notin K\) which is not \(\eta_{\mathrm{cyl}}\)-cylindrical
is contained in the set \(\mathcal T_0(R_{\rm tube})\) defined in (5.1), and it satisfies
\[
 \beta_0(x)<b_0 .
\]
We choose the tube radius such that it contains this radial range with a fixed margin. Moreover, for
every fixed \(b_*>b_0\), all points outside \(K\) and outside the core \(\{\beta_0<b_*\}\)
are \(\eta_{\mathrm{cyl}}\)-cylindrical.
\end{proposition}
\begin{proof}
Let $Q_i=R(x_i)$ be the rescaling factors in
Proposition~\ref{secNN:models}. We also need the comparison at
positive normalized times. By the scale-invariant estimates of
\cite[Remark~A.3 and Corollary~A.9]{CMZ}, we have
$|\nabla^j\Rm_i|\le C_jR_i^{1+j/2}$. Then by the evolution of the scalar
curvature, we get $|\partial_t R_i|\le AR_i^2$. On rescalings of the exact
Bryant product, we know that the same bounds hold by its differentiated asymptotics.

Suppose that the convergence holds through time $u$, and that the complete model
satisfies $R<K$ on $[u,b]$ for a fixed $K\ge1$. Then on each compact subset
of the model we have $R_i(u)\le2K$ when $i$ is large. We integrate at each fixed
spatial point, and we get $R_i\le4K$ through time
$\min\{b,u+(4AK)^{-1}\}$. Note that this time step does not depend on the
compact set. By the evolution of the metric, we have uniform equivalence to the
metric at $u$. Moreover, by the curvature derivative bounds and the evolution of
the connection, we have smooth compactness in the same spatial charts.
We see that a diagonal limit on an exhaustion is complete by metric equivalence, and it
has bounded curvature. Then it follows from the uniqueness of complete solutions with
bounded curvature (Chen--Zhu, \emph{Uniqueness of the Ricci flow on complete noncompact
manifolds}, Theorem~1.1) that this limit coincides with the prescribed model.
We can iterate this argument with the same $K$, and then we get the convergence through $b$.

For a cylinder centered at time zero, we take
$b=1-\eta_{\rm cyl}/2$, which is strictly beyond the last comparison time
$1-\eta_{\rm cyl}$. The scalar curvature of the cylinder is $(1-t)^{-1}$ there,
and the actual flows are eternal. Thus the cylindrical limit gives
a neck in the sense of Definition~3.1. The same continuation also applies
to points far out on the exact Bryant product. By (2.6)--(2.8), the backward limits
of these points at their own curvature scale are cylindrical.
We choose $b_0$ with a strict margin beyond the resulting model neck
threshold. In a Bryant-product limit, the remaining basepoints can only
lie at bounded scaled distance from the tip line. By Propositions~5.7
and~5.9, there is an actual tip at that distance with comparable curvature,
and the basepoint lies in the tube of this tip with $\beta_0<b_0$.
Since the curvatures are comparable, these tips escape together with the basepoints,
as in \cite[Proposition~5.10]{F}.
\end{proof}

In the finite chains of comparisons in \cite[Proposition~5.15]{F},
we need one Bryant-product chart which contains the projected tips and the images
of all the cylindrical charts. To this end, we bound those images in the scale of
the initial tip, and we apply the continuation argument above in order to cover
the full time interval of each chart.

\begin{proposition}[Fang, Proposition~5.15]\label{secNN:fang515}
Let \(\Lambda>0\), let \(J\ge0\) be an integer, and let
\(\mathcal F=\{(L_\ell,T_\ell)\}\) be a finite family of positive pairs. We use
\(\mathcal R^{\rm par}_{\mathcal F}\) from Definition~5.2, and we also allow an edge from
\((x,s)\) to \((\pi_t(y),t)\) whenever \((x,s)\mathcal R^{\rm par}_{\mathcal F}(y,t)\) and
\(y\in\mathcal{C}_t\).
\begin{enumerate}
\item[(i)] There exist \(D_J,T_J<\infty\) and a compact set \(K\) such that every center
which can be reached in at most \(J\) steps from \((p,0)\), \(p\in\Gamma\setminus K\), lies in
\[
 B_{g(0)}(p,D_Jr_p)\times[-T_Jr_p^2,0],\qquad r_p=R(p)^{-1/2}.
\]
\item[(ii)] For any fixed \(0<\eta_{\mathrm{cyl}}<1/2\), one can choose a single
Bryant-product comparison chart with spatial image \(O_p\) and time interval
\([-S_-r_p^2,S_+r_p^2]\) such that
\[
 B_{g(0)}(p,(D_J+1)r_p)\subset O_p,\qquad S_->T_J+1,\qquad S_+>0 .
\]
For each \(\eta_{\mathrm{cyl}}\)-cylindrical center \((x,s)\) in this family, we denote by
\(\Psi_{x,s}\) its map from Definition~3.1. Then the chart also contains
\begin{gather*}
 \Psi_{x,s}\big(B_{\bar g(-1)}(\bar x,\eta_{\mathrm{cyl}}^{-1})\big)\subset O_p,\\
 [s+(1-\eta_{\mathrm{cyl}}^{-1})r_{x,s}^2,\ s+r_{x,s}^2]\subset[-S_-r_p^2,S_+r_p^2],
\end{gather*}
with fixed spatial and time margins. In particular, the cylindrical comparison interval
which ends at \(s+(1-\eta_{\mathrm{cyl}})r_{x,s}^2>s\) is contained in the chart. Here the chart
sizes and the compact set are allowed to depend on \(J,\Lambda,\mathcal F,\eta_{\mathrm{cyl}}\), on
the prescribed comparison accuracy and on the derivative order, but they do not depend on
\(p\).
\end{enumerate}
\end{proposition}
\begin{proof}
We first prove (i). We apply the finite induction (5.4) of
\cite[Proposition~5.15]{F} on the tip-normalized Bryant product.
At each step, we know that the scalar minimum on the previous compact model region
is positive, and so the next neck radius and backward time are bounded.
Moreover, a projection adds at most the fixed radial width of the core.
By smooth convergence on a larger region, the resulting strict
margins carry over to the actual chart. To see this, note that $g(t)\ge g(0)$ for
$t\le0$, so a path of the permitted length cannot leave the larger ball. Also,
the tip projections converge by Propositions~5.9 and~5.13.

Next we prove (ii). Let $m>0$ be the minimum of the model scalar curvature
on a compact neighborhood of these centers. Then we have
$r_{x,s}/r_p\le\sqrt{2/m}$ for all sufficiently distant tips $p$.
By (3.3), the spatial map of each child chart maps radial minimizing paths
to paths whose $g(s)$-length is at most $C\eta_{\rm cyl}^{-1}r_{x,s}$.
The same bound holds for their $g(0)$-lengths. We enlarge the
parent radius in tip units by $C\eta_{\rm cyl}^{-1}\sqrt{2/m}+2$,
and we enlarge its time interval in both directions by
$2(1+\eta_{\rm cyl}^{-1})/m+2$ in parent time units.
On this finite interval, the continuation argument in Proposition~\ref{secNN:fang510}
can be applied to the complete eternal Bryant product.
As a result, we get the required comparison in a single spatial chart through
$s+r_{x,s}^2$.
\end{proof}

In Fang's Proposition~5.21, the cap error is measured on $Q_q$. In the proof
of his Proposition~6.11, Fang fixes larger domains for a further local
improvement. At every step of the cap iteration, we include these spatial, time
and collar enlargements in the error $\widetilde e_C$. In this way, we get
\eqref{secNN:cap-decay-plus} on the domains which are needed for that improvement.

\begin{proposition}[Fang's Proposition~5.21 on enlarged cap domains]\label{secNN:fang521}
Let \(o\in M\) be a fixed point and let \(p>0\) be a decay exponent. Then there exist a core threshold \(2\Lambda\),
fixed neck and cap measurement domains \(P_\rho(x,0)\) and \(Q_q=D_q\times[-Hr_q^2,0]\), together with the
collar \(Q_{q,\mathrm{collar}}\) of (5.6)--(5.7), a compact set \(K_p\subset M\), and a constant \(C_p<\infty\)
such that the following hold.
\begin{enumerate}
\item[(i)] For every \(x\in M\setminus(K_p\cup\mathcal{C}_0)\), we have
\[
e_N(x,0)\le C_p(1+d_g(o,x))^{-p}.
\]
\item[(ii)] For every cap center \(q\in\Gamma\setminus K_p\), we have
\[
e_C(q)\le C_p(1+d_g(o,q))^{-p}.
\]
\end{enumerate}
Here \(\Gamma\) and \(\beta_0\) are the ones given by Proposition~5.7 and (5.2), and we set
\(\mathcal{C}_0=\{x\in\mathcal T_0:\beta_0(x)<2\Lambda\}\).
By the neck error \(e_N(x,0)\) we mean the infimum on \(P_\rho(x,0)\) which is defined in Corollary~3.15.
By the cap error \(e_C(q)\) we mean the infimum (5.8) on \(Q_q=D_q\times[-Hr_q^2,0]\), where the normal and
Gram errors are taken on \(Q_{q,\mathrm{collar}}\) and on whole collar spheres. Both errors use
the common CMC family of Proposition~5.16. We choose all measurement sizes and \(\Lambda\) once for the
given \(p\), and they do not depend on \(x,q\). The constant \(C_p\) and the set \(K_p\) may depend on these
choices and on \(o\).

\par\smallskip\noindent\emph{Enlarged-domain conclusion.}
For each prescribed \(p>0\), we take the fixed enlarged space, time and collar domains which are chosen in
the proof below. Let \(\widetilde e_C\) be the error measured on these domains, and let \(K_p\) be the
corresponding compact set. Then, by the iteration, we have
\begin{equation*}
 \widetilde e_C(q)\le C_p(1+d_g(o,q))^{-p}
 \qquad(q\in\Gamma\setminus K_p).
 \tag{5.21+}\label{secNN:cap-decay-plus}
\end{equation*}
In each cap step, the same domains are used both as input and as output. Moreover, part~(ii) follows from
this bound by restriction.
\end{proposition}
\begin{proof}
We choose \(\sigma>0\) such that \(\log 2/\log(1+2\sigma)>p\). We fix the neck rule of Corollary~3.15 with
gain \(1/2\), and then we fix \(\Lambda,J,B,H\) as in Propositions~4.7 and 5.13. In particular, we have
\(B>4J\Lambda\). Let \(C_{\rm res}\ge1\) be the restriction constant of Proposition~5.17(ii). We remark that
all the enlargements below are fixed before \(E\) is chosen.

\par\smallskip\noindent\emph{A larger error that can be iterated.}\enspace
Let \(q\) be a tip. We write \(r_q=R(q)^{-1/2}\) and use the child time \(\tau=t/r_q^2\), and we freeze the
spatial labels \((\beta_0,z_q)\) at the final time. Here \(z_q\) denotes the signed arclength of the tip
projection from \(q\), divided by \(r_q\). For \(0\le j\le6\), we set
\[
 L_j=H+j,\quad d_j=B+j\Lambda/8,\quad
 b_j=2J\Lambda+j\Lambda/8,\quad a_j=\Lambda/4+\Lambda/2^{j+2},
\]
\[
 D_j=\{\beta_0\le d_j,\ |z_q|\le L_j\},\quad
 A_j=\{a_j\le\beta_0\le b_j,\ |z_q|\le L_j\},\quad I_j=[-L_j,0].
\]
Thus \(D_0,A_0,I_0\) are Fang's original sets, and (5.7) is given by the interiors of the first spatial
enlargements. We enlarge the fixed tube as permitted after (5.2). Let \(\mathcal S_j(\tau)\) be the set of
the actual whole CMC leaves which meet \(A_j\), and let \(\Sigma_j(\tau)\) be their union. In the definition
below, the infimum is taken over time-independent smooth triples \(U=(U_i)_{i=1}^3\) on a neighborhood of
\(D_4\), and we use the normalized metric and potential \((\widehat g,F)\) of (5.8). We define
\begin{align*}
 M(U)&=\sup_{D_3\times I_3}
 \left(\sum_i\left[\sum_{k=0}^2|\widehat\nabla^k\mathcal{L}_{U_i}\widehat g|^2
                  +|U_iF|^2\right]\right)^{1/2},\\
 N(U)&=\sup_{\tau\in I_3,\ x\in\Sigma_3(\tau)}
 \left(\sum_i [|\langle U_i,n_f\rangle|^2+
                   |\langle U_i,\nu\rangle|^2]\right)^{1/2},\\
 G(U)&=\sup_{\tau\in I_3,\ S\in\mathcal S_3(\tau)}
                  |Q_S(U)-I|_{\mathrm F},\qquad
 \widetilde e_C(q)=\inf_U\{M(U)+N(U)+G(U)\}.
\end{align*}
Here \(n_f\) denotes the potential-level normal and \(\nu\) denotes the leaf normal within that level, and
the Gram matrix is taken with the area-squared normalization. Note that in a fixed enlarged Bryant chart the
model leaves are spherical orbits, and each frozen label is constant on each such orbit at every time.
Hence, by the smooth, uniform graph construction in Propositions~2.4 and 5.16, the oscillation of both
frozen labels on each actual whole leaf can be made arbitrarily small. Since the gaps between consecutive
layers are positive, it follows that \(\Sigma_3(\tau)\subset\operatorname{int}D_4\), uniformly for
\(\tau\in I_3\). For the radial lower bound, we note that \(a_j>\Lambda/4\) and that the model radial
potential increases as we go backward in time. The same argument also applies to the CMC leaves which meet
the boundaries of the collars.

By restriction, we have \(e_C\le\widetilde e_C\). Moreover, by Proposition~5.17(ii), we have
\(e_N\le C_{\rm res}\widetilde e_C\) on the specified transition band. Since the construction is natural
under canonical pullback, we also have \(\widetilde e_C(q,t_*)=\widetilde e_C(\Phi_{t_*}q)\). Next we
consider any escaping tip sequence. Along such a sequence, we apply Propositions~5.7, 5.9 and 5.16 on a
neighborhood of \(D_4\times I_3\), and we transport the three area-normalized model rotations by the same
time-independent comparison map. Then their metric, potential, normal and Gram errors tend to zero uniformly
on the compact measurement regions. Applying this to subsequences, we conclude that \(\widetilde e_C(q)\to0\)
at infinity along \(\Gamma\).

\par\smallskip\noindent\emph{Supplying a field on each entire core-box domain.}\enspace
Work in parent tip units and take the geometric comparison on \(\widehat Q_{2E}\) of (4.6).
We glue on an interior mesh covering
\[
 Q_{2E}=\{B_t<2E,\ |Z|<4E^2,\ -4E^2<t\le0\},
 \qquad Q_E^+\subset Q_{2E}.
\]
Use the core row and the whole-leaf coordinates \((v,\zeta)\) of Proposition~4.2 on the annulus.
Here \(v,\zeta\) are the leaf averages of \(B_t,Z\). Choose the comparison accuracy,
for each fixed \(E\) and mesh spacing \(0<h<1\), so that both label errors are
less than \(h/16\), with the derivative control in the proof of that proposition.
Use spacing \(h\) in each grid
direction, with first annular center \(2\Lambda+h/2\) and final time row centered at zero.
Retain a full rectangular array of indices, ending each index range at the first center
at least one mesh step beyond the corresponding radial, axial or past face of \(Q_{2E}\). Keep the boxes and their
enlargements intact there. At time zero use their restrictions to \(t\le0\).
For example, intervals of half-width \(2h/3\) cover each grid direction and admit
cutoffs supported in intervals of half-width \(3h/4\). Adjacent intervals then share a
subinterval of length \(h/3\) even before enlargement. The core row is \(B_t<2\Lambda+h\).
These choices cover the transition to the first annular row and reserve the same overlap at
the outer rows as in the interior.

For a core box, let \(t_\alpha\) be its
final time and let \(Z_\alpha\) be its axial midpoint. At that time, the actual tip line is a small graph in
the chosen fixed chart, and we let \(q_\alpha\) be the point on it with axial coordinate \(Z_\alpha\). We
denote by \(D_j^m,A_j^m\) the layers which are defined with \(B_{t_\alpha}\) and \(Z-Z_\alpha\) on the model
chart. We choose
\[
 W_\alpha=\operatorname{int}D_2^m,\qquad T_\alpha=[t_\alpha-L_2,t_\alpha].
\]
These sets contain the minimum open spatial domain and the closed time interval on p.~47. In parent units \(B_t\le d_2+L_2\) there,
and their axial and past extents are at most \(4E^2+O(h)+L_2\). Thus they are inside
\(\widehat Q_{2E}\), whose axial and past extents are \(8E^2+1\), for all sufficiently large \(E\).
By Proposition~4.7(ii), the enlarged core box, its connecting paths and their
neighborhoods lie in the innermost layer. On core--core overlaps, we use an open family of common spheres
near radial label \(\Lambda\). On core--annulus overlaps, we use the common spheres near \(2\Lambda\) which
are given by the grid. These choices have radial and axial margins in the larger collar, and they include
the common spheres on triple overlaps.

Let \(s_\alpha=r_{q_\alpha,t_\alpha}/r_q\). For fixed \(E\), by the smooth tip and projection convergence,
we have \(s_\alpha\to1\) and \(\beta_{t_\alpha}-B_{t_\alpha}\to0\), and the child axial label converges to
\(Z-Z_\alpha\) on these domains. Note that in parent coordinates the child time is exactly
\((t-t_\alpha)/s_\alpha^2\). We make the label errors smaller than half of the adjacent-layer gaps, and we
require that \(s_\alpha^2>L_2/L_3\). Then we have
\[
 W_\alpha\times T_\alpha\ \subset\
 D_3(q_\alpha,t_\alpha)\times I_3(q_\alpha,t_\alpha)
\]
after identifying the time coordinates. By the same whole-leaf oscillation estimate as above, the common
collar leaves lie in \(\Sigma_3\). Moreover, the final times agree exactly. For fixed \(E\), all these
choices hold simultaneously, since we have smooth control on one enlarged chart.

At corresponding physical times, we have \(g'=s_\alpha^{-2}\bar g\) and
\(F'=s_\alpha^{-1}\bar F+\mathrm{const}(t)\), where a prime denotes child units and a bar denotes parent
units. For the same field \(V\), the same sphere \(S\) and the corresponding unit normal, we have
\begin{align*}
 |\nabla^k\mathcal{L}_V g'|_{g'}&=s_\alpha^k|\nabla^k\mathcal{L}_V\bar g|_{\bar g},
 & |V\bar F|&=s_\alpha|VF'|,\\
 |\langle V,n'\rangle_{g'}|&=s_\alpha^{-1}|\langle V,\bar n\rangle_{\bar g},
 & Q_{S,g'}(U)&=Q_{S,\bar g}(U).
\end{align*}
To see the last identity, we scale the inner product and the area element in the numerator and the area
squared in the denominator. Now suppose that \(C_{\rm res}\widetilde e_C(q_\alpha,t_\alpha)<\varepsilon\).
Then we choose \(U\) such that \(M+N+G<\varepsilon/C_{\rm res}\). By the inclusions and the identities
above, we get all three estimates on p.~77 with \(e_\alpha\le A_{\rm grid}\varepsilon\), where
\(A_{\rm grid}\) is independent of \(E\), since \(1/2<s_\alpha<2\). They also give control of
\(\sum_i|U_i\bar F|\) on the enlarged box. We remark that geometric closeness was used to obtain the
inclusions, and it does not add any term of order \(\eta\) to these error comparisons.

On an enlarged core box we have \(0\le B_t\le2\Lambda+O(h)\), so the weight \((B_t+c_{16})^{-16}\) has a
fixed positive lower bound there. Moreover, the size estimate for a triple in the proof of Proposition~4.1
applies here. In fact, the initial field bounds follow from its Gram and transverse bounds on a common
sphere, and by Killing transport along the bounded paths, three spatial derivatives are controlled by using
the first two derivatives of the metric Lie error. Thus we have the required \(O(\sqrt E)\) size bound.
\par\smallskip\noindent\emph{The backward neck domain for an annular box.}\enspace
We denote by \(T_\alpha\) the expanded time interval of the box
\(\mathcal B_\alpha^+\), which is truncated at zero if necessary, and we
set \(t_\alpha=\max T_\alpha\). On that time slice, we choose a point \(x_\alpha\)
on the whole CMC sphere which has the central radial and axial labels of the box.
Note that the central label of the first annular row is \(v=2\Lambda+h/2\).
We take the chart so accurate that
\(\lvert\beta_{t_\alpha}(x_\alpha)-v\rvert<h/4\) at the centers.
Thus every annular center lies outside the core at its chosen time, while the
first row still overlaps the core row \(B_t<2\Lambda+h\).

Let \(r_\alpha=R(x_\alpha,t_\alpha)^{-1/2}\), measured in the same units as
the grid. In parent tip units, the length of \(T_\alpha\) is \(O(h)\).
On the other hand, by the Bryant asymptotics, we have \(r_\alpha^2\asymp 1+v\).
Hence, by the fixed mesh, we get
\(T_\alpha\subset[t_\alpha-\rho r_\alpha^2,t_\alpha]\).
To get the spatial inclusion, we freeze the model labels at \(t_\alpha\).
It follows from (4.3) that \(B_t\) moves by at most \(t_\alpha-t\) on this interval,
and we note that \(Z\) is fixed. So we see that both frozen labels vary by
\(O(h)\) on the enlarged box. We also know that \(\lvert\nabla B\rvert\) is bounded below on the annulus,
and then the transverse distance to the central sphere is \(O(h)\) in tip units.
Note also that the model sphere diameters, measured in units of \(r_\alpha\), tend to
\(\pi\sqrt{2}\) as \(v\to\infty\). Now if \(\rho\ge100\) and \(\Lambda\)
is large, then we can see that inside \(B_{g(t_\alpha)}(x_\alpha,\rho r_\alpha)\)
there is room for the whole spheres on overlaps, for paths along the spheres and across the
transverse intervals, and also for a fixed positive multiple of \(r_\alpha\) around
those paths, which is what we need in Proposition~4.1. When \(E\) is fixed, it is clear that this
strict spatial margin is preserved by the smooth chart convergence
and by the whole-leaf construction of Proposition~5.16. So we can choose
the open field domain \(W_\alpha\) and the whole-leaf collar
\(\mathcal C_\alpha^+\) of Proposition~4.2 such that
\[
 \mathcal B_\alpha^+\cup\mathcal C_\alpha^+
 \ \subset\ W_\alpha\times T_\alpha
 \ \subset\ P_\rho(x_\alpha,t_\alpha)
 =B_{g(t_\alpha)}(x_\alpha,\rho r_\alpha)
       \times[t_\alpha-\rho r_\alpha^2,t_\alpha].
\]
For any adjacent pair or triple of annular boxes used in the grid, take a time and
transverse labels in their common intervals. The resulting whole sphere lies in each
field domain. On the model the centers differ transversely by \(O(h)\), their time
coordinates differ by \(O(h)\), and the entire sphere has diameter less than \(6r_\alpha\)
for large \(\Lambda\). Hence a collar of transverse width \(c r_\alpha\) about that
sphere and balls of radius \(2c r_\alpha\) about the connecting paths fit simultaneously
in all the spatial neck balls, for a fixed small \(c>0\). One may take \(W_\alpha\)
to contain these neighborhoods and the expanded box inside a ball of radius
\(20r_\alpha\) about \(x_\alpha\). On overlaps between core and annular boxes, use a common sphere near \(2\Lambda\)
in the intersection of the annular label intervals and the core collar. The
incident annular rows have centers within \(O(h)\) of \(2\Lambda\), and there
\(r_\alpha\asymp\sqrt{1+2\Lambda}\) is bounded in parent units, so the fixed gaps between the core
layers permit the same choice after decreasing \(c\). For core-only overlaps use
the collar near \(\Lambda\) and the paths in Proposition~4.7(ii). The common CMC
family and these strict margins persist under smooth comparison for fixed \(E\).
To check the exterior rows, observe that on an annular \(W_\alpha\times T_\alpha\), the extra spatial extent in either
\(B_t\) or \(Z\) is \(O(\sqrt E)+O(h)\), because \(r_\alpha\le C\sqrt E\),
\(B_t\) is spatially Lipschitz and (4.3) controls its time variation.
This fits between the radial bounds \(2E,4E+1\) and the axial bounds \(4E^2,8E^2+1\).
The additional backward time is \(O(h)\). All overlap constants in curvature
units are therefore independent of \(E\). For the extra radial rows,
apply Corollary~5.19 with its size parameter \(2E\). All their center labels
are at most \(2E+O(h)<4E\). The fixed enlargement of the axial and time ranges
is allowed in the proof of Proposition~5.18(ii). Corollary~5.20 includes all
the corresponding direct and descendant centers. This enlargement affects the
required chart and comparison accuracy for fixed \(E\), while the radial
weight constant supplied by Corollary~5.19 is independent of \(E\).
We restrict the triple of Corollary~5.19
with \(N=16\) to this domain. Then its weighted estimates hold through
the final face of the expanded box. For the transitions which are supplied by projected
tips, we use Proposition~5.17(ii).

\par\smallskip\noindent\emph{Weighted gluing on the interior mesh.}\enspace
Put \(w_{16}(b)=(b+c_{16})^{-16}\). We use the alignment argument in the proof of
\cite[Proposition~4.2]{F} on the rectangular index set just constructed.
For each pair of overlapping boxes Proposition~4.1 gives a constant orthogonal
transition. By (4.10) its aligned field difference in unscaled \(C^3\) is at most
\(C\sqrt E\,\varepsilon w_{16}(B_t)\). On a triangle the product of the three
transitions differs from the identity by at most \(C\varepsilon w_{16}(B_t)\),
by evaluating the differences on the common whole sphere and using its Gram lower bound.
On the core row the weight is comparable to a fixed positive constant.

There are \(O(E)\) radial rows and \(O(E^2)\) rows in each of the two other
directions. Choose a tree with inward radial spokes from the outermost radial row,
axial chains on that row, and one time chain at its axial endpoint. Transport the
bases along this tree. For an axial edge, sliding it to the outer row fills its tree loop with
one elementary face per radial step. For a time edge, first
slide it radially to the outer row, then axially to the time chain. These two
strips have at most \(C(E+E^2)\) faces in total. Both fillings stay in the retained
array, including for edges on its boundary. Split each face into two triangles.
The boxes at its corners share an interval in each varying coordinate, so the
preceding collar comparisons apply. Overlaps along a diagonal reduce to coordinate
edges inside one such cell. Every face used has radial label at least that of the
edge up to a fixed mesh error. Thus its weight is at most a fixed multiple of the
edge weight. Summing the triangle errors, using invariance under orthogonal
conjugation, gives a residual transition error at most \(CE^2\varepsilon w_{16}(B_t)\).
The \(C^3\) size bound \(C\sqrt E\) now yields a difference at most
\(CE^{5/2}\varepsilon w_{16}(B_t)\) between the fields in the transported bases.

On \(Q_{2E}\) normalize the fixed mesh cutoffs by their sum. The extra rows,
bounded multiplicity and label derivative bounds give a positive lower bound for the
denominator and uniform bounds through three spatial derivatives and one time
derivative. Use the core cutoff constant in the radial variable near the tip.
In the weighted sum of the aligned fields, each differentiated cutoff multiplies
a field difference after subtracting one active field. This proves (4.8) with
coefficient \(CE^{5/2}\varepsilon\) for the metric Lie derivatives and the time
derivative. The divergence--trace identity gives the vector heat operator bound,
and the potential pairing is the undifferentiated
weighted sum of the local pairings. For (4.9), compare each local triple to the
model rotations by Proposition~4.1 in the model metric, with error \(C_E\eta\)
since \(\varepsilon\le\eta\). Match the model basis at the tree root and propagate
along the finite tree. This gives a single orthogonal choice with \(C^3\) error
\(C_E\eta\) on \(Q_{2E}\). Only this model comparison may lose a constant depending
on the length of the tree. Restrict to \(Q_E^+\) and divide the fields by
\(\sqrt{6\pi}\). Since \(E^{5/2}\le E^4\), the weighted hypothesis of Corollary~4.6
holds with \(C_{\rm gl}\) independent of \(E\).

\par\smallskip\noindent\emph{Making the output control the same error.}\enspace
At the parent tip, we choose \(H_{\rm out}=L_4\), with one additional unit of backward time for the interior
estimates. We also choose a fixed open set \(D_{\rm out}\) containing \(D_4^m\) compactly inside
\(\operatorname{int}D_5^m\), together with a further open enlargement inside \(\operatorname{int}D_6^m\).
We smooth the boundaries between these nested sets if it is necessary. As the output collar, we take the family
of whole leaves which are labeled by the interior of \(A_4^m\), and its outer margin is in layer \(5\). We can also
choose the larger spatial domain so that it contains the closed outer collar with margin. Note that \(b_5<d_5\),
and so the radial margin is available as well. Then we see from the label estimates above that \(D_4(q)\) lies in \(D_{\rm out}\)
and that every leaf which meets \(A_3(q)\) lies in the output collar, and this is true simultaneously for
\(\tau\in I_3\).

We know that the model satisfies \(0\le B_t-B_0\le-t\) on \([-H_{\rm out}-1,0]\). It follows that the radial range of the larger
collar is contained in \([a_5,b_5+H_{\rm out}+1]\), and that \(|Z|\le L_5\). Recall that \(\Lambda\) was chosen
large at the beginning, and then the radial distance corresponding to \(a_5\) is beyond the regular-collar threshold.
Now we choose the parameters \(s_-,s_+,Z_*,\delta\) of Corollary~4.6 around this compact range with positive
margins, and we note that they are preserved by smooth comparison, uniformly in the tilt. Note that these output sets are
fixed in advance, and they lie inside \(Q_E\) when \(E\) is large. By Corollary~4.6, we get one
time-independent, Gram-normalized triple \(Z\) on a neighborhood of \(D_4(q)\). From its \(C^8\) output
bounds, we obtain
\[
 \widetilde e_C(q)\le M(Z)+N(Z)+G(Z)
           \le (CE^{-9}+C_E\eta)\varepsilon .
\]
Here the constant \(C\) depends on the fixed domains and on \(C_{\rm gl}\), but it does not depend on \(E\)
or on the tilt. In particular, layer \(3\) is used on both sides of every subsequent cap improvement.

\par\smallskip\noindent\emph{Iteration.}\enspace
We first choose \(E\) and then \(\eta\) such that \(C_{\rm res}CE^{-9}<1/4\) and
\(C_{\rm res}C_E\eta<1/4\). For fixed \(E\), by Proposition~\ref{secNN:fang515}, all intermediate, terminal,
projected and direct grid centers of one construction lie in a fixed parent chart. By a fixed enlargement,
we include the extra measurement neighborhoods in the chart. This is possible because the model scalar has
a positive minimum on the compact set of centers and each added neighborhood has finite size in its own tip
units. Note that no comparison centers are added to the recurrence in this way.
On \((M\setminus\mathcal{C}_0)\sqcup\Gamma\), we define \(\mathcal E=e_N\)
on the neck copy and \(\mathcal E=C_{\rm res}\widetilde e_C\) on the
cap copy. The assigned scales are \(Q=R\) on the neck copy and \(Q=\rho_*R\) on the cap copy.
By canonical pullback, we take the parent to be \((x,0)\). For each required
center \((y,t)\), where we use the projected tip for a core replacement,
Propositions~5.13--5.14 and Corollary~5.20 give fixed constants \(L_0,T_*>0\) such that
\[
 d_{g(0)}(x,y)\le L_0 Q(x)^{-1/2},\qquad
 -T_*Q(x)^{-1}\le t\le0,\qquad Q(y,t)\le(1+\sigma)Q(x).
\]
We remark that in a cap construction every required center, including the intermediate
vertices, is compared directly with \(x\) in these inequalities.
Since \(\rho_*R(x)\le Q(x)\le R(x)\), we see that the frozen ball and the time interval
above satisfy the hypotheses of Proposition~5.12 with constants
\(L_0/\sqrt{\rho_*}\) and \(T_*/\rho_*\) in scalar-curvature units.
By that proposition, the \(g(t)\)-length of a minimizing
\(g(0)\)-geodesic from \(x\) to \(y\) is bounded. Thus, after enlarging the fixed
distance constant to \(L_*\) and the compact exclusion if necessary, we have
\[
 d_{g(t)}(x,y)\le L_*Q(x)^{-1/2}.
\]
Note that in the edge estimates only the distance constant has been enlarged.
In particular, the factor in the assigned-scale bound is still \(1+\sigma\).
For the original-soliton representative \(z=\Phi_t(y)\), we have
\(Q(z)=Q(y,t)\). As in Lemma~5.1, using the isometry \(\Phi_t\) and
\(\lvert\nabla f\rvert\le1\), we get
\[
 d_g(x,z)\le d_g(x,\Phi_t(x))+d_{g(t)}(x,y)
       \le T_*Q(x)^{-1}+L_*Q(x)^{-1/2}.
\]
We denote by \(Z(x)\) the set of all these representatives, where each representative keeps its neck or cap
type. Note that all errors are uniformly small outside a compact set, and all representatives which are used
by a sufficiently distant parent avoid that set. We write \(S=\sup_{z\in Z(x)}\mathcal E(z)<\eta\). For each
\(S<\varepsilon<\eta\), we choose approximate minimizers with the required strict bounds. For a neck, we
apply Corollary~3.15, and for the neck steps within the cap construction we use Proposition~5.18 and
Corollary~5.19. For a cap, we apply the preceding output estimate. Letting \(\varepsilon\downarrow S\),
where the case \(S=0\) is included, we obtain
\[
 \mathcal E(x)\le\tfrac12\sup_{z\in Z(x)}\mathcal E(z).
\]
Now we apply Corollary~5.3 and Proposition~5.4 in their assigned-scale form with \(\bar C=1+2\sigma\). We
get \(K_j\subset B(o,D_*\bar C^j)\) and \(\mathcal E\le 2^{-j}\varepsilon_*\) outside \(K_j\). For
\(D_*\bar C^j<d_g(o,x)\le D_*\bar C^{j+1}\), this gives decay with exponent \(\log2/\log\bar C>p\). Since
\(\mathcal E=C_{\rm res}\widetilde e_C\) on the cap copy, this estimate proves
\eqref{secNN:cap-decay-plus}. Then (ii) follows by restriction, and (i) follows from the neck copy.
\end{proof}

For the exterior atlas, \eqref{secNN:cap-decay-plus} supplies low-order input
on each enlarged core-box domain, including its earlier time face and whole
collar leaves. We retain the output of one further local improvement at each
atlas center.

\begin{proposition}[Fang, Proposition~6.11]\label{secNN:fang611}
Let $w(x)=1+d_g(o,x)$. Then for each sufficiently large $N$ and each
$q>0$, the exterior atlas of \cite[Proposition~6.9]{F} exists with
$b=1$ and with transitions in $SO(3)$. As in \cite[Proposition~6.11]{F},
the triples of this atlas have metric and potential defects of order
$O(w^{-N})$ through eight derivatives and cyclic bracket defects of order
$O(w^{-q})$ through three derivatives, and their normal and Gram errors
are of order $O(w^{-q})$ on whole collar spheres.
\end{proposition}
\begin{proof}
We fix the output domains with the margins which are required for eight
derivatives and for whole collar leaves. Consider a comparison centered
at $x$. It is shown in \cite[Proposition~6.8]{F} that every intermediate center or
projected tip has a representative $z$ in the original soliton such that
$w(z)\ge c_* w(x)^{1/2}$. Now we apply Proposition~\ref{secNN:fang521}
with exponent $P$. Then for the input errors at the neck centers and at the
cap centers, we have respectively
\[
 e_N(z,0)\le C_P w(x)^{-P/2},\qquad
 C_{\rm res}\widetilde e_C(z)\le C_P w(x)^{-P/2}.
\]
Choose the core-box input triples within $w(x)^{-P/2}$ of the respective
infima in \eqref{secNN:cap-decay-plus}, with the fixed restriction factor
absorbed in the constant. At annular boxes use the weighted neck triples
of Corollary~5.19. At a cap center glue these inputs on the interior mesh
in the proof of Proposition~\ref{secNN:fang521}, restrict to $Q_E^+$,
and apply Corollary~4.6 after the bracket normalization. Choose its
time-independent, Gram-normalized output $Z$ as the cap atlas triple.
At a neck center choose the corresponding output of one cylindrical
improvement. Write $U_i$ for the output triple on chart $i$. On the
reserved output domains the metric Lie derivatives and potential pairings
of $U_i$ are $O(w(x)^{-P/2})$ through eight normalized spatial derivatives.
Claim~3.14 and its cap counterpart give the cyclic bracket bound through
three derivatives for these same triples, with their normal and Gram bounds
on whole collar spheres.

On an overlap, Proposition~4.1 aligns the two output triples on a common
whole leaf by a constant orthogonal matrix and gives a $C^3$ bound for
their aligned difference along the controlled paths. For each component
$D$ of this difference, $h=\mathcal L_D g$ has the eight derivative bounds
by linearity and constancy of the matrix. In curvature units differentiate
the Killing transport identity
\[
 \nabla^2D=\mathrm{Rm}*D+\mathcal A(\nabla h)
\]
through order seven, where $\mathcal A$ is a fixed linear combination of
index permutations. Induction using the bounded curvature derivatives
gives $D=O(w(x)^{-P/2})$ through order nine, hence the corresponding
pairs $(D,\operatorname{skew}\nabla D)$ through order eight. The curvature
length lies between $1$ and $Cw(x)$, so conversion to the original metric
costs at most $w(x)^{b_*}$ for a fixed $b_*$ determined by the derivative
orders. Taking $P/2>\max\{N,q\}+b_*$ gives the required exponents.
Finally, as in \cite[pp.~92--93]{F}, it follows from
the Gram lower bound and the common nonzero cyclic bracket coefficient
that the aligning orthogonal matrices have positive determinant.
\end{proof}

\subsection{Global angular action}\label{secNN:global-angular}

In this subsection, we apply Fang's gluing and drift correction to the
rotation fields of Proposition~\ref{secNN:fang611}. The convex levels
which are needed for the gluing come from Lemma~\ref{secNN:exhaustion},
and we use Lemma~\ref{secNN:vanishing} to show that the corrected fields
are Killing. Then we use the isotropy representation at a fixed point to
find the axis needed in Proposition~\ref{secNN:axis}.

\begin{proposition}[Angular symmetry under nonnegative sectional curvature]\label{app:angular-action}
Let $(M^4,g,f)$ be a complete $\kappa$-noncollapsed  steady gradient Ricci soliton
with $\sec\ge0$ and $\Ric>0$.
Assume that $(M,g)$ is not the four-dimensional Bryant soliton up to
scaling. Then there is an effective isometric $O(3)$ action on $M$ which
preserves $f$. Moreover, in suitable global coordinates
$(z,v)\in\mathbb R\oplus\mathbb R^3$, this action is given by
$A\cdot(z,v)=(z,Av)$, and its fixed set is the axis $v=0$.
\end{proposition}
\begin{proof}
We use the exhaustion $\psi$ of Lemma~\ref{secNN:exhaustion} in Fang's
construction \cite[Section~6]{F}. Since $\operatorname{Hess}\psi>0$, the
inward flow contracts. Integrating the Hessian of $\psi$ along a minimizing
geodesic which starts from the minimum of $\psi$, we get a lower bound for
the outward derivative of $\psi$ outside a fixed ball. By proper Morse
theory, every positive level is a sphere. The second fundamental form of
such a level is $\operatorname{Hess}\psi|_{T\Sigma}/|\nabla\psi|>0$. By
the Gauss equation, the intrinsic sectional curvature of the level is
strictly positive, even when an ambient plane is flat. These facts are the
convexity, level geometry and topology inputs in Propositions 6.1--6.7.
Moreover, volume comparison needs only $\operatorname{Ric}\ge0$, and the
local charts are available from Lemma 2.1. Next we apply Fang's gluing
result \cite[Propositions~6.9 and~6.12]{F} to the atlas of
Proposition~\ref{secNN:fang611}. We obtain three exterior fields $U_i$
which, together with their first derivatives, have at most linear growth.
The metric and potential defects of these fields and the first two
derivatives of these defects are $O(w^{-p})$, where $p>9$ can be chosen
arbitrarily. Moreover, the cyclic bracket defects of these fields and
their two normal components on escaping whole leaves tend to zero, and
their normalized leaf Gram matrices tend to a positive multiple of the
identity. Note that such whole leaves escape along a fixed collar of
escaping cap centers, or in necks if the cap centers are bounded. This is
because $R^{-1/2}=o(d(o,\cdot))$ by Proposition 5.6. The latter estimate
follows directly by dimension reduction. In fact, if it fails, then after
rescaling along an offending sequence with $R(x_i)\to0$, the fixed point
$o$ lies in a bounded smooth comparison ball, while its rescaled scalar
curvature satisfies $R(o)/R(x_i)\to\infty$. This is a contradiction.

\paragraph{Exact fields and isotropy.}
Note that Fang's Proposition 7.2 uses bounded local geometry,
$\operatorname{Ric}\ge0$, $R(o)>0$ and the soliton identities. In
particular, we know that the operator $\Delta+\nabla_X$ is conjugated by $e^{-f/2}$ to
$\Delta-(1+R)/4$. In the Green estimate, we have
$f(x)\ge f(o)-|X(o)|d(o,x)$ with $|X(o)|<1$, which follows from convexity. Now we extend $U_i$ across the
core and subtract the decaying drift corrections of that proposition. Let
$W_i$ be a corrected field. Then we see that the tensor $\mathcal L_{W_i}g$ solves the
equation in Lemma~\ref{secNN:vanishing}. Its decay exponent is
$p-5>4$, and so it is integrable by volume comparison, and it follows that it vanishes.
Note that the gradient of each $W_i f$ is parallel and that $\operatorname{Ric}>0$. Then, as
in Corollary 7.4, we find that this function is constant, and we also know that its limit at infinity is
zero. Now the bracket defects of the $W_i$ are decaying Killing fields, and they
satisfy the homogeneous drift equation. So we get the exact cyclic brackets by the scalar maximum
principle, and the linear independence
follows from the escaping Gram bound. Finally, we integrate these complete Killing
fields and get an isometric $SU(2)$ action.

Next we average $\psi$ over this compact group. It is easy to see that the average is proper and
has positive Hessian, and that its unique minimum $s$ is fixed. As in
Proposition 7.5, we can use equivariant Morse coordinates and then the
outward level flow to identify $M$ equivariantly with $T_sM$. We know that the
nontrivial real four-dimensional representation is either quaternionic or
$\mathbb{R}\oplus\mathbb{R}^3$ with the adjoint action on $\mathbb{R}^3$,
and we need to exclude the former. To do this, we take an escaping sequence
of leaf points $x_j$ and choose unit coefficient vectors $v_j$ in the
kernel of the tangential projection of the approximate triple. Note that the
normal components of the approximate triple vanish in the limit and that
$W_i-U_i\to0$. It follows that $|W_{v_j}(x_j)|\to0$. In the case of the quaternionic action,
we have $X(s)=0$ by invariance, and then we see that $f$ is proper by strict convexity and by
integrating $\operatorname{Ric}$ near $s$.
We also see that the action is free away from $s$. Then the inward gradient flow takes
$x_j$ to a fixed compact regular level. Note that $[X,W_v]=0$ and
$X|W_v|^2=2\operatorname{Ric}(W_v,W_v)\ge0$, and so the norm also tends to zero
there. On the other hand, since the action is free and the level times the
coefficient unit sphere is compact, we get a positive lower bound. This is
a contradiction. Thus the adjoint representation with a trivial axis
occurs. Note that the central element of $SU(2)$ acts trivially. In the global
linear coordinates $(z,v)\in\mathbb{R}\oplus\mathbb{R}^3$, because of the
$SO(2)$ stabilizer, the sphere metric must be round up to a scalar and
there are no mixed angular--horizontal terms. Consequently, the normal
inversion $(z,v)\mapsto(z,-v)$ is an isometry, and we get the standard
angular $O(3)$ action with a fixed axis. Finally, the original potential
is constant on these orbits.
\end{proof}

\end{document}